\documentclass[12pt]{article}

\usepackage[dvips]{graphicx}
\usepackage{textcomp}
\usepackage{amsbsy}
\usepackage{latexsym}
\usepackage[mathscr]{eucal}
\usepackage{amsfonts}
\usepackage{amssymb}
\usepackage[usenames]{color}

\usepackage{fullpage}

\begin{document}
\bibliographystyle{plain}
\title
{
Approximation of high-dimensional parametric PDEs
}
\author{ 
Albert Cohen and Ronald DeVore
\thanks{%
  This research was supported by the ONR contracts  
N00014-11-1-0712 and N00014-12-1-0561, the  NSF grant  DMS 1222715,
the Institut Universitaire de France and the ERC advanced grant BREAD.}  }
\hbadness=10000
\vbadness=10000
\newtheorem{lemma}{Lemma}[section]
\newtheorem{prop}[lemma]{Proposition}
\newtheorem{cor}[lemma]{Corollary}
\newtheorem{theorem}[lemma]{Theorem}
\newtheorem{remark}[lemma]{Remark}
\newtheorem{example}[lemma]{Example}
\newtheorem{definition}[lemma]{Definition}
\newtheorem{proper}[lemma]{Properties}
\newtheorem{assumption}[lemma]{Assumption}
%
\newenvironment{disarray}{\everymath{\displaystyle\everymath{}}\array}{\endarray}

\def\RR{\rm \hbox{I\kern-.2em\hbox{R}}}
\def\NN{\rm \hbox{I\kern-.2em\hbox{N}}}
\def\ZZ{\rm {{\rm Z}\kern-.28em{\rm Z}}}
\def\CC{\rm \hbox{C\kern -.5em {\raise .32ex \hbox{$\scriptscriptstyle
|$}}\kern
-.22em{\raise .6ex \hbox{$\scriptscriptstyle |$}}\kern .4em}}
\def\vp{\varphi}
\def\<{\langle}
\def\>{\rangle}
\def\t{\tilde}
\def\i{\infty}
\def\e{\varepsilon}
\def\sm{\setminus}
\def\nl{\newline}
\def\o{\overline}
\def\wt{\widetilde}
\def\wh{\widehat}
\def\cT{{\cal T}}
\def\cA{{\cal A}}
\def\cI{{\cal I}}
\def\cV{{\cal V}}
\def\cB{{\cal B}}
\def\cF{{\cal F}}

\def\cR{{\cal R}}
\def\cD{{\cal D}}
\def\cP{{\cal P}}
\def\cJ{{\cal J}}
\def\cM{{\cal M}}
\def\cO{{\cal O}}
\def\Chi{\raise .3ex
\hbox{\large $\chi$}} \def\vp{\varphi}
\def\lsima{\hbox{\kern -.6em\raisebox{-1ex}{$~\stackrel{\textstyle<}{\sim}~$}}\kern -.4em}
\def\lsim{\hbox{\kern -.2em\raisebox{-1ex}{$~\stackrel{\textstyle<}{\sim}~$}}\kern -.2em}
\def\[{\Bigl [}
\def\]{\Bigr ]}
\def\({\Bigl (}
\def\){\Bigr )}
\def\[{\Bigl [}
\def\]{\Bigr ]}
\def\({\Bigl (}
\def\){\Bigr )}
\def\L{\pounds}
\def\pr{{\rm Prob}}
\newcommand{\cs}[1]{{\color{magenta}{#1}}}
\def\ds{\displaystyle}
\def\ev#1{\vec{#1}}     
\newcommand{\lt}{\ell^{2}(\nabla)}
\def\Supp#1{{\rm supp\,}{#1}}
\def\R{\mathbb{R}}
\def\E{\mathbb{E}}
\def\nl{\newline}
\def\T{{\relax\ifmmode I\!\!\hspace{-1pt}T\else$I\!\!\hspace{-1pt}T$\fi}}
\def\N{\mathbb{N}}
\def\Z{\mathbb{Z}}
\def\P{\mathbb{P}}
\def\N{\mathbb{N}}
\def\L{\mathbb{L}}
\def\Zd{\Z^d}
\def\Q{\mathbb{Q}}
\def\C{\mathbb{C}}
\def\Rd{\R^d}
\def\sb{\subset}
\def\gsim{\mathrel{\raisebox{-4pt}{$\stackrel{\textstyle>}{\sim}$}}}
\def\sime{\raisebox{0ex}{$~\stackrel{\textstyle\sim}{=}~$}}
\def\lsim{\raisebox{-1ex}{$~\stackrel{\textstyle<}{\sim}~$}}
\def\div{\mbox{ div }}
\def\M{M}  \def\NN{N}                  
\def\Le{{\ell^1}}            
\def\Lz{{\ell^2}}
\def\Let{{\tilde\ell^1}}     
\def\Lzt{{\tilde\ell^2}}
\def\Ltw{\ell^\tau^w(\nabla)}
\def\t#1{\tilde{#1}}
\def\la{\lambda}
\def\La{\Lambda}
\def\ga{\gamma}
\def\BV{{\rm BV}}
\def\Ga{\eta}
\def\al{\alpha}
\def\cZ{{\cal Z}}
\def\cA{{\cal A}}
\def\cU{{\cal U}}
\def\argmin{\mathop{\rm argmin}}
\def\argmax{\mathop{\rm argmax}}
\def\prob{\mathop{\rm prob}}
\def\A{\mathop{\rm Alg}}

\def \bphi{{\bf\phi}}

\def\cO{{\cal O}}
\def\cA{{\cal A}}
\def\cC{{\cal C}}
\def\cS{{\cal F}}
\def\bu{{\bf u}}
\def\bz{{\bf z}}
\def\bZ{{\bf Z}}
\def\bI{{\bf I}}
\def\cE{{\cal E}}
\def\cD{{\cal D}}
\def\cG{{\cal G}}
\def\cI{{\cal I}}
\def\cJ{{\cal J}}
\def\cM{{\cal M}}
\def\cN{{\cal N}}
\def\cT{{\cal T}}
\def\cU{{\cal U}}
\def\cV{{\cal V}}
\def\cW{{\cal W}}
\def\cL{{\cal L}}
\def\cB{{\cal B}}
\def\cG{{\cal G}}
\def\cK{{\cal K}}
\def\cS{{\cal S}}
\def\cP{{\cal P}}
\def\cQ{{\cal Q}}
\def\cR{{\cal R}}
\def\cU{{\cal U}}
\def\bL{{\bf L}}
\def\bl{{\bf l}}
\def\bK{{\bf K}}
\def\bC{{\bf C}}
\def\X{X\in\{L,R\}}
\def\ph{{\varphi}}
\def\D{{\Delta}}
\def\H{{\cal H}}
\def\bM{{\bf M}}
\def\bx{{\bf x}}
\def\bj{{\bf j}}
\def\bG{{\bf G}}
\def\bP{{\bf P}}
\def\bW{{\bf W}}
\def\bT{{\bf T}}
\def\bV{{\bf V}}
\def\bv{{\bf v}}
\def\bt{{\bf t}}
\def\bz{{\bf z}}
\def\bw{{\bf w}}
\def \span{{\rm span}}
\def \meas {{\rm meas}}
\def\rhom{{\rho^m}}
\def\diff{\hbox{\tiny $\Delta$}}
\def\EE{{\rm Exp}}
\def\lll{\langle}
\def\argmin{\mathop{\rm argmin}}
\def\argmax{\mathop{\rm argmax}}
\def\dJ{\nabla}
\newcommand{\ba}{{\bf a}}
\newcommand{\bb}{{\bf b}}
\newcommand{\bc}{{\bf c}}
\newcommand{\bd}{{\bf d}}
\newcommand{\bs}{{\bf s}}
\newcommand{\bff}{{\bf f}}
\newcommand{\bp}{{\bf p}}
\newcommand{\bg}{{\bf g}}
\newcommand{\by}{{\bf y}}
\newcommand{\br}{{\bf r}}
\newcommand{\be}{\begin{equation}}
\newcommand{\ee}{\end{equation}}
\newcommand{\bea}{$$ \begin{array}{lll}}
\newcommand{\eea}{\end{array} $$}
\def \Vol{\mathop{\rm  Vol}}
\def \mes{\mathop{\rm mes}}
\def \Prob{\mathop{\rm  Prob}}
\def \exp{\mathop{\rm    exp}}
\def \sign{\mathop{\rm   sign}}
\def \sp{\mathop{\rm   span}}
\def \vphi{{\varphi}}
\def \csp{\overline \mathop{\rm   span}}
\def \cost{\mathop{\rm   cost}}
 
%
%
\newcommand{\beqn}{\begin{equation}}
\newcommand{\eeqn}{\end{equation}}
\def\beginproof{\noindent{\bf Proof:}~ }
\def\endproof{\hfill\rule{1.5mm}{1.5mm}\\[2mm]}

\newenvironment{Proof}{\noindent{\bf Proof:}\quad}{\endproof}

\renewcommand{\theequation}{\thesection.\arabic{equation}}
\renewcommand{\thefigure}{\thesection.\arabic{figure}}

\makeatletter
\@addtoreset{equation}{section}
\makeatother

\newcommand\abs[1]{\left|#1\right|}
\newcommand\clos{\mathop{\rm clos}\nolimits}
\newcommand\trunc{\mathop{\rm trunc}\nolimits}
\renewcommand\d{d}
\newcommand\dd{d}
\newcommand\diag{\mathop{\rm diag}}
\newcommand\dist{\mathop{\rm dist}}
\newcommand\diam{\mathop{\rm diam}}
\newcommand\cond{\mathop{\rm cond}\nolimits}
\newcommand\eref[1]{{\rm (\ref{#1})}}
\newcommand{\iref}[1]{{\rm (\ref{#1})}}
\newcommand\Hnorm[1]{\norm{#1}_{H^s([0,1])}}
\def\int{\intop\limits}
\renewcommand\labelenumi{(\roman{enumi})}
\newcommand\lnorm[1]{\norm{#1}_{\ell^2(\Z)}}
\newcommand\Lnorm[1]{\norm{#1}_{L_2([0,1])}}
\newcommand\LR{{L_2(\R)}}
\newcommand\LRnorm[1]{\norm{#1}_\LR}
\newcommand\Matrix[2]{\hphantom{#1}_#2#1}
\newcommand\norm[1]{\left\|#1\right\|}
\newcommand\ogauss[1]{\left\lceil#1\right\rceil}
\newcommand{\QED}{\hfill
\raisebox{-2pt}{\rule{5.6pt}{8pt}\rule{4pt}{0pt}}%
  \smallskip\par}
\newcommand\Rscalar[1]{\scalar{#1}_\R}
\newcommand\scalar[1]{\left(#1\right)}
\newcommand\Scalar[1]{\scalar{#1}_{[0,1]}}
\newcommand\Span{\mathop{\rm span}}
\newcommand\supp{\mathop{\rm supp}}
\newcommand\ugauss[1]{\left\lfloor#1\right\rfloor}
\newcommand\with{\, : \,}
\newcommand\Null{{\bf 0}}
\newcommand\bA{{\bf A}}
\newcommand\bB{{\bf B}}
\newcommand\bR{{\bf R}}
\newcommand\bD{{\bf D}}
\newcommand\bE{{\bf E}}
\newcommand\bF{{\bf F}}
\newcommand\bH{{\bf H}}
\newcommand\bU{{\bf U}}
\newcommand\cH{{\cal H}}
\newcommand\sinc{{\rm sinc}}
\def\enorm#1{| \! | \! | #1 | \! | \! |}

\newcommand{\dm}{\frac{d-1}{d}}

\let\bm\bf
\newcommand{\bbeta}{{\mbox{\boldmath$\beta$}}}
\newcommand{\bal}{{\mbox{\boldmath$\alpha$}}}
\newcommand{\bbi}{{\bm i}}

\def\nnew{\color{Red}}
\def\mnew{\color{Blue}}

\newcommand{\dI}{\Delta}
\maketitle
\date{}

\begin{abstract}
Parametrized families of PDEs arise in various contexts such
as inverse problems, control and optimization, risk assessment,
and uncertainty quantification. In most of these applications, the number of parameters is large or perhaps even infinite. Thus, the development of numerical methods for these parametric problems is  faced with the possible curse of dimensionality. 
This article is directed at (i) identifying and understanding which 
properties of parametric equations 
allow one to avoid this curse and (ii) developing and analyzing 
effective numerical methodd which fully exploit 
these properties and, in turn, are immune to the growth in dimensionality.

The first part of this article studies the smoothness and approximability of the
solution map, that is, the map $a\mapsto u(a)$ 
where $a$ is the parameter value and $u(a)$ is
the corresponding solution to the PDE.
It is shown that for many relevant parametric PDEs, the 
parametric smoothness of this map is typically holomorphic and also highly
anisotropic in that the relevant parameters 
are of widely varying importance in describing the solution.  
These two properties are then exploited to establish convergence rates of 
$n$-term approximations to the solution map for which each term is separable in 
the parametric and physical variables. These results reveal that, 
at least on a theoretical level, the solution map
can be well approximated by discretizations of moderate complexity, 
thereby showing how the curse of dimensionality is broken.  
This theoretical analysis is carried out through concepts of approximation theory
such as best $n$-term approximation, sparsity, and  $n$-widths.  
These notions determine a priori the best possible performance of numerical
methods and thus serve as a benchmark for concrete
algorithms.

The second part of this article turns to the development of numerical algorithms based on the theoretically established sparse separable approximations.  
The numerical methods studied fall into two general categories.
The first uses polynomial expansions 
in terms of the parameters to approximate the solution map.  The second one
searches for suitable low dimensional spaces for 
simultaneously approximating all members of the parametric family.
The numerical implementation of these approaches is carried out through adaptive and greedy algorithms.   An a priori analysis of the performance of these algorithms  establishes  how well they  meet the theoretical benchmarks. 
\end{abstract}

\section{Overview}\label{S:intro}

\subsection{Parametric and stochastic PDEs}
\label{SS:parpdes}

  Partial differential equations (PDEs) are   commonly used to model complex systems
in a variety of physical contexts. When solving a given PDE, one typically fixes
certain {\it parameters}: the shape of the physical domain, the diffusion or velocity field, the source term, the flux or reaction law, etc.
We use the terminology  {\it parametric PDEs} when some of these parameters are allowed to vary  over a certain range of interest.  When treating such parametric PDEs, one is interested in finding the solution for all parameters in the range of interest.

To describe such problems in their full generality, we adopt the formulation
\be
\cP(u,a)=0,
\label{genpar}
\ee
where $a$ denotes the parameters, $u$ is the unknown of the problem, and
\be
\cP:  V \times X \to W,
\ee
is a linear or nonlinear partial differential operator, with $(X,V,W)$ a triplet of  Banach spaces.  
We assume that the parameter $a$ ranges  over a compact set $\cA\subset X$,
and that for any $a$ in this range there exists a unique solution $u=u(a)\in V$
to \iref{genpar}. This allows us to define the {\it solution map}
\be
u: a\mapsto u(a),
\label{solmap}
\ee
which acts from $X$ onto $V$ and is well defined over $\cA$. We also define
the {\it solution manifold} as the family
\be
\cM=u(\cA)=\{u(a) \; : \; a\in\cA\},
\ee
which gathers together all solutions as the parameter varies within its range.    

One simple guiding example, which will be often used throughout this article, is the linear elliptic equation
\begin{eqnarray}
-{\rm div}(a\nabla u)&=&f,\quad {\rm on} \ D,\nonumber\\
u&=&0\,\quad {\rm on} \ \partial D,
\label{ellip}
\end{eqnarray}
set on a Lipschitz domain $D\subset \R^m$.
Here, we fix the right  side $f$ as a real valued function and consider   real valued    diffusion coefficients $a$ as the parameter. The corresponding
operator $\cP$ is therefore given by
\be
\cP(u,a)=f+ {\rm div}(a\nabla u).
\label{Pellip}
\ee
A possible choice for the triplet of spaces is then
\be
(X,V,W)=(L^\infty(D),H^1_0(D),H^{-1}(D)).
\ee
Indeed, if $u\in V$, $a\in X$ and $f\in W$, one then defines 
$\cP(u,a)$ as an element of $W$ by requiring that 
\be
\<\cP(u,a),v\>=\<f,v\>+\int_D a\nabla u\cdot\nabla v,\quad v\in V,
\ee
where $\<\cdot,\cdot\>$ is the duality bracket between $V'=W$ and $V$. Lax-Milgram theory ensures the
existence and uniqueness of a solution $u(a)$ to \eref{genpar}  from $V$, if for some $r>0$, the diffusion $a$ satisfies the ellipticity condition
\be
a(x)\geq r,\quad x\in D.
\ee
Therefore, a typical parameter range is a set $\cA\subset \{a\in L^\infty(D)\; : \; a\geq r\}$, which in addition
is assumed to be compact in $L^\infty$.

Although elementary, the above example gathers important features
that are present in other relevant examples of parametric PDEs. In particular,
the solution map $a\mapsto u(a)$ acts from
an infinite dimensional space into another infinite dimensional space.
  Also note that, while the operator 
 $\cP$ of \iref{Pellip} is linear both in $a$ and $u$ (up to the constant
additive term $f$) the solution map is {\it nonlinear}.  Because of the high dimensionality of the parameter space $X$, such problems  represent a significant challenge when trying to capture
this map numerically.   One objective of this article is to understand which properties of this map allow for a successful
numerical treatment.  Concepts such as holomorphy, sparsity, and adaptivity are at the heart of our development.  

The solution map may also be viewed as a function
\be
(x,a)\mapsto u(x,a)
\ee 
of both   the physical variable $x\in D$ and the parametric 
variable $a\in \cA$.
The parametric variable has a particular status because
its different instances are uncoupled: for any fixed instance $a=a_0$, we may solve the PDE 
exactly or approximately and therefore compute $u(x,a_0)$ for all values of $x$,
while ignoring all other values $a\neq a_0$. This plays an 
important role for certain numerical methods 
which are based on solving the parametric PDE for different particular values of $a$, 
since this task can be parallelized.

Parametric PDEs occur in a variety of modeling contexts.    
We draw the following major distinctions in their setting:
\begin{itemize}
\item
{\bf Deterministic modeling:} the parameters are deterministic {\it design or control} variables
which may be tuned by the user so that the solution $u$, or some
quantity of interest $Q(u)$, has   prescribed
properties. For instance, if the elliptic equation \iref{ellip} is used to model the
heat conduction in a component produced by an industrial process,
one may want to design the material in order to minimize the heat flux 
on a certain part on the boundary $\Gamma\subset\partial D$, in which case
the quantity of interest to be optimized is
\be
Q(u)=\int_\Gamma a\nabla u\cdot {\bf n},
\label{qofu}
\ee
where ${\bf n}$ is the outer normal. This amounts to optimizing 
 the function 
\be
a\mapsto F(a):=Q(u(a)).
\label{qofua}
\ee
over $\cA$.
\item
{\bf Stochastic modeling:} the parameters are random variables with prescribed 
probability laws, which account for {\it uncertainties} in the model. Therefore
$a$ has a certain probability distribution $\mu$ supported on $\cA$. One is then typically 
interested in the resulting probabilistic properties of the solution $u$,
which is itself a random variable over $\cA$ with values in $V$, or in the probabilistic properties of a quantity of interest $Q(u)$.
For instance, if the elliptic equation 
\iref{ellip} is used to model oil or ground water diffusion, a common way to
deal with the uncertainty of the underground porous media is to define $a$ as a
random field with some prescribed law. Then one might want to estimate
the mean solution 
\be
\o u:=\E(u),
\ee
or its variance 
\be
V(u):=\E(\|u-\o u\|_V^2),
\ee 
or the average flux through a certain interface $\Gamma$, that is, $\E(Q(u))$ with $Q(u)$ as in \iref{qofu},
or the probability that this flux exceeds a certain quantity, etc.
\end{itemize}

 In both the deterministic and stochastic settings, the given application may   require  evaluating
$u(a^i)$ for a large number $n\gg1$ of instances $\{a^1,\dots,a^n\}$ of the parameter $a$.
This is the case, for instance, when using a descent method for optimizing 
a quantity of interest in the deterministic
framework, or a Monte-Carlo method for evaluating an expectation in
the stochastic framework.  Since each  individual instance
$u(a)$ is  the solution of a PDE, its exact evaluation is typically out of reach. 
Instead, each query for  $u(a^i)$ is approximately 
evaluated by a numerical solver, which may itself be computationally
intensive for high accuracy.

In order to significantly reduce the number of computations 
needed for attaining a prescribed accuracy,  alternate
strategies, commonly referred to as {\it reduced modeling},  have been developed.    
Understanding which of these strategies are effective, in the case where
the parameter has large or infinite dimension, and why, is the subject of this article.

\subsection{Affine representation of the parameters}
\label{SS:affine} 

So far our description of the set $\cA$ of parameters allows it to be  any compact subset of $X$.    An important ingredient
in both our theoretical and numerical developments is to identify any $a\in\cA$ through a sequence of real   numbers.   We are especially interested in  affine representations of $\cA$.
We say that a sequence $(\psi_j)_{j\ge 1}$ of functions $\psi_j\in X$ is an {\it affine representer} for $\cA$, or {\it representer} for short, if  we can write
each $a$ as 
\be
\label{affine}
 a=a(y)=\o a+\sum_{j\geq 1}y_j \psi_j, \quad y:=(y_j)_{j\geq 1},\quad y_j=y_j(a),
\ee
where the $y_j$ are real   numbers, $\o a$ is a fixed function from $X$, and the series converges in
the norm of $X$ for each $a\in\cA$. 
We are making a slight abuse of notation here since we use  $a$ to represent a general 
element of $\cA$  and also use $a$ to
represent the map
\be
a: y \mapsto a(y),
\ee
from $\R^{\N}$ to $X$.  But the meaning will always  be clear from the context.

It is easy to see that for any compact set $\cA$ in a Banach space $X$ affine representers exist.
For example, if $X$ has a Schauder basis then any such basis will be a representer.  Even if
$X$ does not admit a Schauder basis, as is the case for our example $X=L^\infty(D)$, we can still find 
representers as follows.  Choose any $\bar a\in \cA$.   Since $\cK:=\cA-\bar a$ is compact there exist
finite dimensional spaces $(X_n)_{n\geq 0}$, with $\dim(X_n)=n$, such that
\be
\dist(\cK,X_n)_X:=\sup_{a\in\cK} \min_{b\in X_n} \|a-b\|_X \to 0\quad {\rm as}\quad n\to \infty.
\ee
We can also take the spaces $X_n$ to be nested, that is
\be
X_n\subset X_{n+1}, \quad n\geq 0.
\ee
Let  $(\phi_{j,n})_{j=1,\dots,n}$ be any basis for $X_n$ and  define $N_n:=n(n-1)/2$.  The sequence  
\be
\label{defphij}
\psi_j:= \phi_{j-N_n},n,\quad j=N_n,\dots,N_{n+1},
\ee
contains each of the bases $(\phi_{j,n})_{j=1,\dots,n}$ for all $n\geq 1$.  Given any $a\in\cA$,
let $a_n$ be a best approximation in $X$ to $a-\bar a$ from $X_n$,with  $a_{0}:=0$.  Then, we can write
\be
\label{repa}
a=\bar a+ \sum_{n=1}^\infty (a_n-a_{n-1}).
\ee
Each term $a_{n}-a_{n-1}$ is in $X_{n}$ and hence can be written as a linear combination
of the $\psi_j$.  Therefore, $(\psi_j)_{j\geq 1}$ is a representer for $\cA$.

 Affine 
representations \iref{affine} often occur in the natural formulation of the parametric problem. For instance, if the diffusion coefficient $a$
in \iref{ellip} is piecewise constant over a fixed partition $\{D_j\}_{j=1,\dots,d}$ of
the physical domain $D$, then it is natural to set
\be
a(y)=\o a +\sum_{j=1}^d y_j\Chi_{D_j},
\label{pwc}
\ee
where $\o a$ is a constant and the $\Chi_{D_j}$ are the characteristic functions
of the subdomains $D_j$. Similarly, if the parameter $a$ describes 
the shape of the boundary of the physical domain in a computer-aided
design setting, a typical format is
\be
a(y)=\o a+\sum_{j=1}^d y_j B_j,
\label{spline}
\ee
where $\o a$ represents a nominal shape and the $B_j$ are B-spline functions
associated to control points. In these two examples, $d$ is  finite, yet possibly very large.

 In the statistical context, if $a$ is a second order random field over a domain $D$, a frequently used
choice in \eref{affine} is $\o a:=\E(a)$, the average field, and $(\psi_j)_{j\geq 0}$,
the Karhunen-Loeve basis, that is, the eigenfunctions of the covariance
operator
\be
v\mapsto R_a v:=\int_D C_a(\cdot,x)v(x)dx, \quad C_a(z,x):=\E((a(z)-\o a(z))(a(x)-\o a(x)).
\ee
Then, the resulting scalar variables are centered and uncorrelated, 
that is, $\E(y_i)=0$ and $\E(y_iy_j)=0$ when $i\neq j$.

Even if an affine representation of the form \eref{affine} is not given in the formulation of the problem, 
one can be derived by taking any representation system $(\psi_j)_{j\ge 1}$ in the Banach space $X$.  For example, if $X$ admits a Schauder basis,
 then one can take any such   basis $(\psi_j)_{j\ge 1}$  for $X$ and arrive at such an expansion.
In classical spaces $X$, such as $L^p$ or Sobolev spaces,  standard systems of approximation, such as Fourier series, splines, or wavelets can be used.

The advantage of the representation \eref{affine} is that $a$ can now be identified through the sequence $(y_j)_{j\ge 1}$.
When considering all $a\in \cA$, we obtain  a family of  such sequences.   Note that this family can be quite complicated.
In order to simplify matters, we   normalize the $\psi_j$,  so that for any $j$,
\be
\label{renormalize}
\sup_{a\in\cA} |y_j(a)|=1.
\ee
Such a renormalization is usually possible because $\cA$ is compact and $y_j(a)$ depends continuously on $a$.
After this normalization, for each $a\in\cA$, the sequence $(y_j(a))_{j\geq 1}$ belongs to the infinite dimensional cube 
\be
U:=[-1,1]^{\N}.
\ee

Notice  that taking a general sequence $(y_j)_{j\ge 1}$ from this cube, there may not be an $a\in \cA$ with $y_j=y_j(a)$, $j\ge 1$. Also, if  $\{\psi_j\}_{j\ge 1}$ is not a basis , the representation \iref{affine} may not be unique.  We define
\be
\label{defUA}
U_\cA:=\Big \{(y_j)_{j\geq 1}\in U \; : \; \sum_{j\ge 1}y_j\psi_j\in\cA\Big \}.
\ee
We are mainly interested in representers $\bar a$, $(\psi_j)_{j\ge 1}$ for which
  \be
 \label{complete} 
\bar a+  \sum_{j\ge 1}y_j\psi_j 
\ee
converges in $X$ for each $(y_j)\in U$.   We call such representers {\it complete}.  In this case, we  may define
\be
\label{aofU}
a(U):=\{a=a(y)=\o a+\sum_{j\geq 1}y_j \psi_j \; : \; (y_j)_{j\geq 1}\in U\},
\ee
so that
\be
\cA\subset a(U).
\ee
A typical case of a complete representer is when $(\|\psi_j\|_X)$ is a sequence in $\ell^1(\N)$.

Once an affine representation
has been chosen, the  initial solution map
$a\mapsto u(a)$ becomes equivalent to  the map $y\mapsto u(a(y))$
which is defined on $U_\cA$.  With an abuse of notation, we write
this new solution map as
\be
\label{ymap}
y\mapsto u(y):=u(a(y)).
\ee
This is a Banach space valued function of an infinite number of variables.
Note that in the case where the affine representation has a finite number $d$ of terms, 
the range of $y$ is $[-1,1]^d$. However, the infinite dimensional case subsumes the finite dimensional case,
since the latter may be viewed as a particular case with $\psi_j=0$ for $j>d$.


In the case of a complete representer,  $a(y)$ is defined on all of $U$.   However,  we do not know whether  
the solution map $u$ is defined on all of $U$.   
To guarantee this,  the following assumption will be used often.
\nl
\nl
{\bf Assumption A:} {\it  The parameter set $\cA$  has a complete representer 
$(\psi_j)_{j\geq 1}$ and  the solution map $a\mapsto u(a)$ is well defined on the
whole set $a(U)$, or equivalently the solution map $y\mapsto u(y)$ is well defined on the 
whole set $U$.}
\nl

This assumption naturally holds when the set $\cA$ is exactly defined
as $a(U)$.  

\subsection {Smoothness of the solution map}
\label{SS:uprop}

One objective of  this article is to develop efficient numerical approximations
to the solution maps of \iref{solmap} or \iref{ymap}. One of 
the main difficulties is that these maps are high or infinite dimensional,
in the sense that the dimension of the variable $a$ or $y$ is high or infinite.
In order to understand what might be good strategies for
constructing such approximations,
we need first to understand the inherent properties of these maps
that might allow us to circumvent this difficulty.

We initiate such a program in  \S \ref{S:holomorphy}, where we first analyze the smoothness of
the solution map $a\mapsto u(a)$. In the case of the elliptic equation \iref{ellip},
it is easily seen that this map is not only infinitely differentiable, but also admits a holomorphic extension to certain subdomains of the complex valued $X=L^\infty(D)$.  
We  propose two general approaches which allow us to establish
similar holomorphy properties for other relevant instances of linear and nonlinear 
parametric PDEs.  One first approach is based on the Ladyzenskaia-Babushka-Brezzi
theory. It applies to a range of linear PDEs where the operator and the right hand side
have holomorphic dependence in $a$. These include parabolic and saddle-point problems,
such as the heat equations or the Stokes problem, with parameter $a$ in
the diffusion term, similar to \iref{ellip}. One second approach is based
on the implicit function theorem in complex valued Banach spaces.
In contrast to the first approach, it can be applied to certain nonlinear PDEs.

Using the affine representation \iref{affine} of $a$,
we then study the solution map $y\mapsto u(y)$
under {\bf Assumption A}, which means that it
is defined on the whole of $U$. 
In addition to holomorphy, an important property 
of $u$ can be extracted from the affine representation \eref{affine}.  
The functions $\psi_j$ appearing in \eref{affine} have norms 
$\|\psi_j\|_X$ of varying size.  Since the variable $y_j$ is  scaled to be in
$[-1,1]$, when $\|\psi_j\|_X$ is small, this variable has a reduced effect on
the variations of $u(y)$.  Thus the variables $(y_j)_{j\geq 1}$ are not democratic, 
but rather they have varying importance.  
In other words, the map $y\mapsto u(y)$ is highly anisotropic.
More specifically, we derive holomorphic extension results for this map on certain
multivariate complex domains of tensor product type. In particular,
we consider polydiscs of the general form
\be
\cU_{\rho}:=\otimes\{|z_j|\leq \rho_j\}=\{z=(z_j)_{j\geq 1}\in\C^\N\; : \; |z_j|\leq \rho_j\}
\ee
where $\rho=(\rho_j)_{j\geq 1}$ is a positive sequence which serves to
describe the anisotropy of the solution map. We also consider 
polyellipses which deviate less far from the real axis. These holomorphy
domains play a key role in the derivation of approximation results.

\begin{remark}
\label{complexremark}
While we are generally interested in  real valued
solutions $u$ to the parametric PDE \iref{genpar}, corresponding
to real valued parameters $a$ or $y$, our analysis of holomorphic
smoothness leads us naturally to complex valued solutions, corresponding to complex valued parameters.  
For this reason,  the spaces $X,V,W$ are always assumed
to be complex valued Banach spaces throughout this paper.
\end{remark}

\subsection{Approximation of the solution map} 

Reduced modeling methods seek to   take advantage of  the   properties  of the solution maps $a\mapsto u(a)$ or $y\mapsto u(y)$ such as the holomorphy and anisotropy mentioned above.
These properties suggest strategies for approximating these map $u$ by simple functions $u_n$ in which the physical variables $x$ and the parametric variable $a$ or $y$ are separated and hence take the form
\be
(x,a)\mapsto u_n(x,a):=\sum_{i=1}^n v_i(x)\phi_i(a),
\label{separxa}
\ee
or
\be
(x,y)\mapsto u_n(x,y):=\sum_{i=1}^n v_i(x)\phi_i(y),
\label{separxy}
\ee
where $\{v_1,\dots,v_n\}$ are functions
of $x$ living in the solution space $V$ and
$\{\phi_1,\dots,\phi_n\}$ are functions of $a$ or $y$
with values in $\R$ or $\C$.

We may view $u_n$ as a rank $n$ approximation
to $u$, in analogy with low rank approximation of matrices.
We adopt the notations
\be
a\mapsto u_n(a)=u_n(\cdot,a)\quad {\rm and}\quad y\mapsto u_n(y)=u_n(\cdot,y),
\ee
for the above approximations. 

Let us discuss the potential accuracy of separable approximations
of the form \eref{separxa}.   If our objective is to
capture $u(a)$ for all $a\in\cA$ with a prescribed accuracy $\e(n)$, this
means that we search for an error bound in the uniform sense, i.e., of the
form
\be
\|u-u_n\|_{L^\infty(\cA,V)}:=\sup_{a\in\cA} \|u(a)-u_n(a)\|_V \leq \e(n).
\label{unibound}
\ee

For certain applications, in particular in the stochastic framework,
we may instead decide to measure the error on average,
for instance by searching for an error bound in the mean square sense,
\be
\E(\|u-u_n\|_V^2):=\|u-u_n\|_{L^2(\cA,V,\mu)}^2=\int_\cA \|u(a)-u_n(a)\|_V^2d\mu(a) \leq \e(n)^2,
\label{msbound}
\ee
where $\mu$ is the probability measure for the distribution of $a$ over $\cA$.
Since $\mu$ is a probability measure, one has for any $v$,
\be
\|v\|_{L^2(\cA,V,\mu)}\leq \|v\|_{L^\infty(\cA,V)}.
\label{L2uni}
\ee
Therefore the uniform bound is stronger than the average bound,
in the sense that \iref{unibound} implies \iref{msbound}
with the same value of $\e(n)$.

Likewise, for the approximation of the map $y\mapsto u(y)$,
we may search for a uniform bound
\be
\|u-u_n\|_{L^\infty(U_\cA,V)}:=\sup_{y\in U_\cA} \|u(y)-u_n(y)\|_V \leq \e(n)
\label{uniboundy}
\ee
or a mean square bound
\be
\E(\|u-u_n\|_V^2):=\|u-u_n\|_{L^2(U_\cA,V,\mu)}^2=\int_{U_\cA} \|u(y)-u_n(y)\|_V^2d\mu(y) \leq \e(n)^2,
\label{msboundy}
\ee
where $\mu$ is the probability measure for the distribution of $y$ over $U_{\cA}$.

\begin{remark}
We do not indicate the measure $\mu$ in our notation $L^\infty(\cA,V)$ or $L^\infty(U_\cA,V)$ , since in
all relevant examples considered in this article we
always consider the exact supremum over $a$ in $\cA$ or over $y$ in $U_\cA$,
rather than the essential supremum.
\end{remark}

For any $a\in\cA$, the approximation $u_n(a)$ belongs to 
\be
V_n:={\rm span}\{v_1,\dots,v_n\},
\ee
which is a fixed $n$-dimensional subspace of $V$. Ideal
benchmarks for the performance of separable expansions of the form \iref{separxa}
may thus be defined by selecting optimal $n$-dimensional
spaces for the approximation of $u(a)$ in either a uniform or an average sense.

For uniform error bounds, this benchmark is given by the concept
of Kolmogorov's {\it $n$-width}, which is well known in approximation theory.
If $\cK$ is a compact set in a Banach space $V$, we define its
Kolmogorov $n$-width as
\be
d_n(\cK)_V:=\inf_{\dim(V_n)\leq n} \sup_{v\inÊ\cK} \min_{w\in V_n}\|v-w\|_V.
\label{kolnwidth}
\ee
This quantity, first introduced in \cite{Kol}, describes the best achievable accuracy, in the norm of $V$, when approximating
all possible elements of $\cK$ by elements from a linear $n$-dimensional space $V_n$.  Obviously, the
optimal choice of $V_n^*$ for approximation of $u(a)$ in a uniform sense
corresponds to the space, if it exists, that reaches the above infimum 
when $\cK$ is taken to be the solution manifold $\cM$. 
The best achievable error in the uniform sense is thus given by
\be
\e(n):=d_n(\cM)_V.
\ee
There exist many other notions of widths that are used to measure the size of
compact sets. Here we only work with the above one and refer 
to \cite{Pin} for a more general treatment.

For mean square bounds, in the case where $V$ is a Hilbert space, the corresponding benchmark
is related to the concept of principal component analysis, which is of common
use in statistics. By choosing an arbitrary orthonormal basis
$(e_i)_{i\geq 0}$ of $V$, we may expand $u(a)$ according to
\be
u(a)=\sum_{i\geq 0} z_ie_i,\quad z_i=z_i(a):=\<u(a),e_i\>_V.
\ee
Approximation of $u(a)$ from an $n$-dimensional space of $V$  
is then equivalent to the approximation of $\bz(a):=(z_i(a))_{i\geq 0}$
from an $n$-dimensional space of $\ell^2(\N)$. The optimal space is
then obtained through the study of the correlation operator
\be
R=(R_{i,j})_{i,j\geq 0}, \quad R_{i,j}:=\E(z_i(a)z_j(a))=\int_{a\in\cA} z_i(a)z_j(a)d\mu(a).
\label{covop}
\ee
This operator is symmetric, positive, compact and of trace class. It therefore admits
an orthonormal basis of eigenvectors $(\bg_k)_{k\geq 1}$ associated
to a positive, non-increasing and summable sequence $(\lambda_k)_{k\geq 1}$ of eigenvalues.
The space 
\be
G_n:={\rm span}\{\bg_1,\dots,\bg_n\},
\ee
minimizes over all $n$-dimensional spaces $G$ the mean square error between 
$\bz$ and its orthogonal projection $P_G\bz$, with
\be
\E(\|\bz-P_{G_n}\bz\|^2_{\ell^2})=\sum_{k>n}\lambda_k.
\ee
In turn, the optimal space $V_n^*$ is spanned by the functions   
\be
v_k:=\sum_{i\geq 0} g_{k,i}e_i, \quad k=1,\dots,n,
\ee
where $g_{k,i}$ is the $i$-th component of $\bg_k$, that is, $\bg_k=(g_{k,i})_{i\geq 0}$. The
best achievable error in the mean square sense is thus given by
\be
\e(n)^2:=\sum_{k>n}\lambda_k.
\label{tailsing}
\ee
Note that, in view of \iref{L2uni}, one has the comparison
\be
\sum_{k>n}\lambda_k \leq d_n(\cM)_V^2.
\ee

The above described optimal spaces are usually out of reach, 
both from an analytic and computational point of view, and 
it is therefore interesting to consider sub-optimal approximations.
In addition, when considering the solution map $y\mapsto u(y)$,
the tensor product structure of $U$ allows us to consider approximations
based on further separation between the parametric variables, that is,
where each function $\phi_i$ in \iref{separxy} is itself a product of univariate
functions of the different $y_j$. The simplest example of such approximations
are multivariate polynomials, which have the general form
\be
u_n(x,y)=\sum_{\nu\in\Lambda_n} v_\nu(x) y^\nu, \quad\quad  y^\nu:=\prod_{j\geq 1} y_j^{\nu_j},
\label{polynomial}
\ee
where each index $\nu=(\nu_j)_{j\geq 1}\in \Lambda_n$ is a finitely supported sequence of positive integers, 
or equivalently such that $|\nu|=\sum_{j\geq 1} \nu_j<\infty$, and $\Lambda_n$ is a set of $n$ such sequences.

In \S 3, we obtain such polynomial approximations
 by taking finite portions of infinite polynomial expansions
of $u$. Here, we work under {\bf Assumption A}, which means
that $u(y)$ is defined on the whole of $U$.
We consider two types of expansions:
\begin{itemize}
\item
Power series the form
\be
\label{repu}
u(y)=\sum_{\nu\in \cF}t_\nu y^\nu,
\ee
where $\cF$ is the set of all finitely supported sequence of positive integers.
\item
Orthogonal series of the form
\be
\label{leg}
u(y)=\sum_{\nu\in \cF}w_\nu P_\nu(y), \quad P_\nu(y):=\prod_{j\geq 1}P_{\nu_j}(y_j),
\label{legu}
\ee
where $P_k$ is the Legendre polynomial of degree $k$ defined on $[-1,1]$.
\end{itemize}
Using the holomorphy and anisotropy  properties of the solution map 
$y\mapsto u(y)$ established in \S 2 for specific classes of parametric PDEs, we 
derive a priori bounds on the $V$-norms $\|t_\nu\|_V$ and $\|w_\nu\|_V$
of the coefficients which appear in these expansions.
In this way, we are able to establish
algebraic convergence rates $n^{-s}$ for certain truncations
of the above expansions, where $n$ is the
cardinality of the truncation set $\Lambda_n$. 

One critical aspect of the truncation strategy is that we retain the $n$ largest coefficients, 
which is a form of {\it nonlinear approximation}, also known
as {\it sparse} or {\it best $n$-term} approximation. With such a choice for $\Lambda_n$,
the exponent $s$ in the convergence rate is related to
the available $\ell^p$ summability for $p<1$ of the $V$-norms
of the coefficients in the considered infinite expansion. 
In particular, for uniform approximation estimates, that is, in $L^\infty(U,V)$,
one has  
\be
s=\frac 1 p-1,
\label{sp}
\ee
once the $\ell^p$-summability of these $V$-norms has been proven.
The main result from \S 3 shows that, under suitable assumptions, 
the $\ell^p$ summability of the sequence $(\|\psi_j\|_{X})_{j\geq 0}$ 
implies that the norms of the coefficients
in the expansions \eref{repu} or \iref{legu} are also $\ell^p$ summable.

The fact that we obtain the algebraic convergence rate $n^{-s}$
despite the infinite dimensional nature of the variable $(y_j)_{j\geq 1}$
reveals that the curse of dimensionality can be avoided
in the approximation of relevant parametric PDEs. 

\subsection{The $n$-widths of the solution manifold $\cM$}  

In both the uniform or mean square ways 
of measuring error, as described above, the success of   
reduced modeling can be proven if $d_n(\cM)_V$ converges to zero
sufficiently fast as $n\to +\infty$.  Thus, the study of these widths constitutes 
a major subject in the theoretical justification of reduced modeling. 

  A common way to measure the widths of compact sets $\cK$  in classical spaces is to embed $\cK$ into an appropriate smoothness space such as a Sobolev or Besov space.   
  For example, in our model parametric elliptic equation \iref{ellip}, the space $V$ is $H_0^1(D)$, and this approach would  lead us to examine the $H^m$ Sobolev smoothness of the individual functions $u(a)$ for $m>1$.   Classical elliptic regularity theory says that  for smooth domains the smoothness of $u(a)$ can be inferred from the smoothness of the right side $f$.   
 However, for general domains, there are severe limits on this regularity due to the irregularity of the boundary.  
 Therefore, bounding the decay of widths of $\cM$ through such regularity results will generally prove only slow decay for $d_n(\cM)_V$ and therefore is not useful for obtaining the fast decay rates we seek.  
Indeed, recall, that classical smoothness spaces, such as Sobolev or Besov spaces of order $m$ in $d$ variables, have widths that decay, at best, like $n^{-(m-r)/d}$ as $n\to \infty$ when these widths are measured in an $W^{r,p}$ norm
for some $r<m$. For example, it is known that 
if $\cK$ is the unit ball of $C^m([-1,1]^d)$, then its $n$-width in $L^\infty$ satisfies
\be
\label{Sobwidth}
c_dn^{-m/d}\leq d_n(\cK)_{L^\infty}\leq C_dn^{-m/d},\quad n\ge 1.
\ee
Likewise, if $\cK$ is the unit ball of $H^m([-1,1]^d)$, then its $n$-width in $H^1$ satisfies
\be
\label{Sobwidth}
c_dn^{-(m-1)/d}\leq d_n(\cK)_{H^1}\leq C_dn^{-(m-1)/d},\quad n\ge 1.
\ee
Also note that the regularity of $u(a)$, as measured by membership in Sobolev and Besov spaces,  is closely related to the performance of piecewise polynomial approximation such as that used in finite element methods and is the reason why these algorithms have rather slow convergence.  So the fast decay of $d_n(\cM)_V$ to zero cannot be obtained by such an approach.

The widths $d_n(\cM)_V$ go to zero fast not because the individual elements in $\cM$ are smooth in the
physical variable, but rather because $\cM$ is the image of the solution map $u$ which is, 
as previously discussed, smooth and anisotropic in the parametric variable.
However, let us remark that it is by no means trivial 
to deduce fast decay for $d_n(\cM)_V$ from this fact alone, because of what is called {\it the curse of dimensionality}.  Namely, approximation rates for a given target function are generally proved by showing that the target function, in our case the function $u$, has sufficiently high regularity.  But, in high dimensions, regularity by itself is usually not enough.  Indeed, returning to the bounds \eref{Sobwidth}, we see that the large dimension $d$ affects the approximation rate in two detrimental ways.   The exponent in the smoothness rate is divided by $d$ and the constants $C_d$ are known to grow exponentially with increasing $d$.  In our case $d$ can even be infinite and this makes the derivation of approximation rates 
a subtle problem.

In \S 4, we discuss general principles for estimating by above 
 the $n$-widths $d_n(\cM)_V$ of the solution manifold $\cM$.
 One immediate consequence of the approximation results in \S 3 is that,
 for specific classes of parametric PDEs, when {\bf Assumption A}Ê holds and
$(\|\psi_j\|_{X})_{j\geq 0}\inÊ\ell^p$, then these widths decay at least like $n^{-s}$ with $s$ given by \iref{sp}.
We extend this analysis to the general case when {\bf Assumption A}Ê does not
necessarily hold. Since the manifold $\cM$ is not directly accessible, one would like to understand
what properties of $\cA$, which is assumed to be completely known to us, imply
decay rates for $d_n(\cM)$.  We show that the asymptotic decay of the $n$-width of $\cM$
in $V$ is related to that of the $n$-width of $\cA$ in $X$.  This follows from general results on widths  of images of
compact sets under holomorphic maps.  Namely, if $u:X\to V$
 is holomorphic in a neighborhood of a general compact set $\cA\subset X$, then we prove the following result
 on the $n$-width of the image $\cM=u(\cA)$ in $V$:
\be
\sup_{n\geq 1}n^rd_n(\cA)_X<\infty \Rightarrow \sup_{n\geq 1}n^sd_n(\cM)_V<\infty, \quad s<r-1.
\ee
This result shows that from the view point of preserving $n$-widths,  holomorphic maps behave almost as good as linear maps, up a loss of $1$ in the convergence rate. One open problem is to understand if this loss is sharp
or could be improved.

\subsection{Numerical methods for reduced modeling}

The above mentioned results on holomorphic extensions, 
sparse expansions for $u$, and $n$-widths of the manifold $\cM$,
can be thought of as theoretical justifications for the role of reduced modeling in solving
parametric and stochastic equations.  They  provide evidence that reduced modeling numerical methods should yield  significant computational savings over traditional methods such as finite element solvers for parametric problems or
Monte Carlo methods for stochastic problems.  However, they do not constitute actual numerical methods.

The second part of our article turns to the construction of numerical algorithms motivated by
these theoretical results. Such algorithms compute specific separable approximations of the form
\iref{separxa} or \iref{separxy} for a given value of $n$ at an affordable computational cost.
Our objective, in this regard, is not to give an exhaustive description
of  numerical reduced modeling methods and their numerous variants.
Rather, our main focus is to introduce some important representative examples
of these methods for which an a priori analysis quantifies 
the gains in numerical performance of these methods.

One important distinction is between {\it non-adaptive} and {\it adaptive} methods. In the first ones,
the choice of the functions $\{\phi_1,\dots,\phi_n\}$ or of $\{v_1,\dots,v_n\}$ used in \iref{separxa}
or \iref{separxy}, for a given value of $n$ 
is made in an a priori manner, if possible using the available information
on the problem. In the second ones, the computations executed for lower values of $n$
are exploited in order to monitor the choice at stage $n$. One desirable feature
of an adaptive algorithm is that it should be {\it incremental} or {\it greedy}: 
only one new function $\phi_n$ or $v_n$ 
is added at each stage to the $n-1$ previously selected 
functions which are left unchanged. Adaptive algorithms are known to often perform
better than their non-adaptive counterpart, but their convergence analysis is 
usually more delicate.

A second important distinction between the various numerical methods 
is whether they are  {\it non-intrusive}Ê or {\it intrusive}. 
A non-intrusive algorithm builds on an existing exact or approximate solver
for the PDE which may be computationally expensive.
It derives approximations of the form \iref{separxa}
or \iref{separxy} by choosing instances $a^1,\dots,a^n\in \cA$
or $y^1,\dots,y^n\in U_\cA$ and using the values $u(a^i)$ or $u(y^i)$
computed by the solver. Non-intrusive algorithms may be implemented 
even when this solver is a black box, and therefore with a possibly limited knowledge on the 
exact PDE model.  An intrusive algorithm, on the other hand, directly exploits
the precise form of the PDE for computing the approximation 
\iref{separxa} or \iref{separxy}, and therefore requires the full knowledge 
of the PDE model for its implementation.  

It should be noted that instances $u(a^i)$ as well 
as the functions $\{v_1,\dots,v_n\}$ used in \iref{separxa}
or \iref{separxy} can only be computed with a certain level of spatial discretization,
for example resulting from a finite element solver. In such a case they
belong to a finite dimensional space $V_h\subset V$. 
If the same finite element space
$V_h$ is used to discretize all instances $u(a)$
up to a precision that is satisfactory for the user, this means
that we are actually trying to capture the {\it approximate solution
maps}
\be
a\mapsto u_h(a)\in V_h,
\label{appromapa}
\ee
or 
\be
y\mapsto u_h(y)\in V_h.
\label{appromap}
\ee
The analysis of the performance of numerical reduced
modeling methods needs to incorporate the additional error produced by this 
discretization.

\subsection{Sparse polynomial approximation algorithms}

One first class of methods that we analyze consists in finding a numerically 
computable polynomial approximation of the form \eref{polynomial}.
There are two major tasks in constructing
such a numerical approximation:  (i) find good truncation sets $(\Lambda_n)_{n\geq 1}$
and (ii) numerically compute an approximation to the coefficients $v_\nu$ for each $\nu\in\Lambda_n$.  
By far, the most significant issue in numerical methods based on polynomial expansions is  to find a good choice for the sets $(\Lambda_n)_{n\geq 1}$.  If everything was known to us, we would simply take for 
 $\Lambda_n$ the set of indices $\nu$ corresponding to
 the $n$ largest $V$-norms of the coefficients in \eref{repu} or \iref{legu}.  
However, finding such an optimal $\Lambda_n$ would require in principle 
that we compute the coefficients for all values $\nu\in\cF$ which is obviously out of reach.
In addition, the structure of the optimal set $\Lambda_n$ can be quite complicated.  
One saving factor is that our analysis in \S 3 gives a priori bounds on the size of these coefficients.   
These bounds can be used in order to make an a priori selection of the sets $\Lambda_n$.
This is a non-adaptive approach, which generally gives suboptimal performance due to the 
possible lack of sharpness in the a priori bounds. This leads one to try to enhance 
performance by combining the a priori bounds together with an adaptive selection of the
sets $(\Lambda_n)_{n\geq 1}$.

The numerical methods are facilitated by imposing that the 
selected index sets $\Lambda_n$ 
are {\it downward closed} (or {\it lower sets}), i.e.
satisfy the following property:
\be
\nu\in \Lambda_n\;\; {\rm and} \;\; \t \nu \leq \nu \Rightarrow \t \nu\in \Lambda_n,
\ee
where $\t \nu\leq \nu$ means that $\t \nu_j\leq \nu_j$ for all $j$. 

In \S 6, we discuss algorithms which compute the polynomial approximation
by an interpolation process. These algorithms are non-intrusive
and apply to a broad scope of problems.
One key issue is the choice of the interpolation points, which is
facilitated by the following result: if $\{z_k\}_{k\geq 0}$ is a sequence of pairwise distinct points
in $[-1,1]$ and if $z_{\nu}:=(z_{\nu_j})_{j\geq 1}\in U$, then for any 
downward closed set $\Lambda_n$, any polynomial of the form \iref{polynomial}
is uniquely characterized by its value on the grid $\{z_{\nu}\; : \; \nu\in \Lambda_n\}$.
This allows us to construct the interpolation in a hierarchical
manner, by simultaneously incrementing the polynomial space and the
interpolation grid. 
The sets $\Lambda_n$
can either be a priori chosen based on the bounds on the coefficients established in 
\S 3 or adaptively generated. These sets generally differ from the 
ideal sets corresponding to the $n$ largest coefficients (which may not
fullfill the downward closedness property). 
We show that certain choices of the univariate sequence $\{z_k\}_{k\geq 0}$,
known as Leja or $R$-Leja points lead to stable interpolation processes
(in the sense that the Lebesgue constant have moderate growth with the number of points) allowing
us to retrieve by interpolation the same algebraic convergence rate $n^{-s}$ which are proved
for the best polynomial approximations.

In \S 7, we discuss another class of algorithms which recursively compute the exact
coefficients in the Taylor series of the approximate solution map \iref{appromap}.
These algorithms are intrusive and they only apply to problems where $\cP$ is linear 
both in $u$ and $a$ up to a constant term, such as the elliptic equation \iref{ellip}
which serves a guiding example. The recursive computation is facilitated by imposing that
the index sets $\Lambda_n$ in the truncated Taylor expansion are downward closed.
Similar to the interpolation algorithm from \S 6, the sets $\Lambda_n$
can either be a priori chosen based on the available bounds on the coefficients established in 
\S 3, or adaptively generated. One main result shows that 
adaptive algorithms based on a so-called {\it bulk chasing} procedure have the same convergence rate $n^{-s}$ as 
the one which is established when using the index sets corresponding to the $n$ largest
Taylor coefficients.

\subsection{Reduced basis algorithms}

The  second  class of methods that we analyze seeks to find, in an {\it offline} stage,
a set of functions $\{v_1,\dots,v_n\}$  for which the resulting $n$-dimensional space 
$V_n:={\rm span}\{v_1,\dots,v_n\}$ is close to an  optimal linear $n$-dimensional approximation space.   For mean square estimates, one such approach, known as {\it proper orthogonal decomposition},
consists in building the functions $\{v_1,\dots,v_n\}$ based on an approximation
of the exact covariance operator \iref{covop} computed from a sufficiently dense sampling of the random 
solution $u(a)$. Another approach, which targets uniform estimates,
is the {\it reduced basis method}, which consists in generating $V_n$ by a selection 
of $n$ particular solution instances $\{u(a^1),\dots,u(a^n)\}$
chosen from a very large set of potential candidates.   In both cases, the offline stage is 
potentially very costly. 

Once such a space $V_n$ is chosen, one builds an    {\it online} solver, such that  for any given 
$a\inÊ\cA$, the approximate solution $u_n(a)$ is   an element from $V_n$.   There are several possibilities on how to build
this online solver.  The most prominent of these is to take the 
 Galerkin projection of $u(a)$ onto $V_n$ which consist of finding $u_n(a)$ by solving the system of equations
\be
\<\cP(u_n(a),a),w\>=0, \quad w\in V_n
\ee
for a suitable duality product $\<\cdot,\cdot\>$. This online computation determines
for each $a$ the values $\phi_i(a)$ for $i=1,\dots,n$ which appears in \iref{separxa}.  The advantage of using the Galerkin solver, is that, for certain problems such as elliptic ones, it is known to give the 
best error in approximating $u(a)$ by elements of $V_n$ when error is measured
in the norm of $V$.   Its disadvantage its computational 
cost in finding $u_n(a)$ given the query $a$.  For this reason other projections 
onto $V_n$ are also studied, some of these based on interpolation.

The key issue when using reduced basis methods is how to find a good space $V_n$, i.e. how to find 
good basis functions.  In \S 8, we  discuss  an elementary greedy 
strategy for the offline selection of the instances $v_i=u(a^i)$, that consists in
picking the $n$-th instance which deviates the most from the space $V_{n-1}$
generated from the $n-1$ previously selected ones. The approximation error
\be
\sigma_n(\cM):=\sup_{v\in \cM}\inf_{w\in V_n} \|v-w\|,
\ee
produced by such spaces may be significantly
larger than the ideal benchmark of the $n$-width of the solution manifold 
for a given value of $n$. However, a striking result is that 
both are comparable in terms of rate of decay: for any $s>0$, there is 
a constant $C_s$ such that 
\be
\sup_{n\geq 1} n^s \sigma_n(\cM) \leq C_s\sup_{n\geq 1} n^sd_n(\cM)_V.
\ee
Similar results are established for exponential convergence rates.

While both classes of numerical methods aim to construct separable
approximations of the form \iref{separxa} or \iref{separxy}, there is a significant
distinction between them in the way they organize computation. For the first
class of polynomial approximation methods, the offline stage fixes the polynomial functions $\phi_i$
through the selection of the set $\Lambda_n$
and computes the coefficients $v_i$. Then the online stage is in some
sense trivial since it simply computes $u_n(y)$ through the linear combination \iref{polynomial}.
For the second class, the online stage still requires solving PDE approximately
in the chosen reduced space $V_n$. This offline/online splitting makes it difficult to draw a fair comparison
between the different methods from the point of view of computational time vs accuracy.
\nl
\nl
{\bf Notation for constants:} {\it Numerous multiplicative constants appear throughout this paper,
for example in convergence estimates. We use the generic notation
$C$, which may therefore change value between different formulas,
and if necessary we indicate the parameters on which $C$ depends.
We use a more specific notation if we want to express the dependence of
the constant with respect to 
a cetain parameter (for example the dimension $d$ in \iref{Sobwidth})
or if we want to refer to this specific constant later in the paper.}

\subsection{Historical orientation}

Numerical methods for parametric and stochastic PDEs
using polynomials (or other approximation tools) in the parametric variable
have been widely studied since the 1990's. We refer in particular
to \cite{GS,GWZ,KH,KL,X} for general introductions
to these approaches, and to \cite{ABS,GS1,FST,KX,KX1,ST03,ST07}
for related work prior to the results exposed in our paper.

The approximation results presented in \S 3 have been 
obtained by the authors and their co-authors in a series of paper \cite{CDS1,CDS2,CCS1},
and the results in \S 4 on the evaluation of $n$-widths are from \cite{CD}. These results
establish for the first time convergence rates immune to the curse of dimensionality,
in the sense that they hold with infinitely many variables, see also \cite{GitS} for a survey dealing
in particular with these issues. In a similar infinite dimensional framework,
and not covered in our paper,
let us mention the following related works: (i)
similar holomorphy and approximation results are established in \cite{HS,HS2,KS} for specific type of PDEs and control problems,
(ii) approximation of integrals by quadratures is discussed in \cite{KSS1,KSS2}, (iii) inverse problems are discussed in \cite{ScS1,ScS2,SSt}, following the Bayesian perspective from \cite{Stu}, and (iv) diffusion problems with lognormal coefficients are treated in
\cite{HoS,Git1,GKNSSS}.

The sparse interpolation method presented in \S 6 is introduced and studied in \cite{CCS}.
See also \cite{BabNobTem07,BNTT1,NobTemWeb08a,NobTemWeb08b} for related
work on collocation methods.
Other non-intrusive methods are based on least-square regression 
as discussed in \cite{CCMNT,DI,DO,MNTV}, or on pseudo-spectral approximation
as discussed in \cite{Xiu,CEP2}.
The Taylor approximation algorithm presented in \S 7 is introduced and studied in \cite{CCDS}.
Other intrusive methods based on Galerkin projection are discussed in \cite{BNTT1,BNTT2,CDS1,Git}.

Reduced basis methods have been studied since the 1980's
\cite{NP}. The greedy algorithms presented in \S 8
have been introduced and discussed in \cite{VPRP,MMOPR,MPT,MPT1,RHP},
and their convergence analysis was given in \cite{BCDDPW} and \cite{DPW2}.
We refer to \cite{KV} for a general introduction on the related POD method, which is not
discussed in our paper.

\section*{Part I. Smoothness and approximation results}

\section{Holomorphic extensions} 
\label{S:holomorphy}

In this section, we discuss smoothness properties of the solution maps
$a\mapsto u(a)$ and $y\mapsto u(y)$ which are central to the 
development of efficient numerical methods that are immune to the curse of
dimensionality. We show that, under suitable 
assumptions on  the parametric PDE, these maps
admit holomorphic extensions to certain complex domains.  Recall that 
a map $F$ from a complex Banach space
$X$ to a second complex Banach space $Y$ is said to be holomorphic on an open set $\cD\subset X$ if for each $x\in \cD$,
$F$ has a Frechet derivative $dF(x)$ at $x$.  Here $dF(x)$ is a linear operator mapping $X$ to $Y$ such that
\be
\label{Frechet}
\|F(x+h)-F(x)-dF(x)h \|_Y=o(\|h\|_X),\quad h\in X.
\ee
\subsection{Extension of $a\mapsto u(a)$ for the model elliptic equation}

In order to formulate results on holomorphy,
we need to introduce existence-uniqueness theory for solutions to \eref{genpar} 
in the case that $a$ is complex valued.
   
We begin with our guiding example of the elliptic equation \iref{ellip}.
We now consider  
\be
(X,V,W)=(L^\infty(D),H^1_0(D),H^{-1}(D))
\ee
as spaces of complex valued functions,
and extend the standard variational formulation to such spaces in a 
straightforward manner:  for a given $a\in X$, and with $f\in W$, find $u=u(a)\in V$ such that 
\be
\int_D a \nabla u\cdot \nabla v=\<f, v\>_{W,V}, \quad v\in V,
\label{variat}
\ee
where in the left integrand,
\be
\nabla u\cdot \nabla v:=\sum_{i=1}^m\partial_{x_i}u\o{\partial_{x_i}v},
\ee
is the standard hilbertian inner product, and $\<f,v\>$ is the
anti-duality pairing between $W$ and $V$, which when $f\in L^2(D)$ is given by
the hilbertian inner product
\be
\<f,v\>=\int_D f \o v.
\ee
We recall that 
\be
\|v\|_V:=\|\nabla v\|_{L^2(D)},
\ee
and
\be
\|v\|_W:=\sup\{\<v,w\>\; : \; \|w\|_V\leq 1\}.
\ee
This is therefore a particular case of a linear problem with 
the following general variational formulation.  Let $\frak{B}$ denote the set of all sesquilinear forms   defined on $V\times V$ and let $W=V^*$ be the set of all antilinear functionals defined on $V$, i.e., $W$ is the antidual of $V$.  We define the following norm on $\frak{B}$:
\be
\label{Bnorm}
\|B\|:=\sup_{\|v\|_V\le 1,\ \|w\|_V\le 1}|B(v,w)|.
\ee

\vskip .1in
\noindent
{\bf Problem:}  {\it Given $B\in\frak{B}$ and $L\in W$,  find $u\in V$ such that
\be
B(u,v)=L(v),\quad \forall v\in V.
\label{generellip1}
\ee}

The existence-uniqueness theory for such problems
can be proven from  the   
complex version of Lax-Milgram theorem given in Theorem \ref{theolaxmilgram}.  This theorem  is a particular case
of Theorem \ref{theoinfsup} proved below.   To formulate this theorem and for later use, we introduce the notation
$\cL(X,Y)$ for  the space of all linear operators $T$ mapping the Banach space $X$ into the Banach space $Y$ with its usual norm 
\be
\label{defnormT}
\|T\|_{\cL(X,Y)}:=\sup_{x\in X}\frac{\|Tx\|_{Y}}{\|x\|_X}.
\ee

Given $B\in\frak{B}$, one has that $B(u,\cdot)$ is an anti-linear functional and hence for any $u\in V$, there is a $\cB u\in W$
such that 
\be
\label{defcB}
B(u,v)=\langle \cB u,v\rangle_{W,V},\quad v\in V,
\ee
where $\<\cdot,\cdot\>_{W,V}$ is the anti-duality pairing between $W$ and $V$.
Therefore, $\cB$ is a linear operator from $V$ into $W$ and its norm is the same as that of $B$:
\be
\label{normcB}
\|\cB\|_{\cL(V,W)}=\|B\|.
\ee
So the operator $\cB$ is bounded and hence continuous. The problem \iref{generellip1}
is equivalent to the equation
\be
\cB u=L,
\ee
set in $W$. With this notation and remarks in hand, we can now state the complex version of the Lax-Milgram theorem.
\begin{theorem}
\label{theolaxmilgram}
Assume that,   $B\in\frak{B}$ is a sesquilinear form on $V\times V$  such that  
\be
|B(u,u)|\geq \alpha \|u\|_V^2,\quad u\in V,
\label{coerc}
\ee
for some  $\alpha>0$.
Then, $\cB$ is   is invertible and its inverse satisfies 
\be
\|\cB^{-1}\|_{\cL(W,V)}\leq \frac 1 \alpha.
\label{invbound}
\ee
Thus, for each $L\in W$, the problem \eref{generellip1} has a unique solution $u_L=\cB^{-1}(L)$ 
which satisfies the a priori estimate
\be
\|u_L\|_V\leq \frac{ \|L\|_W}{\alpha}.
\ee
\end{theorem}

For the particular $B$ and $L$ given by the left and right side of \iref{variat},
the ellipticity condition \iref{coerc} holds with $\alpha=r$ under the assumption
\be
\Re(a(x))\geq r,\quad x\in D,
\label{ellipre}
\ee
since the latter implies, for all $v\in V$,
\be
|B(v,v)|\geq \Re(B(v,v))=\int_D\Re(a) |\nabla v|^2 \geq r \|\nabla v\|_{L^2}^2=r\|v\|_V^2.
\ee
Therefore, for any $r>0$ we may extend the solution map $a\mapsto u(a)$ 
of the elliptic problem \iref{ellip} to the complex domain
\be
\cD_r:=\{a\in X\; : \; \Re(a)\geq r\},
\label{Dr}
\ee
with the uniform bound 
\be
\|u\|_{L^\infty(\cD_r,V)}=\sup\{ \|u(a)\|_V\; : \; a\in \cD_r\} \leq \frac {\|f\|_{W}}r.
\ee
This extension is therefore defined on the open set $\cD:=\cup_{r>0}\cD_r$.

The fact that this map is holomorphic immediately follows by viewing it
as a chain of holomorphic maps: introducing for 
any $a$ the operator $\cB(a): v\mapsto -{\rm div}(a\nabla v)$ acting 
from $V$ into $W$, we can decompose $a\mapsto u(a)$ into the chain
of maps
\be
a\mapsto \cB(a) \mapsto \cB(a)^{-1} \mapsto \cB(a)^{-1}f = u(a).
\label{chain}
\ee
The first and third maps are continuous linear and therefore holomorphic,
from $X$ into $\cL(V,W)$ and from $\cL(W,V)$ onto $V$ respectively. The second
map is the operator inversion which is holomorphic at any invertible $\cB\in \cL(V,W)$.

For further purposes, it is interesting to compute  the Frechet complex derivative $du(a)\in \cL(X,V)$
for the elliptic problem \iref{ellip}. 
Fix an  $a\in \cD_r$ and let  $h\in X$ be such that $\|h\|_X\leq \frac r 2$.  Then, the solution map
is also defined at $a+h\in \cD_{\frac r 2}$. Substracting the variational formulations \iref{variat}
for $u(a+h)$ and $u(a)$, we find that
\be
\int_D a \nabla (u(a+h)-u(a))\cdot  \nabla v=-\int_D h \nabla u(a+h)\cdot  \nabla v.
\ee 
We first use this identity to obtain a Lipschitz continuity bound: by taking
$v=u(a+h)-u(a)$ and taking the real part of both sides, we find that
$$
\begin{disarray}{ll}
r\|u(a+h)-u(a)\|_V^2  &\leq \|h\|_X \|u(a+h)\|_V \|u(a+h)-u(a)\|_V \\
& \leq \frac {2 \|f\|_W} r \|h\|_X\|u(a+h)-u(a)\|_V,
\end{disarray}
$$
and therefore
\be
\|u(a+h)-u(a)\|_V \leq C \|h\|_X, \quad C:=\frac {2 \|f\|_W} {r^2}.
\label{lip}
\ee
We next show that $du(a) h$ can be defined as the solution $w=w(h)\in V$ to the problem
\be
\int_D a \nabla w\cdot  \nabla v=-\int_D h \nabla u(a)\cdot  \nabla v, \quad v\in V,
\ee 
which is well-posed in the sense of the above Lax-Milgram theorem. Indeed, on the one hand, the dependence
of $w$ on $h$ is linear and  it is continuous because taking $v=w$ we
find that
\be
\|w\|_V \leq C\|h\|_X,\quad C:=\frac { \|u(a)\|_V}{r}.
\ee
On the other hand, the remainder $g=u(a+h)-u(a)-w$ is the solution to 
\be
\int_D a \nabla g\cdot  \nabla v=\int_D h(\nabla u(a)- \nabla u(a+h))\cdot  \nabla v,\quad v\in V,
\ee
which by taking $v=g$ and using \iref{lip} gives the quadratic bound
\be
\|g\|_V \leq \frac 1 r \|h\|_X \|u(a)-u(a+h)\|_V \leq C \|h\|_X^2, \quad  C:=\frac {2 \|f\|_W} {r^3}.
\ee
This confirms that $du(a) h=w(h)$.

\subsection{Extensions by the Ladyzhenskaya-Babushka-Brezzi theory}
Our next goal is to treat more general linear parametric problems that are not necessarily elliptic.
In particular, we have in mind parabolic problems such as
the heat equation, or saddle points problems such as the Stokes equations.
 In order to formulate this general class of problems, we suppose that $V$ and $\t V$ are two complex Hilbert spaces with inner products 
 $\<\cdot,\cdot\>_V$ and $\<\cdot,\cdot\>_{\t V}$, respectively.   
 So, for example,  for every $v\in V$,
 $\<v,\cdot\>_{V}$ is an anti-linear functional on $V$ and $\<\cdot,v\>_{V}$  is a linear functional on $V$, and the same holds for $\t V$.  We let   
 $W:=\t V^*$ denote the space of all anti-linear functionals on $\t V$, 
 i.e.  $W$ is the anti-dual space of $\t V$. 
 
 We now denote by $\frak{B}=\frak{B}(V,\t V)$   the set of all such sesquilinear forms
on $V\times \t V$ and introduce the following topology on $\frak{B}$,
 \be
 \label{defnormB}
 \|B\|= \max_{\|v\|_V=1,\|w\|_{\t V}=1} |B(v,w)|
\ee
As in the previous section, given $B\in\frak{B}$, one has that $B(u,\cdot)$ is an anti-linear functional on $\t V$.  Therefore, as in the previous section, we can define the linear operator $\cB$ that maps $V$ into $W$ by \eref{defcB} 
and we again have $\|\cB\|_{\cL(V,W)}=\|B\|$.

We consider the following general problem.
\vskip .1in
\noindent
{\bf Problem:} {\it  Given   $B\in\frak{B}(V,\t V)$ and  $L\in W$,  find $u\in V$ such that
\be
B(u,v)=L(v),\quad \forall v\in \t V,
\label{generellip}
\ee
}
This problem is again equivalent to the equation $\cB u=L$. In order to establish the existence-uniqueness for complex formulations of such problems, we use the following complex valued version
of Ladyzenskaya-Babushka-Brezzi theorem. 
\begin{theorem}
\label{theoinfsup}

Assume that $B\in \frak{B}(V,\t V)$ satisfies
\be
\inf_{u\in V}\sup_{v\in \t V} \frac{|B(u,v)|}{\|u\|_V\|v\|_{\t V}}\geq \alpha\;\; {\rm and}\;\; \inf_{v\in \t V}\sup_{u\in  V} \frac{|B(u,v)|}{\|u\|_V\|v\|_{\t V}}\geq \alpha,
\label{infsupa}
\ee
for some $\alpha>0$
Then, the operator $\cB$ defined via \eref{defcB}
     is invertible and its inverse satisfies 
\be
\|\cB^{-1}\|_{\cL(W,V)}\leq \frac 1 \alpha.
\label{invbound1}
\ee
Hence, for each $L\in W$, the problem \eref{generellip} has a unique solution 
$u_L=\cB^{-1}(L)$ 
which satisfies the a priori estimate
\be
\|u_L\|_V\leq \frac{ \|L\|_W}{\alpha}.
\ee

\end{theorem}

\noindent
{\bf Proof:}   From the first inf-sup condition in \iref{infsupa} we obtain
that 
\be
\alpha\|u\|_V\leq \|\cB u\|_W,\quad u\inÊV.
\label{inject}
\ee
This shows that $\cB$ is injective and that its range $\cB(V)$ is closed in $W$.
In order to prove that $\cB$ is invertible we need to
show that $\cB(V)$ is all of $W$. We prove it by contradiction: 
if $\cB(V)$ was strictly contained in $W$, we can pick a non-trivial
$w\in W$ which is orthogonal in $W$ to all elements of $\cB(V)$. 
Then, we define $v=v(w)$ in the antidual of $W$, that is, $v\in \t V$ by setting
\be
v: e\mapsto \o{\<e,w\>}_W = \<e,v\>_{W,\t V}.
\ee
It follows that 
\be
B(u,v)=\< \cB u,v \>_{W,\t V}=v(\cB u)=\o{ \<\cB u,w\>}_W=0,
\ee
for all $u\in V$. This contradicts  the second
inf-sup condition. Hence $\cB$ is invertible and
the bound \iref{invbound1} on its inverse follows from \iref{inject}. \hfill $\Box$

\begin{remark}
One particular case of Theorem \ref{theoinfsup} is 
Theorem \ref{theolaxmilgram} since \iref{coerc} implies
\iref{infsupa} in the case $\t V=V$.
\end{remark}

The argument, centering on \iref{chain}, which justified the holomorphy of the solution
map for our canonical elliptic setting, may   be generalized to any linear problem of the form
\be
\cB(a)u=f(a),
\label{genlinear}
\ee
where $a\mapsto \cB(a)$ and $a\mapsto f(a)$ are holomorphic maps from
an open set $\cD\subset X$ into $\cL(V,W)$ and $W$ respectively,
where $(X,V,W)$ are complex Banach spaces.  Namely, if 
$\cB(a)$ is invertible for all $a\in\cD$, we find by
\be
a\mapsto \cB(a)\mapsto \cB(a)^{-1} \mapsto \cB(a)^{-1}f(a) =u(a),
\label{chain2}
\ee
that the solution map is holomorphic over $\cD$.

In particular, we may consider the following parametric linear problems of the form \iref{generellip}
for  a pair of Hilbert spaces $(V,\t V)$: for a given $a\in X$, find $u(a)\in V$ such that
\be
B(u(a),v;a)=L(v;a), \quad v\in \t V,
\label{generalvariatlin}
\ee
where $B(\cdot,\cdot;a)$ and $L(\cdot;a)$ are continuous sesquilinear
and antilinear forms over $V\times \t V$ and $\t V$ respectively, which depend on $a\in X$.

\begin{cor}
\label{corvariat}
If the assumptions of Theorem \ref{theoinfsup} are satisfied for the problem \iref{generalvariatlin}
for each  $a$ in a set $\cD\subset X$, then the operator $\cB(a)$ defined by $B(u,v;a)=\<\cB(a)u,v\>_{W,\t V}$ is well defined
and invertible from $V$ to $W=\t V^*$ for all $a\in \cD$ and the
solution map $a\mapsto u(a)$ is well defined from $\cD$ into $V$. If the constant $\alpha>0$  
can be chosen independent of $a\in\cD$ and if $ \sup_{a\in\cD}\|L(\cdot,a)\|_W=:M<\infty$, then the solution map is uniformly bounded on $\cD$ with
\be
\|u\|_{L^\infty(\cD,V)} =\sup_{a\in\cD} \|u(a)\|_V\leq \frac {M}\alpha.
\ee
In addition, if $\cD$ is an open set
and if the maps $a\mapsto L(\cdot;a)$ and $a\mapsto B(\cdot,\cdot;a)$ are 
holomorphic from $\cD$ into $W$ and $\cD$ into $\frak{B}=\frak{B}(V,\t V)$, respectively,
then\ the solution map $a\mapsto u(a)$ is holomorphic over $\cD$.
\end{cor}

We next give two examples that fall in this general framework. The first one is
a linear parabolic equation with parametrized evolution operator.
As a simple model, we consider the heat equation parametrized by its diffusion coefficient:
\be
\partial_t u={\rm div}(a\nabla u)+f,\quad {\rm in}\; ]0,T[\times D,
\label{parab}
\ee
where $D\subset \R^m$ is a Lipschitz domain, $f\in L^2(]0,T[;H^{-1}(D))$,   and the initial and boundary value conditions are
\be
u_{|t=0}=u_0\in L^2(D)\quad {\rm and}\quad u_{|\partial D}=0.
\ee
It   is well known that a solution space for this PDE is 
\be
V:=L^2(]0,T[;H^1_0(D))\cap H^1(]0,T[;H^{-1}(D)).
\ee
We obtain a space-time variational formulation of the type \iref{generalvariatlin} by introducing the
auxiliary space 
\be
\t V:=L^2(]0,T[;H^1_0(D))\times L^2(D)),
\ee 
and defining for $a\in X:=L^\infty$, $u\in V$ and $v=(v_1,v_2)\in \t V$,
\be
B(u,v;a):=\int_0^T\int_D\(\partial_t u(x,t)\o{v_1(x,t)}+a\nabla u(x,t)\cdot \nabla v_1(x,t)\)dx dt+\int_D u(\cdot,0) \o{v_2(x)}dx,
\ee
and
\be
L(v;a):=\int_0^T\<f(\cdot,t),v_1(\cdot,t)\>dt +\int_D u_0(x)\o{v_2(x)}dx,
\ee
where $\<\cdot,\cdot\>$ is the anti-duality pairing  between $H^{-1}(D)$ and $H^1_0(D)$.

The fact that these are bounded sesquilinear and antilinear forms 
follows readily from the choice of spaces $X$, $V$ and $\t V$. 
By using the general arguments from \cite{SS}, one can show that 
whenever the diffusion coefficient comes from the uniform ellipticity class $\cD_r$ of \iref{ellipre},
then the inf-sup condition \iref{infsupa} holds, with the values of $\alpha$  in \iref{infsupa}
depending on that of $r$.
Therefore, from Theorem \ref{theoinfsup},  the solution map 
$a\mapsto u(a)$ is defined on  $\cD_r$  with a uniform bound 
\be
\|u\|_{L^\infty(\cD_r,V)}=\sup\{ \|u(a)\|_V\; : \; a\in \cD_r\} \leq C_r.
\ee
Since $r>0$ is arbitrary the solution map is therefore defined on the open set $\cD:=\cup_{r>0}\cD_r$,
and its holomorphy follows from Corollary \ref{corvariat} since the sesquilinear form $B(\cdot,\cdot;a)$ depends
on $a$ in an affine manner.

The second example is a linear elliptic PDE parametrized by the shape of
the physical domain. As a simple model we consider the Laplace equation
\be
-\Delta w=1,
\label{delta1} 
\ee
set on a domain $D_a\subset \R^2$ with homogeneous Dirichlet boundary 
conditions $w_{|\partial D_a}=0$. Here $a$ describes the shape of the domain
$D_a$ in polar coordinates, according to
\be
D_a:=\{x=(\rho\cos\theta,\rho\sin \theta)\; : \; 0\leq \rho<a(\theta)\}.
\ee
In order to obtain a Lipschitz domain, 
we take $a\in X$ with
\be
X:=W^{1}_{\rm per}(L^\infty([0,2\pi[),
\ee 
the space of $2\pi$ periodic
Lipschitz continuous functions, which is equiped with the norm
\be
\|a\|_X:=\|a\|_{L^\infty([0,2\pi[)}+\|a'\|_{L^\infty([0,2\pi[)}.
\ee
If in addition, for some $r>0$, we have
\be
a(\theta)\geq r, \quad \theta\inÊ[0,2\pi[.
\label{starshape}
\ee
then $D_a$ is a Lipschitz domain and thus there exists a unique solution $w=w(a)\in H^1_0(D_a)$ to 
\iref{delta1} in the sense of the variational formulation
\be
\int_{D_a} \nabla w\cdot \nabla v=\int_{D_a} v,\quad v\in H^1_0(D_a).
\label{variatDa}
\ee
Note that \iref{starshape}Ê implies that $D_a$ is star-shaped with respect
to a ball of sufficiently small radius centered at the origin. In order to study the parametric smoothness
of $a\mapsto w(a)$, we need to represent the solution in a function space $V$ which does not
change with $a$. One  way to do this is to utilize  the pullback of the solution
to a reference domain $D$ under  a suitable transformation $F_a$ which maps $D$ into $D_a$. A natural choice is to take
for $D$ the unit disc of $\R^2$ centered at the origin, and use the transformation 
\be
F_a(\rho\cos\theta,\rho\sin\theta):= (a(\theta)\rho\cos\theta,a(\theta)\rho\sin\theta).
\ee
For $x=(\rho\cos\theta,\rho\sin\theta)$,  the jacobian matrix of $F_a(x)$ 
is given by
\be
dF_a(x)=\left ( \begin{array}{cc} a(\theta) \cos\theta &   a'(\theta)\cos\theta-a(\theta)\sin\theta\\ a(\theta)\sin\theta & a'(\theta)\sin\theta+a(\theta)\cos\theta\end{array}\right ),
\ee
and its determinant by
\be
J_a(\hat x):=a(\theta)^2 \geq r^2>0.
\label{minJa}
\ee 
We denote by 
\be
u=w\circ F_a\in V:=H^1_0(D),    
\ee 
the pullback solution and study the solution map $a\mapsto u(a)$. Using $F_a$ as a change of variable in \iref{variatDa}, we find
that $u$ satisfies
\be
\int_D  M_a\nabla u\cdot \nabla  v=\int_D J_a  v, \quad v\in V,  
\label{variatpull}
\ee
that is, the variational formulation of the equation
\be
-{\rm div}(M_a\nabla u)=J_a,
\label{eqpull}
\ee
set over the domain $D$ with homogeneous Dirichlet boundary conditions,
where 
\be
M_a(x):=J_a(\hat x)dF_a^{-1}(x)dF_a^{-t}(x)=\left (\begin{array}{cc} 1+b(\theta)^2  &  -b(\theta) \\ -b(\theta) & 1 \end{array}\right ),
\quad b(\theta):=\frac {a'(\theta)}{a(\theta)}.
\ee
We define the complex extension as a solution of \iref{generalvariatlin}, with $V=\t V=H^1_0(D)$ and the forms $B$ and $L$ given by
and 
\be
B(u,v; a):=\int_D  M_a\nabla u\cdot \nabla  v\quad {\rm and}\quad L(v; a):=\int_D J_a  \o v.
\ee
In view of the expressions of $M_a$ and $J_a$, it is readily seen
that $a\mapsto L(\cdot;a)$ is holomorphic from $X$ onto $W=V^*$ and that
$a\mapsto B(\cdot,\cdot;a)$ is holomorphic from the open set of $X$ of nowhere
vanishing functions into $\frak{B} =\frak{B} (V,V)$. It remains to understand for which $a\in X$
the problem has a solution. Introducing the real symmetric matrix
$R_a(x)=\Re(M_a(x))$ and denoting by $\lambda_{\min}(a,x)$ its smallest eigenvalue, 
we have for any $u\in V$,
\be
\Re(B(u,u;a))=\int_D R_a\nabla u\cdot \nabla  u \geq \lambda_{\min}(a) \|u\|_V^2, \quad \lambda_{\min}(a):=\min_{x\in D}\lambda_{\min}(a,x).
\ee
Therefore the coercivity condition \iref{coerc} holds if $\lambda_{\min}(a)>0$. A straightforward computation shows that
\be
\det(R_a(x))=1-\Im(b(\theta))^2\quad {\rm and} \quad {\rm tr}(R_a(x))=1+\Re(b(\theta))^2-\Im(b(\theta))^2.
\ee
We are thus are ensured of the existence of the solution $u(a)\in X$ for those $a\in X$ which are nowhere vanishing
and such that
\be
|\Im(b(\theta))|<1, \quad \theta\inÊ[0,2\pi[, \quad  b(\theta):=\frac {a'(\theta)}{a(\theta)}.
\label{imb}
\ee
By application of Corollary \ref{corvariat}, the solution map has a holomorphic extension onto the open domain $\cD\in X$ consisting of those $a\in X$ which are
nowhere vanishing and such that \iref{imb} holds. 

One important observation for this last example is the following: if the parameter domain $\cA$ is a compact
set of real valued functions in $X$ such that \iref{starshape}Ê holds for all $a\in \cA$, then there exists an
open neighbourhood $\cO$ of $\cA$ in the complex valued $X$ such that the holomorphic extension
of the solution map is uniformly bounded over $\cO$. Indeed, in view of the above remarks,
for every $a\in \cA$, there exists $\e=\e(a)>0$ such that 
\be
\mathring B(a,\e):=\{\t a\; : \; \|\t a-a\|_X<\e\} \subset \cD,
\ee
and such that the assumptions of Theorem \ref{theoinfsup} 
are satisfied for the problem \iref{generalvariatlin} 
with constants $\alpha$ and $C_L$ that are uniform over $\mathring B(a,\e)$. By compactness of $\cA$,
we may define $\cO$ as a finite cover of $\cA$ of the form
\be
\cO=\bigcup_{i=1}^M \mathring B(a_i,\e(a_i)),
\ee
for $\{a_1,\dots,a_M\}\in \cA$, so that $a\mapsto u(a)$ is holomorphic and uniformly bounded over $\cO$.

\subsection{Extensions by the implicit function theorem}

In this section, we consider a further generalization of   problems
of the form \iref{genpar} for which we can prove holomorphy of the solution map on certain subsets of the complex Banach space $X$.    In particular, this generalization can be applied to certain   nonlinear PDEs.  As a simple example, to motivate what follows, we consider the nonlinear
elliptic equation
\be
u^3-{\rm div}(a\nabla u)=f,
\label{nonlin}
\ee
set on a bounded Lipschitz domain $D\subset \R^m$ where $m=2$ or $3$, with homogeneous Dirichlet boundary conditions $u_{|\partial D}=0$,
parametrized by the diffusion coefficient $a$. Similar to the linear equation \iref{ellip}, 
we set 
\be
(X,V,W)=(L^\infty(D),H^1_0(D),H^{-1}(D)),
\ee
and consider for $f\in W$ and $a\in X$ the variational formulation
\be
\int_D u^3 v+\int_D a\nabla u\nabla v=\<f,v\>_{W,V},\quad v\in V.
\label{nonlinvariat}
\ee
By the theory of monotone operators, see for example Theorem 1 in Chapter 6 of \cite{RS},
and using the Sobolev embedding $H^1_0(D)\subset L^4(D)$, one can easily check that for any 
real valued $a\in X$ such that $a\geq r$ for some $r>0$,
there exists a unique solution $u(a)$ to \iref{nonlinvariat} which satisfies the a priori estimate
\be
\|u(a)\|_V\leq \frac {\|f\|_W}{r}.
\ee

 However, we cannot use monotone operator theory to derive a solution $u(a)$  to \eref{nonlinvariat} for complex valued $a$ since in this form the monotonicity is lost.
One could consider an alternative extension of  \iref{nonlinvariat} to complex valued functions given by
\be
|u|^2u-{\rm div}(a\nabla u)=f.
\label{nonlinmod}
\ee
For this equation, one can now apply monotone operator theory to the
real and imaginary part of the equation in order to show that the problem is well
posed when $\Re(a)\geq r$ for some $r>0$. However, this extension
is not holomorphic in the variable $a$ 
due to the presence of the modulus in \iref{nonlinmod}. 
We thus want to adhere to
the original problem \iref{nonlin} for complex valued $a$,
but find an alternative to using monotone operator theory.  This alternative is provided by 
 the following general theorem, which is based on the
holomorphic version
of the implicit function theorem in Banach spaces.

\begin{theorem}
\label{theoimp}
Let $\cP: V \times X \to W$ where $X$, $V$ and $W$ are complex Banach spaces
and let $\cA\subset X$ be a compact set such that
for each $a\in \cA$, there exists a unique solution $u(a)\in V$ to \iref{genpar}.
Assume, in addition, that there exists an open set $\cD$ of $X$ containing  $\cA$
for which
\vskip .1in
\noindent
{\rm (i)} $\cP$ is a holomorphic map from $V\times \cD$ to $W$,

\noindent
{\rm (ii)} for each $a\in \cA$, the partial differential $\partial_u \cP(u(a),a)$ is an isomorphism from $V$ to $W$.
\nl

\noindent
Then, there exists an open set  $\cO\subset X$
containing $\cA$, such that $u$ has a holomorphic extension to $\cO$
which takes values in $V$ and which is uniformly bounded:
\be
\|u\|_{L^\infty(\cO,V)}=\sup_{a\in \cO} \|u(a)\|_V <\infty.
\ee
\end{theorem}

\noindent
{\bf Proof:} Let $a\in \cA$. The assumptions (i) and (ii)
allow us to apply the holomorphic version of the implicit function theorem on complex
Banach spaces, see \cite[Theorem 10.2.1]{Di}, and conclude that 
there exists an $\e=\e(a)>0$, and a unique
holomorphic extension of $u$ from the open ball $\mathring B(a,  \e)$ 
of $X$, with center $a$
and radius $\e$ into $V$, such that $\cP(u(b),b)=0$ 
for any $b\in \mathring B(a,  \e)$. 
In addition, the map $u$ is uniformly bounded and holomorphic on $\mathring B(a, \e)$ 
with
\be
du_b 
= 
-\(\partial_u \cP(u(b),b)\)^{-1} \circ \partial \cP_b(u(b),b),\;\;\; 
b\in \mathring B(a,\e)
\;.
\ee
From the compactness of $\cA$, we may define $\cO$ as a finite cover of $\cA$ of the form
\be
\cO=\bigcup_{i=1}^M \mathring B(a_i,\e(a_i)),
\ee
for $\{a_1,\dots,a_M\}\in \cA$. Therefore $u$ has a uniformly bounded holomorphic extension over $\cO$.
\hfill $\Box$
\nl

There are many settings where Theorem \ref{theoimp} can be applied, including
nonlinear equations.  As an example, we return to   \iref{nonlin} where the 
operator $\cP$ is given by
\be
\cP(u,a)=u^3-{\rm div}(a\nabla u)-f,
\ee
or in variational form by
\be
\<\cP(u,a),v\>_{W,V}=\int_Du^3v+\int_Da\nabla u\nabla v-\<f,v\>_{W,V}.
\ee
Using the fact that $H^1_0(D)$ is continuously embedded into $L^4(D)$,
it is easily seen
that $\cP$ acts as a holomorphic map from $V\times X$ to $W$,
and therefore assumption (i) holds.

We now take for $\cA$ any compact set of $X$
contained in the set of real valued
functions $a\in X$ such that $a\geq r$ where $r>0$ is fixed,
so that there exists a unique solution $u(a)\in V$ for each $a\in \cA$.
We observe that, for $a\in\cA$,
\be
\partial_u \cP (u(a),a)(w) = 3u(a)^2 w-{\rm div}(a\nabla w).
\ee
The operator $\partial_u \cP (u(a),a)$ is associated to the sesquilinear form
\be
A(v,w;a)=\<\partial_u \cP (u(a),a)(v),w\>_{V',V}=\int_D  3u(a)^2  v \o w+\int_Da \nabla v\cdot  \overline {\nabla w}.
\ee
which is continuous over $V\times V$ (by the continuous embedding of  $H^1_0(D)$ into $L^4(D)$)
and satisfies the coercivity condition
\be
|A(v,v;a)|\geq \Re(A(v,v;a)) \geq r \|v\|_V^2, \;\; v\in V,
\ee
By the complex version of Lax-Milgram theorem, $\partial_u \cP (u(a),a)$ is thus an isomorphism
from $V$ onto $W$, and therefore assumption (ii) holds. We therefore conclude
from Theorem \ref{theoimp} that there exists an open set  $\cO\subset X$
containing $\cA$, such that the solution map has a uniformly bounded holomorphic extension over $\cO$.

\begin{remark}
Theorem \ref{theoimp} may also be applied to treat linear parametric problems such as 
the previously discussed elliptic, parabolic or domain dependent elliptic equations
which are already covered by the LBB theory.
Its weakness however is that it does not give an explicit description of the domain
where the holomorphic extension is defined, in contrast to the explicit conditions
on $a$ that can be established for these specific problems using Theorem \ref{theoinfsup}. 
Nevertheless, as it will be seen further, the
sole existence of a holomorphic extension on an open neighbourhood of the compact 
parameter domain $\cA$ turns out to be sufficient for
deriving approximation results for the solution map which are immune to the curse of
dimensionality.
\end{remark}

\subsection{The uniform ellipticity assumption}
\label{subuea}

For the remainder of this section, and all of \S 3, we assume that the parameter space $X$ is the complex Banach space $L^\infty$ and  that the parameter set $\cA\subset X$ has an affine scalar representation of the form \iref{affine}.   This allows us to view the solution map as   $y\mapsto u(y):=u(a(y))$.  
The focus of thie remainder of this section is to show that this map has a holomorphic extension to certain complex domains.

Here and in \S 3, we assume that {\bf Assumption A} holds 
for a suitable affine representer $(\psi_j)_{j\geq 1}$, which means that the solution map
$a\mapsto u(u)$ is well defined over 
\be
a(U):=\Big\{a(y) =\o a+\sum_{j\geq 1} y_j\psi_j\; : \; y=(y_j)_{j\geq 0}\in U\Big\},
\label{rectangle}
\ee
where $U=[-1,1]^{\N}$, so that the map $y\mapsto u(y)$ is well defined from $U$ to $V$. We also assume that
\be
(\|\psi_j\|_{X} )_{j\geq 1}\in \ell^1(\N),
\label{summability}
\ee
which implies that the series in \iref{affine} 
converge absolutely for all $y\in U$. In addition, this assumption guarantees the compactness
of the set $a(U)$ defined by \iref{rectangle}, as shown by the following result.

\begin{lemma}
\label{lemmacomp}
Under the assumption \iref{summability}, 
the set $a(U)$ defined by \iref{rectangle} is compact in $X$.
\end{lemma}

\noindent
{\bf Proof:} Let $(a_n)_{n\geq 1}$ be a sequence in $a(U)$.
Since $(\|\psi_j\|_{X})_{j\geq1}\in \ell^1(\N)$, the sequence
$(a_n)_{n\geq 1}$ is bounded in $X$.
Each $a_n$ is of the form $a_n=\sum_{j\geq 1} y_{n,j} \psi_j$.
Using a Cantor diagonal argument, we infer that 
there exists $y^* = (y_j)_{j\geq 1}\in U$ such that
\be
\lim_{n\to +\infty} y_{\sigma(n),j} = y_j^*, \;\; j\geq 1,
\ee
where $(\sigma(n))_{n\geq 1}$ 
is a monotone sequence of positive integers.
Defining $a^*:=\sum_{j\geq 1} y_{j}^* \psi_j \in a(U)$, 
we may write for any $k\geq 1$,
\be
\|a_{{\sigma(n)}}-a^*\|_{X} 
\leq 
\Big \|\sum_{j= 1}^k (y_{j}^*-y_{\sigma(n),j}) \psi_j\Big \|_X+2 \sum_{j\geq k+1} \|\psi_j\|_X.
\ee
It follows that $a_{\sigma(n)}$ converges towards $a^*$ in $X$ 
and therefore $a(U)$ is compact. \hfill $\Box$
\nl

Let us recall the four previously discussed examples
of parametric PDEs, that is, equations \iref{ellip}, \iref{parab}, \iref{eqpull}, and \iref{nonlin}.
For these problems, we have seen that the  solution map is defined at any 
  real valued $a\in X$  satisfies 
\be
a(x)\geq r, \quad x\in D,
\ee
for some $r>0$.   Here, 
  the physical domain $D$ is replaced by the angular domain $[0,2\pi[$ in
the case of \iref{eqpull}. When $\cA$ is a compact set of the form \iref{rectangle}, 
this condition is met for all $a\in \cA$ if and and only if
\be
a(x,y)=\o a(x)+\sum_{j\geq 1} y_j\psi_j(x)\geq r,  \quad  x\in D, \; y\in U.
\label{uea1}
\ee
By taking the particular choice $y_j=-{\rm sign}(\psi_j(x))$,
we find that the above inequality is equivalent to
\be
\sum_{j\geq 1} |\psi_j(x)| \leq \o a(x)-r, \quad x\in D.
\label{uea}
\ee
We refer to \iref{uea1} or \iref{uea} as the {\it uniform ellipticity assumption} of constant $r$, or ${\rm \bf UEA} (r)$.
In the case of \iref{eqpull},  the physical domain $D$ is replaced by 
the angular domain $[0,2\pi[$, and ${\rm \bf UEA} (r)$ thus
means that, for all $a\in \cA$, the domains $D_a$ are star-shaped with respect to a 
ball of sufficiently small radius centered at the origin.

Our previous analysis showed that  the assumption ${\rm \bf UEA} (r)$ ensures that in all four examples    the solution map
$y\mapsto u(y)$ is well defined from $U$ to $V$. 
We now want to build
an extension $z\mapsto u(z)$ by setting 
\be
a(z)=\o a+\sum_{j\geq 1} z_j \psi_j,
\label{affinez}
\ee
for suitable $z=(z_j)_{j\geq 1}\in \C^N$ and defining
\be
u(z):=u(a(z)).
\ee
This only makes sense for those $z\in\C^\N$ for which $a(z)$ is well defined and
falls inside the domain $\cD\subset X$
where $a\mapsto u(a)$ admits its holomorphic extension. At such a $z$, the chain rule
ensures that the resulting map
$z\mapsto u(z)$ is holomorphic in each variable $z_j$, with partial derivatives given by
\be
\partial_{z_j} u(z)=du(a(z))\psi_j.
\ee
Our next objective is to describe some relevant domains of $\C^\N$ on which the holomorphic extension 
$z\mapsto u(z)$ exists and is uniformly bounded.

\subsection{Holomorphic extensions of $y\mapsto u(y)$ on polydiscs}
\label{secpolydisc}

We first consider the elliptic problem \iref{ellip} and the  parabolic problem \iref{parab}.
For such problems, we have
seen that the solution map $a\mapsto u(a)$ admits a holomorphic extension
on the complex domain $\cD_r$ defined by the condition $\Re(a)\geq r$, with
uniform bound
\be
\|u\|_{L^\infty(\cD_r,V)} \leq C_r.
\ee
If ${\rm \bf UEA} (r)$ holds,
then
\be
\Re(a(x,z))=\o a(x) +\sum_{j\geq 1}\Re(z_j)\psi_j(x) \geq r,\quad x\in D,
\ee
holds for all $z\in\C^N$ such that $|\Re(z_j)|\leq 1$, and 
in particular for all $z \in \cU$, where $\cU$ is the unit {\it polydisc}
\be
\cU:=\{z=(z_j)_{j\geq 1}\; : \; |z_j|\leq 1\}=\otimes_{j\geq 1} \{|z_j|\leq 1\}.
\ee
This shows that the set
\be
a(\cU):=\Big\{a(z)=\o a+\sum_{j\geq 1} z_j\psi_j\; : \; z\in \cU\Big\},
\label{acU}
\ee
is contained in $\cD_r$.
In turn, the map $z\mapsto u(z)$ is holomorphic in each variable $z_j$
over $\cU$ with the  uniform bound 
\be
\label{ubound}
\sup_{z\in \cU}\|u(z)\|_V \leq C_r.
\ee

We next consider general polydiscs of the form
\be
\cU_{\rho}:=\{z=(z_j)_{j\geq 1}\; : \; |z_j|\leq \rho_j\}=\otimes_{j\geq 1} \{|z_j|\leq \rho_j\},
\ee
where $\rho=(\rho_j)_{j\geq 1}$ is a sequence of
positive numbers. Then, for any $t >0$ and any positive sequence $\rho$ that satisfies the constraint
\be
\sum_{j\geq 1} \rho_j |\psi_j(x)| \leq \o a(x)-t,\quad x\in D,
\label{constrho}
\ee
we find that 
\be
z\in\cU_\rho \Rightarrow \Re(a(x,z))\geq t, \quad x\in D,
\ee
which shows that $a(\cU_\rho)\subset \cD_t$.
Therefore, the map $z\mapsto u(z)$ is holomorphic in each variable $z_j$
over $\cU_\rho$ and the uniform bound \eref{ubound} now holds for some   constant $C_{t}>0$. 

If ${\rm \bf UEA}(r)$ holds
and if $0<t<r$, we can find sequences $\rho$ which satisty \iref{constrho}
and such that $\rho_j\geq 1$ for all $j\geq 1$, so that the 
polydisc $\cU_\rho$ contains $\cU$. In particular, let $\e:= r-t>0$ and $\rho=(\rho_j)_{j\geq 1}$ be any sequence
of numbers such that $\rho_j\geq 1$ for all $j\geq 1$ and that satisfies the constraint
\be
\sum_{j\geq 1} (\rho_j-1)\|\psi_j\|_X \leq \e.
\label{constrho1}
\ee
Then, using ${\rm \bf UEA}(r)$, we have
\be
\sum_{j\geq 1} \rho_j |\psi_j(x)| \leq \sum_{j\geq 1}  |\psi_j(x)|+\sum_{j\geq 1} (\rho_j-1)\|\psi_j\|_X \leq \o a(x)-r+\e =\o a(x)-t,\quad x\in D.
\ee
Therefore the map $z\mapsto u(z)$ is holomorphic in each variable $z_j$
over $\cU_\rho$ with again a uniform bound 
\be
\sup_{z\in \cU_\rho}\|u(z)\|_V \leq C_{t}.
\ee

We shall make further use of the following observation: if $\rho$ satisfies one the above constraints
\iref{constrho} or \iref{constrho1}, then   for each $j\geq 1$, there is  an open set $\cO_{\rho_j}\subset \C$
that contains the disc $\{|z_j|\leq \rho_j\}$ and
such that the map $z\mapsto u(z)$ is holomorphic in each variable $z_j$
over the tensorized set
\be
\cO_{\rho}:=\otimes_{j\geq 1}\cO_{\rho_j}.
\label{orho}
\ee
One possible choice is to take for $\cO_{\rho_j}$ the open disc
\be
\cO_{\rho_j}:=\{|z_j|<\t \rho_j\}, \quad \t \rho_j:=\rho_j+ \frac {\t t}{\sum_{j\geq 1}\|\psi_j\|_X},
\ee
for some $0<\t t <t$, since we then have 
\be
\sum_{j\geq 1} \t \rho_j |\psi_j(x)| \leq \sum_{j\geq 1} \rho_j |\psi_j(x)| +\t t\leq \o a(x)-(t-\t t),\quad x\in D,
\ee
which shows that $a(\cO_{\rho})\subset \cD_{t-\t t}$.

In \S 3, we exploit these domains of bounded holomorphy in order
to derive convergence results for polynomial approximations of the
type \iref{polynomial} that are obtained by truncation of the Taylor development
of $u(z)$ on suitable sets $\Lambda_n$. For now, let us observe
that the varying radii $\rho_j$ in each variable of the polydiscs $\cU_\rho$ 
reflect the {\it anisotropy} of the solution map $z\mapsto u(z)$.
Let us also note that the above discussion does not identify 
one particular polydisc $\cU_\rho$. Instead, it shows that bounded 
holomorphy holds on all of the polydiscs $\cU_\rho$
associated to any of the sequences $\rho$ which satisfy the constraint \iref{constrho1} or \iref{constrho}.

Let us next observe that the above procedure of extending the solution map to polydiscs is not restricted to just
the problems of Examples 1 and 2.  More generally, we can obtain bounded holomorphic extensions on similar polydiscs with constrainted radii,
for any parametric PDE that satisfies the assumptions of the following theorem.

\begin{theorem}
\label{theoneiacU}
Consider a parametric problem of the form \iref{genpar}
such that {\bf Assumption A} holds 
for a suitable affine representer $(\psi_j)_{j\geq 1}$. Assume that $(\|\psi_j\|_X)_{j\geq 1}\in \ell^1(\N)$
and that the solution map $a\mapsto u(a)$ admits a
holomorphic extension to an open set $\cO\subset X$ which contains the set $a(\cU)$
defined by \iref{acU}, with uniform bound
\be
\sup_{a\in\cO} \|u(a)\|_V \leq C.
\label{uaC}
\ee
Then, there exists $\e>0$ such that for any sequence
$\rho=(\rho_j)_{j\geq 1}$ of numbers larger than or equal to $1$  which satisfies the
constraint \iref{constrho1}, the following holds: for all $j\geq 1$, there exists
an open set $\cO_{\rho_j}\subset \C$ that contains the disc $\{|z_j|\leq \rho_j\}$ for  which  the map $y\mapsto u(y)$ admits an extension
to the set $\cO_{\rho}$ defined by \iref{orho}, and this extension is holomorphic in each variable $z_j$ 
with uniform bound
\be
\sup_{z\in\cO_\rho} \|u(z)\|_V \leq C,
\ee
with the same value of $C$ as in \iref{uaC}
\end{theorem}

\noindent
{\bf Proof:} We first observe that there exists $\delta>0$ small enough such that
the $\delta$-neigbourhood of $a(\cU)$ is contained in $\cO$, i.e.
\be
\bigcup_{a\in a(\cU)} B(a,\delta) \subset \cO.
\label{deltanei}
\ee
To see this, we observe that, by the same argument
as used for $a(U)$ in the proof of Lemma \ref{lemmacomp}, the set $a(\cU)$ is compact.
The distance function
\be
a\mapsto {\rm dist}(a,\cO^c):=\inf\{\|a-b\|_X \; : \; b\notin \cO\},
\ee
is continuous and strictly positive over $a(\cU)$. By compactness of $a(\cU)$, it
 reaches a strictly positive minimal value, and therefore \iref{deltanei}
holds by taking $\delta >0$ strictly smaller than this minimal value.

Next, we take $\e>0$ strictly smaller than $\delta$ and define
\be
\cO_{\rho_j}:=\{|z_j|<\t \rho_j\}, \quad \t \rho_j:=\rho_j+ \frac {\delta-\e}{\sum_{j\geq 1}\|\psi_j\|_X},
\ee
so that by \iref{constrho1}, we have
\be
\sum_{j\geq 1} (\t \rho_j-1)\|\psi_j\|_X =\sum_{j\geq 1} (\rho_j-1)\|\psi_j\|_X+\delta -\e \leq \delta.
\ee
For any $z\in \cO_{\rho}$, we define  $\t z_j:=z_j\min\{1,|z_j|^{-1}\}$ which gives that $\t z:=(\t z_j)_{j\ge 1}$ is in $\cU$ and
$$
\begin{disarray}{ll}
a(z) &= \o a +\sum_{j\geq 1} z_j \psi_j\\
& = \o a+ \sum_{j\geq 1}\t z_j \psi_j +\sum_{j\geq 1} (z_j-\t z_j)\psi_j \\
&=   a(\t z)+ r(z). \\
\end{disarray}
$$
Since,
\be
\|r(z)\|_X \leq \sum_{j\geq 1} |z_j-\t z_j| \|\psi_j\|_X \leq  \sum_{j\geq 1} (\t\rho_j-1) \|\psi_j\|_X \leq \delta,
\ee
it follows from \iref{deltanei} that $a(z)\in \cO$. Therefore $y\mapsto u(y)$ admits a holomorphic extension
over $\cO_{\rho}$ with at least the same uniform bound. \hfill $\Box$

\subsection{Holomorphic extensions of $y\mapsto u(y)$ on polyellipses}
The reader should notice that the results of the last section on extensions of the solution map to polydiscs
were not applied to two of our main examples: the 
parametrized domain problem
\iref{eqpull} and the nonlinear problem \iref{nonlin}.   For such problems, in contrast to the elliptic
and parabolic problems \iref{ellip} and \iref{parab},
we are not ensured that the solution map $a\mapsto u(a)$ 
admits a holomorphic extension to the whole domain $\cD_r$
defined by the condition $\Re(a)\geq r$.
In turn, while the uniform ellipticity assumption ${\rm \bf UEA}(r)$ 
ensures that the map $y\mapsto u(y)$ is well defined over $U$, 
it does not allow us to define its holomorphic extension
on the polydisc $\cU$, or more generally on polydiscs $\cU_\rho$ for
sequences $\rho$ which satisfy \iref{constrho}. 

On the other hand, for both problems \iref{eqpull} and \iref{nonlin},
we have seen that if $\cA$ is any compact set of real valued functions in $X$ such that
$a\geq r$ for all $a\in \cA$, there exists an open set $\cO\subset X$
which contains $\cA$ and such that the solution map $a\mapsto u(a)$
admits a holomorphic extension on $\cO$. In particular, if ${\rm \bf UEA}(r)$ holds
and $(\|\psi_j\|_X)_{j\geq 1}\in\ell^1(\N)$, we are ensured that such an open set
exists for $\cA=a(U)$ defined by \iref{rectangle}. This allows us to define
a bounded holomorphic extension to $y\mapsto u(y)$ on complex domains
that contain $U$, however with shorter extensions in the imaginary axes
than the polydiscs $\cU_\rho$. In \S 3, we exploit 
these domains of bounded holomorphy in order
to derive convergence results for polynomial approximations of the
type \iref{polynomial} that are obtained by truncation of the development
of $u(z)$ into orthogonal Legendre series.

To formulate the extensions we seek, we introduce some standard concepts from complex analysis.
 For any real number $s>1$,
we define in $\C$ the so-called {\it Bernstein ellipse},
\be
\cE_s:=\left\{\frac {z+z^{-1}} 2\; : \;  \; |z|=s\right \},
\ee
which has semi-axes of length $\frac {s+s^{-1}} 2$ in the real axis and $\frac {s-s^{-1}} 2$ in the imaginary axis.
Note that in the limit $s\to 1$,  we obtain $\cE_1=[-1,1]$.
The convex hull of $\cE_s$ is given by the filled-in ellipse
\be
\cH_s:=\left\{\frac {z+z^{-1}} 2\; : \; \; 1\leq |z|\leq s \right \}.
\ee
Note that
\be
[-1,1] \subset \cH_s \subset \{|z| \leq s\},
\ee
Therefore, defining for any sequence $(\rho_j)_{j\geq 1}$ of numbers strictly larger than $1$
the polyellipse 
\be
\cE_\rho:=\otimes_{j\geq 1} \cE_{\rho_j}
\ee
and the filled-in polyellipse
\be
\cH_\rho:=\otimes_{j\geq 1} \cH_{\rho_j}, 
\ee
we find that
\be
U\subset \cH_\rho \subset \cU_\rho.
\ee
However, the set $\cH_\rho$ has much shorter extension than $\cU_\rho$ in the imaginary axis
for the coordinates $j$ for which $\rho_j$ is close to $1$. This allows us to derive bounded holomorphic
extensions on such domains for the solution map $z\mapsto u(z)$ 
in the case of problems \iref{eqpull} and \iref{nonlin}, and more generally for any parametric 
PDE that fall under the assumptions of the following result. 

\begin{theorem}
\label{theoneiaU} Consider a parametric problem of the form \iref{genpar}
such that {\bf Assumption A} holds 
for a suitable affine representer $(\psi_j)_{j\geq 1}$. 
Assume that $(\|\psi_j\|_X)_{j\geq 1}\in \ell^1(\N)$
and that the solution map $u\mapsto u(a)$ admits a
holomorphic extension to an open set $\cO\subset X$ which contains the set $a(U)$
defined by \iref{rectangle}, with uniform bound
\be
\sup_{a\in\cO} \|u(a)\|_V \leq C.
\label{uaC1}
\ee
Then, there exists $\e>0$ such that for any sequence
$\rho=(\rho_j)_{j\geq 1}$ of numbers strictly larger than $1$  that  satisfies the
constraint \iref{constrho1}, the following holds: for all $j\geq 1$, 
there exists an open set $\cO_{\rho_j}\subset \C$ that contains the 
filled-in ellipse $\cH_{\rho_j}$ and such that 
the map $y\mapsto u(y)$ admits an extension
over the set $\cO_{\rho}$ defined by \iref{orho}, which is holomorphic in each variable $z_j$
with uniform bound
\be
\sup_{z\in\cO_\rho} \|u(z)\|_V \leq C,
\ee
with the same value of $C$ as in \iref{uaC1}
\end{theorem}

\noindent
{\bf Proof:} By the same compactness argument as used for $a(\cU)$ in the proof 
of Theorem \ref{theoneiacU}, there exists
$\delta>0$  sufficiently small that the  $\delta$-neigbourhood of $a(U)$ is contained in $\cO$, i.e.
\be
\bigcup_{a\in a(U)} B(a,\delta) \subset \cO.
\label{deltaneiU}
\ee
We now define $\e=\delta$  and  set $\cO_{\rho_j}$ to be the oval-shaped domain
\be
\cO_{\rho_j}:=\{z\in\C \; : \; {\rm dist}(z,[-1,1]):=\min_{y\in [-1,1]} |z-y| < \rho_j-1\},
\ee
for which is is easily checked that $\cH_{\rho_j}\subset \cO_{\rho_j}$.
For any $z=(z_j)_{j\geq 1}\in \cO_{\rho}$, there exists a $y=(y_j)_{j\geq 1} \in U$ such that
\be
|y_j-z_j| \leq \rho_j-1,\quad j\geq 1.
\ee
We may therefore write
\be
a(z) = \o a +\sum_{j\geq 1} z_j \psi_j = a(y)+ r(z), 
\ee
where 
\be
a(y) = \o a +\sum_{j\geq 1} y_j \psi_j  \in a(U),
\ee
and
\be
\|r(z)\|_X \leq \sum_{j\geq 1} |z_j-y_j | \|\psi_j\|_X \leq  \sum_{j\geq 1} (\rho_j-1) \|\psi_j\|_X \leq \e=\delta.
\ee
It follows from \iref{deltanei} that $a(z)\in \cO$.  Therefore $z\mapsto u(z)$ admits a uniformly bounded extension
over $\cO_{\rho}$. \hfill $\Box$

\begin{remark}
A general setting  in which the existence of $\cO$ in the above result
is satisfied is provided by Theorem \ref{theoimp}.
\end{remark}

\begin{remark}
The assumptions of Theorem \ref{theoneiaU} are obviously weaker than those of Theorem \ref{theoneiacU},
since $a(\cU)$ is replaced by $a(U)$.
\end{remark}

\begin{remark}
\label{remgenu1}
Theorems \ref{theoneiacU}  and \ref{theoneiaU} can be formulated
for a general map $u$ from $\cA$ to $V$ that is not necessarily the solution map
of a parametric PDE. Indeed, the only assumptions on $u$ which are used in the proof of these results is that 
it admits a bounded holomorphic extension in neighborhoods of $a(\cU)$ or $a(U)$. In other words,
the same conclusions in these theorems hold for any map $u$ which admits a bounded holomorphic extension 
on an open set containing $a(\cU)$ or $a(U)$.
On the other hand, as seen in \S 2.1, \S 2.2 and \S 2.3, the fact that $u$ is the solution map to a parametric
PDE can be utilized to prove the validity of these assumptions.
\end{remark}

\section {Best $n$-term polynomial approximations} 

In this section, we place ourself in
the same framework as in \S \ref{subuea}: we consider a
parametric problem \iref{genpar}, and assume that
{\bf Assumption A} holds
for a suitable affine scalar representation \iref{affine}.
The solution map $y\mapsto u(y):=u(a(y))$ 
is then well defined from $U$ to $V$. Our goal is now to establish 
convergence rates for specific separable approximations of this map
which are polynomials in the $y$ variable.
We construct these approximations by suitable truncations
of infinite expansions. 

\subsection{Approximation by $n$-term truncated expansions}

Let us begin with some general remarks concerning
the convergence towards $u$ of separable expansions of the form
\be
\sum_{\nu\in\cF} u_\nu \phi_\nu,
\label{expan}
\ee
where $\phi_\nu$ acts from $U$ to $\R$ and $u_\nu\in V$, for some countable index 
set $\cF$.

\begin{definition}
A sequence $(\Lambda_n)_{n\geq 1}$ of finite subsets of $\cF$ is called an exhaustion
of $\cF$ if and only if for any $\nu\in\cF$ there exists $n_0$ such that $\nu\in\Lambda_n$ for all $n\geq n_0$.
Here we do not impose that $\#(\Lambda_n)=n$.
\end{definition}

\begin{definition}

The series \iref{expan} is said to converge {\it conditionally} with limit $u$ in a given norm $\|\cdot\|$ if and
only if there exists an exhaustion $(\Lambda_n)_{n\geq 1}$ of $\cF$ such that
\be
\lim_{n\to +\infty}Ê\Big\|u-\sum_{\nu\in\Lambda_n} u_\nu\phi_\nu\Big\| =0.
\label{condconv}
\ee
The series \iref{expan} is said to converge {\it unconditionally} towards $u$ in the same norm, if and only if \iref{condconv}
holds for every exhaustion $(\Lambda_n)_{n\geq 1}$ of $\cF$. 
\end{definition}

We are interested both
in establishing unconditional convergence and providing estimates for the approximation error.
One first instance where this is feasible is when $(\phi_\nu)_{\nu\in\cF}$ is an orthonormal basis,
as indicated by the following result which simply gathers well known facts
from Hilbert space theory.

\begin{theorem}
\label{theoconvL2}
Let $(\phi_\nu)_{\nu\in\cF}$ be an orthonormal basis
of $L^2(U,\mu)$ for some given measure $\mu$ on $U$,
and let $u\in L^2(U,V,\mu)$. Then,  the inner products
\be
u_\nu:=\int_{U} u(y)\phi_\nu(y) d\mu(y),\quad \nu\in\cF,
\ee
are elements of $V$,  and the series \iref{expan}
converges unconditionally towards $u$ in $L^2(U,V,\mu)$, with the error given by
\be
\Big \|u-\sum_{\nu\in\Lambda_n} u_\nu\phi_\nu\Big\|_{L^2(U,V,\mu)}= \( \sum_{\nu\notin\Lambda_n} \|u_\nu\|_V^2\)^{1/2},
\label{estimL2}
\ee
for any exhaustion $(\Lambda_n)_{n\geq 1}$.
\end{theorem}

Let us observe that if $\mu$ is any probability measure, 
the $L^\infty(U,V)$ norm controls the $L^2(U,V,\mu)$ norm.
In the previous section, we have given various examples for
which we are ensured that $u$ is uniformly bounded over $U$,
and we may therefore apply the above result whenever
$\mu$ is a probability measure.

We next  give a general result which can be used to establish convergence and give error bounds
in the $L^\infty$ norms for truncating the expansion \eref{expan}.   
\begin{theorem}
\label{theoincond}  Consider an expansion \eref{expan} for which the following hold:

\noindent
{\rm (i)} The functions $\phi_\nu:U\mapsto \R$
  are such that $\|\phi_\nu\|_{L^\infty(U)}=1$,
for all $\nu\in\cF$. 

\noindent
{\rm (ii)} The functions $u_\nu$ are in $V$ and  $(\|u_\nu\|_V)_{\nu\in \cF}\in \ell^1(\cF)$,

\noindent
Then, whenever the 
expansion \iref{expan} converges conditionally to a function $u$
in $L^\infty(U,V)$, it also converges unconditionally to $u$ in $L^\infty(U,V)$, and for any exhaustion $(\Lambda_n)_{n\geq 1}$, we have the error estimate
\be
\Big \|u-\sum_{\nu\in\Lambda_n} u_\nu\phi_\nu\Big\|_{L^\infty(U,V)}\leq \sum_{\nu\notin\Lambda_n} \|u_\nu\|_V.
\label{estimLinf}
\ee

\end{theorem}

\noindent
{\bf Proof:} Let $(\Lambda_n)_{n\ge 1}$ be any given exhaustion of $\cF$ and suppose that $\e>0$ is arbitrary.
We know that there exists an exhaustion $(\Lambda^*_n)_{n\geq 1}$  and an   $n_0$ such that 
\be
\|u-\sum_{\nu\in\Lambda_n^*} u_\nu\phi_\nu\Big\|_{L^\infty(U,V)}\leq \e, \quad n\geq n_0.
\ee
In addition, there exists $m$ larger than $n_0$ such that 
\be
\sum_{\nu\notin\Lambda^*_{m}} \| u_\nu\|_V \leq \e.
\ee
Since  $(\Lambda_n)_{n\geq 1}$ is an   exhaustion, there exists $n_1$ such that
$\Lambda_m^*\subset \Lambda_n$ for all $n\geq n_1$, and therefore  \be
\|u-\sum_{\nu\in\Lambda_n} u_\nu\phi_\nu\Big\|_{L^\infty(U,V)}\leq\|u-\sum_{\nu\in\Lambda_m^*} u_\nu\phi_\nu\Big\|_{L^\infty(U,V)}+\sum_{\nu\notin\Lambda^*_{m}} \| u_\nu\|_V \leq 2\e, \quad n\geq n_1.
\ee
This confirms the unconditional convergence.  The estimate \iref{estimLinf} follows by an application of the  triangle inequality. \hfill $\Box$.
\nl

In the particular case where $(\phi_\nu)_{\nu\in \cF}$ is an orthogonal basis normalized in $L^\infty$,
the next theorem shows the same result holds without the need to assume conditional convergence.

\begin{theorem}
\label{theoLinforth}
Let $(\phi_\nu)_{\nu\in\cF}$be  an orthogonal basis
of $L^2(U,\mu)$ for some given probability measure $\mu$ on $U$, normalized so that $\|\phi_\nu\|_{L^\infty(U)}=1$
for all $\nu\in\cF$.  If  $u\in L^2(U,V,\mu)$ and the inner products
\be
u_\nu:=\frac {1}{\|\phi_\nu\|_{L^2(U,V,\mu)}^2}\int_{U} u(y)\phi_\nu(y) d\mu(y),\quad \nu\in\cF,
\ee
 satisfy  $(\|u_\nu\|_V)_{\nu\in \cF}\in \ell^1(\cF)$, then 
$u\in L^\infty(U,V)$ and   the series \iref{expan}
converges unconditionally towards $u$ in $L^\infty(U,V)$ and  the estimate \iref{estimLinf} holds.
\end{theorem}

\noindent
{\bf Proof:} The summability of $(\|u_\nu\|_V)_{\nu\in \cF}$ ensures that \iref{expan} converges to a limit in
$L^\infty(U,V)$ and in turn in $L^2(U,V,\mu)$. On the other hand, 
we know from Theorem \ref{theoconvL2} that it converges
toward $u\in L^2(U,V,\mu)$. Therefore, its limit in $L^\infty(U,V)$ is also $u$. \hfill $\Box$
\nl

If the  expansion \iref{expan} converges unconditionally towards $u$ in some given norm $\|\cdot\|$, then a
  crucial issue is the choice of sets $\Lambda_n$ that we decide to use 
to truncate \iref{expan} and define an $n$-term approximation. 
Since $n$ measures the complexity of this approximation,
we would like to find the set $\Lambda_n$ which minimizes the truncation
error in some given norm $\|\cdot\|$ among all
sets of cardinality $n$, i.e.
\be
\Lambda_n:={\rm argmin}\left\{Ê\Big\|u-\sum_{\nu\in\Lambda} u_\nu\phi_\nu\Big\|\; :\; \#(\Lambda)=n \right\},
\ee
provided that such a set exists. This is an instance of {\it best $n$-term} approximation, which itself
is an instance of {\it nonlinear approximation}. We refer to \cite{De} for a 
general survey on nonlinear approximation.

In the case where the error is measured in $L^2(U,V,\mu)$, and if
$(\phi_\nu)_{\nu\in\cF}$ is an orthonormal basis of $L^2(U,\mu)$ and $(u_\nu)_{\nu\in\cF}$
are the coefficients of $u$ in this basis, \iref{estimL2} shows that the optimal $\Lambda_n$ is
the set of indices corresponding to the $n$ largest $\|u_\nu\|_V$. Note that such a set
is not necessarily unique, in which case any realization of $\Lambda_n$ is optimal.

In the case where the error is measured in $L^\infty(U,V)$, there is generally
no simple description of the optimal set $\Lambda_n$. However, when the 
functions $\phi_\nu$ are normalized in $L^\infty(U)$, the right-hand side 
in the estimate \iref{estimLinf}  provides a bound for the error of $n$ term approximation.
This upper  bound is   minimized by again defining $\Lambda_n$ as
the set of indices corresponding to the $n$ largest $\|u_\nu\|_V$. The only difference
is that the error is measured by the $\ell^1$ tail of the sequence $(\|u_\nu\|_V)_{\nu\in \cF}$,
in contrast to the $\ell^2$ tail which appears in \iref{estimL2}.  Let us emphasize that
this procedure gives only a  bound for the error of best $n$ term approximation in $L^\infty(U,V)$
but is not guaranteed to be the best error.

There is a good understanding of the properties of a given
sequence $(c_\nu)_{\nu\in\cF}$ of real or complex numbers,  which ensure a certain rate of decay $n^{-s}$ of its
$\ell^q$ tail after one retains its $n$ largest entries.
The following result, originally due to Stechkin in the particular case $q=2$,
show that this rate of decay is related to the $\ell^p$ summability of the sequence
for values of $p$ smaller than $q$.

\begin{lemma}
\label{stechkin}
Let $0<p<q<\infty$ and let $(c_\nu)_{\nu\in\cF}\in \ell^p(\cF)$ be a sequence of positive numbers.
Then, if $\Lambda_n$ is a set of indices which corresponds to the $n$ largest $c_\nu$, one has
\be
\(\sum_{\nu\notin\Lambda_n} c_\nu^q\)^{1/q} \leq C (n+1)^{-s}, \quad C:=\|(c_\nu)_{\nu\in\cF}\|_{\ell^p}, \quad s:=\frac 1 q-\frac 1 p.
\label{bestnrate}
\ee
\end{lemma}

\noindent
{\bf Proof:} Let $(c_k)_{k\geq 1}$ be the decreasing rearrangement of the sequence $(c_\nu)_{\nu\in\cF}$.
From the definition of $\Lambda_n$, we have
\be
\sum_{\nu\notin\Lambda_n} c_\nu^q =\sum_{k\geq n+1} c_k^q \leq c_{n+1}^{q-p}\sum_{k\geq n+1}c_k^{p}  \leq C^p c_{n+1}^{q-p}.
\ee
On the other hand, we also have
\be
(n+1)c_{n+1}^p \leq \sum_{k=0}^{n+1}c_k^{p}\leq C^p.
\ee
Combining both estimates gives
\be
\sum_{\nu\notin\Lambda_n} c_\nu^q \leq C^q(n+1)^{\frac q p-1},
\ee
which is \iref{bestnrate}. \hfill $\Box$
\newline

Combining the above result with either \iref{estimL2} or \iref{estimLinf}  shows that a suitable $\ell^p$ summability of the sequence $(\|u_\nu\|_V)_{\nu\in \cF}$
is  a sufficient condition to guarantee a convergence rate $n^{-s}$   when retaining  the terms 
corresponding to the $n$ largest $\|u_\nu\|_V$ in \iref{expan}.   
For the $L^2(U,V,\mu)$ error, and when $(\phi_\nu)_{\nu \in\cF}$ is an
orthonormal basis of $L^2(U,\mu)$, we obtain the rate $s=\frac 1 p-\frac 1 2$ if $p<2$.
For the $L^\infty(U,V)$ error, and when the $\phi_\nu$ are normalized in $L^\infty(U)$,
we obtain the rate $s=\frac 1 p-1$ if $p<1$.

\begin{remark}
\label{wlp}
Lemma \ref{stechkin} shows that $\ell^p$ summability implies
that the $\ell^q$ tail of $(c_\nu)_{\nu\in \cF}$
after retaining the largest $n$-terms decays with rate $n^{-s}$ where $s:=\frac 1 p-\frac 1 q$.
It is actually possible to {\em exactly} characterize the properties which governs this rate of decay through
weaker summability properties. Let us recall that for $0<p<\infty$ the space $w\ell^p(\cF)$ consists of those sequences
$(c_\nu)_{\nu\in \cF}$ of real or complex numbers such that for a finite constant $C\geq 0$,
\be
\#\{\nu \; :\; |c_\nu| \geq \eta\} \leq C^p  \eta^{-p},\quad \eta >0,
\ee
or equivalently such that for a finite constant $C\geq 0$,
\be
c_k\leq Ck^{-1/p},\quad k\geq 1,
\ee
where $(c_k)_{k\geq 1}$ is the decreasing rearrangement of $(|c_\nu|)_{\nu\in \cF}$. 
The quasi-norm $\|(c_\nu)_{\nu\in \cF}\|_{w \ell^p(\cF)}$ can be defined as the smallest 
$C$ for which either one of these inequalities holds. Then, for $0<p<q\leq \infty$ 
one can check that the $\ell^q$ tail of $(c_\nu)_{\nu\in \cF}$
after retaining the largest $n$-terms decays with rate $n^{-s}$ where $s:=\frac 1 p-\frac 1 q$
if and only if $(c_\nu)_{\nu\in \cF}\in w\ell^p(\cF)$, see \cite{De}. 
\end{remark}

\subsection{Convergence of $n$-term truncated polynomial expansions}

We now restrict our attention to polynomial series.  This corresponds to particular choices of the functions $\phi_\nu$ as polynomials.
For  the remainder of this section, we take  $\cF$ to be the set of all sequences $\nu=(\nu_j)_{j\geq 1}$ of non-negative integers
which are finitely supported. For $\nu\in\cF$, we use the notation
\be
\|\nu\|_0:=\#({\rm supp}(\nu))<\infty,\quad {\rm supp}(\nu):=\{j\geq 1\; : \; \nu_j\neq 0\}.
\ee
as well as
\be
|\nu|:=\|\nu\|_1=\sum_{j\geq 1} \nu_j <\infty.
\ee
For any $z=(z_j)\in\C^\N$, and $\nu\in \cF$, we  define
\be
z^\nu:=\prod_{j\geq 1}z_j^{\nu_j}.
\ee
We consider three type of polynomial series:
\begin{itemize}
\item
Taylor (or power) series of the form
\be
\sum_{\nu\in\cF} t_\nu y^\nu,
\label{taylor}
\ee
where 
\be
t_{\nu}:=\frac 1 {\nu !}\partial^\nu u(y=0), \quad \nu!:=\prod_{j\geq 1} \nu_j!,
\ee
with the convention that $0!=1$.
\item
Legendre series of the form
\be
\sum_{\nu\in\cF} v_\nu L_\nu(y), \quad L_\nu(y)=\prod_{j\geq 1}L_{\nu_j}(y_j),
\label{legendre}
\ee
where $(L_k)_{k\geq 0}$ is the sequence of Legendre polynomials on $[-1,1]$ 
normalized with respect to the uniform measure, i.e. such that
\be
\int_{-1}^1 |L_k(t)|^2 \frac {dt}2=1.
\ee
It follows that $(L_\nu)_{\nu\in \cF}$ is an orthonormal basis of $L^2(U,\mu)$, 
with 
\be
\mu=\otimes_{j\geq 1}\frac {dy_j} 2,
\ee 
the uniform measure over $U$. The coefficients $v_\nu$ are therefore given by
\be
v_\nu:=\int_U u(y)L_\nu(y) d\mu(y).
\label{unu}
\ee
\item
Renormalized Legendre series of the form
\be
\sum_{\nu\in\cF} w_\nu P_\nu(y), \quad P_\nu(y)=\prod_{j\geq 1}P_{\nu_j}(y_j),
\label{renormlegendre}
\ee
where $(P_k)_{k\geq 0}$ is the sequence of Legendre polynomials on $[-1,1]$ with the standard normalization 
\be
\|P_k\|_{L^\infty([-1,1])}=P_k(0)=1.
\ee
One has $L_k=\sqrt{1+2k}P_k$, and therefore the coefficients $w_\nu$ are given by
\be
w_\nu:=\(\prod_{j\geq 1} (1+2\nu_j)\)^{1/2} v_\nu,
\ee
where $v_\nu$ is defined by \iref{unu}.
\end{itemize}

In the case of the Taylor series, the following result shows that 
the assumptions in Theorem \ref{theoneiacU} ensure the conditional
convergence of \iref{taylor} towards $u$ in $L^\infty(U,V)$.

\begin{theorem}
\label{theocond}
Consider a parametric problem of the form \iref{genpar}
such that {\bf Assumption A} holds 
for a suitable affine representer $(\psi_j)_{j\geq 1}$. 
If the assumptions of Theorem \ref{theoneiacU} are satisfied, then
the Taylor expansion \iref{taylor} converges conditionally towards $u$ in $L^\infty(U,V)$.
\end{theorem}

\noindent
{\bf Proof:} Under the assumptions of Theorem \ref{theoneiacU},
the Frechet derivative of the solution map $a\mapsto u(a)$ is uniformly bounded
over $a(\cU)$ and therefore
\be
\label{boundFder}
M:=\max_{a\in a(\cU)} \|du(a)\|_{\cL(X,V)} <\infty.
\ee
From the assumption that  $(\|\psi_j\|_X)_{j\geq 1}\in \ell^1(\cF)$, for any $n\geq 1$, there exists $J=J(n)$ be such that 
\be
\sum_{j\geq J+1} \|\psi_j\|_X \leq \frac {1}{2nM}.
\ee
Increasing the value of $J$ decreases the left side, so we may assume that $J(n)\geq n$.

We know from Theorem \ref{theoneiacU} that the map $y\mapsto u(y)$
admits a holomorphic extension $z\mapsto u(z)$ to  domains $\cU_\rho$
that contain $\cU$. For any $z=(z_j)_{j\geq 1}\in\cU$, we define its truncation
\be
T_Jz:=(z_1,\dots,z_J,0,0,\dots),
\label{Jtrunc}
\ee
and the map 
\be
v(z):=u(T_Jz)=u(a(T_Jz))=u\(\o a+\sum_{j=1}^J z_j \psi_j\).
\ee
Since, for $z\in \cU$, we have
\be
\|a(z)-a(T_Jz)\|_X \leq \sum_{j\geq J+1}\|\psi_J\|_X \leq \frac {1}{2nM},
\ee
it follows from \eref{boundFder} that
\be
\|u-v\|_{L^\infty(\cU,V)} \leq \frac 1 {2n}.
\ee

Now, we can write
\be
v(z)=w(z_1,\dots,z_J),
\label{finitew}
\ee
where the finite dimensional map $w$ is bounded and holomorphic in each variable $z_j$
on an open neighborhood of the unit polydisc $\cU_J:=\otimes_{j=1}^J \{|z_j|\leq 1\}$. 
It follows that $w$ has a Taylor expansion that converges on $\cU_J$. Its Taylor coefficients are given by the 
$t_{\nu}$ for all $\nu$ of the form $(\nu_1,\dots,\nu_J,0,0,\dots)$. Therefore, there exists 
$K=K(n)\ge n$ such that for
\be
\Lambda_n:=\{\nu\in\cF \; :\; {\rm supp}(\nu)\subset \{1,\dots,J\}\; {\rm and} \; |\nu|\leq K\},
\ee
one has
\be
\sup_{z \in \cU} \Big\|v(z)-\sum_{\nu\in \Lambda_n}t_\nu z^\nu\Big\|_V \leq \frac 1 {2n},
\ee
and therefore
\be
\sup_{z \in \cU} \|u(z)-\sum_{\nu\in \Lambda_n}t_\nu z^\nu\|_V \leq \frac 1 {n},
\ee
Since both    $K(n)$ and $J(n)$ tend to infinity with $n$,  the sequence of sets $(\Lambda_n)_{n\geq 0}$ is a exhaustion of $\cF$.
We have thus proved
the conditional convergence of \iref{taylor} towards $u$ in $L^\infty(\cU,V)$,
and in turn in $L^\infty(U,V)$.
 \hfill $\Box$
\nl

We are now in position to state our main result which gives simple
conditions that guarantee the $\ell^p$ summability, $0<p<1$,  of the sequence $(\|u_\nu\|_V)_{\nu\in\cF}$,
where $u_\nu$ is either $t_\nu$, $v_\nu$ or $w_\nu$. 
These conditions are expressed in terms of the $\ell^p$ summability
of the sequence $(\|\psi_j\|_X)_{j\geq 1}$ for the {\it same value of $p$},
and   the assumptions in Theorem \ref{theoneiacU} in the case of $t_\nu$
or in Theorem \ref{theoneiaU} in the case of $v_\nu$ or $w_\nu$.

\begin{theorem}
\label{theomain}
Consider a parametric problem of the form \iref{genpar}
such that {\bf Assumption A} holds 
for a suitable affine representer $(\psi_j)_{j\geq 1}$.  Then, the following summability results hold:
\begin{itemize}
\item[{\rm (i)}]
If the assumptions of Theorem \ref{theoneiacU} are satisfied,
and if in addition $(\|\psi_j\|_X)_{j\geq 1}\in \ell^p(\N)$ for some $p<1$, 
then $(\|t_\nu\|_V)_{\nu\in\cF}\in \ell^p(\cF)$
for the same value of $p$. 
\item[{\rm (ii)}]
If the assumptions of Theorem \ref{theoneiaU} are satisfied,
and if in addition $(\|\psi_j\|_X)_{j\geq 1}\in \ell^p(\N)$ for some $p<1$, 
then $(\|v_\nu\|_V)_{\nu\in\cF}\in \ell^p(\cF)$ and $(\|w_\nu\|_V)_{\nu\in\cF}\in \ell^p(\cF)$
for the same value of $p$.
\end{itemize} 
\end{theorem}

The proof of this result is given in \S \ref{subproof}.
For now, we use this theorem together with  the previous results of this section to  obtain corollaries
on the rate of convergence of $n$-terms approximations
obtained by truncation of Taylor or Legendre series.

\begin{cor}
\label{corratetay}
Consider a parametric problem of the form \iref{genpar}
such that {\bf Assumption A} holds 
for a suitable affine representer $(\psi_j)_{j\geq 1}$.   If the assumptions of Theorem \ref{theoneiacU} are satisfied,
and if in addition $(\|\psi_j\|_X)_{j\geq 1}\in \ell^p(\N)$ for some $p<1$,
then the Taylor series \iref{taylor} converges unconditionally towards $u$ in $L^\infty(U,V)$. Moreover,
for any set $\Lambda_n$  of indices corresponding  to   $n$ largest of  $\|t_\nu\|_{V}$,
we have
\be
\sup_{y\in U} \Big\|u(y)-\sum_{\nu\in\Lambda_n} t_\nu y^\nu\Big\|_V \leq C(n+1)^{-s}, \quad s=\frac 1 p-1,
\label{convtaylor}
\ee
where $C:=\|(\|t_\nu\|_V)_{\nu\in\cF}\|_{\ell^p}<\infty$.
\end{cor}

\noindent
{\bf Proof:} Using the first part of Theorem \ref{theomain}, we are ensured that $(\|t_\nu\|_V)_{\nu\in\cF}\in\ell^p(\cF)$.
Since, by Theorem \ref{theocond}, the   series \iref{taylor} converges conditionally,  by application of Theorem \ref{theoincond},  we find that it also converges unconditionally with the error bound
\be
\sup_{y\in U} \Big\|u(y)-\sum_{\nu\in\Lambda_n} t_\nu y^\nu\Big\|_V \leq \sum_{\nu\notin\Lambda_n}\|t_\nu\|_V.
\ee
We now use  \iref{bestnrate} with $c_\nu=\|t_\nu\|_V$
and $q=1$ to obtain the error bound \iref{convtaylor}. \hfill $\Box$

\begin{cor}
\label{corrateleg} Consider a parametric problem of the form \iref{genpar}
such that {\bf Assumption A} holds 
for a suitable affine representer $(\psi_j)_{j\geq 1}$. 
If the assumptions of Theorem \ref{theoneiaU} are satisfied,
and if in addition $(\|\psi_j\|_X)_{j\geq 1}\in \ell^p(\N)$ for some $p<1$,
then the Legendre series \iref{legendre} and \iref{renormlegendre} converges unconditionally towards $u$ in 
$L^\infty(U,V)$ and in 
$L^2(U,V,\mu)$ where $\mu$ is the uniform probability measure. In addition, we have
the following error bounds:
\begin{itemize}
\item
If $\Lambda_n$ is the set of indices that corresponds to the $n$ largest $\|v_\nu\|_{V}$,
we have
\be
\Big \|u-\sum_{\nu\in\Lambda_n} v_\nu L_\nu \Big\|_{L^2(U,V,\mu)}\leq C(n+1)^{-s}, \quad s=\frac 1 p-\frac 1 2,
\label{convlegendre}
\ee
where $C:=\|(\|v_\nu\|_V)_{\nu\in\cF}\|_{\ell^p}<\infty$.
\item
If $\Lambda_n$ is the set of indices that corresponds to the $n$ largest $\|w_\nu\|_{V}$,
we have
\be
\Big \|u-\sum_{\nu\in\Lambda_n} w_\nu P_\nu \Big\|_{L^\infty(U,V)} \leq C(n+1)^{-s}, \quad s=\frac 1 p-1,
\label{convrenormlegendre}
\ee
where $C:=\|(\|w_\nu\|_V)_{\nu\in\cF}\|_{\ell^p}<\infty$.
\end{itemize}
\end{cor}

\noindent
{\bf Proof:} Using the second part of Theorem
\ref{theomain}, we are ensured that $(\|v_\nu\|_V)_{\nu\in\cF}$ and $(\|w_\nu\|_V)_{\nu\in\cF}$ 
belong to $\ell^p(\cF)$. The unconditional convergence claims in the theorem are ensured by Theorems \ref{theoconvL2} and \ref{theoLinforth}.  These latter two theorems  also give the estimates
\be
\|u-\sum_{\nu\in\Lambda_n} v_\nu L_\nu\|_{L^2(U,V,\mu)}
=\(\sum_{\nu\notin\Lambda_n} \|v_\nu\|_V^2\)^{1/2},
\ee
and
\be
\|u-\sum_{\nu\in\Lambda_n} w_\nu P_\nu\|_{L^\infty(U,V,\mu)}
=\sum_{\nu\notin\Lambda_n} \|w_\nu\|_V.
\ee
The  application of \iref{bestnrate} with $c_\nu=\|v_\nu\|_V$ and $q=2$,  or  with $\|w_\nu\|_V$
and $q=1$, give the error bounds \iref{convlegendre} and \iref{convrenormlegendre}. \hfill $\Box$

\begin{remark}
Note that since we have
\be
v_\nu L_\nu=w_\nu P_\nu, \quad \nu\in \cF,
\ee
the terms in the series \iref{legendre} and \iref{renormlegendre} are actually identical. However
the sets $\Lambda_n$ defined by the $n$ largest $\|v_\nu\|_V$ or the
$n$ largest $\|w_\nu\|_V$, which are used to define the truncations for $L^2$ or $L^\infty$ estimates
in the previous result, generally differ from each other.
\end{remark}

The above corollaries show the curse of dimensionality can be broken
for relevant class of parametric PDEs:
although the solution map $y\mapsto u(y)$ has
infinitely many variables, it can be approximated in various norms 
with an algebraic rate $n^{-s}$, where $n$ is the number of term in the separable expansion. 
The exponent $s$ can be large if $(\|\psi_j\|_X)_{j\geq 1}\in \ell^p(\N)$ for a small value of $p$.
Several critical ingredients have been used in order to reach this conclusion:
\begin{itemize}
\item
The holomorphic extension of the solution map $a\mapsto u(a)$.
\item
The anistrotropy of the solution map with respect to the different variables $y_j$.
\item
The use of best $n$-term polynomial approximations.
\end{itemize}

The fact that anisotropic smoothness may allow certain numerical methods to break the curse of dimensionality,
in the sense that approximation results are immune to the growth in the number of variables,
has also been studied in information based complexity, using certain weight sequences
in oder to quantify anisotropy, see \cite{KSWW}.

\begin{remark}
\label{remgenu2}
Theorem \ref{theomain} and its corollaries can be formulated
for a general map $u$ from $\cA$ to $V$ that is not necessarily the solution map
of a parametric PDE, since as observed in Remark \ref{remgenu1}, 
Theorems \ref{theoneiacU} and \ref{theoneiaU} hold in this more general framework.
\end{remark}

\subsection{Estimates of Taylor coefficients}
\label{subesttay}

In this section, as well as the two that follow, we establish upper estimates for the $V$-norms 
of the Taylor coefficients $t_\nu$ and Legendre coefficients $v_\nu$ and $w_\nu$,
which are instrumental in the proof of Theorem \ref{theomain}. 
These estimates are derived from the results on holomorphic extensions  of  the map $y\mapsto u(y)$
established in Theorems \ref{theoneiacU} and \ref{theoneiaU}.  Namely,  by an application
of the Cauchy integral formula in the different 
complex variables $z_j$.

We recall that the Cauchy formula states that
if $\vp$ is a function from $\C$ to a Banach space $V$ which is 
holomorphic on a simply connected open set $\cO\subset \C$ and if $\Gamma$ is a closed rectifiable
curve contained in $\cO$, then for any $\t z$ contained in the bounded domain
delimited by $\Gamma$,
\be
\vp(\t z):=\frac 1 {2i\pi} \int_{\Gamma}\frac {\vp(z)} {\t z-z}  dz,
\label{cauchy}
\ee
where the fraction in the integrand stands for the scalar multiplication of $\vp(z)\in V$
by   $(\t z -z)^{-1}\in\C$ and the curve $\Gamma$ is positively oriented in the integral,
see for instance Theorem 2.1.2 of \cite{H}. 

We begin with the estimates on Taylor coefficients which are based
on the bounded holomorphic extensions onto polydiscs $\cU_\rho$ that were established in Theorem \ref{theoneiacU}.

\begin{lemma}
\label{lemmaesttaylor}
Consider a parametric problem of the form \iref{genpar}
such that {\bf Assumption A} holds 
for a suitable affine representer $(\psi_j)_{j\geq 1}$. 
If the assumptions of Theorem \ref{theoneiacU} are satisfied,
then there exists an $\e>0$ and a $C>0$ such that the estimates
\be
\|t_\nu\|_V\leq C\rho^{-\nu}=C\prod_{j\geq 1}\rho_j^{-\nu_j}, \quad \nu\in \cF,
\label{esttaylor}
\ee
hold for any sequence 
$\rho=(\rho_j)_{j\geq 0}$ of numbers larger than or equal to $1$ for which
\be
\sum_{j\geq 1}(\rho_j-1)\|\psi_j\|_X\leq \e.
\label{constrhoeps}
\ee
\end{lemma}

\noindent
{\bf Proof:} Let $\e>0$ and $C$ be as in Theorem \ref{theoneiacU},
and let $\rho=(\rho_j)_{j\geq 1}$ of numbers larger than or equal to $1$ satisfying the
constraint \iref{constrhoeps}.  
For each $j\geq 1$, let $\cO_{\rho_j}\subset \C$  be the open set that contains 
the disc $\{|z_j|\leq \rho_j\}$ given in Theorem  \ref{theoneiacU}.     Then, we know that  the map $y\mapsto u(y)$ admits an extension $z\mapsto u(z)$
onto  the set $\cO_{\rho}$ defined by \iref{orho}, which is holomorphic in each variable $z_j$ 
with uniform bound 
\be
\sup_{z\in\cO_{\rho}}\|u(z)\|_V\leq C.
\ee
For any given $\nu\in \cF$, we define
\be
\label{defJ}
J:=J(\nu):=\max\{j \; : \; \nu_j\neq 0\}.
\ee
Similar to \iref{finitew}, we introduce the finite dimensional function $w$ defined by
\be
w(z_1,\dots,z_J)=u(T_Jz), \quad T_Jz=(z_1,\dots,z_J,0,0,\dots),
\ee
so that we have for this particular $\nu$,
\be
\partial^\nu u(0)
=
\frac {\partial^{|\nu|}w}{\partial z_1^{\nu_1}\ldots \partial z_J^{\nu_J}} (0,\dots,0).
\ee
We know that $w$ is holomorphic on the set
\be
\cO_{\rho,J}:=\otimes_{1\leq j\leq J}\ \cO{\rho_j},
\ee
which is an open neighborhood of the $J$-dimensional polydisc
\be
\cU_{\rho,J}:=\otimes_{1\leq j\leq J}\ \cU{\rho_j}.
\ee
In addition, we have
\be
\sup_{(z_1,\dots,z_J)\in\cU_{\rho,J}} \|w(z_1,\dots,z_J)\|_V \leq C
\label{wunibound}
\ee
We may thus apply the Cauchy formula \iref{cauchy}
recursively in each variable $z_j$ and obtain for any $(\t z_1,\dots,\t z_J)$
in the interior of $\cU_{\rho,J}$ a representation of $w(\t z_1,\ldots,\t z_J)$ as a multiple integral
\be
w(\t z_1,\ldots,\t z_J)=(2\pi i)^{-J}
\int_{|z_1|=\rho_1} \dots \int_{|z_J|=\rho_J} \frac {w(z_1,\dots,z_J)}{(\t z_1-z_1)\dots(\t z_J-z_J)}
dz_1\dots dz_J.
\ee
By differentiation, this yields
\be
\frac {\partial^{|\nu|}}{\partial z_1^{\nu_1}\dots \partial z_J^{\nu_J}} w(0,\ldots,0)=\nu ! (2\pi i)^{-J}
\int_{|z_1|=\rho_1} \ldots \int_{|z_J|=\rho_J} \frac {w(z_1,\ldots,z_J)}{z_1^{\nu_1+1}\ldots z_J^{\nu_J+1}}
dz_1\ldots dz_J,
\ee
and therefore, 
using \iref{wunibound}, 
we obtain the estimate
\be
\|\partial^\nu u(0)\|_V=\left \|\frac {\partial^{|\nu|}  w}{\partial z_1^{\nu_1}\ldots \partial z_J^{\nu_J}}(0,\ldots,0)\right \|_V 
\leq C\nu! \prod_{j\leq J} \rho_j^{-\nu_j},
\ee
which is equivalent to \iref{esttaylor}. \hfill $\Box$
\nl

Let us comment on the estimate \iref{esttaylor}. Since we
may take any sequence $\rho$ on the right-hand side,  as long as it
satisfies the constraint \iref{constrhoeps}, we have the estimate
\be
\|t_\nu\|_V\leq C \min\Big \{ \rho^{-\nu} \; : \;  \; \sum_{j\geq 1}(\rho_j-1)\|\psi_j\|_X\leq \e\; {\rm and}\; \rho_j\geq 1,\; j\geq 1\Big\}.
\label{esttaylor1}
\ee
It is possible to characterize the sequence $\rho^*$ for which the minimum in the above right-hand side
is attained. An important observation is that this minimizing sequence 
{\it depends on $\nu$}.

To find $\rho^*$, we observe
that for a given $\nu$, this minimization problem is in fact finite dimensional since $\rho^{-\nu}$
is not influenced by the values of $\rho_j$ for those $j$ such that $\nu_j=0$. 
Since $\rho^{-\nu}$ is monotone non-increasing with $\rho_j$ for the other values of $j$,
and in view of the constraint \iref{constrhoeps}, we should thus set
\be
\rho_j^*=1,\quad j\notin {\rm supp}(\nu).
\ee
It remains to solve the finite dimensional problem
\be
\min\Big\{ \prod_{j\in  {\rm supp}(\nu)}\rho_j^{-\nu_j} \; : \; 
\sum_{j\in  {\rm supp}(\nu)}(\rho_j-1)\|\psi_j\|_X\leq \e\; {\rm and}\; \rho_j\geq 1,\; j\in  {\rm supp}(\nu)\Big\},
\ee
or equivalently
\be
\max\Big\{ \sum_{j\in  {\rm supp}(\nu)}\nu_j\log(\rho_j) \; : \; 
\sum_{j\in  {\rm supp}(\nu)}(\rho_j-1)\|\psi_j\|_X\leq \e\; {\rm and}\; \rho_j\geq 1,\; j\in  {\rm supp}(\nu)\Big\},
\ee
which admits a unique solution since we minimize a strictly concave function over a convex set.
The solution necessarily satisfies the equality constraint 
\be
\sum_{j\in  {\rm supp}(\nu)}(\rho_j^*-1)\|\psi_j\|_X= \e.
\label{equality}
\ee
For the optimal solution $\rho^*$, if $E\subset  {\rm supp}(\nu)$ is the subset of those $j\in  {\rm supp}(\nu)$
such that $\rho_j^*>1$, there exists a Lagrange multiplier $\lambda\in \R$ such that
\be
\frac {\nu_j}{\rho_j^*}= \lambda\|\psi_j\|_X, \quad j\in E.
\label{lagrange}
\ee
For any index $\nu=(\nu_j)_{j\geq 1}\in\cF$ and $E\subset \N$, we use the notation
\be
\nu_E:=(\t \nu_j)_{j\geq 1},\quadÊ\t \nu_j=\nu_j \quad {\rm if}\quad  j\in E,\quad  \t\nu_j=0\quad {\rm otherwise}.
\label{nuE}
\ee
Combining \iref{lagrange} and \iref{equality}, we thus find that 
\be
\label{deflambda}
\lambda=\frac {|\nu_E|}{\sigma_E+\e}  \quad {\rm where}\quad \sigma_E:=\sum_{j\in E} \|\psi_j\|_X.
\ee
Therefore, the solution $\rho^*=\rho^*(\nu)=(\rho_j^*)_{j\geq 1}$, has the form
\be
\rho_j^*=\frac {\nu_j (\sigma_E+\e )}{ |\nu_E|\, \|\psi_j\|_X}\quad {\rm if}\; j\in E, \quad \rho_j^*=1\; {\rm  if}\; j\notin E.
\label{solrhoj}
\ee

This characterization is not satisfactory since the set $E$ is not explicitly given.  However,  given any set $E$,
we can define $\lambda$  by \eref{deflambda} and define a coresponding sequence $(\rho_j)$ as in \eref{solrhoj}.   Therefore, the minimum we seek is the same as
\be
\min_{E}\(\frac {|\nu_E|}{\sigma_E+\e}\)^{|\nu_E|} \prod_{j\geq 1} \(\frac {\|\psi_j\|_X}{\nu_j}\)^{\nu_j},
\ee
over all sets $E\subset   {\rm supp}(\nu)$ for which the corresponding
$\rho_j$ given as in \iref{solrhoj} are strictly larger than $1$ for all $j\in E$. The optimal 
set $E$ is the one for which this minimum is reached. This is a combinatorial problem
which is not easy to solve except for those $\nu\in \cF$ of small support. For this reason,
we do not make further use of the above optimal sequence $\rho^*(\nu)$ in   bounding  $\|t_\nu\|_V$.
Instead, we use in \S \ref{secsummability} certain suboptimal choices $\rho(\nu)$ 
which have an explicit expression inspired by \iref{solrhoj}.

\subsection{Refined estimates for elliptic and parabolic PDEs}
\label{secrefined}

The estimate \iref{esttaylor} can be refined in the particular case of
the elliptic and parabolic problems \iref{ellip} and \iref{parab}.
We recall that for each of these problems, the parameter $a$ is taken in
\be
X=L^\infty(D),
\ee
and that the uniform
boundedness and holomorphy of the solution map is ensured 
under a condition of the form $\Re(a)\geq t$ for some $t>0$. 
In such a case, we have seen in \S \ref{secpolydisc} that when the sequence $\rho=(\rho_j)_{j\geq 1}$   fulfills the constraints
\be
\sum_{j\geq 1} \rho_j |\psi_j(x)| \leq \o a(x)-t,\quad x\in D,
\label{constrhox}
\ee
for some $t>0$, the holomorphic extension
is defined over the polydisc $\cU_\rho$ with uniform bound 
\be
\sup_{z\in \cU_\rho}\|u(z)\|_V \leq C_{t}.
\ee
By a recursive application of Cauchy's  formula, as in the proof of Lemma \ref{lemmaesttaylor},
we now obtain, for any fixed $t>0$, the estimate
\be
\|t_\nu\|_V \leq C_t \min \Big\{ \rho^{-\nu} \; : \; \sum_{j\geq 1} \rho_j |\psi_j(x)| \leq \o a(x)-t,\quad x\in D\Big\}.
\label{esttaylor2}
\ee
It is not clear how to give a simple characterization of the above minimization problem, due to the
form of the constraints \iref{constrhox} which need to be fullfilled for every $x\in D$. There are however
two particular instances of affine decompositions where such a simple characterizations exists.

The first of these is when the $\psi_j$ have {\it disjoint supports}, by which we mean that
\be
|{\rm supp}(\psi_i)\cap {\rm supp}(\psi_j)|=0,\quad i\neq j.
\ee
In this case, the uniform ellipticity assumption ${\rm \bf UEA}(r)$ holds if and only if
\be
|\psi_j(x)| \leq \o a(x)-r, \quad x\in D, \quad j\geq 1.
\ee
This instance is sometimes referred to as the model of {\it disjoint inclusions}.
One particular example is when $a$ is piecewise constant over a finite or infinite partition $(D_j)_{j\geq 1}$
of $D$,
in which case $\o a$ is a strictly positive constant and $\psi_j=c_j\Chi_{D_j}$ for some positive numbers $c_j$
each of them smaller than $\o a -r$. 

In the general case of  disjoint supports of the $\psi_j$, the constraint
\iref{constrhox} can be decoupled, so that the minimization problem on the right-hand side
of \iref{esttaylor2} is equivalent to
\be
\min \Big\{ \rho^{-\nu} \; : \;  \rho_j |\psi_j(x)| \leq \o a(x)-t,\quad x\in D, \quad j\geq 1\Big\}.
\label{minesttaylor}
\ee
The optimal solution $\rho^*$ to this problem is obviously given by
\be
\rho_j^*=\inf_{x\in D} \frac {\o a(x)-t}{|\psi_j(x)|}.
\label{rhojdisjoint}
\ee
Let us note that in that case $\rho^*$ does not depend on $\nu$.
This leads us to the estimate
\be
\|t_\nu\|_V \leq C_t \prod_{j\geq 1} \( \sup_{x\in D}\frac {|\psi_j(x)|}{\o a(x)-t}\)^{\nu_j}.
\label{estimdisjoint}
\ee
If ${\rm \bf UEA}(r)$ holds, we see that we can take each $\rho_j^*$ strictly larger than $1$ if
we take $0<t<r$, for example by setting $t=\frac r 2$. In such a case, we have indeed
\be
\frac {\o a(x)-t}{|\psi_j(x)|} \geq \frac {\o a(x)-t}{\o a(x)-r} \geq 1+\frac {r-t}{\o a(x)-r} \geq 1+\frac {r}{2\|\o a\|_{X}},
\label{rhojbound}
\ee
which shows that $\rho_j^*>1$. Note that the values $\rho_j^*$ increase as
$t$ decrease, which in principle results in a better bound for $\|t_\nu\|_V$. However
the constant $C_t$ tends to $+\infty$ as $t\to 0$. One may in principle
search for an optimal value of $t$, however we do not enter this discussion.

The second instance is when the $\psi_j$ are functions of constant moduli,
such as complex exponentials. In this case, the uniform ellipticity assumption ${\rm \bf UEA}(r)$ holds if and only if
\be
\sum_{j\geq 1} \|\psi_j\|_{X}\leq \o a_{\min}-r, \quad \o a_{\min}:=\min_{x\in D} \o a(x),
\ee
and the minimization problem on the right-hand side
of \iref{esttaylor2} is equivalent to
\be
\min \Big\{ \rho^{-\nu} \; : \;  \sum_{j\geq 1}\rho_j \|\psi_j\|_{X} \leq \o a_{\min}-t \Big\}.
\ee
By the same Lagrange multiplier approach which we used above for the
characterization of the minimizer in \iref{esttaylor1}, we find that the above minimum
is attained for $\rho^*=\rho^*(\nu)=(\rho_j^*)_{j\geq 1}$ given by
\be
\rho^*_j=\frac {\nu_j  (\o a_{\min}-t)}{|\nuÊ| \|\psi_j\|_{X}}.
\ee
This leads us to the estimate
\be
\|t_\nu\|_V \leq C_t \(\frac {|\nu|}{ \o a_{\min}-t }\)^{|\nu|} 
\prod_{j\geq 1}\(\frac {\|\psi_j\|_{X}}{ \nu_j }\)^{\nu_j}.
\ee

\subsection{Estimates of Legendre coefficients}
\label{subestleg}

Returning to general parametric PDEs of the form \iref{genpar}
with an affine representation \iref{affine}, our next objective
is to establish similar estimates for the Legendre coefficients
$\|v_\nu\|_V$ and $\|w_\nu\|_V$. Let us recall that these
coefficients are given by
\be
v_\nu=\int_U u(y)L_\nu(y)d\mu(y),
\ee
and
\be
w_\nu=\prod_{j\geq 1} (2\nu_j+1) \int_U u(y)P_\nu(y)d\mu(y),
\label{wnu}
\ee
They are linked by the relation
\be
w_\nu=\(\prod_{j\geq 1} (1+2\nu_j)\)^{1/2} v_\nu.
\label{propvw}
\ee
We introduce the function
\be
t\mapsto \theta(t):=\frac{\pi t}{2(t-1)},
\ee
which is monotone non-increasing over $]1,+\infty[$.

The following result establishes estimates on the Legendre coefficients, based
on the bounded holomorphic extension of $u$  onto
the polyellipses $\cH_\rho$ and their neighborhood $\cO_\rho$ established in Theorem \ref{theoneiaU}.

\begin{lemma}
\label{lemmaestlegendre}
Consider a parametric problem of the form \iref{genpar}
such that {\bf Assumption A} holds 
for a suitable affine representer $(\psi_j)_{j\geq 1}$. 
If the assumptions of Theorem \ref{theoneiaU} are satisfied,
then there exists $\e>0$ and $C>0$ such that the estimates
\be
\|v_\nu\|_V\leq C\prod_{j\in {\rm supp}(\nu)} \theta(\rho_j)(1+2\nu_j)^{1/2}\rho_j^{-\nu_j},
\label{estlegendrew}
\ee
and
\be
\|w_\nu\|_V\leq C\prod_{j\in {\rm supp}(\nu)} \theta(\rho_j)(1+2\nu_j)\rho_j^{-\nu_j},
\label{estlegendrev}
\ee
hold for any sequence  $\rho=(\rho_j)_{j\geq 0}$ of numbers strictly larger than $1$,
which satisfies the constraint \iref{constrhoeps}.
\end{lemma}

\noindent
{\bf Proof:}  We only need to prove \iref{estlegendrew}, since \iref{estlegendrev} then 
follows by \iref{propvw}. Let $\e>0$ and $C$ be as in Theorem \ref{theoneiaU},
and let $\rho=(\rho_j)_{j\geq 1}$ be a sequence of numbers strictly larger than $1$, which satisfies the
constraint \iref{constrhoeps}.
We know that for each  $j\geq 1$ there exists an open set $\cO_{\rho_j}\subset \C$ that contains 
the filled-in ellipse $\cH_{\rho_j}$ and such that 
the map $y\mapsto u(y)$ admits an extension $z\mapsto u(z)$
over the set $\cO_{\rho}$ defined by \iref{orho}, which is holomorphic in each variable $z_j$ 
with uniform bound 
\be
\sup_{z\in\cO_{\rho}}\|u(z)\|_V\leq C.
\label{uniorho}
\ee
We observe that $U\subset \cO_\rho$.

In the case $\nu=0$, the estimate \iref{estlegendrew}
is immediate since
\be
\|w_0\|_V=\Big\|\int_U u(y)d\mu(y)\Big\|_V 
\leq 
\sup_{y\in U}\|u(y)\|_V 
\leq 
C,
\ee
where we have used the fact that $\mu$ is a probability measure. We now assume that $\nu\neq 0$. 
Up to a reordering of $(\psi_j)_{j\geq 1}$, 
we may assume without loss of generality 
that  $\nu_j\neq 0$ for $j\leq J$ and $\nu_j=0$ for $j>J$ for $J=|{\rm supp}(\nu)|\geq 1$.
We partition the variable $y$ into
\be
y=(y_1,\ldots,y_J,y'),\;\;  y':=(y_{J+1},y_{J+2},\dots)\in [-1,1]^\N=U,
\ee
and rewrite \iref{wnu} as
\be
w_\nu=\prod_{j=1}^J (2\nu_j+1) \int_U v(y') d\mu(y'), 
\label{defunu1}
\ee
where
\be
v(y'):=\int_{[-1,1]^J} u(y_1,\dots,y_J,y') 
\left(\prod_{j=1}^J P_{\nu_j}(y_j)\right) \frac {dy_1}2\dots \frac{dy_J}2.
\ee
For a fixed $y'\in U$, we use the holomorphy 
of the finite dimensional map $(z_1,\dots,z_J)\mapsto u(z_1,\dots,z_J,y')$ 
in order to evaluate $\|v(y')\|_V$.
For this purpose, we introduce for any integer $n\geq 1$ 
the following function of a single complex variable $z$ 
\beqn
\label{defQ}
Q_n(z):=  \int_{-1}^1 \frac {P_n(y)}{z-y} dy,
\eeqn
and the corresponding multivariate functions
 \beqn
 \label{defQ1}
Q_\nu(z_1,\dots,z_J):=  \prod_{j=1}^JQ_{\nu_j}(z_j),  
\eeqn
which are well defined as long as $|z_j|>1$ for $j=1,\dots,J$.
For our given $\rho$, we introduce the $J$-dimensional polyellipse
\be
\cE_{\rho,J}:=\displaystyle{\otimes_{1\leq j\leq J} \cE_{\rho_j}}.
\ee
Since   $\rho_j>1$, $1\leq j\leq J$,  the unit interval 
$[-1,1]$ is contained in the interior of each filled-in ellipse $\cH_{\rho_j}$.
Therefore, we may recursively 
apply Cauchy's integral formula on each ellipse $\cE_{\rho_j}$ for each of the variables $z_j$ , $j=1,...,J$ , and obtain 
\be
\label{CF}
u(y_1,\ldots,y_J,y')=\frac{1}{(2\pi i)^J}\int_{\cE_{\rho_1}}\cdots\int_{\cE_{\rho_J}}\frac{u(z_1,\ldots,z_J,y')}{(y_1-z_1)\ldots (y_J-z_J)} dz_1
\ldots dz_J,
\ee
for any $(y_1,\ldots,y_J)\in [-1,1]^J$ and any $y'\in U$.
Multiplying by $\prod_{j=1}^J P_{\nu_j}(y_j)$ and integrating over $[-1,1]^J$ with respect
to $\frac {dy_1} 2 \ldots \frac{ dy_J} 2$, we therefore obtain
\be
v(y')=\(\frac {i}{4\pi}\)^J\int_{\cE_{\rho_1}}\cdots\int_{\cE_{\rho_J}}u(z_1,\ldots,z_J,y')Q_\nu(z_1,\dots,z_J)dz_1\ldots dz_J.
\ee
From the uniform bound \iref{uniorho} we know that
\be
(z_1,\ldots,z_J)\in \cE_{\rho,J} \mbox{ and } y'\in U \Rightarrow (z_1,\ldots,z_J,y')\in \cO_\rho
\Rightarrow \|u(z_1,\ldots,z_J,y')\|_V\leq C.
\ee
Injecting this bound in the above integral yields
\be
\|v(y')\|_V \leq C\(\prod_{j=1}^J \frac {\rho_j}2\)  \max_{(z_1,\dots,z_J)\in \cE_{\rho,J}}|Q_\nu(z_1,\dots,z_{J})|, \quad y'\in U,
\ee
where we have used the fact the perimeter of $\cE_{\rho_j}$ has length smaller than $2\pi \rho_j$.
We now use the following estimate (see page 313 of \cite{Dav})
\be
\max_{z\in \cE_t} |Q_n(z)|\leq \frac {\pi\ t^{-n}}{t-1},
\ee
which yields 
\be
\max_{(z_1,\dots,z_J)\in \cE_{\rho,J}}|Q_\nu(z_1,\dots,z_{{J}})| 
\leq 
\prod_{j=1}^J \frac {\pi\ \rho_j^{-\nu_j}}{\rho_j-1},
\ee
and therefore
\be
\|v(y')\|_V \leq C\prod_{j=1}^J \theta(\rho_j)\rho_j^{-\nu_j}, \quad y'\in U.
\ee
Combining this estimate with \iref{defunu1}, we obtain 
\be
\|w_\nu\|_V \leq \prod_{j=1}^J (1+2\nu_j) \sup_{y'\in U} \|v(y')\|_V 
\leq C\prod_{j\in {\rm supp}(\nu)} \theta(\rho_j)(1+2\nu_j)\rho_j^{-\nu_j},
\ee
which is \iref{estlegendrew}. \hfill $\Box$
\nl

The estimates \iref{estlegendrev} and \iref{estlegendrew} 
are very similar to the estimate \iref{esttaylor}Ê obtained 
in Lemma \ref{lemmaesttaylor}Ê for the Taylor coefficients,
however with two noticable differences:
\begin{itemize}
\item
On the one hand, the estimates for the Legendre coefficients
are a bit more pessimistic,Ê due to the presence of the additional factors $\theta(\rho_j)$ and $(1+2\nu_j)$.
Intuitively, these factors are absorbed by the decay of the factor $\rho_j^{-\nu_j}$
when $\rho_j$ or $\nu_j$ become large. The analysis in the next section
confirms that they do not affect the $\ell^p$ summability properties
of the estimate.
\item
On the other hand, these estimates are obtained under much weaker
conditions than those of Lemma \ref{lemmaesttaylor}. Indeed 
Theorem \ref{theoneiaU} only requires 
the existence of a holomorphic extension of the solution map $a\mapsto u(a)$
in a neigborhood of $a(U)$, in contrast to Theorem \ref{theoneiacU}
which requires a neighborhood of $a(\cU)$. In particular, for problems such
as \iref{eqpull} or \iref{nonlin}, we know that the conditions of Theorem \ref{theoneiaU}
are met but not those of Theorem \ref{theoneiacU}.
\end{itemize}

Similar to the estimate for Taylor coefficients, we can 
use the fact that \iref{estlegendrev} and \iref{estlegendrew} 
hold for any sequence $\rho$ satisfying the prescribed constraints,
in order to obtain the estimates
\be
\|v_\nu\|_V\leq C\inf\Big\{\prod_{j\in {\rm supp}(\nu)} \theta(\rho_j)(1+2\nu_j)^{1/2}\rho_j^{-\nu_j}\Big\},
\label{estlegendrev1}
\ee
and 
\be
\|w_\nu\|_V\leq C\inf\Big\{\prod_{j\in {\rm supp}(\nu)} \theta(\rho_j)(1+2\nu_j)\rho_j^{-\nu_j}\Big \},
\label{estlegendrew1}
\ee
where the infima are taken over all sequences $\rho$ of numbers strictly larger than $1$, such that
$\sum_{j\geq 1}Ê(\rho_j-1)\|\psi_j\|_X\leq \e$.

\begin{remark}
The values of $\rho_j$ enter the above estimates only for $j\in {\rm supp}(\nu)$.
This implies that we can consider the above infimas over all sequences 
$\rho$ of numbers larger or equal to $1$ with
$\rho_j>1$ if $j\in {\rm supp}(\nu)$ and such that
$\sum_{j\in {\rm supp}(\nu)}Ê(\rho_j-1)\|\psi_j\|_X\leq \e$,
which amounts in taking $\rho_j=1$ if $j\notin {\rm supp}(\nu)$.
\label{remconstrho}
\end{remark}

\subsection{Summability of multi-indexed sequences}
\label{secsummability}

We want to use the upper estimates obtained for
$\|t_\nu\|_V$, $\|v_\nu\|_V$ and $\|w_\nu\|_V$ derived in the previous sections
in order to prove Theorem \ref{theomain}.
As a preliminary step, we establish in this section several results concerning
the $\ell^p$ summability of certain type of multi-indexed sequences,
which appear in the proof of Theorem \ref{theomain} that follows.

 We begin by considering sequences of the form $(b^\nu)_{\nu\in\cF}$ where
$b=(b_j)_{j\geq 1}$ is a given sequence of positive numbers. For such sequences
we have the following elementary result.

\begin{lemma}
For any $0<p<\infty$, the sequence $(b^\nu)_{\nu\in\cF}$
belongs to $\ell^p(\cF)$ if and only if $b\in \ell^p(\N)$ and $\|b\|_{\ell^\infty}<1$. Moreover
\be
\|(b^\nu)_{\nu\in\cF} \|_{\ell^p} \leq {\rm exp} \(c_p \frac {\|b\|_{\ell^p}^p}p\), \quad c_p:=\frac 1 {1-\|b\|_{\ell^\infty}^p}.
\label{bnubound}
\ee
\label{lemmabnu}
\end{lemma}

\noindent
{\bf Proof:}  For any positive integer $J$, let  $\cF_J$ denote the set of those  $\nu\in \cF$ such that ${\rm supp}(\nu)\subset \{1,\dots, J\}$.  Now, if  $\|b\|_{\ell^\infty}<1$,
we can write  
\be
\sum_{\nu\in\cF_J } b^{p\nu}=\prod_{1\leq j\leq J} \sum_{n\geq 0} b_j^{pn}=\prod_{1\leq j\leq J} \frac 1 {1-b_j^p}, \quad J=1,2,\dots.
\ee
If    $b\in \ell^p(\N)$, we can let $J$ tend to $+\infty$ and obtain
\be
\sum_{\nu\in\cF } b^{p\nu}=\prod_{j\geq 1} \frac 1 {1-b_j^p}<\infty.
\ee
This proves the one implication of the theorem.
Since,
\be
\prod_{j\geq 1} \frac 1 {1-b_j^p} =\prod_{j\geq 1}\(1+ \frac {b_j^p} {1-b_j^p}\) \leq \prod_{j\geq 1}{\rm exp}\( \frac {b_j^p} {1-b_j^p}\)\leq \prod_{j\geq 1}{\rm exp}(c_p b_j^p)={\rm exp}(c_p \|b\|_{\ell^p}^p).
\ee
we also have  the bound \iref{bnubound}. 

For the other implication, we observe that
the sequences $(b_j)_{j\geq 1}$ and $(b_j^n)_{n\geq 0}$ for any $j\geq 1$,
are subsequences of $(b^\nu)_{\nu\in\cF}$ corresponding to particular selections
of indices $\nu$. This shows that the $\ell^p$ summability of $(b^\nu)_{\nu\in\cF}$
implies both that $b\in \ell^p(\N)$ and $\|b\|_{\ell^\infty}<1$. \hfill $\Box$
\nl

One immediate application of the above lemma
concerns the $\ell^p$ summability of the Taylor coefficients
for the elliptic or parabolic problems in the model of disjoint inclusions
discussed in \S \ref{secrefined}. In this case, the estimate
\iref{estimdisjoint}Ê has the form
\be
\|t_\nu\|_V \leq C_tb^{\nu},Ê\quad {\rm where} \ b=(b_j)_{j\ge 1} \ {\rm with} \ b_j:=\sup_{x\in D}\frac {|\psi_j(x)|}{\o a(x)-t}
\ee
Working under ${\rm \bf UEA}(r)$ and taking $t=\frac r 2$, we know from \iref{rhojbound} that for $X:=L^\infty(D)$,
\be
\|b\|_{\ell^\infty}Ê\leq \frac {2\|\o a\|_{X}}{ 2\|\o a\|_{X}+r} <1.
\ee
  Since in addition 
\be
b_j \leq \frac {2\|\psi_j\|_{X}}{r},
\ee
this shows that $(\|\psi_j\|_{X})_{j\geq 1} \in \ell^p(\N)$ implies $b\in \ell^p(\N)$.
Combining these observations with Lemma \ref{lemmabnu}, we thus find that
if ${\rm \bf UEA}(r)$ holds and if $(\|\psi_j\|_{X})_{j\geq 1} \in \ell^p(\N)$, then
the sequence $(\|t_\nu\|_V)_{\nu\in\cF}$
belongs to $\ell^p(\cF)$, which is a particular case of Theorem \ref{theomain}.

\begin{remark}
We have mentioned in Remark \ref{wlp} that the convergence 
rate $n^{-s}$ of best $n$-term approximation in $\ell^q$ spaces 
is equivalent to the property of weak $\ell^p$ summability with $s=\frac 1 p-\frac 1 q$. 
Therefore, a  relevant question is whether the above Lemma \ref{lemmabnu}
is valid with $\ell^p$ replaced by $w\ell^p$. Surprisingly, the answer
is negative, and closely related to classical results in number theory. 
Indeed, fix any   $0<p<1$  and consider the prototype sequence $b\in w\ell^p(\N)$ given by
\be
b_j=(j+1)^{-1/p}.
\ee
This sequence also satisfies $\|b\|_{\ell^\infty}<1$.  If we were to have 
$(b^\nu)_{\nu\in\cF}\in w\ell^p(\cF)$  then there would be a constant 
$C$ such that for any $\eta>0$, we have
\be
\#\{\nu \in\cF \;  : \; b^\nu\geq \eta\} \leq C\eta^{-p},
\ee
or equivalently, such that for any $A\geq 2$, 
\be
t(A):=\#\Big\{\nu \in\cF \; : \; \prod_{j\geq 2} j^{\nu_j} \leq A\Big \} \leq CA.
\label{mA}
\ee
The left side can be rewritten as
\be
t(A)=\sum_{n=2}^{\lfloor A\rfloor} f(n),
\ee
where $f(n)$ is the number of possible multiplicative partitions of $n$.
The problem of counting multiplicative partitions of natural numbers,
sometimes refered to as {\em factorisatio numerorum}, has been extensively studied in
number theory, see in particular \cite{CEP} which gives a sharp
asymptotic bound for $f(n)$. In \cite{LMS}, it is proved that the total number
of multiplicative partitions $t(A)$ has the asymptotic behaviour
\be
\frac {t(A)}{A}\sim {\rm exp}\left \{\frac {4\sqrt{\log(A)}}{\sqrt{2e}\log(\log(A))}(1+o(1))\right \}\to +\infty
\ee
as $A\to +\infty$.  This shows   that \iref{mA} does not hold, and
thus that $(b^\nu)_{\nu\in\cF}$ does not belong to $w\ell^p(\cF)$.
\end{remark}

We make further use of a slightly more general version of 
Lemma \ref{lemmabnu} where we incorporate additional algebraic factors
into the sequence $b^\nu$.

\begin{lemma}
\label{lemmabnualg}
For a given sequence $b=(b_j)_{j\geq 1}$ of positive numbers,
and for non-negative numbers $c$ and $r$, let $(b_\nu)_{\nu\in\cF}$ be defined by
\be
b_\nu:=b^\nu\prod_{j\geq 1} (1+c\nu_j^r)=\prod_{j\geq 1} (1+c\nu_j^r) b_j^{\nu_j}.
\ee
For any $0<p<\infty$, this sequence
belongs to $\ell^p(\cF)$ if and only if $b\in \ell^p(\N)$ and $\|b\|_{\ell^\infty}<1$. 
\end{lemma}

\noindent
{\bf Proof:} Since $b_\nu\geq b^\nu$, the ``only if'' part follows
from Lemma \ref{lemmabnu} and therefore we only need to prove the if part.
With $\cF_J$ as in the proof of Lemma \ref{lemmabnu}, we write
\be
\label{prod}
\sum_{\nu\in\cF_J } b_\nu^{p}=\prod_{1\leq j\leq J} \sum_{n\geq 0} (1+cn^r)^pb_j^{pn},
\ee
Since $\|b\|_{\ell^\infty} \leq 1$ we find that
\be
\sum_{n\geq 0} (1+cn^r)^p b_j^{pn} \leq 1+Cb_j^p,
\ee
where the constant $C$ depends on $c$, $r$, $p$ and $\|b\|_{\ell^\infty}$.
Since $b\in \ell^p(\N)$, this shows that the   product on the right side of \eref{prod} converges as $J \to \infty$.
Therefore $(b_\nu)_{\nu\in\cF}\in \ell^p(\cF)$. \hfill $\Box$
\nl

The estimates obtained for $\|t_\nu\|_V$, $\|v_\nu\|_V$ and $\|w_\nu\|_V$ also 
involve quantities of the form
\be
\frac {|\nu|^{|\nu|}}{\prod_{j\geq 1} \nu_j^{\nu_j}}Êd^\nu,
\ee
for sequences $d=(d_j)_{j\geq 1}$ of positive numbers. In view of the Stirling inequalities
\be
n! \leq  n^n \leq  n! e^n ,
\label{stir1}
\ee
we may write 
\be
\frac {|\nu|^{|\nu|}}{\prod_{j\geq 1} \nu_j^{\nu_j}}Êd^\nuÊ\leq e^{|\nu|} \frac {|\nu|!}{\nu !}Êd^\nu= \frac {|\nu|!}{\nu !}Êb^\nu,
\label{stir2}
\ee
where
\be
\label{defb}
b=(b_j)_{j\geq 0}, \quad b_j=ed_j.
\ee
This suggest studying the $\ell^p$ summability
of sequences of the form $\(\frac {|\nu|!}{\nu !}Êb^\nu\)_{\nu\in\cF}$.
Due to the presence that the multinomial factor $\frac {|\nu|!}{\nu !}$ which 
can be much larger than $1$, we expect that the conditions for  $\ell^p$ summability 
are more stringent than for the sequence $(b^\nu)_{\nu\in \cF}$. This is confirmed
by the following result.

\begin{lemma}
\label{lemmabnumulti}
For any $0<p<1$, a  sequence  $\(\frac {|\nu|!}{\nu !}Êb^\nu\)_{\nu\in\cF}$
belongs to $\ell^p(\cF)$ if and only if $b\in \ell^p(\N)$ and $\|b\|_{\ell^1}<1$.
\end{lemma}

\noindent
{\bf Proof:} We first observe that whenever $b\in \ell^1(\N)$,  
  the multinomial formula gives  
\be
\sum_{|\nu|=k}\frac {|\nu|!}{\nu !}Êb^{\nu}=\(\sum_{j\geq 1} b_j\)^k.
\ee
Summing over $k$  we see that $\(\frac {|\nu|!}{\nu !}Êb^\nu\)_{\nu\in\cF}$ is in $\ell^1(\cF)$ if and only if $b\in \ell^1(\N)$ and $\|b\|_{\ell^1(\N)}<1$.  Moreover,
\be
\label{moreover}
\Big \|\(\frac {|\nu|!}{\nu !}Êb^\nu\)_{\nu\in\cF} \Big\|_{\ell^1(\cF)}= \sum_{\nu\in \cF}\frac {|\nu|!}{\nu !}Êb^{\nu}=\frac 1 {1-\|b\|_{\ell^1}},
\ee

Now suppose that $\(\frac {|\nu|!}{\nu !}Êb^\nu\)_{\nu\in\cF}\in\ell^p(\cF)$ for some $p\le 1$.  Then, $b$ is  in $\ell^p(\N)$ since it is a subsequence of $\bar b$ corresponding to a particular selection
of indices $\nu$.  Also $\bar b$ is in $\ell^1(\cF)$ so   $b$ must be in  $\ell^1(\cN)$ with norm smaller than one.
  
Conversely, assume that $b\in \ell^p(\N)$ and $\|b\|_{\ell^1}<1$. We claim that there exists
two positive sequences $c=(c_j)_{j\geq 1}$ and $d=(d_j)_{j\geq 1}$ with the following
properties:
\begin{enumerate}
\item[(i)]
$b_j=c_jd_j$ for all $j\geq 1$.
\item[(ii)]
$c\in \ell^1(\N)$ with $\|c\|_{\ell^1}<1$.
\item[(iii)]
$d\in \ell^q(\N)$ with $\frac 1 q=\frac 1 p-1$, or equivalently $q=\frac {p}{1-p}$, and $\|d\|_{\ell^\infty}<1$.
\end{enumerate}
Before proving this claim, let us show that it implies the $\ell^p$ summability of $\(\frac {|\nu|!}{\nu !}Êb^\nu\)_{\nu\in\cF}$.
Indeed, from H\"older's inequality, we have
$$
\begin{disarray}{ll}
\sum_{\nu\in\cF} \(\frac {|\nu|!}{\nu !}Êb^\nu\)^p & =\sum_{\nu\in\cF} \(\frac {|\nu|!}{\nu !}Êc^\nu\)^p d^{p\nu} \\
& \leq \(\sum_{\nu\in\cF} \frac {|\nu|!}{\nu !}Êc^\nu\)^p \(\sum_{\nu\in\cF} d^{q\nu} \)^{1-p}.
\end{disarray}
$$
 As observed previously in \eref{moreover}, the first factor
is finite due to the fact that $\|c\|_{\ell^1}<1$. The second factor is finite by application of Lemma \ref{lemmabnu}.

It remains to prove the claim by constructing specific sequences $c$ and $d$ having the prescribed
properties. With $\delta:=1-\|b\|_{\ell^1}>0$, we define
\be
\eta:=\frac \delta 3,
\ee
and take $J$ large enough   that
\be
\sum_{j>J} b_j^p\leq \frac {\delta} 3.
\ee
We then define $c$ and $d$ by
\be
c_j=(1+\eta)b_j\;\; {\rm and}\;\; d_j=\frac 1 {1+\eta},\quad j \leq J,
\ee
and 
\be
c_j=b_j^p\;\;Ê{\rm and} \;\; d_j= b_j^{1-p}, \quad j>J.
\ee
By construction, we have $c_jd_j=b_j$ for all $j\geq 1$.
For the sequence $c$, we have
\be
\|c\|_{\ell^1} \leq (1+\eta)\|b\|_{\ell^1}+\sum_{j>J} b_j^p \leq \(1+\frac \delta 3\)(1-\delta) +\frac \delta 3\leq 1-\frac \delta 3,
\ee
We next bound $\|d\|_{\ell^\infty}$.  For $1\le j\le J$, we have $d_j=\frac 1 {1+\eta}<1$ and for $j>J$, we have
\be
d_j=\(b_j^p\)^{\frac {1-p}{p}} \leq \(\frac {\delta} 3\)^{\frac {1-p}p}<1.
\ee
Therefore, we have $\|d\|_{\ell^\infty} <1$.   Finally, since $d_j^q=b_j^p$ for $j>J$, we find that $d\in \ell^q(\N)$,
which completes the confirmation of the claim. \hfill $\Box$
\nl

Similar to Lemma \ref{lemmabnualg}, the following result shows
that $\ell^p$ summability is maintained if
we incorporate additional algebraic factors.

\begin{lemma}
\label{lemmabnumultalg}
For a given sequence $b=(b_j)_{j\geq 1}$ of positive numbers,
and for non-negative numbers $c$ and $r$, let $(b_\nu)_{\nu\in\cF}$ be defined by
\be
b_\nu:=\frac {|\nu|!}{\nu !} b^\nu\prod_{j\geq 1} (1+c\nu_j^r).
\ee
For any $0<p<1$, this sequence
belongs to $\ell^p(\cF)$ if and only if $b\in \ell^p(\N)$ and $\|b\|_{\ell^1}<1$. 
\end{lemma}

\noindent
{\bf Proof:} Since $b_\nu\geq \frac {|\nu|!}{\nu !} b^\nu$, the ``only if'' part follows
from Lemma \ref{lemmabnumulti} and we only need to prove the if part.

Using the same sequences $c$ and $d$ as in the proof of Lemma \ref{lemmabnumulti},
and introducing
\be
d_\nu=d^\nu \prod_{j\geq 1} (1+c\nu_j^r),
\ee
we write
$$
\begin{disarray}{ll}
\sum_{\nu\in\cF} b_\nu^p & =\sum_{\nu\in\cF} \(\frac {|\nu|!}{\nu !}Êc^\nu\)^p d_\nu^{p} \\
& \leq \(\sum_{\nu\in\cF} \frac {|\nu|!}{\nu !}Êc^\nu\)^p \(\sum_{\nu\in\cF} d_\nu^q \)^{1-p},
\end{disarray}
$$
and conclude in a similar manner that both factors are finite, using Lemma \ref{lemmabnualg}
for the second factor. \hfill $\Box$

\subsection{Proof of Theorem \ref{theomain}}
\label{subproof}

In order to prove Theorem \ref{theomain}, we use the estimates
\iref{esttaylor}, \iref{estlegendrev} and \iref{estlegendrew} 
for the $\|t_\nu\|_V$, $\|v_\nu\|_V$ and $\|w_\nu\|_V$, respectively.
The right-side of these estimates has a general form $C r(\nu,\rho)$
for any sequence  $\rho$ of numbers larger than $1$ that satisfy the constraint \iref{constrhoeps}.
Our objective is to build for each $\nu$ such a sequence $\rho=\rho(\nu)$,
and show that, for $0<p<1$ the resulting quantities
\be
r_\nu:=r(\nu,\rho(\nu)),
\ee
are $\ell^p$ summable provided that $(\|\psi_j\|_X)_{j\geq 1}\in \ell^p(\N)$.
Obviously, it is sufficient to treat the case when
\be
r(\nu,\rho):=\prod_{j\in {\rm supp}(\nu)} \theta(\rho_j)(1+2\nu_j) \rho_j^{-\nu_j},
\label{rnurho}
\ee
which appears in the right of \iref{estlegendrew}, since it is the largest estimate.

We fix an arbitrary $\nu\in\cF$ and describe our choice for the sequence $\rho$ 
that we insert into the above expression. In what follows, we use the notation
\be
b=(b_j)_{j\ge 1},\quad {\rm where} \quad b_j:=\|\psi_j\|_X,\quad j\ge 1.
\ee
For $J\geq 1$ to be fixed further, we split 
$\N$ into 
\be
E:=\{1,\dots,J\}\;\; {\rm and} \;\; F:=\{J+1,J+2,\dots\},
\ee
and use the notation $\nu_E=(\nu_1,\dots,\nu_J)\in \N^J$ and $\nu_F=(\nu_{J+1},\nu_{J+2},\dots)\in\cF$.
In view of Remark \ref{remconstrho},
we may take 
\be
\rho_j=1, \quad j\notin {\rm supp}(\nu).
\ee
With $\e$ the right side of the constraint \iref{constrhoeps}, we then take
\be
\rho_j=\kappa:=1+\frac {\e} {2\|b\|_{\ell^1}},\quad j\in E\cap {\rm supp}(\nu),
\label{defrhojE}
\ee
and
\be
\rho_j=\kappa+\frac {\e \nu_j}{2b_j |\nu_F|}, \quad j\in F \cap {\rm supp}(\nu).
\label{defrhojF}
\ee
Therefore $\rho_j>1$ when $j\in {\rm supp}(\nu)$, and in addition
\be
\sum_{j\geq 1}(\rho_j-1)b_j \leq \frac {\e \sum_{j\leq J} b_j}{2\|b\|_{\ell^1}}+Ê
\sum_{j>J} \(\frac {\e b_j} {2\|b\|_{\ell^1}}+\frac {\e \nu_j}{2 |\nu_F|}\)
\leq \e.
\ee
which shows that the constraint \iref{constrhoeps} is satisfied.

When using this choice for the sequence $\rho$,  the resulting estimate
may be written
\be
r_\nu=r(\nu,\rho(\nu))=r_E(\nu)r_F(\nu), 
\ee
where
\be
r_E(\nu):=\theta(\kappa)^J\prod_{j\in E}(1+2\nu_j)Ê\kappa^{-\nu_j}\quad {\rm and}\quad r_F(\nu):=\prod_{j\in F\cap{\rm supp}(\nu)}\theta(\rho_j)(1+2\nu_j)\rho_j^{-\nu_j}.
\ee
Denoting by $\cF_E$ and $\cF_F$ the multi-indices in $\cF$ supported on $E$ and $F$, respectively, we
may then write
\be
\sum_{\nu\in \cF} r_\nu^p =\Sigma_E \Sigma_F,
\ee
where
\be
\Sigma_E:=\sum_{\nu\in \cF_E} r_E(\nu)^p\quad {\rm and} \quad \Sigma_F:=\sum_{\nu\in \cF_F} r_F(\nu)^p,
\ee
provided that both sums converge.

The first sum $\Sigma_E$ is estimated by
$$
\begin{disarray}{ll}
\Sigma_E & =\theta(\kappa)^{pJ}\sum_{\nu\in \N^J} \prod_{j=1}^J(1+2\nu_j)^pÊ\kappa^{-p\nu_j}\\
&=\theta(\kappa)^{pJ}\(\sum_{n\geq 0} (1+2n)^p\kappa^{-pn}\)^J <\infty,
\end{disarray}
$$
For the second sum $\Sigma_F$, we first notice that for each $\nu\in \cF_F$,
\be
r_F(\nu)\leq \prod_{j\in F\cap{\rm supp}(\nu)}\theta(\kappa)(1+2\nu_j)\(\frac {\e \nu_j}{2b_j |\nu_F|}\)^{-\nu_j},
\ee
where we have used the fact that $\theta(\kappa)=\max_{t\geq \kappa} \theta(t)\geq \theta(\rho_j)$ for $j\in F$.
Therefore, with $c:=3\theta(\kappa)$, we find that
$$
\begin{disarray}{ll}
r_F(\nu) 
& \leq |\nu_F|^{|\nu_F|}
\prod_{j\in F}\frac {(1+c\nu_j)\(\frac {2 b_j} \e\)^{\nu_j}}  {\nu_j^{\nu_j}} \\
&  \leq \frac {|\nu_F|!}{\nu_F !}
\prod_{j\in F}(1+c\nu_j)\(\frac {2e b_j} \e\)^{\nu_j},
\end{disarray}
$$
where we have used \iref{stir1}. Introducing the sequence $d=(d_j)_{j\geq 1}$ defined by
\be
d_j=\frac {2e b_{j+J}} \e,
\ee
we thus find that 
\be
 \Sigma_F\leq \sum_{\nu\in\cF}Êd_\nu^p \quad {\rm where } \quad d_\nu:=\frac {|\nu |!}{\nu!}
d^\nu \prod_{j\geq 1}(1+c\nu_j).
\ee
We now choose $J$ sufficiently large so that
\be
\|d\|_{\ell^1} =\sum_{j>J}\frac {2e b_{j}} \e <1.
\ee
Since our assumption   $b\in \ell^p(\N)$ implies
that $d\in \ell^p(\N)$, we may   apply Lemma \ref{lemmabnumultalg}
to conclude that $\Sigma_F$ is finite. The proof of Theorem \ref{theomain}Ê is complete. 

\begin{remark}
One defect in the proof the Theorem \ref{theomain} is that,
while it establishes the $\ell^p$ summability of the
sequences $(\|t_\nu\|_V)_{\nu\in\cF}$, 
$(\|v_\nu\|_V)_{\nu\in\cF}$ and $(\|w_\nu\|_V)_{\nu\in\cF}$,
it does not provide us with a simple bound of the $\ell^p$ norms of
these sequences in terms of the $\ell^p$ norm of the 
sequence $(\|\psi_j\|_X)_{j\geq 1}$.
\end{remark}

\subsection{Approximation using downward closed sets}
\label{sublower}

Theorem \ref{theomain} has implications on 
the rate convergence of polynomial approximations obtained by
retaining the terms 
 in Taylor and Legendre series corresponding the $n$ largest coefficients.
Corollaries \ref{corratetay} and \ref{corrateleg}
show that these approximations converge 
with the rates $n^{-s}$,
where $s=\frac 1 p-1$ for uniform convergence
and $s=\frac 1p -Ê\frac 1 2$ for convergence in $L^2(U,V,\mu)$.

These results should be viewed as a theoretical
justification that reduced modeling methods based on
polynomial approximations may perform well for 
parametric PDEs which satisfy the assumptions of Theorem \ref{theomain}.
However, they constitute, by   no means,  a numerical algorithm
since finding the optimal sets $(\Lambda_n)_{n\geq 1}$ are, in
practice,  out of reach, and so is the exact computation of 
the Taylor and Legendre coefficients.

Practical algorithms for the computation of 
polynomial approximations are discussed later in this paper in 
\S 5 and \S 6. The implementation and analysis
of the algorithms presented there benefit from 
imposing additional structure on  the index sets $\Lambda_n$   used to define
the polynomial approximation.  To define this structure, we first recall that $\cF$ has a partial ordering:
for $\nu,\t \nu\in\cF$, we write
$\t \nu \leq \nu$ if and only if $\t \nu_j\leq \nu_j$ for all $j\geq 1$.
We also write $\t \nu<\nu$ if and only if $\t\nu\leq \nu$ and 
$\t\nu_j< \nu_j$ for at least one value of $j$.

\begin{definition}
A set $\Lambda\subset \cF$ is called {\em downward closed} or a {\em lower set}
if and only if 
\be
\nu\in \Lambda \;\; {\rm and} \;\; \t \nu \leq \nu  \ {\rm implies} \  \t \nu\in \Lambda.
\ee
\end{definition}

When considering polynomial spaces 
\be
\P_\Lambda:={\rm span} \{y \mapsto y^\nu \; : \; \nu\in \Lambda\},
\ee
it is quite natural to make the assumption that $\Lambda$ is a downward closed set. 
In particular, this assumption allows us to describe $\P_\Lambda$ in
terms of any tensorized polynomial basis of the form 
\be
\phi_\nu(y)=\prod_{j\geq 1} \phi_{\nu_j}(y_j),
\ee
where $(\phi_k)_{k\geq 0}$ is any family of univariate polynomials such that
$\phi_0=1$ and $\phi_k$ has degree exactly $k$. This includes in particular
the tensorized Legendre polynomials $L_\nu$. By expressing each monomial
$y\mapsto y^k$ as a linear combination of the $\phi_l$ for $0\leq l\leq k$, we find
that $\P_\Lambda$ is equivalently defined by
\be
\P_\Lambda:={\rm span} \{\phi_\nu \; : \; \nu\in \Lambda\},
\ee
Polynomial spaces associated to 
downward closed sets have been introduced in \cite{K}, 
in dimension $d=2$ and refered to as {\em polyn\^omes pleins}.
Later, these notions were studied in general dimension $d$,
in \cite{DR} and \cite{LL}. Note that in dimension $d=1$, 
a downward closed set is simply of the form $\{0,1,\dots,n\}$.

The sets index sets $\Lambda_n$ corresponding to the 
$n$ largest $\|t_\nu\|_V$, $\|v_\nu\|_V$ or $\|w_\nu\|_V$
are generally {\it not} downward closed sets. A legitimate question is
therefore: does there exists nested sequences $(\Lambda_n)_{n\geq 0}$ of
downward closed sets such that the truncated Taylor or Legendre series using
such sets have the same convergence rates
as those obtained in Corollaries \ref{corratetay} and \ref{corrateleg},
using the $n$ largest $\|t_\nu\|_V$, $\|v_\nu\|_V$ or $\|w_\nu\|_V$?
The results of the present section give a positive result to this question.

Let us begin by observing  that if a sequence $(c_\nu)_{\nu\in\cF}$ of positive numbers is
monotone non-increasing, that is, if
\be
\nu\leq \t \nu \Rightarrow  c_{\t \nu} \leq c_{\nu},
\ee
then the set $\Lambda_n$ corresponding to the $n$ largest values of $c_\nu$
is downward closed, provided that it is unique. In case of non-uniqueness, 
there is at least one realization of such a set which is downward closed.
In addition, there exists a sequence $(\Lambda_n)_{n\geq 1}$ of such
realizations which is nested. Note that in such a realization, we necessarily have $\Lambda_0=\{0\}$.

For an arbitrary sequence $c=(c_\nu)_{\nu\in\cF}\in \ell^\infty(\cF)$ we 
introduce its {\it monotone majorant} which is the sequence $\hat c=(\hat c_{\nu})_{\nu\in\cF}$
defined by
\be
\hat c_\nu:=\sup_{\t \nu\geq \nu} |c_{\t \nu}|.
\ee
This is  the smallest monotone non-increasing sequence that dominates $c$. In order
to study best $n$-term approximations using downward closed sets, we introduce
the following sequence spaces.

\begin{definition}
\label{deflpm}
For $0<p<\infty$, we say that a sequence $c\in \ell^\infty(\cF)$ 
belongs to $\ell^p_m(\cF)$ if and only its monotone majorant $\hat c$ 
belongs to $\ell^p(\cF)$ and we define
\be
\|c\|_{\ell^p_m}:=\|\hat c\|_{\ell^p}.
\ee
\end{definition}

Combining this definition with Lemma \ref{stechkin} shows that
if $0<p<q\leq \infty$ and if $(c_\nu)_{\nu\in\cF}$ is a positive
sequence which belongs to $\ell^p_m(\cF)$, then one has the tail bound
\be
\(\sum_{\nu\notinÊ\Lambda_n} c_\nu^q\)^{1/q} \leq Cn^{-s}, \quad 
C=\|(c_\nu)_{\nu\in\cF}\|_{\ell^p_m}, \quad s:=\frac 1p-\frac 1q,
\label{stechlower}
\ee
where $\Lambda_n$ is any downward closed set of indices corresponding to
the $n$ largest terms of the monotone majorant $\hat c$ of $c$.
We may therefore obtain the same rate $n^{-s}$ as in 
Lemma \ref{stechkin} now using downward closed sets. 

We would therefore  like to know under which circumstances the sequences
$(\|t_\nu\|_V)_{\nu\in\cF}$, $(\|v_\nu\|_V)_{\nu\in\cF}$ and $(\|w_\nu\|_V)_{\cF}$
belong to $\ell^p_m(\cF)$. The following result, originally proved
in \cite{CCDS} in the case of elliptic parametric PDEs and in \cite{CCS1} for other
models, shows that this
holds under the exact same assumptions as in Theorem \ref{theomain}.

\begin{theorem}
\label{theomainlower}
Consider a parametric problem of the form \iref{genpar}
such that {\bf Assumption A} holds 
for a suitable affine representer $(\psi_j)_{j\geq 1}$.   Then, the following summability results hold:
\begin{itemize}
\item
If the assumptions of Theorem \ref{theoneiacU} are satisfied,
and if in addition $(\|\psi_j\|_X)_{j\geq 1}\in \ell^p(\N)$ for some $p<1$, 
then $(\|t_\nu\|_V)_{\nu\in\cF}\in \ell^p_m(\cF)$
for the same value of $p$. 
\item
If the assumptions of Theorem \ref{theoneiaU} are satisfied,
and if in addition $(\|\psi_j\|_X)_{j\geq 1}\in \ell^p(\N)$ for some $p<1$, 
then $(\|v_\nu\|_V)_{\nu\in\cF}\in \ell^p_m(\cF)$ and $(\|w_\nu\|_V)_{\nu\in\cF}\in \ell^p_m(\cF)$
for the same value of $p$.
\end{itemize} 
\end{theorem}

\noindent
{\bf Proof:}  Similar to the proof of Theorem \ref{theomain}, we use the estimates
\iref{esttaylor}, \iref{estlegendrev} and \iref{estlegendrew} 
for the $\|t_\nu\|_V$, $\|v_\nu\|_V$ and $\|w_\nu\|_V$.

In the case of $\|t_\nu\|_V$, the estimate has the form
\be
\|t_\nu\|_{\nu\in\cF} \leq e_\nu:= C\inf \rho^{-\nu},
\ee
where the infimum is taken over the set of sequences $\rho$
of numbers larger than $1$ that satisfy the constraint \iref{constrhoeps}.
Since for any such $\rho$, the sequence $(\rho^{-\nu})_{\nu\in \cF}$ 
is monotone non-increasing, it follows that the 
sequence $(e_\nu)_{\nu\in\cF}$ is also monotone non-increasing.
On the other hand, the proof of Theorem \ref{theomain} shows that
$(e_\nu)_{\nu\in\cF}\in \ell^p(\cF)$. This implies that $(\|t_\nu\|_V)_{\nu\in\cF}\in \ell^p_m(\cF)$.

We cannot proceed in the same way for the Legendre coefficients
$\|v_\nu\|_V$ and $\|w_\nu\|_W$ since the right side $C r(\nu,\rho)$ in the estimates
\iref{estlegendrev} and \iref{estlegendrew} do not have the monotone
non-increasing property due to the presence of the factors $\theta(\rho_j)$ and $(1+2\nu_j)$.
Instead we slightly modify the construction of the sequence $\rho=\rho(\nu)$
in the proof of Theorem \ref{theomain}, and show that the
resulting sequence of estimates
\be
r_\nu=r(\nu,\rho(\nu)),
\ee
has a monotone majorant which is $\ell^p$ summable over $\cF$.
Here again, it suffices to work with the estimate \iref{estlegendrew} which is the largest one.
 
We use the same notation as in Theorem 3.9, in particular   $b_j:=\|\psi_j\|_X$.    For a constant $\beta>0$ to be fixed later, we take
 $J\geq 1$ large enough such that
\be
\sum_{j>J} b_j\leq \frac \e {3\beta},
\ee
where $\e$ is the right side of the constraint \iref{constrhoeps}. 

We now let $\nu\in\cF$ and fix $\nu$ and proceed to define an appropriate sequence $\rho=\rho(\nu)$ for this $\nu$.  Namely, using  the same splitting of
$\N$ into $E$ and $F$, we   take
\be
\rho_j=\kappa:=1+\frac {\e} {3\|b\|_{\ell^1}},\quad j\in E\cap {\rm supp}(\nu),
\ee
where $b=(b_j)_{j\geq 1}$ and
\be
\rho_j=\kappa+\beta+\frac {\e \nu_j}{3b_j |\nu_F|}, \quad j\in F \cap {\rm supp}(\nu).
\ee
We again take $\rho_j=1$ if $j\notin{\rm supp}(\nu)$.
Therefore $\rho_j>1$ when $j\in {\rm supp}(\nu)$, and in addition
\be
\sum_{j\geq 1}(\rho_j-1)b_j \leq \frac {\e \sum_{j\leq J} b_j}{ 3\|b\|_{\ell^1}}+Ê\sum_{j>J} \(\frac {\e b_j} {3\|b\|_{\ell^1}}+\beta b_j+\frac {\e \nu_j}{3 |\nu_F|}\)
\leq \e.
\ee
which shows that the constraint \iref{constrhoeps} is satisfied.

For this choice of $\rho$, the   estimate   \iref{estlegendrew} 
may be written 
\be
\|w_\nu\|_V\le r_\nu=r(\nu,\rho(\nu))=r_E(\nu)r_F(\nu),
\ee
with $r_E(\nu)$ as in the proof of Theorem \ref{theomain}, and a slightly modified $r_F(\nu)$ that  incorporates the new form of $\rho_j$ for $j\in F$.   This new $r_F(\nu)$ satisfies
$$
\begin{disarray}{ll}
r_F(\nu) &\leq \t r_F(\nu)Ê:=\prod_{j\in F\cap{\rm supp}(\nu)} \theta(\kappa)(1+2\nu_j)\(\beta+\frac {\e \nu_j}{3b_j |\nu_F|}\)^{-\nu_j} \\
&
\leq \prod_{j\in F\cap{\rm supp}(\nu)} \theta(\kappa)(1+2\nu_j)\(\frac {\e \nu_j}{3b_j |\nu_F|}\)^{-\nu_j}.
\end{disarray}
$$
Since $\kappa>1$, we there exists $C_0=C_0(\kappa) > 0$ such that
$(1+2n) \leq C_0 (\frac{1+\kappa}{2})^n$ for any $n\geq 1$ and so we can write
\be
r_E(\nu) \leq \t r_E(\nu):= C \prod_{j\in E} \eta^{\nu_j}, \quad \eta:=\frac {1+\kappa}{2\kappa}<1\;\;{\rm and}\;\; 
C=(C_0\theta(\kappa))^J.
\label{tre}
\ee
The same argument as in the proof of Theorem \ref{theomain} shows that, up to choosing a larger $J$, the estimates
\be
\t r_\nu:=\t r_E(\nu)\t r_F(\nu),
\label{trnu}
\ee
are $\ell^p$ summable over $\cF$.

We conclude by showing that $(\t r_\nu)_{\nu\in\cF}$ is monotone non-increasing
if $B$ has been chosen large enough. 
On the one hand, since $\eta<1$, it is readily seen that
\be
\nu\leq \t \nu \Rightarrow \t r_E(\t\nu) \leq \t r_E(\nu).
\ee
For proving a similar monotonicity property for the second factor $\t r_F$, it suffices 
to show that $\t r_F(\nu)$ is reduced if we increase 
$\nu_j$ by $1$ for any $j>J$, that is
\be
\t r_F(\nu+e_j) \leq \t r_F(\nu),
\ee
where $e_j=(0,\dots,0,1,0,\dots)$
is the Kroenecker sequence with $1$ at position $j>J$. 
In the case where $\nu_j\neq 0$, we may write
$$
\begin{disarray}{ll}
\frac {\t r_F(\nu+e_j)}{\t r_F(\nu)} & =\frac {1+2\nu_j+2}{1+2\nu_j}
\frac{ \(\beta+\frac {\e \nu_j}{3b_j |\nu_F|}\)^{\nu_j}}
{ \(\beta+\frac {\e (\nu_j+1)}{3b_j( |\nu_F|+1)}\)^{\nu_j+1}}
\prod_{k\in F\cap{\rm supp} (\nu)-\{j\}}\(\frac {\beta+\frac {\e \nu_k}{3b_k |\nu_F|}}{\beta+\frac {\e \nu_k}{3b_k (|\nu_F|+1)}}\)^{\nu_k}
\\
& \leq \frac {2}{\beta+\frac {\e (\nu_j+1)}{3b_j( |\nu_F|+1)}}
\prod_{k\in F\cap{\rm supp} (\nu)}\(\frac {\beta+\frac {\e \nu_k}{3b_k |\nu_F|}}{\beta+\frac {\e \nu_k}{3b_k (|\nu_F|+1)}}\)^{\nu_k}
 \leq \frac 2 \beta \(\frac {1+|\nu_F|} {|\nu_F|}\)^{|\nu_F|},
\end{disarray} 
$$
and therefore 
\be
\frac {\t r_F(\nu+e_j)}{\t r_F(\nu)}  \leqÊ\frac {2e}{\beta}
\ee
In the case where $\nu_j=0$, we have
\be
\frac {\t r_F(\nu+e_j)}{\t r_F(\nu)} =
\frac {3c_\kappa}{\beta+\frac {\e (\nu_j+1)}{3b_j( |\nu_F|+1)}}
\prod_{k\in F\cap{\rm supp} (\nu)-\{j\}}\(\frac {\beta+\frac {\e \nu_k}{3b_k |\nu_F|}}{\beta+\frac {\e \nu_k}{3b_k (|\nu_F|+1)}}\)^{\nu_k}
 \leqÊ\frac {3c_\kappa e}{\beta},
\ee
We thus find that $(\t r_\nu)_{\nu\in\cF}$ is monotone non-increasing
provided that $\beta\geq \max\{2e,3c_\kappa e\}$. \hfill $\Box$
\nl

Combining the above Theorem with \iref{stechlower}, we obtain the following result.

\begin{cor}
\label{corratemonot}
Corollaries \ref{corratetay} and \ref{corrateleg} remain valid, 
with the sets $\Lambda_n$ of corresponding to $n$ largest 
terms in the sequences $(\|t_\nu\|_V)_{\nu\in \cF}$,  $(\|v_\nu\|_V)_{\nu\in \cF}$ or  $(\|w_\nu\|_V)_{\nu\in \cF}$,
replaced by downward closed sets $\Lambda_n$ corresponding to the $n$ largest
terms in the monotone majorants of each of these sequences.
\end{cor}

\subsection{Exponential approximation rates}
\label{subexp}

The rates of convergence $n^{-s}$ that are established for polynomial approximations
in Corollaries \ref{corratetay} and \ref{corrateleg} are of algebraic type.
We conclude this study of polynomial approximation by a brief 
discusssion on the circumstances where faster rates of exponential type
can be established.
For this, we focus on the finite dimensional case, that is, when 
finitely many $\psi_j$ are non-zero in the affine representation \iref{affine}.
In such a case, one first obvious observation is that since
$(\|\psi_j\|_X)_{j\geq1}\in \ell^p(\N)$ for all values of $p>0$,
Corollaries \ref{corratetay} and \ref{corrateleg} give 
convergence rates $n^{-s}$ for all $s>0$. However, a more detailed inspection
shows that the multiplicative constant $C_s$ obtained in front of 
this rate grows very fast to $+\infty$ as $s\to +\infty$. Instead of trying
to search for a fast rate by optimizing $C_sn^{-s}$ over $s$ for a given $n$, we return 
to the estimates on the polynomial coefficients and use them to obtain exponential
convergence rates for the truncated series \iref{taylor}, \iref{legendre} or \iref{renormlegendre}.

Without loss of generality, we assume that only $\{\psi_1,\dots,\psi_d\}$ are
non-zero, meaning that the scalar parameter vector
is now
\be
y=(y_1,\dots,y_d)\in U:=[-1,1]^d,
\ee 
and that the solution map
$y \mapsto u(y)$ from $U$ to $V$ is finite dimensional. Polynomial approximations
are again based on truncation of the series \iref{taylor}, \iref{legendre} or \iref{renormlegendre},
now with
\be
\cF=\N^d.
\ee
For the sake of simplicity, we focus our attention on Taylor series and make some further
remarks on the case of Legendre series.

A particularly simple case for estimates of Taylor coefficients is
that of the disjoint inclusion model for the elliptic
and parabolic PDEs \iref{ellip}Ê and \iref{parab} discussed in \S \ref{secrefined}.
In this case, working under ${\rm \bf UEA}(r)$, we explicitly solved \iref{minesttaylor}
for any given  $0<t<r$ and obtained the estimate
\be
\|t_\nu\|_V \leq C\rho^{-\nu}= C\prod_{j=1}^d \rho_j^{-\nu_j}
\label{esttaylorfinite}
\ee
with $C=C_t$ and
\be
\rho_j=\rho_j^*=\inf_{x\in D} \frac {\o a(x)-t}{|\psi_j(x)|}>1,\quad j=1,\dots,d.
\label{rhojdisjoint}
\ee
More generally, if we work under the assumptions of Theorem \ref{theoneiacU},
we know from Lemma \ref{lemmaesttaylor} that we have an estimate
of the form \iref{esttaylorfinite} for any choice of $\rho_j\geq 1$ such that
\be
\sum_{j=1}^d(\rho_j-1)\|\psi_j\|_X\leq \e.
\ee
We may, for instance, take
\be
\rho_j:=1+\frac {\e}{\|\psi_j\|_X} >1,\quad j=1,\dots,J.
\ee
We thus again reach the estimate \iref{esttaylorfinite}
with a fixed finite vector $(\rho_1,\dots,\rho_d)$ independent of $\nu$ and whose oordinates are strictly larger than $1$.

Based on such an estimate, a natural choice for the sets $\Lambda_n$ is
to pick the indices $\nu$ corresponding to the $n$ largest values of $\rho^{-\nu}$. Equivalently,
for any  given threshold $\eta>0$ we define
\be
\Lambda_n:=\{\nu\in\cF \; : \; \rho^{-\nu} \geq \eta\},\quad {\rm where}\quad n=n(\eta):=\#\{\nu\in\cF \; : \; \rho^{-\nu} \geq \eta\}.
\ee
Notice that as we vary $\eta>0$, it may be that not all values of $n$ 
arise because of possible ties in the values of $\rho^{-\nu}$.

Let us now focus on the
particular thresholds
\be
\eta=2^{-k}, \quad k\geq 0,
\ee
we may write, with $n:=n(k)$ growing with $k$,
\be
\Lambda_n=S_k:=\{\nu\in\cF \; : \; \sum_{j=1}^d \lambda_j \nu_j \leq k\}, \quad \lambda_j:= \log_2(\rho_j)> 0.
\ee
Sets of this type consist of all integer lattice points inside the simplex with bounding 
  hyperplanes given by the coordinate hyperplanes to gether with the hyperplane   $\sum_{j=1}^d t_j\lambda_j=k$. Note that
these sets  are downward closed. 

The cardinality of the above $\Lambda_n$ is bounded
from above by the volume of the continuous
simplex
\begin{eqnarray*}
\label{growthsk}
T_k:&=&\{(t_1,\dots,t_d)\in\R^d: t_jÊ\geq -1, \ j=1,\dots,d, \ {\rm and} \ \; :\; \sum_{j=1}^d \lambda_j t_j \leq k\}\\
&=&\{(t_1,\dots,t_d)\in\R^d: t_jÊ\geq -1, \ j=1,\dots,d, \ {\rm and} \ \; :\; \sum_{j=1}^d \lambda_j t_j \leq k\}.
\end{eqnarray*}
This gives the crude cardinality bound
\be
 \#(\Lambda_n)=\#(S_k) \leq |T_k| =
\frac 1{d!} \prod_{j=1}^d\(\frac {k+\sum_{j=1}^d \lambda_j} {\lambda_j}\)\le C k^d
\ee
where $C$  depends on $d$ and on $(\lambda_1,\dots,\lambda_d)$.

Likewise, we may estimate the approximation error when retaining the $n$ terms whose indices are in $\Lambda_n$  by
$$
\begin{disarray}{ll}
\sup_{y\in U} \Big \|u(y)-\sum_{\nu\in S_k}t_\nu y^\nu\Big \| &
\leq \sum_{\nu\notin S_k} \|t_\nu\|_VÊ\\
&\leq C\sum_{\nu\notin S_k} \rho^{-\nu}Ê\\
& \leq C\sum_{l\geq k} 2^{-l} \#\{\nu \: : \: 2^{-l-1}\leq  \rho^{-\nu} <2^{-k}\}\\
& \leq C\sum_{l\geq k} 2^{-l}\#(S_{l+1}).
\end{disarray}
$$
Using the estimate \iref{growthsk} on the asymptotic growth of $\#(S_k)$,
we this find that
\be
\sup_{y\in U} \Big \|u(y)-\sum_{\nu\in S_k}t_\nu y^\nu\Big \| \leq C \sum_{l\geq k} 2^{-l}(l+1)^d. 
\ee
Combining this estimate with \iref{growthsk}, we obtain
\be
\sup_{y\in U} \Big \|u(y)-\sum_{\nu\in S_k}t_\nu y^\nu\Big \| \leq  C {\rm exp}(-ck),
\ee
which is equivalent to the exponential rate
\be
\sup_{y\in U} \Big \|u(y)-\sum_{\nu\in \Lambda_n}t_\nu y^\nu\Big \| \leq  C {\rm exp}(-cn^{-1/d}),
\ee
with multiplicative constants $c$ and $C$ that depend on $d$ and on $(\lambda_1,\dots,\lambda_d)$.
Since this rate is valid for all $n$ of the form $\#(S_k)$ which grow like $k^d$, it is easily seen
that it is also valid for all values of $n\geq 1$, up to a change in the multiplicative constants.

\begin{remark}
We notice that this exponential rate deteriorates as $d$ grows, due to the power $1/d$,
as well as to the hidden dependence on $d$ in the constants $c$ and $C$. However
in the case where the $\ell^p$ norm of $(\|\psi_j\|_X)_{j=1,\dots,d}$ remains uniformly bounded
for some $0<p<1$ as we raise $d$, our analysis of the infinite dimensional case
always ensures the algebraic rate $n^{-s}$ with $s:=\frac 1 p-1$.
\end{remark}

\begin{remark}
A similar analysis leads to the same exponential rates for approximation by truncated Legendre series,
now under the assumptions of Theorem \ref{theoneiaU}, based on the estimates
\iref{estlegendrew} and \iref{estlegendrev}, up to a proper treatment of the algebraic
factors $\theta(\nu_j)$ and $(1+2\nu_j)$ appearing in these estimates.
\end{remark}

\section{Estimating the $n$-widths of solution manifolds} 

We have already
noted that, when approximating the solution map by separable
expansions of the form \iref{separxa} or \iref{separxy}, the best achievable error in $L^\infty(\cA,V)$ or in $L^\infty(U_\cA,V)$
is described by the $n$-width of the solution manifold $\cM=u(\cA)$ in $V$, that is,
\be
d_n(\cM)_V:=\inf_{\dim(V_n)=n} \sup_{v\in\cM} \min_{w\in V_n} \|v-w\|_V
\ee
In this section, we use the polynomial approximation results established in the
previous section  to derive a priori estimates for the decay of $d_n(\cM)_V$.

\subsection{Estimates of $n$-width by polynomial approximation}

In the case where {\bf Assumption A}Ê holds, we may use the polynomial approximation results 
of \S 3  to estimate $d_n(\cM)_V$ from above. 
Indeed, if $u_n(y)=\sum_{\nu\in\Lambda_n} c_\nu y^\nu$ is
a polynomial approximation to the map $y\mapsto u(y)$ for some
set $\Lambda_n\subset \cF$ of cardinality $n$, we  
define the $n$ dimensional space
\be
V_n:={\rm span}\{c_\nu \; : \; \nu\in \Lambda_n\}Ê\subset V,
\ee
and observe that
\be
d_n(\cM)_V \leq \sup_{v\in\cM} \min_{w\in V_n} \|v-w\|_V
= \sup_{y\in U_{\cA}} \min_{w\in V_n} \|u(y)-w\|_V \leq \|u-u_n\|_{L^\infty(U,V)}.
\label{widthpol}
\ee
Therefore a polynomial approximation bound in $L^\infty(U,V)$ induces
an estimate on the $n$-width of $\cM$ in $V$. Combining this observation
with Corollary \ref{corrateleg}, we obtain the following result.

\begin{cor}
\label{corwidth}
Consider a parametric problem of the form \iref{genpar}
such that {\bf Assumption A} holds 
for a suitable affine representer \iref{affine}. 
Assume that the solution map $u\mapsto u(a)$ admits a
holomorphic extension over an open set $\cO\subset X$ which contains the 
compact set $a(U)$ and this extension satisfies the  uniform bound
\be
\sup_{a\in\cO} \|u(a)\|_V \leq C.
\label{uaC1}
\ee
If $(\|\psi_j\|_X)_{j\geq 1}\in \ell^p(\N)$ for some $0<p<1$.
then 
\be
d_n(\cM)_V \leq C(n+1)^{-s}, \quad n\geq 1, \quad  s:=\frac 1 p-1,
\label{nwidthest1}
\ee
for a suitable constant $C$.
\end{cor}

\noindent
{\bf Proof:} 
We consider the truncated Legendre expansion
\be
u_n=\sum_{\nu\in \Lambda_n} w_\nu P_\nu,
\ee
where $\Lambda_n$ is the set of indices corresponding to the $n$ largest $\|w_\nu\|_V$.
Since the assumptions of Corollary \ref{corrateleg} are satisfied, we obtain \iref{nwidthest1} with $C:=\|(\|w_\nu\|_V)_{\nu\in\cF}\|_{\ell^p}$, 
by combining  \iref{widthpol} and \iref{convrenormlegendre}. \hfill $\Box$
\nl

One drawback of the above result is that it requires
  {\bf Assumption A}.  For 
some natural examples of  compact sets $\cA$ of $X$, this assumption
may not hold. For instance, in the case of the elliptic equation \iref{ellip}, 
we know that the standard compact sets of $X=L^\infty(D)$
are   described by a smoothness assumption.  A typical  example
for $\cA$ of this type is 
\be
\cA:=\{a\in X\; : \; a>r, \;\; \|a\|_{C^\beta} \leq M\},
\ee
for some $M,\beta,r>0$, where $C^\beta:=C^\beta(D)$ is the H\"older space with smoothness $\beta>0$,
equiped with its usual norm 
\be
\|a\|_{C^\beta}:=\sup_{|\alpha|< m} \|\partial^\alpha a\|_{L^\infty}+ \sup_{|\alpha|=m} \sup_{x,x'\in D} |x-x'|^{-(\beta-m)}
|\partial^\alpha a(x)-\partial^\alpha a(x')|, \quad m:=\lfloor \beta \rfloor.
\ee
For such $\cA$, there are many ways to choose an $\o a\in \cA$ and  a properly normalized basis $(\psi_j)_{j\geq 1}$ such that expanding $a-\o a$ in this basis allows us to write  
\be
\cA\subset a(U),
\ee
with $a(U)$ of the form \iref{aofU}.    However, 
it will generally not follow that there is an $r'>0$ such that  for each $a$ in  $a(U)$, we have $ a>r'$.  Therefore, we are not guaranteed to have well posedness of the PDE for all $u(y)\in U$ and so {\bf Asssumption A} will  not hold for this affine representation.
We fix this defect in the next section by a different
approach based on local polynomial approximations.

\subsection{Estimates of $n$-width by local polynomial approximation}

In this section,  we treat   parameter sets $\cA\in X$ 
which may not have {\bf Assumption A} .    We assume that $(\psi_j)_{j\geq 1}$ is a complete representer for $\cA$
and in addition that $(\|\psi_j\|_X)\in \ell^1(\N)$. 
It follows that for each $(z_j)_{j\ge 1}\in \cU$, the series $\sum_{j\ge 1}z_j\psi_j$ converges in $X$ and so the set
\be
\cR:=\Big \{\sum_{j\geq 1} z_j\psi_j \; : \; z=(z_j)_{j\geq 1} \in \cU\Big\}.
\label{cR}
\ee
is well defined.
We replace {\bf Assumption A} by the requirement
\be
\cA\subset \cR.
\ee
Notice that, in contrast to $\cA$, the set $\cR$
might not be contained in the open set $\cO$ over which the solution map
admits a bounded holomorphic extension.  However, we will  remedy this problem by using  the 
following covering result.

\begin{lemma}
\label{lemmacover}
Let $\cA$ be a compact set in a complex Banach space $X$, and 
assume that $\cA\subset \cR$ where $\cR$ is of the form \iref{cR}
for a family of functions $(\psi_j)_{j\geq 1}$ such that $(\|\psi_j\|_X)_{j\geq 1} \in \ell^1(\N)$.
Let $\cO$ be any open set of $X$ which contains $\cA$.
Then, there exists $\eta,\e>0$, an integer $J\geq 1$ and a finite collection
$\{\o a_1,\dots, \o a_M\} \subset X$ such that defining
\be
\t \psi_j:=\eta \psi_j, \quad j=1,\dots,J, \quad \t \psi_j:=\psi_j, \quad j>J,
\ee
and for any  sequence $z=(z_j)_{j\geq 1} \in\C^\N$
\be
a_i(z):=\o a_i+\sum_{j\geq 1} z_j\t \psi_j,\quad i=1,\dots,M,
\ee
whenever the series on the right converges, the following holds:

\noindent
{\rm (i)}
The compact set $\cA$ admits the following cover
\be
\cA\subset \displaystyle{ \bigcup_{i=1}^M \cA_i}, \quad \cA_i:=a_i(\cU)=\{a_i(z)\; : \; z\in \cU\}.
\ee

\noindent
{\rm (ii)}
The compact sets $\cA_i$, $i=1,\dots,M$, are all contained in $\cO$.

\noindent
{\rm (iii)}
For any sequence $\rho=(\rho_j)_{j\geq 1}$ of numbers, each  larger than $1$,
which satisfies the constraint $\sum_{j\geq 1}(\rho_j-1)\|\t \psi_j\|_X\leq \e$, there exists, 
for each $j\geq 1$,  an open set $\cO_{\rho_j}\subset \C$ which contains the disc $\{|z_j|\leq \rho_j\}$ and  for which  the set 
  $\cO_\rho:=\otimes_{j\geq 1} \cO_{\rho_j}$ satisfies
  \be
a_i(\cO_\rho):=\{a_i(z)\; : \; z\in \cU_\rho\}\subset \cO.
\ee

\end{lemma}

\noindent
{\bf Proof:} Similar to the proof Theorem \ref{theoneiacU}, we first observe that since $\cA$ is compact,
there is an $\e>0$ sufficiently small, such that the $3\e$ neighborhood of $\cA$ is contained in $\cO$, that is,
\be
\bigcup_{a\in\cA}B(a,3\e)\subset \cO.
\ee
For this $\e$, we  next choose $J\geq 1$ large enough so  that 
\be
\label{tailbound}
\sum_{j>J} \|\psi_j\|_X\leq \frac {\e}{4}.
\ee
We then  define 
\be
\eta:=\frac{\e}{4\sum_{j=1}^J\|\psi_j\|_X}.
\ee
This fixes the $\e, \eta$ and $J$ claimed in the theorem.
In going further, we use the notation
\be
\label{XJ}
\cU_J:=\{z\in \cU:\ z_j=0, \;j>J\}.
\ee

Since $\cA\subset \cR$, for any $a\in \cA$ there exists a $z\in \cU$ such that
\be
a=\sum_{j=1}^J z_j\psi_j +\sum_{j>J} z_j \psi_j =:a_J+(a-a_J).
\ee
Note that this decomposition may not be unique - since the $\psi_j$ are not assumed to be
linearly independent - but, for each $a\in \cA$, we assign one such decomposition. 
We can find a finite set $F\subset \cU_J$, such that, for each $z\in \cU_J$, there is a $z'\in F$ such that
\be
\label{B}
 \|z-z'\|_{\ell^\infty(\N)}\le \eta.
\ee
 We let $\{\o a_1,\dots,\o a_M\}$ be the finite set consisting of all elements in   $X$ of the form
 \be
 \o a_i=\sum_{j=1}^Jz_j'\psi_j,
 \ee
 where  $z'\in F$ and in addition there is an 
 $a=\sum_{j=1}^\infty z_j\psi_j\in \cA$,  such that
 \be
 \label{sty}
|z_j-z_j'|\le \eta,\quad j=1,\dots, J.
 \ee

Let us now  show (i). 
If $a\in \cA$ and $a=\sum_{j=1}^\infty z_j\psi_j$,    then according to \eref{B} and \eref{sty}, there is a $\o a_i$ such that 
\be
a_J-\o a_i= \sum_{j=1}^J  c_j\psi_j, \quad |c_j|\le \eta,
\ee
which implies that $a\in \cA_i$.

Next, note that (iii) implies (ii).  Indeed, take  any $\rho$ satisfying the assumptions of (iii), then $\cU\subset\cO_\rho$ and hence  the validity of (iii) implies
 $\cA_i\subset  a_i(\cO_{\rho})\subset  \cO$  for each $i=1,\dots,M$.

We are left to prove (iii).   For this,  let $\rho$ be any sequence satisfying the constraint in (iii) and define   for each $j\ge 1$, the sets
\be
\label{defO}
\cO_{\rho_j}:=\{|z_j| < \t \rho_j\}, \quad \t \rho_j:=\rho_j+ \frac{\e}{\sum_{j\geq 1}\|\t\psi_j\|_X}.
\ee
We need to check that $a_i(\cO_\rho)\subset\cO$, $i=1.2.\dots,M$.   For this, we fix any value of $i$. We know that $\o a_i=\sum_{j=1}^Jz_j'\psi_j$,  and from \eref{sty},  there is an $a^*=\sum_{j=1}^\infty z_j^*\psi_j\in\cA$ 
for which    $|z'_j-z^*_j|\le \eta$ for $j=1,\dots,J$.  In view of \eref{tailbound} and the definition of $\eta$,  we have    
\be
\label{bi}
\|\o a_i-a^*\|_X\le   \frac {\e}{2}.
\ee
Now take any $a\in a_i(\cO_\rho)$, that is
\be
a=\o a_i+\sum_{j\geq 1} z_j\t\psi_j,
\ee
with $z=(z_j)_{j\geq 1}\in \cU_\rho$. We define
\be
\t z_j=z_j\min\{1,|z_j|^{-1}\},\quad j\ge 1,
\ee
so that $(\t z_j)_{j\ge 1}$  is a point in $\cU$.   We can now estimate 

$$
\begin{disarray}{ll}
\|a-a^*\|_X & \leq \|a-\o a_i\|_X+ \|\o a_i-a^*\|_X \\
&\leq  \Big \|\sum_{j\geq 1} z_j\t \psi_j\Big \|_X+\frac {\e}{2}\\
& \leq  \Big \|\sum_{j=1}^{J} \t z_j\t \psi_j\Big \|_X+\Big \|\sum_{j>J} \t z_j\t \psi_j\Big \|_X+ \Big \|\sum_{j\geq 1} (\t z_j-z_j)\t \psi_j\Big \|_X+\frac {\e}{2}\\
&   \leq \eta \sum_{j=1}^{J}\|\psi_j\|_X+\sum_{j>J}\|\psi_J\|_X+ \Big \|\sum_{j\geq 1} (\t z_j-z_j)\t \psi_j\Big \|_X+\frac {\e}{2}\\
& \leq \frac{\e}{4} +\frac{\e}{4}+\Big \|\sum_{j\geq 1} (\t z_j-z_j)\t \psi_j\Big \|_X
+\frac {\e}{2}.
\end{disarray}
$$
Since 
$$|z_j-\t z_j|\le (\t \rho_j-1)\le  \t \rho_j-\rho_j +\rho_j-1, \quad j\ge 1,
$$
 we obtain
$$
 \|a-a^*\|_X \leq \e+\sum_{j\geq 1}(\t \rho_j-\rho_j)\|\t\psi_j\|_X +\sum_{j\geq 1}(\rho_j-1)\|\t\psi_j\|_X \le \e+\e+\e= 3\e,
 $$
 where we have used \eref{defO} to bound the first sum and the assumption on $(\rho_j)_{j\ge 1}$ in estimating the second sum.
This shows that $a$ belongs to the $3\e$-neighborhood of $\cA$ which is contained in $\cO$. 
Therefore $a_i(\cU_\rho)\subset \cO$.
\hfill $\Box$
\nl

With the help of the above lemma, we now establish a result 
which shows that the conclusion of Corollary \ref{corwidth} remains valid
without the assumption that $\cA$ is of the exact form $a(U)$.

\begin{theorem}
\label{theowidth1}
For a parametric problem of the form \iref{genpar}, 
assume that the solution map $u\mapsto u(a)$ admits a
holomorphic extension over an open set $\cO$ of the complex Banach space $X$
which contains $\cA$,
 with uniform bound
\be
\sup_{a\in\cO} \|u(a)\|_V \leq C.
\label{uaC2}
\ee
Assume in addition that there exists functions $(\psi_j)_{j\geq 1}$ in $X$
such that $(\|\psi_j\|_X)_{j\geq 1}\in \ell^p(\N)$ for some $0<p<1$, and such that
$\cA\subset \cR$, where $\cR$ is of the form \iref{cR}. Then, there exists $C>0$
such that one has
\be
d_n(\cM)_V \leq Cn^{-s}, \quad n\geq 1, \quad  s:=\frac 1 p-1.
\label{nwidthest2}
\ee
\end{theorem}

\noindent
{\bf Proof:} Applying Lemma \ref{lemmacover}, we write 
\be
\cA \subset \bigcup_{i=1}^M\cA_i,
\ee
and therefore
\be
\cM\subset \bigcup_{i=1}^M \cM_i, \quad  \cM_i:=u(\cA_i).
\ee
It now sufficient to prove that the estimate \iref{nwidthest2}
holds for each $\cM_i$ in place of $\cM$, that is
\be
d_n(\cM_i)_V \leq Cn^{-s}, \quad n\geq 1, \quad  s:=\frac 1 p-1.
\label{nwidthest3}
\ee
Indeed, if for each $i=1,\dots,M$ one can approximate all elements of $\cM_i$ with accuracy
$\delta$ by elements from an $n$ dimensional space $V_{n,i}$, then
one can approximate all elements of $\cM$ with the same accuracy by elements
from the space $V_{n,1}\oplus\cdots \oplus V_{n,M}$ which has at most dimension $nM$.
This shows that 
\be
d_{Mn}(\cM)_V\leq \max_{i=1,\dots,M} d_n(\cM_i)_V,
\ee
and therefore \iref{nwidthest3} implies \iref{nwidthest2} up to a change in the constant $C$.

The proof of \eref{nwidthest3} follows from Corollary \ref{corratetay}.  We fix  $i\in \{1,\dots,M\}$, we know that
\be
\cA_i\subset a_i(\cU),
\ee
where $a_i(z):=\o a_i+\sum_{j\geq 1} z_j\t\psi_j$. It follows that the assumptions of this corollary
are satisfied for $\cA_i$ and  $a_i$ in place of $\cA$ and $a$. This confirms the estimate \iref{nwidthest3} and therefore concludes the proof.
\hfill $\Box$

\begin{remark}
\label{remgenu3}
The above theorem can be formulated for
a general map $u$ from $\cA$ to $V$ that is not necessarily the solution map
of parametric PDE,  following the arguments from Remarks \ref{remgenu1} and \ref{remgenu2}.
\end{remark}

\subsection{$n$-widths under holomorphic maps: a general result}

In this section, we let $u$ be any map from $\cA$ to $V$, not necessarily
the solution map to a parametric PDE.
In view of Remark \ref{remgenu3},
Theorem \ref{theowidth1} gives an estimate  for the $n$-widths of 
$\cM=u(\cA)$ whenever  $u$ has a bounded   holomorphic 
extension to a neighborhood of $\cA$, provided that  $\cA$ is contained
in a set $\cR$ of the form \iref{cR} with   $(\|\psi_j\|_X)_{j\ge 1}\in \ell^p$, with  $0<p<1$.
Notice that the containment and $\ell^p$ summability assumptions on $\cA$ imply a decay on the $n$-width of $\cA$ in $X$.
Namely, for $s:=\frac 1 p-1$,  we can  write
\be
d_n(\cA)_X \leq \sup_{a\in\cA}\min_{z\in\cU}\Big\|a-\sum_{j=1}^n z_j \psi_j\Big\|_X
\leq \sum_{j>n} \|\psi_j\|_X \leq C(n+1)^{-s},  n\ge 0,
\ee
where we have used Lemma \ref{stechkin}.  Thus, in a certain sense, the results of the previous section can be interpreted as saying the Kolmogorov widths of $ \cM$ inherit the decay rate of the widths of $\cA$.   Since, the sets $\cA$ are generally more accessible and their $n$-widths are more readily computed,  it is natural to ask whether there is a general principle in effect here.   That is,  do the $n$-widths $(d_n(\cM)_V)_{n\geq 1}$ of  an image $\cM=u(\cA)$ of a compact set $\cA$ under a general holomorphic map $u$ have the same decay
as that of $(d_n(\cA)_X)_{n\geq 1}$.  The main goal of this section is to show that there is indeed such a general principle in effect, however with a slight loss in the decay rate of  the widths of $\cM$ when compared
with those of $\cA$.

 Our first step in deriving such comparison results  is to show that whenever a compact set $\cA$ of  a Banach space $X$ has  widths  
 $(d_n(\cA)_X$) with some prescribed decay, then $\cA$ is  contained in a set $\cR$ 
of the form \iref{cR} where the $X$-norms of
the $\psi_j$ defining $\cR$ are $\ell^p$ summable for certain values of $p$.
For this, we need the following classical result due to
Auerbach, the proof of which is given below for completeness.

\begin{lemma}
Let $E$ be an $n$-dimensional subspace of  a complex Banach space $X$. 
Then, there exists a basis $\{\vp_1,\dots,\vp_n\}$ for  $E$ and a dual basis $\{\t \vp_1,\dots, \t\vp_n\}$ in $X'$ such that 
\be
\<\t \vp_i,\vp_j\>_{X',X}=\delta_{i,j}, \quad i,j=1,\dots n,
\ee
and 
\be
\|\vp_i\|_X=\|\t \vp_i\|_{X'}=1,\quad i=1,\dots,n.
\ee
\end{lemma}

\noindent
{\bf Proof:} We start with an arbitrary basis $\psi_1,\dots, \psi_n$  of $E$ and let 
$\t \psi_1,\dots, \t \psi_n$ be its dual basis in $E'$, that is 
\be
\label{dual basis}
\< \t \psi_i,\psi_j\>=\delta_{i,j},
\ee
where  $\delta $ is the Kronecker delta  and $\< \cdot,\cdot\>: = \< \cdot,\cdot\> _{X',X}$ throughout this proof.
Then, any $f\in E$ can be uniquely written as
\be
\label{unique}
f=\sum_{i=1}^n \langle \t \psi_i,f\rangle  \psi_i.  
\ee

Given any   $ (g_1,\dots,g_n)\in E^n$, we define
\be
\label{defJ}
J(g_1,\dots,g_n)=|\det(M(g_1,\dots,g_n))|, \quad M:=(\<\t \psi_i,g_j\>)_{i,j=1,\dots,n}.
\ee
We now take $(\vp_1,\dots,\vp_n)\in E^n$ such that
\be
(\vp_1,\dots,\vp_n):= \argmax J(g_1,\dots,g_n),\ee
where the maximum is taken over all $(g_1,\dots,g_n)\in E^n$ such that $\|g_i\|_X=1$ , for $i=1,\dots,n$.   This  maximum is attained since the function $J$ is continuous and we are maximizing over a compact set.   The   functions
$\vp_1,\dots,\vp_n$  are linearly independent since this  maximum is positive.  Hence, they form a basis for  $E$ and any 
$f\in E$ can be written uniquely as 
\be
\label{unique1}
f=\sum_{i=1}^n \<\t \vp_i,f\>\vp_i
\ee
where $\t \vp_i$, $i=1,\dots,n$, is it dual basis.  Applying the functional $\t\psi_i$ to both sides of \eref{unique1}, we obtain
 \be
\sum_{j=1}^n \<\t\vp_j,f\>\<\t \psi_i,\vp_j \>= \<\t\psi_i,f \>, \quad i=1,2,\dots,n
\ee
From Cramer's rule, it follows that, for any $j\in \{1,\dots,n\}$ and any $f\in E$,
\be
|\<\t \vp_j,f\>|= \frac {J(\vp_1,\dots,\vp_{j-1},f,\vp_{j+1},\dots,\vp_n)}{J(\vp_1,\dots,\vp_n)}\le 1.
\ee
This proves that  for each $j$, we have $\|\t\vp_j\|_{E'}=1$.  By application of the Hahn-Banach theorem, we can extend $\t \vp_j$ over all of $X$ with $\|\t \vp_j\|_{X'}=1$. \hfill $\Box$
\nl

Using  the Auerbach Lemma, we now show that whenever  $\cA$ is a compact set of a complex Banach space $X$ and $d_n(\cA)_X$ has some prescribed
rate of decay, then $\cA$ is  contained in a set $\cR$ 
of the form \iref{cR} and   the $X$-norms of
the $\psi_j$ defining $\cR$ have the same rate of decay as  $d_n(\cA)_X$ .

\begin{lemma}
\label{lemmaloss}
Let $X$ be a complex Banach space and $\cA\subset X$ be a compact space such that 
 \be
 \sup_{n\geq 1} n^sd_n(\cA)_{X} <\infty.
\label{AT1}  
 \ee
Then, there exists a family $(\psi_j)_{j\geq 1}$ of functions from $X$ such that
\be
\sup_{j\geq 1}j^s\|\psi_j\|_{X} <\infty,
\ee
and 
\be
\cA\subset \cR:=\Big \{\sum_{j\geq 1} z_j \psi_j \; : \; z=(z_j)_{j\geq 1} \in \cU \Big \}.
\ee
\end{lemma}

\noindent
{\bf Proof:} From \iref{AT1}, we know that 
here exists a constant $C>0$ and a sequence of spaces 
$(V_k)_{k\geq 0}$ with $V_k\subset X$ and  $\dim(V_k)=2^k$, such that
\be
\max_{a\in \cA} \min_{g\in X_{k} } \|a-g\|_X \leq C2^{-sk}, \;\; k\geq 0.
\ee
By replacing $V_k$ by $V_0+V_1+\dots + V_{k-1}$ and possibly changing the constant $C$, we may assume
that the spaces $V_k$ are nested: $V_{k-1}\subset V_k$, for all $k\ge 1$. 

For any $a\in \cA$, we denote by $a_k$ a best approximation to $a$ from $V_k$ for $k\geq 0$
and set $a_{-1}:=0$.  
Then, $g_k:=a_k-a_{k-1}$ is in $V_k$,  and we have 
\be
a=\sum_{k\geq 0} g_k.
\ee
In addition, there exists a constant $C>0$, such that
\be
\|g_k\|_{X}\leq C 2^{-sk}, \;\; k\geq 0.
\ee
By Auerbach's lemma, 
for every $k\geq 0$, there exists a basis $\{\vp_{k,l}\}_{l=1,\dots,2^k}$ 
of the space $V_{k}$, and a dual basis $\{\t\vp_{k,l}\}_{l=1,\dots,2^k}\subset X'$
such that $\|\vp_{k,l}\|_{X}=\|\t\vp_{k,l}\|_{X'}=1$. It follows that any $a\in \cA$ can be written as
\be
a=\sum_{k\geq 0} \sum_{l=1}^{2^k} z_{k,l}\vp_{k,l}, \;\; |z_{k,l}|\leq C 2^{-sk}.
\ee
Each integer $j\ge 1$ can be written uniquely as $j=2^k+l-1$ with $l\in \{1,\dots,2^k\}$.  We use this to define 
\be
\psi_j:=C 2^{-sk}\vp_{k,l},\quad j=2^k+l-1.
\ee
This gives that any $a\in \cA$ is of the form
\be
a=\sum_{j\geq 1} z_{j}\psi_j, \quad |z_{j}|\leq 1,
\ee
that is, $\cA\subset \cR$. In addition, we have
\be
\|\psi_j\|_{X} \leq Cj^{-s},
\ee
up to a change in the constant $C$.
\hfill $\Box$
\nl

An immediate consequence of the above lemma is that if 
$d_n(\cA)_X$ has the rate of decay $n^{-s}$, then $(\|\psi_j\|_X)_{j\geq 1}\in \ell^p(\N)$
for any $p$ such that $sp>1$. Combining this observation with Theorem \ref{theowidth1},
leads to the following result which shows that the rate of decay of $n$-width
is almost preserved under holomorphic maps, up to a loss of $1$ in the rate.

\begin{theorem}
\label{theowidth}
For a pair of complex Banach spaces $X$ and $V$, assume that
$u$ is a general holomorphic map from an open set $\cO\subset X$
into $V$ with uniform bound
\be
\sup_{a\in \cO} \|u(a)\|_V \leq C.
\ee
If $\cA \subset \cO$ is a compact subset of $X$ and $\cM=u(\cA)$,  then for any $s>1$ and $t<s-1$,
\be
\sup_{n\geq 1} n^sd_n(\cA)_{X}<\infty\;  \Rightarrow \;\sup_{n\geq 1} n^td_n(\cM)_V <\infty.
\label{widthmap}
\ee
\end{theorem}

Some comments on this result are in order. If $u$ was a linear map, could write
for any subspace $X_n\subset X$ of dimension $n$ and any $a\in \cA$,
\be
\min_{v\in V_n} \|u(a)-v\|_V \leq C \min_{\t a\in X_n}\|a-\t a\|_X, \quad C:=\|u\|_{\cL(X,V)}, 
\ee
with $V_n:=u(X_n)\subset V$ also of dimension $n$. Therefore, we would obtain
\be
d_n(\cM)_V\leq Cd_n(\cA)_X,
\ee
which implies that $d_n(\cM)_V$ has at least the same rate of decay as $d_n(\cA)_X$.
Theorem \ref{theowidth} shows that holomorphic maps behave almost as good as linear maps,
except for the loss of $1$ in the rate expressed by the inequality $t<s-1$.
This loss occurs due to a lack of sharpness in Lemma \ref{lemmaloss}:
if we start from the conclusion $\|\psi_j\|_X \leq C j^{-s}$ of this lemma,
we may only retrieve that 
\be
d_n(\cA)_X\leq d_n(\cR)_X \leq \sum_{j>n}\|\psi_j\|_X \leq C n^{1-s}.
\ee
An open question is if the implication \iref{widthmap} in
Theorem \ref{theowidth} remains valid with $t=s$.

\subsection{Towards faster low rank approximations}
\label{subrank}

We close the first part of this article with some remarks concerning
our approximation results. As explained in the introduction, 
our general interest is in the  accuracy of separable approximations of the form \eref{separxa} and \eref{separxy}. These approximations can 
be thought of  as the analog of low rank approximations for finite dimensional matrices.

Optimal approximations are provided by  best
optimal $n$-dimensional spaces $V_n$ either in the sense of $n$-widths for uniform approximation
or Karhunen-Loeve decompositions for approximation in the mean square sense.
Since these spaces are out of reach, both from a theoretical and computational point 
of view, we build sub-optimal approximations $y\mapsto u_n(y)$ based on
best $n$-term truncations of polynomial expansions.
This approach leads us to quantitative convergence results
such as in Corollaries \ref{corratetay}Ê and \ref{corrateleg}, and in turn
to estimates for the decay of the $n$-widths $d_n(\cM)_V$ of solution manifolds as 
discussed in \S 4, for example 
by using the estimate 
\be
d_n(\cM)_V \leq \|u-u_n\|_{L^\infty(U,V)},
\label{widthpol1}
\ee
in the case when {\bf Assumption A}Ê holds.

A legitimate question is to evaluate the possible lack of optimality
of the convergence rates, obtained by our polynomial approximation approach, in comparison
to the rates which could be achieved by using the optimal $n$-dimensional spaces $V_n$.
Equivalently, we would like to know if the rate of decay of the $n$-width 
$d_n(\cM)_V$ could sometimes be much faster than the rate of decay of
the polynomial approximation error on the right of \iref{widthpol1}.

We can give simple examples which reveal this lack of optimality
in the case of the elliptic equation \iref{ellip}. Here, we consider the finite dimensional
setting where
\be
a(y)=\o a+\sum_{j=1}^d y_j\psi_j.
\ee
In this setting, convergence rates of exponential type 
\be
\|u-u_n\|_{L^\infty(U,V)} \leq C{\rm exp}(-cn^{1/d}),
\label{exprate}
\ee
are established in \S \ref{subexp} using for $u_n$ the Taylor series truncated
with the index set corresponding to the $n$ largest values of $\|t_\nu\|_V$. 

Let us here consider the particular case of a piecewise constant diffusion coefficient
of the form
\be
a(y)=\sum_{j=1}^d (1+ \theta y_j)\Chi_{D_j},
\ee
where $\{D_1,\dots,D_j\}$ is a partition of $D$, that is, $\o a=1$ and $\psi_j=\theta  \Chi_{D_j}$. We assume
that $0<\theta=1-r<1$ so that ${\rm \bf UEA}(r)$ holds.

We first examine the case where $D$ is a one dimensional interval partitioned into 
sub-interval $D_j$. The problem now reads 
\be
-(a u')'=f,
\ee
with homogeneous Dirichlet boundary conditions at the endpoints of $D$. Since $a(y)$ is constant on 
each interval $D_j$, we find that the restriction of $u(y)$ to this interval is always the sum of an
affine function and of a scalar multiple of $F$ such that $F''=f$. It follows that, for any $y\in [-1,1]^d$,
the solution $u(y)$ belongs to the finite dimensional space
\be
V_{3d}={\rm span}\{ \Chi_{D_i},\,x \Chi_{D_i},\,F\Chi_{D_i} \; : \; i=1,\dots d\},
\ee
where $x$ stands for the identity function $x \mapsto x$. Using the fact that $u(y)$ is $0$
at the endpoints of $D$ and continuous at the breakpoint between the $D_j$, we find
that it belongs to an even smaller subspace of $V_{3d}$ that has smaller dimension $2d-1$.
This implies that
\be
d_n(\cM)_V=0,
\ee
for $n\geq 2d-1$, therefore showing that the rate in the right-hand side of \iref{exprate} is not sharp 
for $d_n(\cM)_V$.

Let us now examine a less trivial case where $D$ is a domain in higher dimension $m\geq 2$.
In such case, it is not true that $\cM$ belongs to a finite dimensional space, however
we can still show that the rate in the right-hand side of \iref{exprate} is not sharp 
for $d_n(\cM)_V$. For simplicity, consider the case of a two domains partition, that is, $d=2$.
Since $\|\psi_1\|_X=\|\psi_2\|_X=\theta$, the sets $\Lambda_n$ that are used in \S \ref{subexp} to obtain the rate
\be
\|u-u_n\|_{L^\infty(U,V)} \leq C{\rm exp}(-cn^{1/2}),
\label{exprate2}
\ee
have the simple structure
\be
\Lambda_n=\{|\nu|=\nu_1+\nu_2\leq k\},
\ee
for integers $k\geq 0$. Therefore, we use polynomial approximations of total degree $k$ of the form
\be
u_n(y)=\sum_{|\nu|\leq k} t_\nu y^\nu,
\label{unk}
\ee
which have accuracy
\be
\|u-u_n\|_{L^\infty(U,V)}\leq C{\rm exp}(-ck),
\ee
with $n=\frac {k(k+1)}2\sim k^2$. 

This trunctated power series can be interpreted in a different way by writing the elliptic 
equation in operator form
\be
\cB(y) u(y)= f, \quad \cB(y)=\o\cB+y_1\cB_1+y_2\cB_2,
\ee
where 
\be
\o \cB u:=-\Delta u\quad {\rm and}\quad \cB_ju:=-{\rm div}(\theta \Chi_{D_j}\nabla u).
\ee
With $\o\cB^{-1}$ the inverse of $-\Delta$ on $D$ with homogeneous Dirichlet boundary condition, we may 
then rewrite the equation as
\be
(I+y_1\t\cB_1+y_2\t\cB_2)u(y) =g,\quad g:=\o\cB^{-1}f,\quad \t\cB_j:=\o\cB^{-1}\cB_j, \ j=1,2.
\ee
It is easily seen that the Taylor series of $u(y)$ coincides with the Neumann series
\be
u(y)=\sum_{l\geq 0}(-1)^l(y_1\t\cB_1+y_2\t\cB_2)^l g.
\ee
The convergence of this series can be directly checked by observing that
\be
\|y_1\t\cB_1+y_2\t\cB_2\|_{\cL(V,V)} \leq \theta, \quad (y_1,y_2)\in [-1,1]^2.
\ee
In particular this confirms the exponential rate
\be
\|u-u_n\|_{L^\infty(U,V)} \leq C\theta^{k}=C{\rm exp}(-ck).
\ee
We now observe that, due to the fact that $\Chi_{D_1}+\Chi_{D_2}=\Chi_D$, we have the
identity
\be
\t\cB_1+\t\cB_2=\theta I.
\ee
We may therefore rewrite each term in the Neumann series as  
\begin{eqnarray}
\label{terms}
(-1)^l(y_1\t\cB_1+y_2\t\cB_2)^l g&=&(-1)^l(y_2\theta I+(y_1-y_2)\t\cB_1)^l g\nonumber \\
&=&(-1)^l \sum_{j=0}^l (y_2\theta)^{l-j} {l\choose j}(y_1-y_2)^j\t\cB_1^j g.
\end{eqnarray}

Therefore, summing the terms in \eref{terms}  from $l=0$ up to $k$, we may rewrite $u_n(y)$ as
\be
u_n(y)= \sum_{j=0}^k v_j \phi_j(y),
\ee
with 
\be
v_j:=\t\cB_1^j g\in V,
\ee
and
\be
\phi_j(y):=(y_1-y_2)^j\sum_{l=j}^k (-1)^l (y_2\theta)^{l-j}  {l\choose j}.
\ee
This new representation of $u_n(y)$ shows that it belongs to the 
$k+1$ dimensional space
\be
V_k={\rm span}\{v_0,\dots,v_k\}.
\ee
 We may thus conclude that
\be
d_{k+1} (\cM)_V \leq C\theta^{k}=C{\rm exp}(-ck),
\ee
Since $k\sim \sqrt{n}$, this shows  that the rate in the right-hand side of \iref{exprate2} is not sharp 
for $d_n(\cM)_V$.

These examples reveal that in certain relevant cases, 
polynomial approximations based on best $n$-term
truncations may be highly sub-optimal in comparison
to the $n$-width spaces. Note, however, that the rank reduction is made
possible due to fine properties of the affine representation \iref{affine}, such as the
fact that the $\psi_j$ are characteristic functions with disjoint supports. For other affine
representations with general functions $\psi_j$ which have overlapping
support, numerical computations show that polynomial approximation rates
are sometimes close to the optimal rates to be expected from arbitrary separable approximations.
The development of alternate strategies for a sharper convergence analysis of 
separable approximations is thus desirable, and it inevitably requires exploiting
the detailed structure of the affine representation.

\section*{Part II. Algorithms for parametric PDEs}

\section{Towards concrete algorithms}  The results exposed in the first part of this 
paper show that relevant instances of parametric PDEs
admit separable approximations $u_n$ of the form
\iref{separxa} or \iref{separxy} with error bounds 
that reflect a certain rate of convergence in terms of the
number $n$ of terms that are retained. 

However, these approximations
are obtained by mathematical techniques which, as such, cannot
be implemented through a computational algorithm. 
For example, in order to compute the best $n$-term 
truncation of the Legendre series we need in principle
to be able to compute exactly all Legendre coefficients $w_\nu$
and to search for the $n$ largest values of $\|w_\nu\|_V$.
This is unfeasible for two reasons: (i) we can only compute
the $w_\nu$ with limited precision due to spatial discretization,
for example through a finite element space $V_h$ of $V$, 
and (ii) we cannot perform an exhaustive search through the infinite
set $\cF$ of multi-indices.

In this second part of the paper, we discuss concrete numerical methods which 
compute separable approximations $u_n$, still of the form
\iref{separxa} or \iref{separxy}, however at
an affordable computational cost. 

\subsection{Space discretization and computational cost}

Our approach to the computation
of such approximations can be viewed as follows:
\begin{itemize}
\item[(i)] 
We develop and analyze strategies for computing separable expansions 
first based on a few instances of the exact solution 
maps $a\mapsto u(a)$ and $y\mapsto u(y)$, or quantities related
to these maps such as the Taylor coefficients $t_\nu$.
\item[(ii)]
We then instead apply
these strategies to the {\it approximate solution maps}
\be
a\mapsto u_h(a) \in V_h \quad {\rm and}Ê\quad y \mapsto u_h(y)\in V_h,
\ee
which correspond to a certain space discretization process for each instance
of the solution map in a fixed discretization space $V_h$.
\end{itemize}

Ideally, we would like to obtain error bounds for these approximations
which meet the benchmark established in the first part of 
the paper in terms of their decay as $n$ grows, up to an additional term
that reflects the space discretization error. 

We assume that
space discretization can be performed within a certain finite element
space $V_h$ of dimension $N_h$, through a numerical
solver which we may apply for  each individual instance of 
$a\in \cA$ or $y\in U_{\cA}$ to compute approximate solutions
$u_h(a)$ or $u_h(y)$ from $V_h$. For simplicity we assume 
\begin{itemize}
\item[{\rm (i)}]
A cost $C_h$ for computing $u_h(a)$ or $u_h(y)$ that is independent of $a$ or $y$.
\item[{\rm (ii)}]
An error bound 
\be
\sup_{a\in\cA} \|u(a)-u_h(a)\|_V=\sup_{y\in U_\cA} \|u(y)-u_h(y)\|_V\leq \e(h),
\label{errorh}
\ee
therefore also independent of $a$ or $y$.
\end{itemize}

Recall that making $\e(h)$ small requires to make $N_h$ large and $C_h$ even larger,
which is one of the motivations for reduced modeling.  

As an example of such a space discretization, consider the elliptic equation \iref{ellip}. 
We may then define the discrete solution
by the standard Galerkin method on $V_h$, that is, 
$u_h(a)\in V_h$ is defined by
\be
\int_D a\nabla u_h(a)\nabla v_h=\int_D fv_h, \quad v_h\in V_h.
\label{galerellip}
\ee
We may then use classical techniques of finite element approximation
of elliptic PDEs, see \cite{Ciarlet} or \cite{BrennerScott}, in order to obtain 
an error bound of the form \iref{errorh}. First, assuming
that $0<r\leq a\leq R$ for all $a\in \cA$, Cea's Lemma ensures that
\be
\|u(a)-u_h(a)\|_V \leq \sqrt{\frac R r}\min_{v_h\in V_h} \|u(a)-v_h\|_V.
\ee
Then, if $(V_h)_{h>0}$ are Lagrange finite elements spaces of polynomial degree $m\geq 1$
associated to a regular family of conforming simplicial partitions $(\cT_h)_{h>0}$
with mesh size $h>0$, we have for $1<r\leq m+1$ the classical approximation bound
\be
\min_{v_h\in V_h} \|u(a)-v_h\|_V \leq C h^{r-1} \|u(a)\|_{H^{r}(D)}.
\ee
We therefore obtain an error bound \iref{errorh} with $\e(h)\sim h^{r-1}$ provided
that $u(a)$ is bounded in $H^r(D)$ independently of $a\in \cA$.

\begin{remark}
Our approach to space discretization means in particular that, when 
computing polynomial approximations by trunctated expansions, the
Taylor or Legendre coefficients are discretized in the {\em same} finite element space
$V_h$, independently of their index $\nu$. An alternate approach, which
we do not embark in here, is to
search for space discretizations of these coefficients which {\em vary} with $\nu$,
with the objective of optimizing the total number of degrees of freedom
required to reach a given accuracy. This approach is analyzed 
in \cite{CDS1,CDS2} for Legendre and Taylor series. See also \cite{Git} for
computational approaches based on a global adaptivity both in
the parameter and space variable.
\end{remark}

When evaluating the total computational cost for computing
the separable approximation $u_n$, we make the distinction
between two types of cost:
\begin{itemize}
\item[{\rm (i)}]
The {\it offline cost} refers to the computation
of the functions $v_1,\dots,v_n$ which are used in \iref{separxa} or \iref{separxy},
or equivalently of the space $V_n:={\rm span}\{v_1,\dots,v_n\}$ which is used to simultaneously approximate
all members of the solution manifold $\cM$.
\item[{\rm (ii)}]
The {\it online cost} which refers to the computation of the approximate
solution $u_n(a)$ or $u_n(y)$ from $V_n$ for 
any given query  $a\in\cA$ or $y\in U$ .
\end{itemize}
One can view the offline cost as a ``one time only'' fixed cost, while the 
online cost could be repeated many times in certain applications of
reduced modeling.

\subsection{Polynomial approximation algorithms}

The first class of numerical methods that we study searches for
computable polynomial approximations of the 
general form \iref{polynomial}. For any finite set $\Lambda$, we define the space
\be
V_\Lambda:=V \otimes \P_\Lambda,
\ee
of $V$-valued polynomials associated to $\Lambda$, where 
\be
\P_\Lambda:={\rm span}\{ y\mapsto y^\nu \; :Ê\; \nu\in\Lambda\},
\ee
is the corresponding space of real valued polynomials. Therefore
a polynomial approximation of the form \iref{polynomial} belongs to $V_{\Lambda_n}$.

There are two main issues in 
the design of these methods : 
\begin{itemize}
\item[{\rm (i)}]
Given 
an index set $\Lambda_n$, how do we construct the polynomial approximation
\iref{polynomial}.
\item[{\rm (ii)}]
How do we select the index sets $\Lambda_n$.
\end{itemize}
For treating both of these issues, it is very useful to impose that the considered sets 
$\Lambda_n$ are {\it downward closed}, which we assume in going further.

Concerning the first issue, we present two different strategies 
which illustrate the important distinction between {\it non-intrusive}
and {\it intrusive} methods mentionned in the introduction. 

The first strategy, discussed in \S 5,  is non intusive.
It computes a polynomial approximation of the form
\iref{polynomial} by {\it interpolation}
of the solution map 
at well chosen points $y^1,\dots,y^n\in U$
by a method introduced in \cite{CCS}, in the line of \cite{NobTemWeb08a,NobTemWeb08b}. In particular
it could even be applied in a context where the exact model is not known, but only
the solver is given. Other important representatives 
of non-intrusive methods, which we do not discuss in this paper, include {\it least-square}Ê projection methods
as developed in \cite{CCMNT,DI,DO}, and {\it pseudo-spectral} methods
as developped in \cite{CEP2,Xiu}.

The second strategy, discussed in \S 6, performs
an explicit computation of the truncated Taylor series, up to the spatial
discretization of the coefficients $t_\nu$, by a recursive method 
introduced in \cite{CCDS}. 
In contrast to the previous one, this approach is intrusive. It strongly
exploits the particular form of the parametric PDE, and actually 
it can only be easily implemented for parametric problems \iref{genpar} where $\cP$ is linear 
both in $u$  and $a$. Other important representatives of intrusive methods, which we do not discuss in this paper, include
Galerkin projection methods as developed in \cite{BNTT1,BNTT2,CDS1,Git}.

Concerning the second issue, an important distinction should be made between
{\it non-adaptive}  and {\it adaptive}Ê methods. 
In non-adaptive methods, the selection of the
set $\Lambda_n$  for a given value of $n$ is done in an a priori manner, 
based on available information on the problem. Ideally we would like
to use the set $\Lambda_n$ associated to the $n$ largest 
coefficients in the Taylor or Legendre expansion,
however this set cannot be easily identified. Instead,
we consider the set $\Lambda_n$ 
associated to the $n$ largest a priori estimates obtained in \S 3
for the $V$-norms of these coefficients. We detail further in \S \ref{secapriori} the 
algorithmic construction of 
the sequence $(\Lambda_n)_{n\geq 1}$ by this approach.

In adaptive methods, the selection of $\Lambda_n$ is made in an a posteriori
manner, based on the computation for downward closed values of $n$, for instance
using the knowledge of both the previous choice $\Lambda_{n-1}$ and 
the computed approximation polynomial $u_{n-1}$ for this choice.
One reason why adaptive methods might perform significantly better
than their above described non-adaptive counterpart in the present context 
is because the a priori bound $e_\nu$ may lack sharpness
and therefore only gives a limited indication on the real set of the $n$ largest
coefficients. In particular, the guaranteed rate $n^{-s}$ based on these
a priori bounds may be too pessimistic, and a better rate could be obtained
using an adaptive method. However, the convergence analysis of adaptive methods is 
usually much more delicate than that of their non-adaptive counterparts. We give examples of adaptive strategies
both for interpolation in \S 5 and Taylor approximations in \S 6,
convergence analysis being available only for the latter.

\subsection{Non-adaptive constructions of the sets $\Lambda_n$}
\label{secapriori}

We recall that the a priori estimates obtained in \S 3 for the Taylor or Legendre have the 
following general form:
\begin{itemize}
\item
For the Taylor coefficients, under the assumptions of Theorem \ref{theoneiacU},
\be
\|t_\nu\|_V\leq Ce_\nu, \quad   e_\nu:=\prod_{j\in {\rm supp}(\nu)} \rho_j^{-\nu_j},
\label{tayestimate}
\ee
for any given sequence $\rho=(\rho_j)_{j\geq 1}$ of number larger than $1$ that
satisfies the constraint \iref{constrhoeps}.
\item
For the Legendre coefficients, under the assumptions of Theorem \ref{theoneiaU},
\be
\|w_\nu\|_V\leq Ce_\nu, \quad e_\nu:=C\prod_{j\in {\rm supp}(\nu)} \theta (\rho_j) (1+2\nu_j)\rho_j^{-\nu_j},
\label{legestimate}
\ee
for any given sequence $\rho=(\rho_j)_{j\geq 1}$ of number strictly larger than $1$ that
satisfies the constraint \iref{constrhoeps}.
\end{itemize}
Also recall that for certain specific
problems, we can sharpen these estimates by improving on the constraint \iref{constrhoeps} imposed on $\rho$,
see \S \ref{secrefined}. Once an admissible sequence $\rho=\rho(\nu)$ has been fixed for each $\nu$,
each resulting estimate $e_\nu$ is computable as a product of $\|\nu\|_0$ numbers.
In the proof of Theorem \ref{theomain}, we use particular choices
of admissible sequences $\rho=\rho(\nu)$ which ensures the $\ell^p$ summability of the resulting $e_\nu$
provided that $(\|\psi_j\|_X)_{j\geq 1}$ is $\ell^p$ summable, for some $p<1$. 
However, one may hope to further improve the estimate $e_\nu$ by using other sequences.

An important observation is that the above general definition of $e_\nu$
does not guarantee that the set $\Lambda_n$ corresponding 
to the $n$ largest $e_\nu$ is downward closed. Indeed,
we are not ensured that the sequence $(e_\nu)_{\nu\in\cF}$ defined in 
\iref{tayestimate} or \iref{legestimate} is monotone non-increasing, in particular 
due to the fact that the sequence $\rho$ is allowed to vary with $\nu$.
One may try to construct the sequences $\rho(\nu)$ such that
the sequence $(e_\nu)_{\nu\in\cF}$ is monotone non-increasing. However, a 
simpler possibility is to search instead for a surrogate $s_\nu$, with
\be
e_\nu \leq  s_\nu:=\prod_{j\in {\rm supp}(\nu)} s_j(\nu),
\ee
where the $s_j(\nu)$ are again explicitly given, and in addition
$(s_\nu)_{\nu\in\cF}$ is a monotone non-increasing sequence. Then, we know that
at least one of the sets
$\Lambda_n$ corresponding to the $n$ largest $s_\nu$ is downward closed.
One example of such a surrogate
in the case of Legendre coefficients is given by $s_\nu:=\t r_\nu$ defined in \iref{trnu},
for which $\ell^p$ summability is also established
provided that $(\|\psi_j\|_X)_{j\geq 1}$ is $\ell^p$ summable, for some $p<1$. 

We now discuss the complexity of identifying the downward closed set $\Lambda_n$
associated to the $n$ largest $s_\nu$. In addition to the monotonicity of $(s_\nu)_{\nu\in\cF}$, 
the following property is useful for limiting this complexity.

\begin{definition}
\label{defanchoredseq}
A monotone non-increasing positive sequence $(s_\nu)_{\nu\in\cF}$ is said to be {\em anchored} if and only if
\be
l\leq j \Rightarrow s_{e_j} \leq s_{e_l},
\ee
where $e_l$ and $e_j$ are the Kroenecker sequences with $1$ at position $l$ and $j$, respectively.
\end{definition}

This property implies that at least one of the sets
$\Lambda_n$ corresponding to the $n$ largest $s_\nu$ has the following property.

\begin{definition}
\label{defanchoredset}
A finite downward closed set $\Lambda$ is said to be {\em anchored} if and only if
\be
e_j\in \Lambda \quad {\rm and}\quad l\leq j \quad \Rightarrow \quad e_l\in \Lambda.
\ee
where $e_l$ and $e_j$ are the Kroenecker sequences with $1$ at position $l$ and $j$, respectively.
\end{definition}

We now show that for an anchored sequence $(s_\nu)_{\nu\in\cF}$ the identification of the set $\Lambda_n$
can be executed in at most $n^2/2$ evaluations of $s_\nu$. For this purpose, we introduce for any downward closed set
$\Lambda$ its set of {\it neighbors} defined by
\be
N(\Lambda):=\{\nu\notin \Lambda \mbox{ such that } \Lambda\cup\{\nu\} \mbox { is downward closed}Ê\},
\label{neighbors}
\ee
We also  intoduce the set of its {\it anchored neighbors} defined by
\be
\t N(\Lambda):=\{\nu\in N(\Lambda)\; : \; 
\nu_j=0\;\; {\rm if}\;\; j >j(\Lambda)+1\},
\label{redneighbors}
\ee
where 
\be
j(\Lambda):=\max\{j\; : \; \nu_j>0\; \mbox{for some}\; \nu\in \Lambda\}.
\ee 
If $(s_\nu)_{\nu\in\cF}$ is an
anchored sequence, we may define the sets $\Lambda_n=\{\nu^1,\dots,\nu^n\}$ by the following induction:
\begin{itemize}
\item
Take $\nu^1=0$ the null multi-index.
\item
Given $\Lambda_k=\{\nu^1,\dots,\nu^k\}$, pick a $\nu^{k+1}$ maximizing
$s_\nu$ over $\nu\in \t N(\Lambda_k)$ and such that the new set $\Lambda_{k+1}$ is anchored.
\end{itemize}
We observe that $\t N(\Lambda_k)$ is contained in the union of $\t N(\Lambda_{k-1})$ and
of the set consisting of the indices
\be
e_{j(\Lambda_{k})+1} \quad {\rm and}Ê\quad \nu^k+e_j, \quad j\leq j(\Lambda_{k}).
\ee
Therefore, since the values of the $s_\nu$ have already been computed for $\nu\in \t N(\Lambda_{k-1})$, 
the step $k$ of the induction requires at most $j(\Lambda_{k})+1$ evaluations of $s_\nu$. In addition,
the fact that $\Lambda_{k}$ is anchored implies that $j(\Lambda_{k}) \leq k-1$. Therefore, the total number
of evaluations of $s_\nu$ in order to reach $\Lambda_n$ is at most
\be
N_n=1+2+\dots +(n-1)\leq n^2/2.
\ee
Finally, let us observe that the computation of a single $s_\nu$ costs
$\|\nu\|_0$ multiplications, and on the other hand, for all $\nu\in\Lambda_n$,
\be
2^{\|\nu\|_0} \leq \prod_{j\in{\rm supp}(\nu)} (1+\nu_j)
\leq \#\{\t \nu\; : \; \t \nu\leq \nu\} \leq \#(\Lambda_n)=n,
\ee
since $\Lambda_n$ is downward closed.
The total cost of identifying $\Lambda_n$ is therefore at most of the order $n^2\log (n)$
which is generally negligible compared to the computation of the 
approximation polynomial. Indeed, the latter involves $n$ elements from the space $V_h$,
and has therefore complexity at least $nN_h$ which, in the practice of reduced modeling, is much larger than 
$n^2\log(n)$ since $N_h\gg n$.

\subsection{Reduced basis methods}

A second class of numerical methods is not based on polynomial 
approximations. Instead, it directly seeks  choices of functions $v_1,\dots,v_n$ 
for which the approximation of the parametric PDE
in the resulting $n$-dimensional space $V_n:={\rm span}\{v_1,\dots,v_n\}$ 
performs almost as good as the optimal benchmarks for separable approximations.
Recall that these benchmarks are
measured by $n$-width $d_n(\cM)_V$ for the uniform error, or by
the tail of the singular values \iref{tailsing} for the mean-square error.

We discuss in \S 7 the {\it reduced basis} method which targets 
uniform error estimates, and which consists in generating $V_n$ by
a selection of $n$ particular solution instances $u(a^i)$ for $i=1,\dots,n$,
chosen from a very large set of potential candidates.  The selection
process is critical for the success of this algorithm, and one main
result is that a certain greedy strategy meets
the benchmark of the $n$-width in the sense that it results in similar convergence rates.

Another representative of this second class of methods, which we do not discuss 
in this paper, is known as the {\it proper orthogonal decomposition} method
and targets mean square estimates. It builds the functions $\{v_1,\dots,v_n\}$ 
based on an empirical approximation
of the exact covariance operator \iref{covop} using a sufficiently 
dense sampling of the random 
solution $u(a)$.  

One main disadvantage of both 
reduced basis and proper orthogonal decomposition methods, compared to the
first class of methods based on polynomial approximation, is that their offline stage is 
potentially very costly, especially in high parameter dimension. However,
their potential gain is in that they can get significantly
closer to the optimal benchmarks for separable approximations.
This is due to the fact that the best $n$-term polynomial approximation error
may in some cases decay substantially slower than the $n$-width, as discussed in
\S \ref{subrank}.
 
\section{Sparse polynomial interpolation}  

In this section, we discuss the construction of polynomial approximations
to the solution map $y\mapsto u(y)$
by {\it interpolation}. We place ourselves in
the same framework as in \S 3: 
we consider a parametric problem of the form \iref{genpar}
such that {\bf Assumption A} holds 
for a suitable affine representer $(\psi_j)_{j\geq 1}$,
so that the solution map $y\mapsto u(y):=u(a(y))$ 
is then well defined from $U$ to $V$.

Given $\Lambda \subset \cF$  with  $\#(\Lambda)=n$,
we say that a discrete set 
\be
\Gamma\subset  U, \quad \#(\Gamma)=n
\ee
is {\it unisolvent} for $\P_\Lambda$ if and only if for any values $(v_\gamma)_{\gamma\in\Gamma}\in\R^\Gamma$,
there exists a unique polynomial $\pi\in \P_\Lambda$ such that
\be
\pi(\gamma)=v_\gamma,\quad \gamma\in\Gamma.
\ee
In such a case, to any real valued function $v$ defined over $U$, we associate its
interpolation polynomial $I_\Lambda v\in \P_\Lambda$ which satisfies
\be
I_\Lambda v(\gamma)=v(\gamma), \quad \gamma\in\Gamma.
\ee
The interpolation operator $I_\Lambda$ is a linear map from the space of 
real valued functions defined over $U$ onto $\P_\Lambda$. It may be written in
the usual Lagrange form
\be
I_\Lambda v=\sum_{\gamma\in \Gamma} v(\gamma) \ell_{\Lambda,\gamma},
\label{lagrange}
\ee
where the $\ell_{\Lambda,\gamma}\in \P_\Lambda$ are
uniquely defined by $\ell_{\Lambda,\gamma}(\t\gamma)=\delta_{\gamma,\t\gamma}$ for $\gamma,\t\gamma\in \Gamma$.

By a standard vectorization procedure, we may define a similar 
interpolation process that maps the space of $V$-valued functions defined 
on $U$ onto the space $V_\Lambda$. This amounts in now using
the $V$-valued $v(\gamma)$ in the definition of the interpolant by \iref{lagrange}.
With a slight abuse of notation, we again denote by $I_\Lambda$ this operator.
From exactly or approximately computed instances
\be
u_\gamma=u(\gamma), \quad \gamma\in\Gamma,
\ee
of the solution map, we may thus compute $I_\Lambda u\in V_\Lambda$ such that
\be
I_\Lambda u(\gamma)=u_\gamma,\quad \gamma\in\Gamma.
\ee
One of the main attractions of interpolation, also sometimes refered to as
{\it collocation}Ê in the context of parametric PDEs
\cite{BabNobTem07,NobTemWeb08a,NobTemWeb08b} is that it is 
a non-intrusive process.   

In addition to the existence and uniqueness of the interpolation polynomial,  we point out two
other properties of the interpolation process that are of interest to us:
\begin{enumerate}
\item
Stability: one typical way of quantifying the stability of the interpolation process
is through its  {\it Lebesgue constant}.  If $\Gamma
\subset U$  is a set of unisolvent interpolation points
for $\P_\Lambda$ with  Lagrange basis elements $\ell_{\Lambda,\gamma}$, the  Lebesgue constant is defined as
\be
\label{Lebesgue}
\L_\Lambda:=\sup \frac {\|I_\Lambda u\|_{L^\infty(U)}}{\|u\|_{L^\infty(U)}}= \max_{y\in U} \sum_{\gamma\in\Gamma} |\ell_{\Lambda,\gamma}(y)|.
\ee
where the first supremum is taken over all non-zero real valued functions $u$
which are everywhere defined and uniformly bounded over $U$. It is easily seen that we obtain the same quantity if we instead
take the supremum over the set of $V$-valued functions, using the $L^\infty(U,V)$ norm in the quotient.
The Lebesgue constants typically grow
with the number $n$ of interpolation points, however it is well known
that this growth strongly depends on the selection of points. For instance,
on the univariate interval $[-1,1]$, the Lebesgue constant for interpolation
by polynomials of degree $n-1$ at $n$ points grows exponentially 
with $n$ for uniformly spaced points and logarithmically for Chebychev
or Gauss-Lobatto points.

\item
Progressivity: we would like to use sequences  $(\Lambda_n)_{n\geq 1}$ of index sets 
which have the nestedness property $\Lambda_n\subset \Lambda_{n+1}$
in order to define polynomial spaces with increasing accuracy. The sets
$\Lambda_n$ may be defined  a priori, based on the analysis of best $n$-term polynomial 
approximations presented in \S 3, or adaptively generated. In both cases,
it is desirable that the polynomial interpolation operators $I_{\Lambda_{n+1}}$
can be derived in a simple way from $I_{\Lambda_{n+1}}$. This requires
in particular that the associated unisolvent sets of points $(\Gamma_n)_{n\geq 0}$
are themselves nested.
\end{enumerate}

It was shown in \cite{CCS} that such progressive interpolation
processes can be derived provided that the sets $\Lambda_n$ are downward closed.
We present this approach in \S \ref{sparsint} and discuss its stability properties in \S \ref{stabint}.
We finally discuss in \S \ref{compint} the computational cost of such interpolation schemes, taking into account the space discretization
for the computation of the instances $u(\gamma)$, for example using a finite element method.

\subsection{Sparse interpolation using downward closed sets}
\label{sparsint}

We describe the construction of the interpolation operator for real valued functions,
since, as previously explained, it induces a similar interpolation
operator for $V$-valued functions.

We begin by discussing progressive constructions in the case of univariate polynomial interpolation.   The starting point is any sequence 
\be
T=(t_k)_{k\geq 0},
\ee
of   distinct points from $[-1,1]$.
We introduce the abbreviated notation
\be
I_k:=I_{\{t_0,\dots,t_k\}},
\ee
for the univariate interpolation operator associated with the $k$-section $\{t_0,\dots,t_k\}$ of this 
sequence: for any function $u$ defined everywhere over $[-1,1]$, the polynomial $I_ku\in\P_k$
satisfies
\be
I_ku(t_i)=u(t_i),\quad i=0,\dots,k.
\ee
We can express $I_k$ in a hierarchical form
\be
I_ku=I_0u+\sum_{l=1}^k\Delta_l u,\quad \Delta_l:=I_l-I_{l-1},
\label{newton}
\ee
also commonly known as the Newton form. We set $I_{-1}=0$ so that we
can also write
\be
I_ku=\sum_{l=0}^k \Delta_l u.
\ee
Since $I_ku$ and $I_{k-1}u$ agree at the points $\{t_0,\dots,t_{k-1}\}$, it is readily seen that, for $k>0$,
\be
\Delta_k u(t)=\alpha_k h_k(t),
\ee
where
\be
\alpha_k=\alpha_k(u):=u(t_k)-I_{k-1}u(t_k),
\ee
is the error at $t_k$ of interpolation by $I_{k-1}$, and
\be
h_{k}(t):=\prod_{l=0}^{k-1}\frac {t-t_l}{t_k-t_l}.
\ee
We also set
\be
h_0(t):=1.
\ee
For all $k\geq 0$, the system $\{h_0,\dots,h_k\}$ is a basis for $\P_k$,
sometimes called a {\it hierarchical basis}.

Although the sequence $T$ could be arbitrary,  the stability of the resulting interpolation
scheme, as reflected through the growth of it Lebesgue constants, depends very much
on the choice of $T$.
One interesting
choice is the sequence of the so-called {\it Leja points}, which is initiated
from an arbitrary $t_0$ (usually taken to be be $1$ or $0$) and recursively
defined by
\be
t_k:={\rm argmax}\Big\{ \prod_{l=0}^{k-1}|t-t_l| \; : \; t\in [-1,1]\Big \}.
\label{leja}
\ee
With this particular choice, we note that that the hierarchical basis functions satisfy
\be
\|h_k\|_{L^\infty([-1,1])}=1, \quad k\geq 0.
\ee
The Leja points may be viewed as an incremental variant to the classical Fekete points 
\be
\{t_{0,k},\dots,t_{k,k}\}:={\rm argmax}\Big\{ \prod_{i\neq j} |t_i-t_j| \; : \; \{t_0,\dots,t_k\}\in [-1,1]^{k+1} \Big \},
\ee
which, in contrast to the Leja points, are {\it not} $k$-sections of a single sequence.

We turn now to  the multivariate setting.   Starting again with the univariate sequence $T$,  we now define the points 
\be
y_\nu:=(t_{\nu_j})_{j\geq 1}\in U, \quad \nu\in\cF,
\label{multipoint}
\ee
which are therefore extracted from the tensorized grid $T^{\N}$. 
We also define the tensorized operators
\be
I_\nu:=\otimes_{j\geq 1} I_{\nu_j} 
\ee
Recall that the application of a tensorized operator $\otimes_{j\geq 1} A_j$
to a multivariate function amounts in applying each univariate operator $A_j$ by 
freezing all variables except the $j$-th and then applying $A_j$ to the non-frozen variable.
We may define $I_\nu$
by induction. For this, let us introduce $\cF_k$ the set of all $\nu$ such that $\nu_j=0$ for $j\geq k$.
\begin{itemize}
\item 
For $k=1$, there is only $\nu=0$ the null multi-index contained in $\cF_0$. Then
$I_0u$ is the constant function with value $u(y_0)$, where $y_0=(t_0,t_0,\dots)$.
\item
For $k>1$, assuming that $I_{\t \nu}$ has been defined for any $\t \nu\in \cF_{k-1}$,
and taking $\nu\in\cF_k$, we write
\be
\nu=(\nu_1,\t \nu),\quad \t\nu=(\nu_2,\nu_3,\dots)\in\cF_{k-1},
\ee
and for any $y\in U$,
\be
y=(y_1,\t y), \quad \t y=(y_2,y_3,\dots).
\ee
We then define $I_\nu:=I_{\nu_1}\otimes I_{\t \nu}$, that is,
\be
I_\nu u (y)= I_{\hat\nu} v_{y_1}(\t y),
\ee
where $v_{y_1}(\t y):=I_{\nu_1} u_{\t y}(y_1)$ with $u_{\t y}$ 
the univariate function defined on $[-1,1]$ by 
$u_{\t y}(t) = u(y)$.
\end{itemize}
Note that in the finite dimensional
case $U=[-1,1]^d$, this induction terminates after at most $d$ steps. 

It is easily seen that $I_\nu$ is the interpolation operator on
the tensor product polynomial space 
\be
\P_\nu=\otimes_{j\geq 1} \P_{\nu_j},
\ee
for the grid of points
\be
\Gamma_{\nu}=\otimes_{j\geq 1} \{t_{0},\dots,t_{\nu_j}\},
\ee
which is unisolvent for this space. This polynomial space
corresponds to a particular set $\Lambda$ 
which has rectangular shape. Namely $\Lambda=R_\nu$, where, for any $\nu\in\cF$, we define
the {\it shadow} of $\nu$ as
\be
R_\nu:=\{ \t \nu \; : \; \t \nu \leq \nu \}.
\label{shadow}
\ee
We thus have $\P_\nu=\P_{R_\nu}$.

We next define in a similar manner the tensorized 
difference operators
\be
\Delta_\nu:=\otimes_{j\geq 1}\Delta_{\nu_j}.
\ee
It follows that the range of $\Delta_\nu$ is the one dimensional subspace of $\P_\nu$ spanned by 
\be
\label{hnu}
H_\nu(y):= \prod_{\nu_j\neq 0}h_{\nu_j}(y_j),\quad \nu\in\cF.
\ee
To a general finite set $\Lambda\subset \cF$, we associate the operator
\be
I_\Lambda:=\sum_{\nu\in \Lambda} \Delta_\nu,
\label{newtongen}
\ee
and the grid
\be
\Gamma_\Lambda:=\{y_\nu \; : \; \nu\in \Lambda\}.
\ee
In the case where $\Lambda=R_\nu$, we find that $\Gamma_\Lambda=\Gamma_\nu$.
It is thus unisolvent for $\P_\Lambda$.
In addition, we then have
\be
I_\nu=\otimes_{j\geq 1} \(\sum_{l=0}^{\nu_j}\Delta_{l}\)=\sum_{\t \nu\leq \nu} \Delta_{\t\nu}=I_\Lambda,
\ee
which shows that $I_\Lambda$ is the interpolation operator onto $\P_\Lambda$ for this grid.

Let us remark that for a general set $\Lambda$,  the set $\Gamma_\Lambda$ is not  unisolvent on $\P_\Lambda$ and $I_\Lambda$ is not an interpolation operator.  However, an important observation is that this is the case  whenever $\Lambda$ is an arbitrary downward closed set.
This fact was first noticed in \cite{K} for bivariate functions, and then used in higher dimensions 
for particular cases of downward closed sets in \cite{Sm}.

\begin{theorem}
\label{Thmlowerint}
Let $\Lambda\subset \cF$ be a finite downward closed set. Then, the grid $\Gamma_\Lambda$
is unisolvent for $\P_\Lambda$ and $I_\Lambda$ is the interpolation operator onto $\P_\Lambda$ for this grid.
\end{theorem}

\noindent
{\bf Proof:} Because of the downward closed set property,  $P_\nu\subset P_\Lambda$ for all $\nu\in\Lambda$. Hence the image of $I_\Lambda$  is contained in $\P_\Lambda$.
In order to prove that it is the interpolation operator for the grid $\Gamma_\Lambda$, we
need to show that, for any function $u$ defined over $U$,
\be
I_\Lambda u(y_\nu) =u(y_\nu),\quad \nu\in \Lambda.
\label{ilambda}
\ee
Since $\#(\Gamma_\Lambda)=\dim(\P_\Lambda)$ this also ensures the unisolvence of $\Gamma_\Lambda$
for $P_\Lambda$. 

For any $\nu\in\Lambda$, we may write
\be
I_\Lambda u=I_{\nu}u+\sum_{\t\nu\in\Lambda, \t\nu \nleqslant \nu} \Delta_{\t \nu} u\;.
\ee
Since $I_\nu$ is the interpolant on the tensor product grid $\Gamma_\nu$, and this grid contains $y_\nu$,
it follows that 
\be
I_\nu u(y_\nu)=u(y_\nu).
\ee
On the other hand, if $\t\nu\in\Lambda$ is such that $\t\nu \nleqslant \nu$, this means that
there exists a $j\geq 0$ such that $\t \nu_j>\nu_j$. For this $j$ we thus have
$\Delta_{\t \nu} u(y)=0$ for all $y\in U$ with $j$-th coordinate  equal to $t_{\nu_j}$ due to the application
of $\Delta_{\nu_j}$ in the $j$-th variable.    Therefore 
\be
\Delta_{\t \nu} u(y_\nu)=0.
\ee
It follows that $I_\Lambda u(y_\nu)=u(y_\nu)$ which concludes the proof.
\hfill $\Box$
\nl

The decomposition \iref{newtongen}
of $I_\Lambda$ as a sum of  the various $\Delta_\nu$ may be viewed
as a generalization of the Newton form \iref{newton}. 
This decomposition also  yields a simple strategy for the fast computation
of $I_\Lambda u$ that we now describe. 

We first observe that if $\Lambda$ is a downward closed set of 
cardinality $n>0$, we can find at least one 
$\nu\in\Lambda$ which is maximal in $\Lambda$, that is, such that
\be
\t\nu\geq \nu\quad {\rm and}\quad \t\nu\neq \nu \quad \Rightarrow\quad \t\nu\notin \Lambda.
\ee
We may then write
\be
\Lambda=\t \Lambda\cup\{\nu\}.
\ee
where $\t \Lambda$ is a downward closed set of cardinality $n-1$.
Writing
\be
I_\Lambda u=I_{\t\Lambda} u+\Delta_\nu u,
\ee
we observe that $\Delta_\nu$ is characterized by the fact
that it belongs to $\P_\nu$ and is characterized by 
\be
\Delta_\nu u(y_{\t \nu})=0, \quad \t\nu\in \Gamma_\nu-\{\nu\},
\ee
and 
\be
\Delta_\nu u(y_\nu)=I_\Lambda u(y_{\nu})-I_{\t\Lambda} u(y_\nu)=u(y_\nu)-I_{\t\Lambda} u(y_\nu).
\ee
Using the tensorized hierarchical basis function $H_\nu$, 
it follows that
\be
\Delta_\nu u=\alpha_\nu H_\nu, \quad \alpha_\nu=\alpha_\nu(u):=u(y_\nu)-I_{\t\Lambda} u(y_\nu).
\ee
By iteration, we may write
$\Lambda=\{\nu^1,\dots,\nu^n\}$,
where the $\nu^i$ are enumerated in such way that for each $i$, 
$\Lambda_i=\{\nu^1,\dots, \nu^i\}$
is a downward closed set. This allows us to compute $I_\Lambda$ by
$n$ recursive applications of
\be
I_{\Lambda_{i}}u=I_{\Lambda_{i-1}}u +\alpha_{\nu^i} H_{\nu^i}.
\label{recurs}
\ee
Note that $(H_\nu)_{\nu\in \Lambda}$ is a basis of $\P_\Lambda$ and that
any $v\in \P_\Lambda$ has the unique decomposition
\be
v=\sum_{\nu\in \Lambda} \alpha_\nu H_\nu,
\ee
where the coefficients $\alpha_\nu=\alpha_\nu(v)$ are defined by the above procedure applied to $v$. Therefore, although the enumeration 
$\{\nu^1,\dots,\nu^n\}$ is not unique, the coefficients $\alpha_\nu=\alpha_\nu(u)$ in the expression 
\be
I_\Lambda u=\sum_{\nu\in \Lambda} \alpha_\nu H_\nu,
\ee
are unique. Also note that $\alpha_\nu(u)$ does not depend on the choice of $\Lambda$
but only on $\nu$ and $u$.
This computation is exactly the same in the case of $V$-valued functions,
now with uniquely defined coefficients $\alpha_\nu\in V$.

The recursive computation of the interpolation operator by \iref{recurs} can be used
in two different contexts:
\begin{itemize}
\item
Non-adaptive methods: a nested sequence $(\Lambda_n)_{n\geq 0}$
of downward closed sets is prescribed in advance, and we use \iref{recurs} to compute
$I_{\Lambda_n}u$ for increasing values of $n$.
\item
Adaptive methods: the sequence $(\Lambda_n)_{n\geq 0}$ is not prescribed in advance, and we use
the computation of $I_{\Lambda_n}u$ to define $\Lambda_{n+1}$.
\end{itemize}

We next give a typical example of an adaptive interpolation algorithm. In order to present this 
algorithm we begin by an analogy: since we have
\be
I_\Lambda u= \sum_{\nu\in \Lambda}\alpha_\nu H_\nu,
\label{hierexp}
\ee
we may view the interpolant as a truncation of the formal infinite
expansion of $u$ in the hierarchical basis
\be
\sum_{\nu\in \cF}\alpha_\nu H_\nu,
\ee
From elementary results on polynomial interpolation, we know that
this series does not converge for a general function defined everywhere over $U$.
Even for the various models of parametric PDEs discussed
in this paper, we don't know natural conditions that would ensure
the unconditional convergence of this expansion towards $u$,
in contrast to the Taylor and Legendre series discussed in \S 3.
In particular we do not know estimates for the coefficients $\alpha_\nu$ which would allow us to 
establish convergence rates for the best $n$-term truncations.

Nevertheless, we may still take the same view as in \S 3, and 
use for $\Lambda_n$ the set of 
indices corresponding to the $n$ largest terms of \eref{hierexp} measured in some given metric 
$L^p(U,V,\mu)$. 
We take $p=\infty$ if we search for uniform approximation estimates or $p=2$
if we search for mean-square approximation estimates.
This amounts to choosing the indices of the
$n$ largest $c_\nu \|\alpha_\nu\|_V$, where 
$c_\nu$ is given by
\be
c_\nu:=\|H_\nu\|_{L^p(U,\mu)}.
\label{cnu}
\ee
In the case where $\mu$ is the uniform measure, we
also have
\be
c_\nu:=\prod_{j\geq 1} \|h_{\nu_j}\|_{L^p([-1,1],\frac {dt} 2)},
\ee
Note that in the case where $p=\infty$ and if we use the Leja sequence,
we are ensured that $\|H_\nu\|_{L^\infty(U)}=1$ and therefore this amounts to choosing
the largest $\|\alpha_\nu\|_V$. The defect of this strategy is
that the sets $\Lambda_n$ are not ensured to be downward closed. In addition,
we generally cannot afford an exhaustive search for the
$n$ largest contributions in \eref{hierexp}.

In order to build a feasible adaptive algorithm, we need to limitate this
search. In what follows, we describe a greedy algorithm proposed in \cite{CCS} for the selection of the 
sequence $(\Lambda_n)_{n\geq 1}$, which uses the
set of neighbors $N(\Lambda)$ defined by \iref{neighbors}. We first give 
an idealized version of this algorithm which cannot be applied
as such.
\nl
\nl
{\bf Greedy Interpolation Algorithm:} We start with $\Lambda_1:=\{0\}$ the null multi-index.
Assuming that $\Lambda_{n-1}$ has been selected and that
the $(\alpha_\nu)_{\nu\in \Lambda_{n-1}}$ have been computed,
we compute the $\alpha_\nu$ for $\nu\in N(\Lambda_{n-1})$.   
We then set
 \be
\nu^n:={\rm argmax} \{ c_\nu\|\alpha_\nu\|_V\; : \; \nu\in N(\Lambda_{n-1})\},
\label{greedy}
\ee
and define $\Lambda_{n}=\Lambda_{n-1} \cup \{\nu^n\}$.
\nl
\nl
Note that when  $p=\infty$ and $T$ is  the Leja sequence, this strategy 
amounts to choosing the $\nu\in N(\Lambda_{n-1})$ 
that maximizes the interpolation error at the new grid point which would be added by adjoining $\nu$, that is, setting
\be
\nu^n:={\rm argmax}\{\|u(y_{\nu}) - I_{\Lambda_{n-1}}u(y_{\nu})\|_V
\; : \; \nu\in N(\Lambda_{n-1})\}.
\label{nextindex}
\ee
The above greedy algorithm is not computationally feasible since we are 
in principle working with infinitely many variable
$(y_j)_{j\geq 1}$, in which case the set of neighbours $N(\Lambda)$ to be
explored has infinite cardinality. One way to circumvent this defect
is to replace in the algorithm
the infinite set $N(\Lambda_n)$ by the finite 
set of anchored neighbors $\t N(\Lambda_n)$
defined by \iref{redneighbors}.

One more serious defect of this algorithm is that it may fail to converge, even
if there exist sequences $(\Lambda_n)_{n\geq 0}$ such that
$I_{\Lambda_n} u$ fastly converges towards $u$. Indeed, if it happens
that $\Delta_\nu u=0$ for a certain $\nu$, then no index $\t \nu\geq \nu$ will
ever be selected by the algorithm. As an example,
consider a two dimensional function of the form
\be
u(y)=u_1(y_1)u_2(y_2),
\ee 
where $u_1$ and $u_2$ are non-polynomial smooth functions 
such that $u_2(t_0)=u_2(t_1)$. Then the sets $\Lambda_n$ selected
by the algorithms will consists of the indices $\nu=(k,0)$ for $k=0,\dots,n-1$, since the
interpolation error at the point  $(t_k,t_1)$ always vanishes. One way to avoid
this problem is to change the strategy by alternating the selection
of $\nu^n$ using \iref{nextindex} with a selection rule ensuring that all indices
are eventually picked. For example, when $n$ is even, we define 
$\nu^{n}$ according to \iref{nextindex}, and when $n$ is odd we pick
for $\nu^n$ the multi-index $\nu\in\t N(\Lambda_n)$ which has appears at the earliest stage
in the neighbors of the previous sets $\Lambda_k$. In summary, this results in the following algorithm.
\nl
\nl
{\bf Alternating Greedy Interpolation Algorithm:} We start with $\Lambda_1:=\{0\}$ the null multi-index.
Assuming that $\Lambda_{n-1}$ has been selected and that
the $(\alpha_\nu)_{\nu\in \Lambda_{n-1}}$ have been computed,
we compute the $\alpha_\nu$ for $\nu\in \t N(\Lambda_{n-1})$. 
We set, if $n$ is even, 
 \be
\nu^n:={\rm argmax} \{ c_\nu\|\alpha_\nu\|_V\; : \; \nu\in \t N(\Lambda_{n-1})\},
\label{greedy}
\ee
and, if $n$ is odd,
\be
\nu^{n}:={\rm argmin} \{k(\nu)\; : \; \nu\in \t N(\Lambda_{n-1}\}, \quad k(\nu):=\min\{k\; : \; \nu\in \t N(\Lambda_k)\}.
\label{nun}
\ee
We then define $\Lambda_{n}=\Lambda_{n-1} \cup \{\nu^n\}$.
\nl

Even with such modifications, although the adaptive algorithm seem to 
behave well in many practical instances, the convergence of the 
interpolation produced by this algorithm
is still not guaranteed. It is an open problem to
understand which additional assumptions on $u$ ensure convergence,
and more importantly a convergence rate that is comparable to
that which is proved for best $n$-term approximations based on
Taylor and Legendre series. Note that the solution to this problem 
need to involve the initial choice of the univariate sequence $T$, which,
as discussed in the next section, strongly affects the stability and convergence
properties of the interpolation process. In the next section, using this stability
analysis, we establish convergence
rates for the interpolation algorithm, however based on non-adaptive choices of the
sequence $(\Lambda_n)_{n\geq 1}$.

\begin{remark}
A very similar greedy algorithm was proposed in {\rm \cite{GG}}
in the slightly different context of adaptive quadrature, that is,
when we want to approximate the integral of $u$ over the
domain $U$ rather than $u$ itself. In that case, one natural choice
is to pick the new neigbor $\nu$ that maximizes the absolute value
of the integral of $\Delta_\nu u$.
\end{remark}

Let us conclude this section by mentioning that there exists several
natural generalizations to the above described construction of the sparse multivariate
interpolation process. 

The first obvious one is that we could work on more general tensor product
domains of the form
\be
U=\otimes_{j\geq 0} U_j,
\ee
where the $U_j$ are univariate intervals or other bounded domains in $\R$ or $\C$, and define points $y_\nu$
by tensorization of sequences
\be
T_j=(t_{j,k})_{k\geq 0},
\ee
of pairwise distinct points, each of them picked from $U_j$. 

The second generalization is that we could
start with   univariate systems other  than polynomials that still having a hierarchical interpolation
structure. We consider a general index set $\cS$ equiped with a partial order $\leq$ 
and assume that there exists a root index $0$ such that $0\leq \gamma$ for all $\gamma\in \cS$.
Given a grid of pairwise
distinct points $G=(t_\gamma)_{\gamma\in\cS}$, we say that a family of functions
$(h_\gamma)_{\gamma\in\Gamma}$ defined over $[-1,1]$ is a hierarchical basis associated to
the grid $G$ if and only if $h_0(t)=1$ and
\be
h_\gamma(t_\gamma)=1\; \;{\rm and}\;\; h_\gamma(t_{\t\gamma})=0 \; {\rm if}\; \t \gamma \leq \gamma \; {\rm and} \; \t \gamma \neq \gamma.
\ee
By tensorization, we obtain an index set $\cF\subset G^{\N}$ of finitely supported sequences,
equiped with a partial order $\leq$ induced by its univariate counterpart. This allows us to 
define downward closed sets in $\cF$ in a
the  same way that we have  for the particular case $G=\N$. For $\nu=(\nu_j)_{j\geq 1}\in\cF$, we also define the points $y_\nu\in U$ and the tensorized hierarchical functions $H_\nu$ in the same
way as in \iref{multipoint} and \iref{hnu}. Then,
if $\Lambda$ is a downward closed set, we may inuctively define an interpolation operator $I_\Lambda$
onto the space 
\be
H_\Lambda:={\rm span}\{H_\nu\; : \; \nu\in\Lambda\},
\ee
associated to the grid $\Gamma_\Lambda$, using the same recursion
\be
I_\Lambda u=I_{\t \Lambda} u+\alpha_{\nu} H_\nu,\quad \alpha_\nu:=\alpha_\nu(u)=u(y_\nu)-I_{\t \Lambda} u(y_\nu),
\ee
where $\nu\notin\t \Lambda$ and $\t\Lambda$ is a downward closed set such that $\Lambda=\t\Lambda\cup\{\nu\}$. We initialize
this computation for $\Lambda=\{0\}$, where $0$ is the null multi-index, by defining $I_{\{0\}}u$ 
as the constant function with value $u(y_0)$. Examples of relevant hierarchical systems 
include the classical piecewise linear, or more generally piecewise polynomial,
hierarchical basis functions. With such choices the spaces $H_\Lambda$ include as
particular cases the well-studied piecewise polynomial {\it sparse grid} spaces, see \cite{BG} for a survey on this
topic.

\subsection{Stability}
\label{stabint}

We now turn to the stability analysis of the interpolation operator. We recall the Lebesgue constant
defined in \iref{Lebesgue}.
One principal interest of the Lebesgue constant is that it allows us to estimate the error
of interpolation in terms of the best polynomial approximation error in the $L^\infty$ norm.
Indeed, for any $u\in L^\infty(U,V)$ and any $v\in V_\Lambda$ we may write
\be
\|u-I_\Lambda u\|_{L^\infty(U,V)} \leq \|u-v\|_{L^\infty(U,V)}+\|I_\Lambda v-I_\Lambda u\|_{L^\infty(U,V)},
\ee
which yields
\be
\|u-I_\Lambda u\|_{L^\infty(U,V)} \leq (1+\L_\Lambda)\inf_{v\in V_\Lambda}\|u-v\|_{L^\infty(U,V)},
\ee
by taking the infimum over $V_\Lambda$.

We know from the results in \S \ref{sublower}, in particular Corollary
\ref{corratemonot}, that for relevant classes of parametric PDEs, we can find nested sequences
of downward closed sets $(\Lambda_n)_{n\geq 1}$ with $\#(\Lambda_n)=n$, such that
\be
\inf_{v\in V_{\Lambda_n}}\|u-v\|_{L^\infty(U,V)} \leq Cn^{-s},
\ee
where $s>0$ is some given rate.  This holds in particular
with $s:=\frac 1 p-1$ if the assumptions of Theorem \ref{theoneiaU} hold and 
if in addition $(\|\psi_j\|_X)_{j\geq 1}\in \ell^p(\cF)$.
Therefore we have the error bound
\be
\|u-I_{\Lambda_n} u\|_{L^\infty(U,V)} \leq C(1+\L_{\Lambda_n})n^{-s}.
\label{Intern}
\ee
if we use such sequences in the interpolation process. This motivates estimating the
growth of $\L_{\Lambda_n}$ with $n$.

In order to estimate $\L_\Lambda$, we introduce the
univariate Lebesgue constants
\be
\lambda_k:=\sup\frac {\|I_k u\|_{L^\infty([-1,1])}} {\|u\|_{L^\infty([-1,1])}},
\label{lebk}
\ee
where the supremum is taken over all non-zero real valued functions $u$ that are everywhere defined and uniformly bounded on $[-1,1]$.
We define an analogous quantity for the difference operator $\Delta_k$, namely
\be
\delta_k:=\sup\frac {\|\Delta_k u\|_{L^\infty([-1,1])}} {\|u\|_{L^\infty([-1,1])}},
\ee
and observe that 
\be
\delta_k\leq \lambda_{k-1}+\lambda_k, \quad k\geq 0, 
\ee
where we have set $\lambda_{-1}=0$. We introduce for each $\nu\in\cF$ the quantities
\be
\delta_\nu:=\sup \frac {\|\Delta_\nu u\|_{L^\infty(U)}}{\|u\|_{L^\infty(U)}},
\ee
so that we have, on the one hand 
\be
\delta_\nuÊ\leq \prod_{j\geq1} \delta_{\nu_j} \leq \prod_{j\geq 1} (\lambda_{\nu_j-1}+\lambda_{\nu_j}),
\label{deltanu}
\ee
and on the other hand
\be
\L_\Lambda\leq \sum_{\nu\in\Lambda} \delta_\nu.
\label{sumdelta}
\ee
The following result from \cite{CCS} gives an estimate on the
growth of $\L_{\Lambda}$ in terms of $\#(\Lambda)$, provided
that a similar estimate holds for the univariate Lebesgue constant $\lambda_k$
or for the quantity $\delta_k$.

\begin{theorem} 
\label {theoLeb}
If either one of the estimates
\be
\lambda_k \leq (k+1)^\theta ,\quad k\geq 0,
\label{lambdak}
\ee
or 
\be
\delta_k \leq (k+1)^\theta ,\quad k\geq 0,
\label{deltak}
\ee
holds for some $\theta \geq 1$, then the Lebesgue constant $\L_\Lambda$ satisfies
\be
\L_{\Lambda} \leq (\#(\Lambda))^{\theta+1}
\label{theta1}
\ee
for any downward closed set $\Lambda$.
\end{theorem}
{\bf Proof:} The case where \iref{deltak} holds is elementary since the first inequality in \iref{deltanu} yields
\be
\delta_\nu \leq \prod_{j\geq 1} (\nu_j+1)^\theta =\(\prod_{j\geq 1} (\nu_j+1)\)^{\theta} =(\#(R_\nu))^{\theta}\leq
(\#(\Lambda))^{\theta},
\ee
where we have used the fact that $R_\nu\subset \Lambda$ since $\Lambda$ is downward closed. 
Using \iref{sumdelta} we thus obtain \iref{theta1}.

For the case where \iref{lambdak} holds, we observe that 
\be
\lambda_{k}+\lambda_{k-1}
\leq 
(k+1)^\theta +k^\theta
\leq 
(2k+1)(k+1)^{\theta-1}.
\ee
The second inequality in \iref{deltanu} yields
$$
\begin{disarray}{ll}
\delta_\nu
&\leq 
\(\prod_{j\geq1} (\nu_j+1)\)^{\theta - 1}
\prod_{j\geq1} (2\nu_j+1)
\\
&=
(\#(R_\nu))^{\theta - 1} \prod_{j\geq1} (2\nu_j+1)
\\
&
\leq
(\#(\Lambda))^{\theta - 1} \prod_{j\geq1} (2\nu_j+1),
\\
\end{disarray}
$$
In order to establish \iref{theta1}, it thus suffices to prove that
$\sigma(\Lambda) \leq (\#\Lambda)^2$, 
where
\be
\sigma(\Lambda):= \sum_{\nu\in\Lambda} \prod_{j\geq1} (2\nu_j+1)\;.
\ee 
For this, we use induction on $n:=\#(\Lambda)$. 
For $n=1$ and $\Lambda=\{0\}$ the result obviously holds. 
Assuming that it holds for some $n\geq 1$, we 
consider a downward closed set $\Lambda$ of cardinality $n+1$.
We may assume without loss of generality that $\nu_1\neq 0$ for some $\nu\in\Lambda$, 
and denote by $K\geq 1$ the maximal
value attained by the coordinate $\nu_1$ when $\nu\in\Lambda$. 
For $0\leq k \leq K$, we define
\be
\Lambda_k := \{\hat\nu=(\nu_2,\nu_3,\dots) : (k,\hat\nu) \in \Lambda\}
\ee
Each of the set $\Lambda_k$ is downward closed and, since $K\geq1$, we have
$\#(\Lambda_k)< \#(\Lambda)$ for all $k=0,\dots,K$. The induction hypothesis implies
\be
\sigma(\Lambda) 
=
 \sum_{k=0}^K \sum_{\nu\in\Lambda_k} \prod_{j\geq1} (2\nu_j+1)= \sum_{k=0}^K (2k+1) \sigma(\Lambda_k)
\leq
 \sum_{k=0}^K (2k+1) (\#(\Lambda_k))^2.
\ee
Also, we have 
\be
 \Lambda_K
\subset
\cdots
 \subset\Lambda_1 
\subset
  \Lambda_0,
\ee
since for $k\geq1$,
$ \nu\in \Lambda_k \Rightarrow (k,\nu) \in \Lambda \Rightarrow (k-1,\nu) \in \Lambda \Rightarrow \nu \in \Lambda_{k-1}$. 
This implies
\be
k(\#(\Lambda_k))^2 
\leq 
\#(\Lambda_k)
\sum_{j=0}^{k-1}\#(\Lambda_j),
\ee
and therefore
\be
\sigma(\Lambda) 
\leq
 \sum_{k=0}^K (\#(\Lambda_k))^2
 +
2\sum_{k=0}^K 
\#(\Lambda_k)
\sum_{j=0}^{k-1}\#(\Lambda_j) =
\(\sum_{k=0}^K \#(\Lambda_k)\)^2
=
(\#(\Lambda))^2,
\ee
which concludes the proof. 
\hfill $\Box$

\begin{remark}
One noticable feature of the above result is that the bound on $\L_\Lambda$
only depends on the cardinality of $\Lambda$. In particular, it is 
independent of the number of variables, which can be infinite, as well as of the
shape of $\Lambda$.
\end{remark}

In view of the above result, we are therefore interested in choosing univariate sequences
$T=(t_k)_{k\geq 0}$ such that the Lebesgue constant $\lambda_k$
or the quantity $\delta_k$ have moderate algebraic growth
with $k$. It is well known that for particular sets of points such as the Chebychev points
\be
C_k:=\Big\{\cos\(\frac {2 l+1}{2k+2} \pi\)\; : \; l=0,\dots, k\Big\},
\ee
or the Gauss-Lobatto (or Clemshaw-Curtis) points
\be
G_k:=\Big\{\cos\(\frac {l}{k} \pi\)\; : \; l=0,\dots, k\Big\},
\ee
the Lebesgue constant has
logarithmic grows $\lambda_k\sim \log(k)$, therefore slower than algebraic.
However these points are not adapted to our construction since 
the sets $C_k$ and $G_k$ are not nested as $k$ grows, and therefore
are not the $k$-sections of a single sequence.

For the Leja points defined by \iref{leja},
numerical computations of $\lambda_k$ for the first $200$ values of $k$ indicates that the linear
bound 
\be
\lambda_k \leq (1+k),
\label{lineark}
\ee
seems to hold and that this bound could be sharp. However there is currently
no rigorous proof supporting this evidence or establishing 
another algebraic rate. Nevertheless, Leja points seem to be a good
choice for the construction of our multivariate interpolation process.

Leja points have also been considered
on the complex unit disc $\{|z|\leq 1\}$, taking for example $t_0=1$ and using again
the recursion \iref{leja}, now with $|\cdot|$ standing for the modulus.
These points have the property of accumulating in a regular manner on the 
unit circle according to the so-called Van der Corput enumeration \cite{CP1}. 
Their projections on the real axis are called the $\Re$-Leja points, and coincide 
with the Gauss-Lobatto points for values of $k$ of the form $2^n+1$ for $n\geq 0$.
The growth of the Lebesgue constant $\lambda_k$ has been studied 
in \cite{CP1,CP2,C,CC} for these two families of points.
In the case of the complex Leja points, this constant 
is defined as in \iref{lebk}, however taking the supremum over
functions defined everywhere and bounded over the complex unit disc.
It is proved in \cite{C} that the linear bound \iref{lineark}
holds for the complex Leja points. For the $\Re$-Leja points, quadratic
bounds of the type
\be
\lambda_k \leq C(1+k)^2 \quad {\rm and}\quad \delta_k \leq (1+k)^2,
\ee
with $C>1$ are established in \cite{CC}. 

With such estimates, application of Theorem \ref{theoLeb}
gives us bounds of the form
\be
\L_\Lambda\leq (\#(\Lambda))^{1+\theta},
\ee
for example with $\theta=2$ when using the $\Re$-Leja points. 
If we combine this bound with \iref{Intern}, we obtain the convergence
estimate
\be
\|u-I_{\Lambda_n} u\|_{L^\infty(U,V)} \leq C n^{-(s-1-\theta)},
\ee
which expresses a deterioration of the convergence rate 
when using the interpolation process instead of the truncated 
expansions studied in \S 3.

We now present a sharper analysis, introduced in \cite{CCS}, 
which reveals that this deterioration actually 
does not occur for the models of parametric PDEs
which are of interest to us. This analysis is based on the following Lemma
which gives an estimate of the interpolation error in terms
of the tail of the Legendre coefficients of $u$ multiplied by
algebraic factors.

\begin{lemma} 
\label{lemmaresid}
Assume that the Legendre
expansion \iref{renormlegendre} of $u$ is
  unconditionally  convergent in $L^\infty(U,V)$.
If the univariate sequence $T=(t_k)_{k\geq 0}$ is such that
that \iref{lambdak} or \iref{deltak} holds for some $\theta\geq 1$, 
then, for any downward closed set $\Lambda$,
\be
\| u-I_\Lambda u \|_{L^\infty(U,V)} 
\leq 
2\sum_{\nu\notin \Lambda} p_\nu(b) \|w_\nu\|_V\;,
\label{residual}
\ee
where $b:=\theta+1$ and 
\be\label{eq:defpnu}
p_\nu(b) := \prod_{j\geq 1} (1+\nu_j)^{b}\;.
\ee
\end{lemma}

\noindent
{\bf Proof:}
The unconditional convergence of the Legendre series allows us to write
\be
I_\Lambda u 
=
I_\Lambda\(\sum_{\nu\in\cF} w_\nu P_\nu\)
=
 \sum_{\nu\in\cF} w_\nu I_\Lambda P_\nu
=
\sum_{\nu\in\Lambda} w_\nu P_\nu
+
\sum_{\nu\notin\Lambda} w_\nu I_{\Lambda}P_\nu,
\ee
where we have used that $P_\nu\in \P_\Lambda$ because $\Lambda$ is downward closed and hence
$I_\Lambda P_\nu =P_\nu$ for every $\nu\in\Lambda$.
For the second term, we observe that for each $\nu\notin\Lambda$,
\be
I_{\Lambda}P_\nu=\sum_{\t\nu\in \Lambda}\Delta_{\t\nu} P_\nu
=\sum_{\t\nu\in \Lambda\cap R_\nu}\Delta_{\t\nu} P_\nu=I_{\Lambda\cap R_\nu} P_\nu,
\ee
since $\Delta_{\t\nu} v = 0$ whenever $\t\nu \not\leq \nu$ and $v\in \P_\nu$.
Therefore 
\be
u-I_\Lambda u  
= \sum_{\nu\not\in\Lambda} w_\nu (I- I_{\Lambda\cap\cR_\nu})P_\nu,
\ee
where $I$ stands for the identity operator. 
This implies
\be
\|u-I_\Lambda u \|_{L^\infty(U,V)}
\leq
\sum_{\nu\not\in\Lambda} (1+\L_{\Lambda\cap R_\nu})\|w_\nu\|_V 
\leq 
2\sum_{\nu\not\in\Lambda} \L_{\Lambda\cap R_\nu}\|w_\nu\|_V \;.
\label{sumbound}
\ee
Since \iref{lambdak} or \iref{deltak} holds, we obtain from Theorem \ref{theoLeb} that
\be
\L_{\Lambda\cap R_\nu} 
\leq \#(\Lambda\cap R_\nu)^{\theta+1} 
\leq \#(R_\nu)^{\theta+1} 
= p_\nu(b),
\ee
which completes the proof.
\hfill $\Box$
\nl

The estimate for the interpolation error in Lemma \ref{lemmaresid}
is very similar to that of the trunctated Legendre expansion, up to the
presence of the factor $p_\nu(b)$. We may therefore use the same
techniques as those used in \S 3 for this expansion in 
order to establish convergence rates for the interpolation error.
We first establish a summability result which is 
analogous to Theorem \ref{theomainlower}.

\begin{theorem}
\label{theointerlower}
Consider a parametric problem of the form \iref{genpar}
such that {\bf Assumption A} holds 
for a suitable affine representer $(\psi_j)_{j\geq 1}$.
If the assumptions of Theorem \ref{theoneiaU} are satisfied,
and if in addition $(\|\psi_j\|_X)_{j\geq 1}\in \ell^p(\N)$ for some $p<1$, 
then $(p_\nu(b)\|w_\nu\|_V)_{\nu\in\cF}\in \ell^p_m(\cF)$
for the same value of $p$, and for any $b\geq 0$.
\end{theorem}

\noindent
{\bf Proof:} This proof is very similar to the proof of Theorem \ref{theomainlower} and so 
we only sketch it. We obtain a similar estimate
\be
p_\nu(b) \|w_\nu\|_V \leq \t r_\nu:=\t r_E(\nu)\t r_F(\nu),
\ee
where $\t r_E$ has exactly the same form as in \iref{tre} up to a change in the multiplicative
constant $C_0$, and $\t r_F$ is now given by
\be
\t r_F(\nu)Ê:=\prod_{j\in F\cap{\rm supp}(\nu)}c_\kappa(1+2\nu_j)^{1+b}\(\beta+\frac {\e \nu_j}{3b_j |\nu_F|}\)^{-\nu_j}
\ee
By a similar reasoning, up to a modification in the choice of $J$ and $\beta$,
one then shows that the sequence $(\t r_\nu)_{\nu\in\cF}$ belongs to $\ell^p(\cF)$
and that it is monotone non-increasing. \hfill $\Box$
\nl

Combining this result with Lemma \ref{lemmaresid}, we obtain the 
following corollary which shows that the interpolation process converges
without deterioration of the rate established for the Legendre series.

\begin{cor}
\label{corinterlower}
Consider a parametric problem of the form \iref{genpar}
such that {\bf Assumption A} holds 
for a suitable affine representation \iref{affine}.
Assume that the univariate sequence $T=(t_k)_{k\geq 0}$ is such that
that \iref{lambdak} or \iref{deltak} holds for some $\theta\geq 1$.
Then, if the assumptions of Theorem \ref{theoneiaU} are satisfied,
and if in addition $(\|\psi_j\|_X)_{j\geq 1}\in \ell^p(\N)$ for some $p<1$, 
there exists a sequence of 
nested downward closed sets $(\Lambda_n)_{n\geq 1}$ such that $\#(\Lambda_n)=n$
and such that
\be
\|u-I_{\Lambda_n}u\|_{L^\infty(U,V)} \leq Cn^{-s}, \quad n\geq 1, \quad s:=\frac 1 p-1.
\ee
\end{cor}

\subsection{Space discretization and computational cost}
\label{compint}

In the practice of numerical computation to parametric PDEs, as explained in \S 5, we replace
the instances $u(y_\nu)$ by their approximation
$u_h(y_\nu)\in V_h$ computed by a numerical solver. We denote by 
$I_{\Lambda,h} u$ the resulting interpolation polynomial, which belongs to the space
\be
V_{\Lambda,h}:=\P_\Lambda\otimes V_h.
\ee
In other words, we have
\be
I_{\Lambda,h} u:=I_\Lambda u_h,
\ee
where 
\be
u_h:y\mapsto u_h(y),
\ee
is the approximate solution map acting from $U$ into $V_h$.

We begin by discussing the accuracy of this polynomial approximation.
For a given set $\Lambda_n$, one way to estimate the total interpolation error 
is by writing
\be
\|u-I_{\Lambda_n,h} u\|_{L^\infty(U,V)} \leq  \|u-I_{\Lambda_n} u\|_{L^\infty(U,V)}
+\|I_{\Lambda_n} (u-u_h)\|_{L^\infty(U,V)}.
\ee
The first term is estimated by the results in the previous section, such as Corollary
\ref{corinterlower} and for the second term we may write
\be
\|I_{\Lambda_n} (u-u_h)\|_{L^\infty(U,V)}\leq \L_{\Lambda_n} \e(h),
\ee
where $\e(h)$ is the acuracy of the numerical solver. Under the assumptions of 
Corollary \ref{corinterlower}, this result in an error estimate of the form
\be
\|u-I_{\Lambda_n,h} u\|_{L^\infty(U,V)} \leq  Cn^{-s}+ n^{\theta+1} \e(h),
\label{neps}
\ee
where we have also used Theorem \ref{theoLeb} for bounding $\L_{\Lambda_n}$.

There is  a more efficient way to estimate the error in the case
where the approximate solution map $u_h$ may be viewed as
the solution map of a 
discrete parametrized problem of the form
\be
\cP_h(u_h,a)=0,
\label{genparh}
\ee
where $\cP_h: V_h\times X \to W$ and with similar properties as the original parametric problem \iref{genpar}.
Consider for example the case of the elliptic equation \iref{ellip} 
and its Galerkin discretization on $V_h$ defined by \iref{galerellip}.
It is then readily seen that the discrete solution map $a\mapsto u_h(a)$ has the same
boundedness and holomorphy properties
as the original map $a\mapsto u(a)$. In turn, we obtain the same
estimates for the Taylor or Legendre coefficients
of the solution map $y\mapsto u_h(y):=u(a(y))$. This type of problem 
is therefore covered by the following discrete counterpart to Corollary \ref{corinterlower}. 

\begin{cor}
\label{corinterlowerh}
Consider a parametric problem of the form \iref{genparh}
such that {\bf Assumption A} holds 
for a suitable affine representer $(\psi_j)_{j\geq 1}$.
Assume that the univariate sequence $T=(t_k)_{k\geq 0}$ is such that
that \iref{lambdak} or \iref{deltak} holds for some $\theta\geq 1$.
Then, if the assumptions of Theorem \ref{theoneiaU} are satisfied by
the map $a\mapsto u_h(a)$, with the open set $\cO$ and the bound in \iref{uaC1}
independent of $h$,
and if in addition $(\|\psi_j\|_X)_{j\geq 1}\in \ell^p(\N)$ for some $p<1$, 
there exists a sequence of 
nested downward closed sets $(\Lambda_n)_{n\geq 1}$ such that $\#(\Lambda_n)=n$
and such that
\be
\|u_h-I_{\Lambda_n}u_h\|_{L^\infty(U,V_h)} \leq Cn^{-s}, \quad n\geq 1, \quad s:=\frac 1 p-1,
\ee
where the constant $C$ is independent of $h$.
\end{cor}

Under the assumptions of the above corollary, we may now estimate the total interpolation error 
 by writing
\be
\|u-I_{\Lambda_n,h} u\|_{L^\infty(U,V)} \leq  
\|u_h-I_{\Lambda_n} u_h\|_{L^\infty(U,V)} +\|u- u_h\|_{L^\infty(U,V)},
\label{secondsplit}
\ee
which yields the bound
\be
\|u-I_{\Lambda_n,h} u\|_{L^\infty(U,V)}\leq Cn^{-s} +\e(h).
\label{nepsh}
\ee
This bound is clearly more favorable than \iref{neps}.

We next turn to the estimate of the computational cost, starting with the 
offline cost. Assuming that $\Lambda_n$ is given, we want to pre-compute the
coefficients $\alpha_{\nu,h}:=\alpha_\nu(u_h)\in V_h$ in the expression
\be
I_{\Lambda_n,h} u=\sum_{\nu\in \Lambda_n} \alpha_{\nu,h}H_\nu.
\ee
We begin by
computing the discretized instances $u_h(y_\nu)$ for $\nu\in \Lambda_n$, as vectors of size $N_h$ of their
coefficients in a given basis of $V_h$.
This has a cost of order $nC_h$ where $C_h$ is the individual cost of one application of
the discrete solver. The coefficients $\alpha_{\nu,h}$ are then computed by 
recursive application of the formula \iref{recurs} based on these discretized instances.
At the stage $i$ of this recursion, the coefficient $\alpha_ {\nu^i,h}:=\alpha_{\nu^i}(u_h)$ 
is computed by a linear combination of the $i-1$ previous ones, according to
\be
\alpha_ {\nu^i,h}=u_h(y_{\nu^i})-I_{\Lambda_{i-1}}u_h(y_{\nu^i})
=u_h(y_{\nu^i})-\sum_{l=1}^{i-1} \alpha_{\nu^l,h}H_{\nu^l}(y_{\nu^i}).
\ee
Note that in view of the definition
of the $H_\nu$, we can compute $H_\nu(y)$ for any $y\in U$ in $|\nu|=\sum_{j\geq 1}\nu_j$
multiplications. Therefore the total cost in this stage is bounded by
\be
i N_h + \sum_{l=1}^i |\nu^l|=i N_h + \sum_{l=1}^i (l-1) =i N_h + i(i-1)/2.
\ee
Summing over $i=1,\dots,n$, we thus find that the total offline cost is of order
\be
C_{\rm off}\sim nC_h+n^2N_h+n^3 \sim  nC_h+n^2N_h.
\ee
Here we neglect, the cost $n^3$ relative to the computation of the $H_{\nu^l}(y_{\nu^i})$ as well as
the cost of $n^2\log (n)$ needed for the non-adaptive
selection of the sets $\Lambda_n$, as derived in \S \ref{secapriori}, since $N_h$ is typically much larger than $n$.

We finally evaluate the online cost. Since the online stage simply amounts
in the combination of the $\alpha_ {\nu^i,h}$ for computing the interpolant, 
we find that this cost is of the order
\be
C_{\rm on} \sim nN_h,
\ee
where we have neglected the cost $n^2$ relative to the computation of the $H_{\nu}(y)$
for the given $y\in U$.

If $\e$ is a targeted order of accuracy, and if we have the error bound
\iref{nepsh}, then one way to reach this accuracy is to take both
$Cn^{-s}$ and $\e(h)$ of the order of $\e$. Denoting by $h(\e)$
the inverse function of $\e(h)$, that is,
\be
h(\e_0)=h_0\quad \Leftrightarrow \quad \e(h_0)=\e_0,
\label{heps}
\ee
we thus find that the interpolation algorithm reaches the order of accuracy $\e$ at costs
\be
C_{\rm off}(\e)\sim \e^{-1/s} C_{h(\e)}+\e^{-2/s}N_{h(\e)}\quad {\rm and}Ê\quad C_{\rm on}(\e)\sim \e^{-1/s}N_{h(\e)}.
\ee
It should be noticed that this algorithm is immune to the curse of dimensionality
since the above trade-off between accuracy and complexity 
is obtained with infinitely many variables.

\section {Taylor approximation}
\label{sectaylor}

The results established in \S 3 show that effective polynomial approximations
\iref{polynomial} to the solution map $y\mapsto u(y)$ can be obtained by best $n$-term truncations of the Taylor series \iref{taylor}, for relevant
classes of parametric PDEs.

In this section, we discuss a strategy, proposed in \cite{CCDS}, 
for  numerically finding a good $n$ term Taylor approximation to $u$.   This numerical method rests
in part on the effective computation of the
Taylor coefficients $t_\nu$. In contrast to the interpolation method discussed
in \S 5, the strategy for $n$-term Taylor approximations is intrusive and strongly exploits the specific structure
of the parametric PDE. In fact, it only applies to the limited, yet relevant,
range of problem where the parametric problem has the form of a linear operator equation
\be
\cB(a)u=f,
\label{genlin}
\ee
where $f\in W$ and $\cB(a)\in \cL(V,W)$ for a pair of Hilbert spaces $(V,W)$, and where
\be
a\mapsto \cB(a),
\ee
is a continuous linear map from $X$ to $\cL(V,W)$.  
In other words, the problem map
\be
\cP: (u,a) \mapsto f-\cB(a)u,
\ee
is linear both in $a$ and $u$, up to the constant term $f$.  Recalling our four examples
of parametric PDEs discussed in \S 2, that is, equations \iref{ellip}, \iref{parab}, \iref{eqpull}, and \iref{nonlin},
the first two fall in this category while the last two do not.

Any linear parametric problem between Hilbert spaces $V$ and $W$
can be expressed through a variational formulation \iref{generalvariatlin}
for a pair of Hilbert spaces $(V,\t V)$ where $\t V$ is the antidual of $W$. 
In our present setting, this 
formulation is: find $u(a)\in V$ such that
\be
B(u(a),v;a)=L(v), \quad v\in \t V,
\label{generalvariatlin1}
\ee
where $B(\cdot,\cdot;a)$ and $L$ are continuous sesquilinear
and antilinear forms over $V\times \t V$ and $\t V$ respectively, 
and where, throughout this section, we make the additional assumption that
\be
a\mapsto B(\cdot,\cdot;a),
\ee 
is a continuous {\it linear} map $X$ to ${\frak B}$ the set of continuous sesquilinear forms over $V\times \t V$. We work under the following 
assumption.
\nl
\nl
{\bf Assumption AL:} {\it The parameter set $\cA$ has a complete affine representer 
$(\psi_j)_{j\geq 1}$ and the sesquilinear form $B(\cdot,\cdot;a)$ satisfies
the inf-sup conditions \iref{infsupa} for all $a\in a(U)$.}
\nl

Notice that this assumption requires that 
the problem \iref{genlin} has a solution for any $a\in a(U)$ 
and for all $f\in W$, in contrast to {\bf Assumption A} which requires that
that it has a solution for all $a\in a(U)$ but only for the given $f\in W$. 
Under this assumption, the solution map 
\be
y\mapsto u(y)=\cB(a(y))^{-1}f,
\ee 
is well defined
over $U$. Examples of problems falling in this category are the elliptic equation \iref{ellip} 
and the parabolic equation \iref{parab} under the uniform elliptic assumption ${\bf UEA}(r)$
as discussed in \S \ref{subuea}. Since all maps in the chain
\be
a \mapsto \cB(a) \mapsto \cB(a)^{-1} \mapsto u(a) =\cB(a)^{-1} f,
\ee
are infinitely Frechet differentiable at $a\in a(U)$, and since $y\mapsto a(y)$ is affine, 
we are ensured of the existence of the partial derivatives 
\be
\partial^\nu u:=\(\prod_{j\geq 1} \frac{\partial}{\partial^{\nu_j} y_j} \)u,
\ee
at every $y\in U$ and for all $\nu\in \cF$,
and therefore of the Taylor coefficients $t_\nu=\frac 1 {\nu!}Ê\partial^\nu u(0)\in V$
for all $\nu\in\cF$.

We first show in \S \ref{tayrec} that these coefficients can be computed 
by a simple recursive procedure which takes advantage of
the linear structure of the problem. When the truncation sets $\Lambda_n$
are downward closed, this procedure computes the coefficients
$(t_\nu)_{\nu\in \Lambda_n}$ at the cost of solving $n$ times
a linear problem with operator $\o \cB:=\cB(\o a)$.
Similar to the interpolation method
in \S 5, the downward closed sets $(\Lambda_n)_{n\geq 1}$ 
for which we compute the coefficients $t_\nu$
can either be chosen in an a priori manner, based on the a priori estimates
for the coefficients $\|t_\nu\|_V$, or adaptively built. 
Various adaptive selection strategies for finding the sets $\Lambda_n$ are proposed in \S \ref{tayadapt}, 
and a convergence analysis is given in \S \ref{tayconvadapt} for one of them 
in the particular case of the elliptic equation \iref{ellip}: we
show that if the sequence
$(\|t_\nu\|_V)_{\nu\in\cF}$ belongs to the space $\ell^p_m(\cF)$
defined in \S \ref{sublower} for some $p<1$, then the adaptive strategy
generates downward closed sets $(\Lambda_n)_{n\geq 1}$ such that the tructated
Taylor series converges in $L^\infty(U,V)$
with the expected rate $n^{-s}$ where $s:=\frac 1 p-1$.
Recall that Theorem \ref{theomainlower} shows that 
$(\|t_\nu\|_V)_{\nu\in\cF}\in \ell^p_m(\cF)$ whenever 
the assumptions of Theorem \ref{theoneiacU} hold and $(\|\psi_j\|_X)_{j\geq 1}\in\ell^p(\N)$.
Space discretization and computational cost of the
adaptive and non-adaptive strategies are 
discussed in \S \ref{taycost}.

\subsection{Recursive computations of the Taylor coefficients}
\label{tayrec}

Since $B(\cdot,\cdot,a)$ is defined for all $a\in X$, we may introduce the sesquilinear forms
\be
\o B(\cdot,\cdot):=B(\cdot,\cdot;\o a)\quad {\rm and}\quad B_j(\cdot,\cdot):=B(\cdot,\cdot;\psi_j).
\ee
Then, the solution map $y\mapsto u(y)=u(a(y))$ is defined by
\be
B(u(y),v;y)=L(v), \quad v\in \t V,
\label{generalvariatliny}
\ee
where 
\be
B(\cdot,\cdot;y):=\o B(\cdot,\cdot)+\sum_{j\geq 1}y_j B_j(\cdot,\cdot).
\label{affineB}
\ee
Note that $\o B$ as well as each individual $B_j$ is bounded, that is, belong to $\frak B$.
In addition, we know from {\bf Assumption AL} that $\o B$ satisfies the inf-sup conditions
\iref{infsupa}, and so does $B(\cdot, \cdot;y)$ for any $y\in U$.
The following result shows that the Taylor coefficients of the solution
map $y\mapsto u(y)$ satisfy simple equations which allow us to compute
them in a recursive way.

\begin{lemma}
\label{lemmatay}
Consider a parametric problem of the form \iref{genlin}
such that {\bf Assumption AL} holds 
for a suitable affine representer $(\psi_j)_{j\geq 1}$. 
Then the Taylor coefficients of the solution map $y\mapsto u(y)$ satisfy the equations
\be
\o B(t_\nu,v)=L_\nu (v), \quad v\in \t V,
\label{recurstnu}
\ee
where $L_\nu(v):=L(v)$ if $\nu=0$ is the null multi-index,
and
\be
L_\nu(v):=-\sum_{j\in {\rm supp}(\nu)}B_j(t_{\nu-e_j},v), 
\ee
when $\nu\in\cF-\{0\}$, where $e_j:=(0,\dots,0,1,0,\dots)$ is the Kroenecker sequence with $1$ at position $j$.
\end{lemma}

\noindent
{\bf Proof:} The case $\nu=0$ is immediate since $t_0=u(0)$ and $B(0)=\o B$. For the other
values of $\nu$, we apply the operator $\partial^\nu$  to the
equation 
\be
L(v)=B(u(y),v;y).
\ee
Since $L$ does not depends on $y$, and due to the affine form \iref{affineB} of $B(\cdot,\cdot,y)$,  
we obtain by the multivariate Leibniz rule  
$$
\begin{disarray}{ll}
0 & =\partial^\nu (B(u(y),v;y)) \\
& = \partial^\nu \(\o B(u(y),v)+\sum_{j\geq 1} y_j B_j (u(y),v)\)\\
& =\o B(\partial^\nu u(y),v)+\sum_{j\geq 1}y_j B_j(\partial^\nu u(y),v)+\sum_{j\in {\rm supp}(\nu)}\nu_j B_j(\partial^{\nu-e_j} u(y),v).
\end{disarray}
$$
At $y=0$, this gives
\be
\o B(\partial^\nu u(0),v)=-\sum_{j\in {\rm supp}(\nu)} \nu_j B_j(\partial^{\nu-e_j} u(0),v).
\ee
Dividing by $\nu !=(\nu-e_j)! \nu_j$, this gives \iref{recurstnu}. \hfill $\Box$
\nl

For further purposes, we give the particular expression of the equations \iref{recurstnu}
in the case of the elliptic equation \iref{ellip}. In this case, $V=\t V=H^1_0(D)$ and the sesquilinear forms
are given by
\be
\o B(u,v):=\int_D \o a \nabla u\cdot \nabla v,
\ee
and
\be
B_j(u,v):=\int_D \psi_j\nabla u\cdot \nabla v.
\ee
Therefore, for $\nu\in\cF-\{0\}$, the coefficient $t_\nu$ is the solution of the boundary value problem
\be
\int_D \o a \nabla t_\nu\cdot \nabla v=-\sum_{j\in {\rm supp}(\nu)}\int_D \psi_j \nabla t_{\nu-e_j} \cdot \nabla v, \quad v\in V.
\label{elliprecurs}
\ee 
Note that, in the particular case where $\nu=e_j$, that is, 
when $t_\nu=\partial_{y_j} u(0)$, this relation has the form
\be
\int_D \o a \nabla \partial_{y_j} u(0)\cdot \nabla v=-\sum_{j\in {\rm supp}(\nu)}\int_D \psi_j \nabla u(0) \cdot \nabla v, \quad v\in V,
\ee 
and it can be derived from the general expression of the Frechet derivative $du(a)$ obtained
in \S 2.1, applied to $h=\psi_j$.

Lemma \ref{lemmatay} shows that if $\Lambda_n$ is a downward closed
set, then it is possible to compute the $n$ Taylor coefficients
$(t_\nu)_{\nu\in \Lambda_n}$ by solving exactly $n$ linear problems
of the form \iref{recurstnu}. This is done by writing $\Lambda_n=\{\nu^1,\dots,\nu^n\}$ 
where the order is such that all sections $\{\nu^1,\dots,\nu^k\}$ are downward closed sets for $k=1,\dots,n$
and by recursively computing the $t_\nu$ in this order.
Then, the right side of the problem \iref{recurstnu} for computing $t_{\nu^k}$
only depends on the $t_{\nu^i}$ for $i<k$ which have already been computed.

In practice, these linear problems can only be solved approximately, through a space discretization, for example
using the finite element method. This induces an error in the computation of the Taylor coefficients. We deal with
this issue in \S \ref{taycost} and assume for the moment that these problems are solved exactly.

As already explained in \S \ref{secapriori}, one non-adaptive strategy consist in
defining $\Lambda_n$ 
as the set of indices corresponding to the $n$ largest a priori estimates $e_\nu$ 
for the Taylor coefficients or the $n$ largest surrogate $s_\nu$ such that the sequence  
$(s_\nu)_{\nu\in\cF}$ is monotone non-increasing. We have observed that 
the total cost of identifying $\Lambda_n$ is then of the order $n^2\log(n)$
which is negligible compared to the computation of the 
approximation polynomial. We then know from the results in \S \ref{sublower}
that if the assumptions of Theorem \ref{theoneiacU} hold for the solution map $a\mapsto u(a)$ and 
if in addition $(\|\psi_j\|_X)_{j\geq 1}\in \ell^p(\cF)$ for some $p<1$, we have
the convergence rate 
\be
\sup_{y\in U} \Big\|u(y)-\sum_{\nu\in \Lambda_n} t_\nu y^\nu\Big\|_V \leq Cn^{-s}, \quad s:=\frac 1 p-1,
\ee
using these a priori selected $(\Lambda_n)_{n\geq 1}$.

In the next section, we discuss adaptive strategies for the selection of the sets $(\Lambda_n)_{n\geq 1}$.  Given the above
result on a priori choices for the $\Lambda_n$, one might ask why should we even consider adaptive strategies?  The answer is that
it may be that the best $n$-term Taylor approximations of $u$ actually perform much better than the proven rate $O(n^{-s})$ established
using the a priori chosen $\Lambda_n$.  This is possible, since
we do not have any results that say the $O(n^{-s})$ rate is best possible under the assumption $(\|\psi_j\|)_{j\ge 1}\in\ell^p(\N)$
and since the estimates of $\|t_\nu\|_V$ by the computable surrogate $s_\nu$ could be too pessimistic. 

\subsection{Adaptive algorithms}
\label{tayadapt}

For any downward closed set $\Lambda$, 
given the Taylor coefficients $(t_\nu)_{\nu\in\Lambda}$,
the equations \iref{recurstnu} allow us to compute the
Taylor coefficients $t_\nu$ for $\nu\in N(\Lambda)$, where $N(\Lambda)$ is the set of neighbors
of $\Lambda$ defined in \iref{neighbors}. This suggests an adaptive algorithm in the same line
as the greedy interpolation algorithms proposed in \S 5.2.
Once again, we start with 
an idealized version of this algorithm which cannot be applied
as such.  We discuss more practical versions later.
\nl
\nl
{\bf Greedy Taylor Algorithm:}
We start with $\Lambda_1:=\{0\}$ the null multi-index.
Assuming that $\Lambda_{n-1}$ has been selected and that
the $(t_\nu)_{\nu\in \Lambda_{n-1}}$ have been computed,
we compute the $t_\nu$ for $\nu\in N(\Lambda_{n-1})$. We then set
\be
\nu^n:={\rm argmax} \{ \|t_\nu\|_V\; : \; \nu\in N(\Lambda_{n-1})\},
\label{greedytay}
\ee
and define $\Lambda_{n}=\Lambda_{n-1} \cup \{\nu^n\}$.
\nl
\nl
The rationale behind this procedure is that if the sequence $(\|t_\nu\|_V)_{\nu\in\cF}$ were monotone
non-increasing, it would automatically select the $n$ largest terms of this sequence.
Similar to the greedy interpolation algorithm, the 
above algorithm needs to be modified in the following two directions:
\begin{itemize}
\item[(i)] In order to guarantee finite complexity  of the search in the case of infinitely many variables, the infinite set $N(\Lambda)$
should be replaced by a finite one such as the set of anchored neighbors $\t N(\Lambda)$ 
defined by \iref{redneighbors}.
\item[(ii)] In order to guarantee convergence, one should alternate the selection
of $\nu^n$ using \iref{nextindex} with a selection rule ensuring that all indices
are eventually picked.
\end{itemize}
This results in the following algorithm, which is
similar to the alternating greedy interpolation algorithm.
\nl
\nl
{\bf Alternating Greedy Taylor Algorithm:} We start with $\Lambda_1:=\{0\}$ the null multi-index.
Assuming that $\Lambda_{n-1}$ has been selected and that
the $(t_\nu)_{\nu\in \Lambda_{n-1}}$ have been computed,
we compute the $t_\nu$ for $\nu\in \t N(\Lambda_{n-1})$.  If $n$ is even, we set
 \be
\nu^n:={\rm argmax} \{ \|t_\nu\|_V\; : \; \nu\in \t N(\Lambda_{n-1})\},
\label{greedy}
\ee
and, if $n$ is odd, we set $\nu^n$ as in \iref{nun}. We then define $\Lambda_{n}=\Lambda_{n-1} \cup \{\nu^n\}$.
\nl
\nl
It is easily checked that the selection of $\nu^n$ by \iref{nun} ensures that the sequence $(\Lambda_n)_{n\geq 1}$
is an exhausion of $\cF$. Therefore, we are ensured that if the Taylor series \iref{taylor} converges
unconditionally towards $u$ in $L^\infty(U,V)$ or in any other norm, the approximations produced
by the greedy algorithm converge towards $u$ in the same norm. However, we do not have much
information on the rate of convergence. 

In view of the results obtained in \S 3, a legitimate objective is to build adaptive algorithms
that can be proven to converge at a rate that is comparable to
that which is proved for best $n$-term Taylor approximations, for relevant classes
of parametric PDEs. Let us point out that a similar objective
can be attained when considering the spatial
discretization of a {\it single} PDE, by either 
adaptive wavelet methods \cite{CDD1,CDD2,GHS} or by 
adaptive finite element methods \cite{MNS,BDD,St}.
More precisely, these paper show that specific refinement
strategies based on a posteriori analysis generate adaptive wavelet sets or adaptive meshes
such that the approximate solution converges with the optimal  algebraic rate 
allowed by the exact solution, as the number of wavelets or elements grows.
One key tool in these algorithms, is the use of a refinement procedure
which guarantees that the error decreases by a fixed amount after 
the refinement is performed. This procedure is called {\it bulk chasing}, and requires that,  in general,
more than  one wavelet/element is added/refined at each step. 

In the present context of $n$-term Taylor approximation, it is possible to 
introduce similar bulk chasing procedures, which, in general, require the selection of  more than
one term   from the Taylor expansion at each step. Iterating 
these bulk chasing  procedure produces a nested sequence $(\Lambda^k)_{k\geq 1}$
of downward closed sets with $\Lambda^1=\{0\}$. Here, we use the notation $\Lambda^k$
instead of $\Lambda_k$ in order to stress that $\#(\Lambda^k)$ is in general larger than $k$.
When we want to index these sets by their cardinality, we may define the sets
\be
\Lambda_n:=\Lambda^k\quad {\rm for}\quad n=n(k):=\#(\Lambda^k),
\label{cardin}
\ee
which are indexed by the integers $n\in \{n(k)\; : \; k\geq 1\}$.  Of course, some indices $n$ are missed in the $\Lambda_n$ notation but these can be
filled in by simply repeating the sets $\Lambda^k$.   The resulting sequence $(\Lambda_n)_{n\ge 1}$ then satisfies $\#(\Lambda_n)\le n$.

We now discuss one specific procedure of this type for which it is proved,
in the particular case of the elliptic problem \iref{ellip}, that the resulting 
adaptive approximation converges with an algebraic convergence rate that
matches the rate that is established when keeping the largest $n$-terms in the
Taylor expansion.
For this purpose, we introduced for any finite downward closed set $\Lambda$
its {\it margin} $M(\Lambda)$, defined by
\be
M(\Lambda): =\{\nu\notin \Lambda\; : \; \exists j\in {\rm supp}(\nu) \; : \; \nu-e_j\in\Lambda\},
\ee
Note that the set of neighbors $N(\Lambda)$ can be defined by
\be
N(\Lambda): =\{\nu\notin \Lambda\; : \; \forall j\in {\rm supp}(\nu) \; : \; \nu-e_j\in\Lambda\}.
\ee
Therefore, we have $N(\Lambda)\subset M(\Lambda)$ and this inclusion
is generally strict. Still, for any downward closed set $\Lambda$, 
given the Taylor coefficients $(t_\nu)_{\nu\in\Lambda}$,
the relations \iref{recurstnu} allow us to compute the
Taylor coefficients $t_\nu$ for $\nu\in M(\Lambda)$.

For the rest of this section, we assume that the Taylor coefficients of
the solution map are $\ell^2$ summable, that is, 
\be
\sum_{\nu\in\cF} \|t_\nu\|_V^2 <\infty.
\ee
Note that, according to Theorem \ref{theomain}, this holds if the assumptions of
Theorem \ref{theoneiacU} are satisfied and if in addition $(\|\psi_j\|_X)_{j\geq 1}\in \ell^p(\cF)$ for some $p<1$.
For any set $S\subset \cF$, we introduce the quadratic energy
\be
e(S):=\sum_{\nu \in S} \|t_\nu\|_V^2.
\ee
We also define for any finite downward closed set $\Lambda$ the quadratic error
\be
\sigma(\Lambda):=\sum_{\nu\notin \Lambda} \| t_\nu\|_V^2.
\ee
Note that since the functions $y\mapsto y^\nu$ do not form an orthonormal
system, this quantity differs from the mean-square error between $u$ and its truncated
Taylor series.

The bulk chasing procedure consist in building the new set $\Lambda^k$
by adding to $\Lambda^{k-1}$ a subset $S^{k-1}$ of the margin $M^{k-1}:=M(\Lambda^{k-1})$
which captures a prescribed portion of its energy, that is, such that
\be
e(S^{k-1}) \geq \theta e(M^{k-1}),
\label{bulk}
\ee
for some fixed $0<\theta <1$. 

The objective of this procedure is to reduce the quadratic error $\sigma(\Lambda^{k-1})$
by a fixed amount. This  will be achieved provided that the
Taylor expansion \iref{taylor} satisfies a so-called ``saturation property'',
which is analogous to that which is sometimes established in order to prove convergence of adaptive finite element methods,
see \cite{Do,MNS,BDD}.  In \S \ref{tayconvadapt}, we establish the validity of this property in the case of the parametric elliptic problem
\iref{ellip}, provided that the uniform ellipticity property ${\bf UEA}(r)$ holds for some $r>0$.  
For now, we start from this property viewed as a general assumption in order to analyze the convergence
of adaptive algorithms based on bulk chasing.
\nl
\nl
{\bf Saturation property:} {\it there exists a fixed constant $\delta>1$ such that for
any finite downward closed set $\Lambda$,
\be
\sigma(\Lambda)\leq \delta e(M),
\label{satur}
\ee
where $M=M(\Lambda)$.}
\nl
\nl
If this property holds and if we use the bulk chasing procedure
to construct $(\Lambda^k)_{k\geq 1}$ we may write 
\be
\sigma(\Lambda^k) = \sigma(\Lambda^{k-1})-e(S^{k-1})\leq  \sigma(\Lambda^{k-1})-\theta e(M^{k-1}) \leq \kappa  \sigma(\Lambda^{k-1}),
\label{reduce}
\ee
where
\be
\kappa=1-\frac{\theta}{\delta} <1.
\ee
Therefore, the quadratic error decreases by a fixed amount at every step of the algorithm.
and in particular, at step $k$ we have
\be
\sigma(\Lambda^k) \leq C\kappa^k,
\ee
where  $C:=\kappa^{-1} \sum_{\nu\neq 0}\|t_\nu\|_V^2$.

The set $S^{k-1}$ should be chosen as small as possible,
however we need to ensure that the new set $\Lambda^k:=\Lambda^{k-1}\cup S^{k-1}$ is still downward closed.
This is executed by the following algorithm, which we first present in an idealized form
that cannot be applied as such.
\nl
\nl
{\bf Bulk Chasing Taylor Algorithm:} Having fixed $0<\theta <1$, we start with $\Lambda^1:=\{0\}$ the null multi-index.
Assume that $\Lambda^{k-1}$ has been selected
and that the coefficients $t_\nu$ have been computed for $\nu\in\Lambda^{k-1}$.
For all $\nu\in M^{k-1}$ we compute the $t_\nu$ and the quantities
\be
m_\nu:=\max\{\|t_{\t \nu}\|_V \; : \; \t \nu\geq \nu\; {\rm and}\; \t \nu\in M^{k-1}\}.
\label{mnu}
\ee
We define $S^{k-1}$ as the set of indices
$\nu\in M^{k-1}$ corresponding to the $l$ largest $m_\nu$, for the smallest
value of $l$ such that the bulk condition \iref{bulk} is met. We then define 
\be
\Lambda^k:=\Lambda^{k-1}\cup S^{k-1}.
\ee

\begin{remark} 
The quantity $m_\nu$ is introduced in order to guarantee that the new set $\Lambda^k$ is downward closed.
This can always be ensured due to the monotonicity property 
\be
 \nu,\t\nu\in M^{k-1}\quad {\rm and}
\quad \t \nu\geq \nu \quad \Rightarrow \quad m_{\t \nu}\leq m_\nu.
\ee
Another option would be to define $S^{k-1}$ as the smallest
subset of $M^k$ such that the bulk condition \iref{bulk}Ê is met and such
that $\Lambda^k:=\Lambda^{k-1}\cup S^{k-1}$ is monotone, however
it is not clear if there is a simple algorithm for determining such a set.
\end{remark}

The above bulk chasing algorithm is not computationally feasible due to the fact that
the margin $M(\Lambda)$ of a finite downward closed set $\Lambda$ has infinite cardinality in the case of
countably many variable $y_j$. We want to modify it by restricting
the computation of the $t_\nu$ and the bulk search to a finite subset of $M^k$.     In order to accomplish this,
we take the usual view in numerical computation, where we are given a target accuracy $\e>0$ and we want the algorithm
  to achieve this accuracy as efficiently as possible.  So, 
in our modified algorithm, we   design the procedure so  that the algorithm
terminates when $\sigma(\Lambda^k)\leq C \e$ for some fixed constant $C$
to be specified later.     We begin by introducing a procedure that
computes from the Taylor coefficients $(t_\nu)_{\nu\in\Lambda}$ a finite 
version of the margin $M=M(\Lambda)$ which captures its energy
up to accuracy $\e$: if $\Lambda$ is a finite downward closed set with margin $M$, then
\be
\t M={\rm SPARSE}(\Lambda,(t_\nu)_{\nu\in\Lambda},\e),
\ee
is a finite subset of $M$ such that $\Lambda\cup \t M$ is, by definition, a downward closed set and such that
\be
e(M\sm \t M)Ê\leq \e.
\label{sparseaccuracy}
\ee
We present in the end of \S \ref{tayconvadapt} one practical realization of such a procedure in 
the particular case of the elliptic equation \iref{ellip}. Our modified algorithm is the following. 
\nl
\nl
{\bf Bulk Chasing Taylor Algorithm with $\e$-Accuracy:} Having fixed $0<\theta <1$, we start with $\Lambda^1:=\{0\}$ the null multi-index.
Assume that $\Lambda^{k-1}$ has been selected
and that the coefficients $t_\nu$ have been computed for $\nu\in\Lambda^{k-1}$.
We define
\be
\t M^{k-1}={\rm SPARSE}(\Lambda^{k-1},(t_\nu)_{\nu\in\Lambda^{k-1}},\e),
\ee
For all $\nu\in\t M^{k-1}$ we compute the $t_\nu$ and the quantities
\be
m_\nu:=\max\{\|t_{\t \nu}\|_V \; : \; \t \nu\geq \nu\; {\rm and}\; \t \nu\in \t M^{k-1}\}.
\ee
We define $S^{k-1}$ as the set of indices
$\nu\in\t M^{k-1}$ corresponding to the $l$ largest $m_\nu$, for the smallest
value of $l$ such that the bulk condition
\be
e(S^{k-1})\geq \theta e(\t M^{k-1}),
\label{bulk2}
\ee
is met. We then define 
\be
\Lambda^k:=\Lambda^{k-1}\cup S^{k-1}.
\ee
The algorithm is stopped when $e(\t M^k)\leq 2\theta \e$.
\nl

The same computation as in \iref{reduce}, now using \iref{bulk2} and \iref{sparseaccuracy} together with the saturation property
\iref{satur}, leads to the reduction inequality
\be
\sigma(\Lambda^k) \leq  \sigma(\Lambda^{k-1})-\theta e(\t M^{k-1}) \leq \kappa  \sigma(\Lambda^{k-1})+\theta \e,
\label{reduceps}
\ee
Therefore we are ensured that after sufficiently many steps $k$, we have 
\be
e(\t M^k)\leq e(M^k)\leq \sigma(\Lambda^k)\leq 2\theta \e,
\ee
and thus the algorithm terminates. The final step may occur before $\sigma(\Lambda^k)\leq 2\theta \e$, but 
we are still ensured by the saturation property \iref{satur} that
\be
\sigma(\Lambda^k)\leq \delta e(M^k) \leq \delta (e(\t M^k)+\e) \leq C\e, \quad \quad C:=\delta(2\theta+1),
\ee
and thus we have reached the announced order of accuracy for the quadratic error.

However, this does not settle the convergence analysis of the algorithm.
First, we want to relate the accuracy $\sigma(\Lambda^k)$ 
with the number of terms $\#(\Lambda^k)$ in the Taylor approximation
through a quantitative convergence estimate.
Secondly, we also want to retrieve convergence estimates for the error between $u$ and its 
truncated Taylor approximation measured in the $L^\infty(U,V)$ metric. This is the purpose
of the next section. 

\subsection{Convergence analysis of adaptive algorithms}
\label{tayconvadapt}

We know that if the Taylor series converges conditionally towards $u$ in $L^\infty(U,V)$
and if the sequence $(\|t_\nu\|_V)_{\nu\in\cF}$ belongs to the 
sequence space $\ell^p_m(\cF)$ for some $0<p<1$, then there
exists a sequence of downward closed sets $(\Lambda_n)_{n\geq 1}$ such that $\#(\Lambda_n)=n$
and such that the uniform error bound
\be
\sup_{y\in U} \Big \|u(y)-\sum_{\nu\in \Lambda_n} t_\nu y^\nu\Big \|_V \leq Cn^{-s},\quad s:=\frac 1 p-1,
\label{unifgoal}
\ee
holds for all $n\geq 1$. This holds in particular if the assumptions of Theorem \ref{theoneiacU} hold and 
if in addition $(\|\psi_j\|_X)_{j\geq 1}\in \ell^p(\cF)$.
In this section, we show that this benchmark rate is met by the two previously described bulk chasing Taylor
algorithms, provided that the saturation assumption hold. We begin with a result that describes
the rate of decay of the quadratic error $\sigma(\Lambda_n)$.

\begin{theorem}
\label{theoadapttaylor}
Consider a parametric problem of the form \iref{genlin}
such that {\bf Assumption AL} holds 
for a suitable affine representer $(\psi_j)_{j\geq 1}$. Assume that $(\|t_\nu\|_V)_{\nu\in\cF}$ belongs to the 
sequence space $\ell^p_m(\cF)$ for some $0<p<1$,
and that the saturation property holds. Then for the bulk chasing Taylor algorithm,
the convergence estimate
\be
\sigma(\Lambda_n) \leq C \|(\|t_\nu\|_V)_{\nu\in\cF}\|_{\ell^p_m}^2 n^{-2r},\quad r:=\frac 1 p-\frac 1 2,
\label{estimsigma}
\ee
holds for all $n\in \{n(k)\; : \; k\geq 1\}$ where we have used the convention \iref{cardin}.
The constant $C$ depends on $p$, $\theta$ and $\delta$.
The same estimate holds for the bulk chasing Taylor algorithm with $\e$-accuracy 
for all $n\in \{n(k)\; : \; 1\leq k\leq k(\e)\}$, where
$k(\e)$ is the step where the algorithm terminates. 
\end{theorem}

\noindent
{\bf Proof:} We begin by considering the bulk chasing Taylor algorithm. We first control the cardinality of
the set $S^{k-1}$ which is added to $\Lambda^{k-1}$, for any $k>1$. Recall that this set is obtained by picking indices $\nu\in M^{k-1}$
corresponding to the $l$ largest $m_\nu$ defined by \iref{mnu}, for the smallest value of $l$ such that the bulk condition
\iref{bulk} is met. Let $\t S^{k-1}\subset S^{k-1}$ denote the set corresponding to the $l-1$ largest $m_\nu$
for this value of $l$. Since the bulk condition is not met by $\t S^{k-1}$, we have
\be
(1-\theta) e(M^{k-1})\leq e(M^{k-1})-e(\t S^{k-1})=\sum_{\nu\in M^{k-1}\sm \t S^{k-1}} \|t_\nu\|_V^2.
\label{nonbulk}
\ee
On the one hand, using Stechkin's Lemma \ref{stechkin}, we find that
\be
\sum_{\nu\in M^{k-1}\sm \t S^{k-1}} \|t_\nu\|_V^2   \leq 
\sum_{\nu\in M^{k-1}\sm \t S^{k-1}} m_\nu^2 
 \leq \|(m_\nu)_{\nu\in M^{k-1}}\|_{\ell^p}^2 l^{-2r}
 \ee 
and therefore, using the fact that $m_\nu$ is dominated by the monotone majorant of $\|t_\nu\|_V$, and that $l=\#(S^{k-1})$, we obtain
\be
\sum_{\nu\in M^{k-1}\sm \t S^{k-1}} \|t_\nu\|_V^2  \leq 
 \|(\|t_\nu\|_V)_{\nu\in\cF}\|_{\ell^p_m}^2 (\#(S^{k-1}))^{-2r}.
\ee
Combining this with \iref{nonbulk} we find that
\be
(1-\theta) e(M^{k-1}) \leq \|(\|t_\nu\|_V)_{\nu\in\cF}\|_{\ell^p_m}^2 (\#(S^{k-1}))^{-2r}.
\label{cmn1}
\ee
Using the saturation property, it follows that
\be
\frac {1-\theta}{\delta} \sigma(\Lambda^{k-1})\leq \|(\|t_\nu\|_V)_{\nu\in\cF}\|_{\ell^p_m}^2 (\#(S^{k-1}))^{-2r},
\ee
or equivalently
\be
\#(S^{k-1})\leq \(\frac {\delta}{1-\theta}\)^{1/2r}\|(\|t_\nu\|_V)_{\nu\in\cF}\|_{\ell^p_m}^{1/r}\sigma(\Lambda^{k-1})^{-1/2r}
\ee
For any $k>1$, we may thus control the cardinality of $\Lambda^k$ by writing
\be
\#(\Lambda^k) = 1+\sum_{l=1}^{k-1} \#(S^l)
\leq 1+ \(\frac {\delta}{1-\theta}\)^{1/2r}\|(\|t_\nu\|_V)_{\nu\in\cF}\|_{\ell^p_m}^{1/r}\sum_{l=1}^{k-1} \sigma(\Lambda^{l})^{-1/2r}.
\ee
On the other hand, we know from \iref{reduce} that $\sigma(\Lambda^l) \geq \kappa^{l-k} \sigma(\Lambda^k)$ with $\kappa:=1 -\frac\delta \theta$, and therefore
\be
\#(\Lambda^k)  \leq 1+\(\frac {\delta}{1-\theta}\)^{1/2r}\frac{ \kappa^{1/2r}}{1-\kappa^{1/2r}}\|(\|t_\nu\|_V)_{\nu\in\cF}\|_{\ell^p_m}^{1/r}\sigma(\Lambda^k)^{-1/2r}.
\ee
This can be rewritten as
\be
\sigma(\Lambda^k) \leq \frac {\delta}{1-\theta}\frac{ \kappa}{(1-\kappa^{1/2r})^{2r}}\|(\|t_\nu\|_V)_{\nu\in\cF}\|_{\ell^p_m}^{2}(\#(\Lambda^k) - 1)^{-2r}.
\ee
Using the inequality $\#(\Lambda^k) - 1\geq \frac 1 2 \#(\Lambda^k)$ in the case where $k>1$, we have thus established
\iref{estimsigma} with constant $C:=2^{2r}\frac {\delta}{1-\theta}\frac{ \kappa}{(1-\kappa^{1/2r})^{2r}}>1$. If $k=1$,
we simply write
\be
\sigma(\Lambda^{1})\leq \|(\|t_\nu\|_V)_{\nu\in\cF}\|_{\ell^2}^{2}\leq \|(\|t_\nu\|_V)_{\nu\in\cF}\|_{\ell^p_m}^{2},
\ee
which shows that \iref{estimsigma} also holds in this case.

We next consider the bulk chasing Taylor algorithm with $\e$-accuracy and explain how
\iref{estimsigma} can be established for $n=n(k)$ with $1\leq k\leq k(\e)$ up to an inflation in the constant $C$
by a similar argument. First, with the exact
same reasoning which led to \iref{cmn1}, we obtain the estimate
\be
(1-\theta) e(\t M^{k-1}) \leq \|(\|t_\nu\|_V)_{\nu\in\cF}\|_{\ell^p_m}^2 (\#(S^{k-1}))^{-2r}.
\ee
We next observe that since $e(\t M^{k-1})\geq 2\theta \e$ for $k\leq k(\e)$, we have the modified saturation property
\be
\sigma(\Lambda^{k-1})\leq \delta e(M^{k-1}) \leq \delta (e(\t M^{k-1})+\e) \leq \t \delta e(\t M^{k-1}),
\ee
with $\t \delta:=\delta (1+\frac 1 {2\theta})$. By the same reasoning, for 
any $n>1$, we obtain the bound 
\be
\#(\Lambda^k) = 1+\sum_{l=1}^{k-1} \#(S^l)
\leq 1+ \(\frac {\t \delta}{1-\theta}\)^{1/2r}\|(\|t_\nu\|_V)_{\nu\in\cF}\|_{\ell^p_m}^{1/r}\sum_{l=1}^{k-1} \sigma(\Lambda^l)^{-1/2r}.
\ee
The saturation property implies that $\sigma(\Lambda^l) \geq \t \kappa^{l-l} \sigma(\Lambda^k)$ with $\t \kappa:=1-\frac {\t \delta} \theta$.
Therefore, by the same reasoning, we reach \iref{estimsigma} with 
the larger constant $C:=2^{2r}\frac {\t \delta}{1-\theta}\frac{ \t \kappa}{(1-\t \kappa^{1/2r})^{2r}}>1$.
\hfill $\Box$
\nl

Our next result shows that the benchmark rate \iref{unifgoal} is met under the same assumptions 
as those of the above theorem.

\begin{theorem}
\label{theoadapttaylor2}
Consider a parametric problem of the form \iref{genlin} such that the assumptions of Theorem \ref{theoadapttaylor}
hold and such that in addition the 
Taylor series converges conditionally towards $u$ in $L^\infty(U,V)$. Then, we have for all $n\geq 1$ the uniform convergence estimate
\be
\sup_{y\in U} \Big \|u(y)-\sum_{\nu\in \Lambda_n} t_\nu y^\nu\Big\|_V \leq C \|(\|t_\nu\|_V)_{\nu\in\cF}\|_{\ell^p_m} n^{-s},\quad s:=\frac 1 p-1.
\label{estimunif}
\ee
holds for all $n\in \{n(k)\; : \; k\geq 1\}$. The constant $C$ depends on $p$, $\theta$ and $\delta$.
The same estimate holds for the bulk chasing Taylor algorithm with $\e$-accuracy 
for all $n\in \{n(k)\; : \; 1\leq k\leq k(\e)\}$, where
$k(\e)$ is the step where the algorithm terminates. 
\end{theorem}

\noindent
{\bf Proof:} It suffices to prove that \iref{estimsigma} implies \iref{estimunif} for the same value of $n$, up to a change
in the constant $C$. Since the Taylor series converges conditionally and since $(\|t_\nu\|_V)_{\nu\in\cF}$ belongs
to $\ell^1(\cF)$, this series also converges unconditionally. We thus have
\be
\sup_{y\in U} \Big\|u(y)-\sum_{\nu\in \Lambda^k} t_\nu y^\nu\Big \|_V  \leq \sum_{\nu\notin\Lambda^k} \|t_\nu\|_V.
\ee
For $n=n(k)=\#(\Lambda_n)=\#(\Lambda^k)$, we consider the 
set $\Lambda_n^*$ of the indices corresponding to the $n$ largest $\|t_\nu\|_V$. Using Lemma \ref{stechkin}, Cauchy-Schwarz inequality
and \iref{estimsigma}, we write
$$
\begin{disarray}{ll}
\sum_{\nu\notin\Lambda_n} \|t_\nu\|_V &\leq \sum_{\nu\notin\Lambda_n^*} \|t_\nu\|_V
+\sum_{\nu\in\Lambda_n^*\sm \Lambda_n} \|t_\nu\|_V \\
& \leq  \|(\|t_\nu\|_V)_{\nu\in\cF}\|_{\ell^p}(n+1)^{-s}+ n^{1/2}e(\Lambda_n^*\sm \Lambda_n)^{1/2} \\
& \leq \|(\|t_\nu\|_V)_{\nu\in\cF}\|_{\ell^p}(n+1)^{-s}+  n^{1/2}\sigma( \Lambda_n)^{1/2}\\
& \leq \|(\|t_\nu\|_V)_{\nu\in\cF}\|_{\ell^p}(n+1)^{-s}+   n^{1/2}C^{1/2}\|(\|t_\nu\|_V)_{\nu\in\cF}\|_{\ell^p_m}n^{-r}\\
& \leq (1+C^{1/2})\|(\|t_\nu\|_V)_{\nu\in\cF}\|_{\ell^p_m}n^{-s},
\end{disarray}
$$
which confirms \iref{estimunif}.\hfill $\Box$
\nl

The   above theorem assumes  that the saturation property is valid.  We next prove that this property always holds in the particular case of the 
elliptic problem \iref{ellip}.  Recall that we then have
$V=H^1_0(D)$ and $X=L^\infty(D)$. The saturation 
property turns out to be a consequence of the
uniform ellipticity assumption ${\rm \bf UEA}(r)$. In order to see this, we introduce the norm
\be
\|v\|_{\o a}:=\( \int_D \o a |\nabla v|^2 \)^2,
\ee
which is equivalent to the $V$-norm under ${\rm \bf UEA}(r)$, since we then have
\be
r \|v\|_V^2 \leq \o a_{\min} \|v\|_V^2 \leq \|v\|_{\o a}^2 \leq \o a_{\max}\|v\|_V^2,
\label{equivnorm}
\ee
where $\o a_{\min}:=\min_{x\in D} \o a(x)$ and $\o a_{\max}:=\max_{x\in D} \o a(x)=\|\o a\|_{X}$. 
For $\nu\in \cF$ and $j\geq 1$, we use the notation
\be
d_\nu:=\|t_\nu\|_{\o a}^2,
\ee
and 
\be
d_{\nu,j}:=\int_D |\psi_j| |\nabla t_\nu|^2.
\ee
The proof of the saturation property uses the following lemma which relates the above quantities.

\begin{lemma}
\label{lemmatnu}
Consider a parametric problem of the type \iref{ellip}, with affine representer $(\psi_j)_{j\geq 1}$ such that ${\rm \bf UEA}(r)$ holds for some $r>0$.
Then, we have
\be
\sum_{j\geq 1}d_{\nu,j} \leq \gamma d_\nu, \quad \quad \gamma:=1-\frac r{\o a_{\max}}<1,
\label{gammanu}
\ee
and
\be
d_\nu
\le 
\alpha \sum_{j\in {\rm supp}(\nu)} d_{\nu-e_j,j}, \quad \quad \alpha:=\frac {\o a_{\max}}{r+\o a_{\max}} <1,
\label{TL1}
\ee
where $e_j$ is the Kroenecker sequence with $1$ at position $j$.
\end{lemma}

\noindent
{\bf Proof:} The uniform ellipticity assumption implies that, for all $x\in D$, 
\be
\sum_{j\geq 1}|\psi_j(x)| \leq \o a(x)-r \leq  \gamma \o a(x),
\label{gamma}
\ee
which implies \iref{gammanu}. On the other hand, we take $v=t_\nu$
in \iref{elliprecurs}, which gives
\be
d_\nu=-\sum_{j\in {\rm supp}(\nu)} \int_D \psi_j\nabla t_{\nu-e_j}\nabla t_\nu,
\label{recurs1}
\ee
and therefore
\be
\label{firstest}
d_\nu \leq \frac 1 2 \sum_{j\in {\rm supp}(\nu)} \int_D |\psi_j|\, |\nabla t_{\nu-e_j}|^2
+ 
\frac 1 2 \sum_{j\in {\rm supp}(\nu)} \int_D |\psi_j|\, |\nabla t_{\nu}|^2
\;.
\ee
Using \eref{gamma} in  the second term of \eref{firstest} gives 
\be
\(1-\frac \gamma 2\)d_\nu \le  \frac 1 2 
\sum_{j\in {\rm supp}(\nu)} \int_D |\psi_j|\, |\nabla t_{\nu-e_j}|^2,
\ee
from which we derive \eref{TL1}. \hfill $\Box$
\nl

We are now in position to establish the saturation property for the elliptic problem \iref{ellip} under
the uniform ellipticity assumption. For any downward closed set $\Lambda$ and any set $S$, we
introduce the modified quadratic error and energy
\be
\o \sigma(\Lambda):=\sum_{\nu\notin \Lambda}d_\nu \quad{\rm and} \quad 
\o e(S):=\sum_{\nu\in S}d_\nu.
\ee

\begin{theorem}
\label{theosatur}
Consider a parametric problem of the type \iref{ellip} such that ${\rm \bf UEA}(r)$ holds for some $r>0$.
Then the saturation property \iref{satur} holds with $\delta$ depending on $r$ and $\o a_{\max}$.
\end{theorem} 

\noindent
{\bf Proof:} We consider an arbitrary downward closed set $\Lambda$ and its margin $M:=M(\Lambda)$.
We first observe that
\be
\o \sigma(\Lambda)=\o e(M)+\o \sigma(\t \Lambda), \quad \quad \t \Lambda:=\Lambda\cup M.
\label{lambdaMtlambda}
\ee
Using \iref{TL1}, we may write
\be
\o \sigma(\t \Lambda) \leq \sum_{\nu\notin \t \Lambda} d_\nu \leq \alpha \sum_{\nu\notin \t\Lambda} \(\sum_{j\in {\rm supp}(\nu)} d_{\nu-e_j,j}\) \leq A+B,
\label{aplusb}
\ee
where 
\be
A:= \alpha \sum_{\nu\notin \t \Lambda} \(\sum_{j\,{\rm s.t.} \, \nu-e_j\notin\t \Lambda} d_{\nu-e_j,j}\)
=\alpha \sum_{\nu\notin \t\Lambda} \(\sum_{j\,{\rm s.t.} \, \nu+e_j\notin \t \Lambda} d_{\nu,j}\),
\ee
and
\be
B:=\alpha \sum_{\nu\notinÊ\t \Lambda} \( \sum_{j\,{\rm s.t.} \, 
 \nu-e_j\in \t\Lambda} d_{\nu-e_j,j}\)
=\alpha \sum_{\nu\in M} \(\sum_{j\,{\rm s.t.} \, \nu+e_j\notin\t \Lambda} d_{\nu,j}\).
\ee
In this splitting, we have used the fact that if $\nu\notin\t \Lambda$
and $\nu_j\neq 0$, we have either $\nu-e_j\notin\t \Lambda$
or $\nu-e_j\in M$. Using \iref{gammanu}, we control the first term $A$ by
\be
A\leq \alpha\gamma \sum_{\nu\notin \t\Lambda} d_\nu =\alpha\gamma \o \sigma(\t\Lambda),
\label{estima}
\ee
and by the same argument we obtain
\be
B\leq \alpha\gamma \o e(M).
\label{estimb}
\ee
Combining these estimates with \iref{aplusb}, it follows that
\be
(1-\alpha\gamma)\o\sigma(\t \Lambda)\leq \alpha\gamma \o e(M),
\ee
and thus by \iref{lambdaMtlambda}
\be
\o\sigma(\Lambda)\leq \(1+ \frac {\alpha\gamma}{1-\alpha\gamma}\)\o e(M).
\ee
Finally, using the norm equivalence \iref{equivnorm}, we obtain the saturation property \iref{satur}
with constant $\delta:=\frac {\o a_{\max}}{r}\(1+ \frac {\alpha\gamma}{1-\alpha\gamma}\)$. \hfill $\Box$
\nl
\nl
We conclude this section by presenting a concrete realization of the procedure
SPARSE which is used in the bulk chasing Taylor algorithm with $\e$-accuracy,
in the particular case of the 
elliptic problem \iref{ellip}. We again work under ${\rm \bf UEA}(r)$.
We define
$$
\o \psi_j:=\frac {\psi_j}{\o a},
$$
and choose an integer $J>0$ large enough such that
\be
\Big \|\sum_{j>J} |\bar \psi_j| \Big \|_{X}
\leq 
\(\frac {\alpha \o e(\Lambda)}{1-\alpha\gamma}\)^{-1}r\e,
\label{tailJ}
\ee
where $\alpha$ and $\gamma$ are defined as in Lemma \ref{lemmatnu},
and we define 
\be
\t M:={\rm SPARSE}(\Lambda, (t_\nu)_{\nu\in\Lambda},\e):=\{ \nu\in M
\; ; \; \nu-e_j\in\Lambda \Rightarrow j\leq J\}.
\ee
Clearly $\t M$ is finite with $\#(\t M)\leq J\#(\Lambda)$.  
  
\begin{theorem}
\label{taylor11}
With the above definition of $\t M$, one has
\be
e(M\setminus\t M)=\sum_{\nu\in M\setminus \t M} \|t_\nu\|^2_V\le \e.
\label{TL11}
\ee
\end{theorem}
 \noindent
{\bf Proof:} We proceed in a similar way to the proof
of Theorem \ref{theosatur}, by first writing
\be
\o e(M\setminus\t M)  \leq \alpha \sum_{\nu \in M\setminus\t M} \(\sum_{j\in {\rm supp}(\nu)} d_{\nu-e_j,j}\) \leq A+B,
\label{aplusb2}
\ee
where now
\be
A:= \alpha \sum_{\nu \in M\setminus\t M} \(\sum_{j\,{\rm s.t.} \, \nu-e_j\in M\setminus\t M} d_{\nu-e_j,j}\)
=\alpha \sum_{\nu \in M\setminus\t M} \(\sum_{j\,{\rm s.t.} \, \nu+e_j \in M\setminus\t M} d_{\nu,j}\),
\ee
and
\be
B:=\alpha \sum_{\nu \in M\setminus\t M} \( \sum_{j\,{\rm s.t.} \, 
 \nu-e_j\notin M\setminus\t M} d_{\nu-e_j,j}\)
=\alpha \sum_{\nu\in \Lambda\cup \t M} \(\sum_{j\,{\rm s.t.} \, \nu+e_j \in M\setminus\t M} d_{\nu,j}\).
\ee
In this splitting, we have used the fact that if $\nu\in M\setminus\t M$
and $\nu_j\neq 0$, we have either $\nu-e_j\in M\setminus\t M$
or $\nu-e_j\in\Lambda\cup \t M$. Using \iref{gammanu}, we can bound $A$ by
\be
A\leq \alpha\gamma \sum_{\nu\in M\setminus\t M} d_\nu =\alpha\gamma \o e(M\setminus\t M).
\ee

To bound $B$, we  first claim that for any $\nu\in\Lambda\cup\t M$ such that $\nu+e_j\in M\setminus\t M$,
we must have $\nu \in\Lambda$ and $j> J$.  Indeed, since $\nu+e_j\in M\setminus\t M$, the definition of $\t M$ guarantees that $\nu +e_j=\t\nu  +e_k$ for some $\t\nu\in\Lambda$ and $k> J$.  
If $j=k$ we have our claim.  If $j\neq k$ then necessarily $\t\nu-e_j\in\Lambda$ since $\Lambda$ is a 
downward closed set, and therefore $\nu$ can be written as the sum of $\t \nu-e_j\in\Lambda$ and $e_k$, 
which means that $\nu$ is not in $\t M$. Thus, we have verified our claim.  From the claim, it follows that 
the only $j$'s that may contribute in the
summation inside $B$ are such that $j>J$ and 
$\nu-e_j\in\Lambda$. Hence, 
$$
\begin{disarray}{ll}
B &\leq \alpha \sum_{\nu\in \Lambda} \sum_{j>J} d_{\nu,j}\\
& =\alpha\sum_{\nu\in \Lambda}\int_D\(\sum_{j>J} |\psi_j|\) |\nabla t_\nu|^2 \\
& =\alpha\sum_{\nu\in \Lambda}\int_D\(\sum_{j>J} |\o\psi_j|\)\o a |\nabla t_\nu|^2 \\
& \leq \alpha \Big \|\sum_{j>J} |\o\psi_j| \Big \|_{X} e(\Lambda)\le (1-\alpha\gamma) r\e.
\end{disarray}
$$
Combining the bounds for $A$ and $B$ with \iref{aplusb2}, we obtain
\be
\o e(M\setminus\t M)\leq \frac {B}{1-\alpha\gamma}\leq  r\e,
\ee
which by \iref{equivnorm} Êimplies \iref{TL11}.
\hfill $\Box$

\subsection{Space discretization and computational cost}
\label{taycost}

In numerical computation, we need to take into account the
additional space discretization of the solution map
in the space $V_h\subset V$. In the case of variational problems
of the form \iref{generalvariatlin1},
one typical such discretization is by the Petrov-Galerkin method:
we define $u_h(a)\in V_h$ such that
\be
B(u_h(a),v_h;a)= L(v_h), \quad v_h\in \t V_h,
\label{generalvariatlin1h}
\ee
where $\t V_h\subset \t V$ is an auxiliary finite element space such that ${\rm dim}(\t V_h)={\rm dim}(V_h)$. 
For elliptic problems such as \iref{ellip}, we have $\t V=V$ and we may take $\t V_h=V_h$,
which is the standard Galerkin method expressed in \iref{galerellip}.
We make the assumption that the discrete problem is well posed
for all $a\in a(U)$, that is, {\bf Assumption AL} also holds for the 
discrete problem. 

Defining $u_h(y)=u_h(a(y))$ for a given affine representation, we thus have
\be
B(u_h(y),v_h;y)= L(v_h), \quad v_h\in \t V_h,
\label{generalvariatlinyh}
\ee
and the same computation as in Lemma \ref{lemmatay} shows that 
the Taylor coefficients $t_{\nu,h}\in V_h$ of $y\mapsto u_h(y)$ 
are computed by solving 
\be
\o B(t_{\nu,h},v_h)=L_\nu(v_h), \quad v_h\in \t V_h,
\ee
where $L_\nu=L$ when $\nu=0$ is the null multi-index and
\be
L_\nu(v_h):=-\sum_{j\in {\rm supp}(\nu)}B_j(t_{\nu-e_j,h},v_h), 
\ee
when $\nu\in\cF-\{0\}$. Note that these relations amount 
in applying the Petrov-Galerkin approximation in the 
recursive computation of the Taylor coefficients $t_\nu$.

Non-adaptive and adaptive strategies
may therefore be applied in order to compute truncated Taylor expansions of the form
\be
u_{n,h}(y):=\sum_{\nu\in \Lambda_n} t_{\nu,h} y^\nu,
\ee
with a similar convergence analysis as
for the continuous problem \iref{generalvariatliny}.
In particular, if the assumptions of Theorem \ref{theoneiacU} hold 
for the solution map $a\mapsto u_h(a)$ and 
if in addition $(\|\psi_j\|_X)_{j\geq 1}\in \ell^p(\cF)$, 
both non-adaptive methods based on a priori bounds
for the $\|t_{\nu,h}\|_V$ or adaptive methods based on bulk chasing
have convergence rate
\be
\|u_h-u_{n,h}\|_{L^\infty(U,V)} \leq Cn^{-s}, \quad s:=\frac 1 p-1.
\ee
The constant $C$ is independent of $h$ if in the assumptions
of Theorem \ref{theoneiacU} the open set $\cO$ and
the bound in \iref{uaC1} can be fixed independently of $h$.

Similar to the splitting \iref{secondsplit} that was used for the interpolation method, we may split the
resulting error into
\be
\|u-u_{n,h}\|_{L^\infty(U,V)} \leq \|u_h-u_{n,h}\|_{L^\infty(U,V)}+ \|u-u_h\|_{L^\infty(U,V)}.
\ee
The second term is bounded by the error $\e(h)$ of the numerical solver. Therefore
we obtain an global error bound of the form
\be
\|u-u_{n,h}\|_{L^\infty(U,V)} \leq Cn^{-s}+\e(h),
\ee
similar to the bound \iref{nepsh} obtained for the interpolation method after space discretization.

We next turn to the estimate of the computational cost, starting with the 
offline cost. The computation of each individual $t_{\nu,h}$, stored as vectors of dimension $N_h$
of their coordinates in the nodal finite element basis of $V_h$, requires
to solve a system. 
The cost of solving this system is of order $C_h$
where $C_h$ is the individual cost of one application often
the discrete solver. Indeeds it amounts in solving the a discrete problem where we invert
the exact same stiffness matrix as for the computation of the particular
instance $u_h(0)$. Assembling this system requires to compute the right
hand side which necessitates $\|\nu\|_0$ applications of
the stiffness matrices associated to the sesquilinear forms $B_j$.
Since $B_j$ is associated to a partial differential operator,
its stiffness matrices in the nodal basis is sparse
and therefore each such application has 
cost smaller of order $N_h$. We have already observed in \S \ref{secapriori} that 
since $\Lambda_n$ is a downward closed set, we have $2^{\|\nu\|_0}\leq n$ for
each $\nu\in\Lambda_n$. The cost of computing an individual $t_{\nu,h}$
is thus at most of the order
\be
C_{\rm off}(\nu)\sim C_h+\log(n) \,N_h.
\ee
In the non-adaptive algorithm, we compute the $n$ values of $t_{\nu,h}$
for $\nu\in \Lambda_n$, and therefore the total offline cost is
at most of order
\be
C_{\rm off}\sim n C_h +n\log (n) \,N_h.
\ee
In adaptive algorithms, we need to take into account the 
additional computation of the $t_{\nu,h}$ for $\nu$ in the margin of $\Lambda_n$. 
For the bulk chasing Taylor algorithm with $\e$ accuracy, the individual
cost $C_h+\log(n) \,N_h$ is thus multiplied by $n+\#(\t M_n)$
where 
\be
\#(\t M_n):={\rm SPARSE}(\Lambda_n,(t_{\nu,h})_{\nu\in \Lambda_n},\e).
\ee
For the SPARSE procedure that we have proposed in the case of the
elliptic problem \iref{ellip}, we have $\#(\t M_n)\leq Jn$ where
$J=J(\e)$ is such that
\be
\Big \|\sum_{j\geq J} |\psi_j| \Big\|_X \lsim \e.
\ee
Having assumed that $(\|\psi_j\|_X)_{j\in\N}$ is $\ell^p$ summable,
and organizing them such that this sequence is non-increasing, we find
by Lemma \ref{stechkin} that $J(\e)\lsim \e^{-1/s}$ where $s:=\frac 1 p-1$.
It follows that the total offline cost for this algorithm is at most of order
\be
C_{\rm off}\sim \e^{-1/s}n C_h +\e^{1/s}n\log (n) \,N_h.
\ee

As to the online cost, since the online stage simply amounts
in the combination of the $t_{\nu,h}$ for computing $u_{n,h}$, 
we find that this cost is of the order
\be
C_{\rm on} \sim nN_h,
\ee
similar to the sparse polynomial interpolation algorithms.

If $\e$ is a targeted order of accuracy, and if we have the error bound
\iref{nepsh}, then one way to reach this accuracy is to take both
$Cn^{-s}$ and $\e(h)$ of the order of $\e$. With $h(\e)$
the inverse function of $\e(h)$, as in \iref{heps},
we thus find that the non-adaptive Taylor algorithm reaches the order 
of accuracy $\e$ at cost at most of order
\be
C_{\rm off}(\e)\sim \e^{-1/s} C_{h(\e)}+\e^{-1/s}|\log(\e)|N_{h(\e)}\quad {\rm and}Ê\quad C_{\rm on}(\e)\sim \e^{-1/s}N_{h(\e)}.
\ee
For the bulk chasing Taylor algorithm with $\e$-accuracy, we have the more pessimistic estimate
\be
C_{\rm off}(\e)\sim \e^{-2/s} C_{h(\e)}+\e^{-2/s}|\log(\e)|N_{h(\e)}\quad {\rm and}Ê\quad C_{\rm on}(\e)\sim \e^{-1/s}N_{h(\e)},
\ee
due to the inflation by $J(\e)$.

Similar to the interpolation algorithm discussed in \S 6, 
both algorithms are immune to the curse of dimensionality
since these trade-off between accuracy and complexity
are obtained with infinitely many variables.

\section{Reduced basis methods}  
We turn next to the   class of numerical techniques for solving parametric PDEs known as 
reduced basis methods.  These method aim at finding a good 
subspace $V_n\subset V$, of small dimension $n$, to be used for approximating the elements of the solution manifold $\cM$.   We know that, for any fixed value of $n$, the best choice of $V_n$ is one which gives achieves the infimum in the definition \eref{kolnwidth} of the Kolmogorov $n$-width with $\cK=\cM$, however such a space, if it exists, is generally out of reach
from a computational point of view. The reduced basis method uses a  space $V_n$, which may be suboptimal, spanned by $n$ {\it snapshots}
$u(a^1),\dots,u(a^n)$ from the solution manifold $\cM$. While these snapshots can be chosen in various ways,
a particularly interesting strategy proceeds with a recursive greedy selection. We present this strategy
in \S \ref{subgreedy}. In \S \ref{subconvrbh} we prove, in the where $V$ is a Hilbert space,
that, in a certain sense, the resulting spaces $V_n$ perform  almost as well as the optimal $n$-width spaces
in terms of convergence rates. A similar analysis is given in \S \ref{subconvrbb} in the case of a general Banach space.
The effect of space discretization on the convergence of the algorithm
is discussed in \S \ref{subrbnum}, and computational cost is analyzed in \S \ref{subcost}.

\subsection{Greedy selection algorithms}  
\label{subgreedy}

The solution manifold $\cM$ is a compact set in the Banach space $V$.  While in most applications $V$ is a Hilbert space, we describe the greedy algorithm for any compact set $\cK$ in any Banach space   $V$. We then analyze
its performance, first in the case  $V$ is a Hilbert space, and then later in  the case of a general  Banach space.    We  describe two versions of a greedy algorithm for generating approximation spaces for $\cK$.  The first version, called the {\it pure greedy algorithm} is rather ideal, while the second version, called the
{\it weak greedy algorithm} is more amenable to numerical implementation.
\vskip .1in
\noindent
{\bf Pure Greedy Algorithm:} We first choose a function $g_0\in \cK$ such that
\be
\label{first}
{\displaystyle \|g_0\|_V =\max_{g\in\cK}\|g\|_V.}
 \ee
 Since $\cK$ is compact, such a $g_0$ always exists but of course may not be unique.
  Assuming $\{g_0,\dots,g_{n-1}\}$ have been selected, we set $V_n:=\span\{g_0,\dots,g_{n-1}\}$ and
  we then take $g_n\in\cK$ such that   
\be
\label{second}
 \dist( g_n,V_n)_V=\displaystyle{\max_{g\in \cK} \dist(g,V_n)_V},
\ee
where 
\be
\dist( g,V_n)_V:=\min_{h\in V_n} \|g-h\|_V.
\ee
We define $\sigma_0:=\sigma_ 0(\cK)_V=\displaystyle{\max_{g\in\cK}\|g\|_V}$ and
  \be
  \label{second1}
   \sigma_n:=\sigma_ n(\cK)_V:=\sup_{g\in\cK}\inf_{v\in V_n}\|g-v\|_V, \quad n\geq 1,
  \ee
so that
 \be
 \sigma_n :=\dist( \cK,V_n)_V= \dist( g_n,V_n)_V.
 \ee
 This greedy algorithm was introduced, for  the case when  $V$ is  a Hilbert space in \cite{VPRP} and subsequently extensively studied  in \cite {BMPPT, MPT, MPT1}.  
 \nl

In the setting of parametric PDEs, it not possible compute for a given $g\in \cK$ the 
distance $\dist( g,V_n)_V$, so 
that one cannot exactly perform the maximization in \iref{second}. However, it is possible to introduce 
a computable {\it error indicator}
$d(g,V_n)_V$ which satisfies 
\be
c d(g,V_n)_V\leq \dist( g,V_n)_V\leq Cd(g,V_n)_V, \quad g\in\cK,
\ee
for fixed constants $c,C>0$. Performing the maximization \iref{second} on $\d(g,V_n)_V$
is equivalent to the application, with $\gamma:=\frac c C$, of the following weaker form of the greedy algorithm which matches 
better its application.  
  \vskip .1in
 \noindent
 {\bf Weak Greedy Algorithm:}   We fix a constant $0<\gamma\le 1$.
  At the first step of the algorithm, 
one  chooses a function $g_0\in\cK$ such that
\beqn
\|g_0\|_V\ge \gamma \max_{g\in\cK}\|g\|_V.
\eeqn
 At the general step, if $g_0,\dots,g_{n-1}$ have been chosen,  we set
 $  V_{n}:=\span\{g_0,\dots,g_{n-1}\}$,  and we now
 choose $g_{n}\in \cK$ such that
\beqn
\label{gae1}
\dist(g_n,V_n)_V\ge \gamma \max_{g\in\cK} \dist(g,V_n)_V,
\eeqn
 to be the next element in the greedy selection. As in the  pure greedy case, we introduce
 \be
 \label{wgerror}
 \sigma_n:=\sigma_n(\cK)_V:=\dist(\cK,V_n)_V,\quad n\ge 0,
 \ee
 which now measures the performance of the weak greedy algorithm.
 Note that if $\gamma=1$, then the weak greedy algorithm 
reduces to the pure greedy algorithm that we have introduced above. With the same definition
as above for $\sigma_n:=\sigma_n(\cK)_V$, we thus have
\be
 \dist( g_n,V_n)_V \geq \gamma \sigma_n.
\ee

For both of these  algorithms,  the sequence $(\sigma_n)_{n\ge 0}$ is    monotone non-increasing.
It is also important to note that neither the pure greedy algorithm or the weak greedy algorithm   
give a unique sequence $(g_n)_{n\geq 0}$,
nor is the sequence $(\sigma_n)_{n\ge 0}$ unique.  In all that follows, the notation reflects any sequences which can arise in the implementation of the weak greedy selection for the fixed value of $\gamma$.

  \subsection{Convergence analysis of greedy algorithms in a Hilbert space}
  \label{subconvrbh}
  
  We are interested in how well the space $V_n$, generated by the weak greedy algorithm,  
  approximates  the elements of $\cK$.  For this purpose we would like to compare its performance measured by $\sigma_n$
   with the best possible performance which is given by the Kolmogorov width 
   \be
   d_n:=d_n(\cK)_ V.
   \ee 
If $(\sigma_n)_{ n\geq 0}$  were bounded by $(d_n)_{n\geq 0}$ up to a fixed multiplicative constant, this would mean that  the greedy selection provides essentially the best possible accuracy attainable by $n$-dimensional subspaces.   However, such
a general comparison is not to be expected.

Various  comparisons  between $\sigma_n$ and $d_n$ have been proven in the literature.    A first result in this direction, in the case  of the pure greedy algorithm applied to a Hilbert space  $V$,   was   given in \cite{BMPPT}  where it was proved that 
\beqn
\label{BMPPT}
\sigma_n(\cK)_V\le Cn2^nd_n(\cK)_V, \quad n\ge 1,
\eeqn
with $C$ an absolute constant. The
same result holds with $C$ depending on $\gamma$ for the weak greedy algorithm.
While this is an interesting comparison, it is only useful if $d_n(\cK)_V$ decays 
to zero faster than $n^{-1}2^{-n}$ which may be  a severe assumption. Unfortunately, the
above result is sharp in the following sense: it was proved in \cite{BCDDPW}
that for all $n\geq 1$ and $\e>0$ there exists a compact set $\cK$ such that 
\be
\sigma_n(\cK)_V\ge (1-\e)2^nd_n(\cK)_V.
\ee
This reveals that a direct comparison between $\sigma_n(\cK)_V$ and $d_n(\cK)_V$ is doomed to fail. 

Significant improvements on \eref{BMPPT} were given in \cite{BCDDPW}, again in the Hilbert space setting,
by changing the way of comparing $\sigma_n(\cK)_V$ and $d_n(\cK)_V$.  Perhaps the most interesting comparison is the following:   if for some constant $C>0$ and some $s>0$, the compact set $\cK$ satisfies $d_n(\cK)_V\le C(\max\{1,n\})^{-s}$ for all $n\geq 0$, then there is a constant $\t C$ which depends only on $C$ and $s$ such that
 \be
 \label{poly1}
 \sigma_n(\cK)_V\le \t C(\max\{1,n\})^{-s},\quad n\ge 0.
 \ee
 In other words,  for  the scale of polynomial decay,  the greedy algorithm performs with the same decay rates as $n$-widths.
 These results were improved upon in \cite{DPW2} and extended to the case of a general Banach space $V$.

 The analysis of the two greedy algorithms above is quite simple and executed with elementary results from linear algebra.  
We consider the case when $V$ is a Hilbert space and show that the action of the weak greedy algorithm is captured by a certain lower triangular matrix.
Note that in general, the weak greedy algorithm does not terminate and we obtain an infinite sequence $(g_n)_{n\geq 0}$.
In order to have a consistent notation in what follows, we define $g_n:=0$, $n>m$, if the algorithm terminates at $n=m$, i.e. if $\sigma_{m}(\cK)_V=0$.

By $(g_n^*)_{n\geq 0}$ we denote the orthonormal system
obtained from  $(g_n)_{n\geq 0}$ by Gram-Schmidt orthogonalization executed in the natural order.  It follows that  the orthogonal projector $P_n$ from $V$ onto $V_n$ is given by
$$
P_n g=\sum_{i=0}^{n-1}\langle g, g_i^*\rangle g_i^*,
$$
where $\<\cdot,\cdot\>$ denotes the inner product of $V$, and, in particular,
\beqn
\nonumber
g_n = P_{ n+1} g_n =\sum_{j=0}^n a_{n,j} g^*_j,\quad a_{n,j}=\langle g_n,g^*_j\rangle, \,\, 0\le j\leq n.
\eeqn

We consider the infinite lower triangular matrix
$$
A:= (a_{i,j})_{i,j=0}^\infty,\quad a_{i,j}:=0,\, j>i.
$$
This matrix incorporates all the   information about the  weak greedy algorithm on $\cK$.  For example, the $n$-th row of $A$
gives the $n$-th element $g_n$ in the greedy selection. 
The following two properties characterize any lower triangular matrix $A$ generated by the weak greedy algorithm with constant $\gamma$.   With the notation $\sigma_n:= \sigma_n(\cK)_V$, we have:
\vskip .1in
\noindent
 {\bf P1:}  The diagonal elements of $A$  satisfy
$\gamma \sigma_n \leq |a_{n,n}|\leq \sigma_n$.
\vskip .1in
\noindent
  {\bf P2:} For every $m\ge n$, one has
$\sum_{j=n}^m a_{m,j}^2\leq \sigma_n^2$.
\vskip .1in

\noindent
Indeed, {\bf P1} follows from
$$
a_{n,n}^2 = \|g_n-P_{n}g_n\|^2_V=\dist(g_n,V_n)^2_V.
$$
This shows the upper bound in {\bf P1} because each element of $\cK$ is approximated to error $\sigma_n$.  It also shows the lower bound  because of   the weak greedy selection property \eref{gae1}.
To see {\bf P2}, we note that for  $m\ge n$,
$$
\sum_{j=n}^m a_{m,j}^2=\|g_m-P_{n}g_m\|^2_V\leq \max_{g\in\cK} \|g-P_{n}g\|^2_V =\sigma_n^2.
$$

\begin{remark}
\label{remP1P2}
 If $A$ is any infinite  matrix satisfying {\bf P1} and {\bf P2} with $ (\sigma_n)_{n\geq 0}$ a  non-increasing sequence that converges to $0$, then
 the rows of $A$ form a compact subset of $\ell^2(\N_0)$ where $\N_0:=\N\cup\{0\}$. 
  If $\cK$ is the set consisting of these rows, 
 then one of the possible realizations of the weak greedy algorithm on this set $\cK$ with constant $\gamma$
 will choose the rows in  that order and $A$ will be the resulting  matrix.  In this sense, the action of the greedy algorithm on the original set $\cK$ is completely described by the matrix $A$.
 \end{remark}

It follows from the above remark that there is no loss of generality in assuming that the infinite dimensional Hilbert space $V$ is $\ell^2(\N_0)$
and that $g_j^*=e_j$, where $e_j$ is the vector with a one in the coordinate indexed by $j$ and  is zero in all other coordinates, i.e. $(e_j)_i =\delta_{j,i}$.

  With this matrix description of the weak greedy algorithm in hand, estimates for the convergence rate of the algorithm rely
  on an analysis of $A$ or a corresponding matrix when $V$ is not necessarily Hilbertian.   Notice that the diagonal elements of $A$ give the errors $\sigma_n$ and hence we want a general way to bound the diagonal elements of matrices $A$ with the above properties.  The following lemma from \cite{DPW2} gives a general way to bound diagonal elements of a general lower triangular matrix $G$.   It is applied later to the sections of $A$ to obtain convergence results for the weak greedy algorithm.

 \begin{lemma}
 \label{L1}
 Let $G=(g_{i,j})$ be a $K\times K$  lower triangular matrix with rows ${\bf g}_1, \ldots,{\bf g}_K$.  If $W$ is any $m$ dimensional 
 subspace of $\R^K$ for some $0<m\leq K$,  and $P$ is the orthogonal projection from $\R^K$ onto $W$, then
 \be
 \label{L11}
\det(G)^2=  \prod_{i=1}^Kg_{i,i}^2\le  \(\frac{1}{m}\sum_{i=1}^K \|P{\bf g}_i\|_{\ell^2}^2\)^m\( \frac{1}{K-m}\sum_{i=1}^K \|{\bf g}_i-P{\bf g}_i\|_{\ell^2}^2\)^{K-m},
 \ee
 where $\|\cdot\|_{\ell^2}$ is the euclidean norm of a vector in $\R^K$.
 \end{lemma}
 
 \noindent
 {\bf Proof:}  Let ${\bf \varphi}_1,\dots,{\bf \varphi}_m$ be any orthonormal basis for the  space $W$ and complete it into an orthonormal  basis 
 ${\bf \varphi}_1,\dots,{\bf \varphi}_K$ for $\R^K$.  If we denote by $\Phi$ the $K\times K$ orthogonal matrix whose $j$-th column is ${\bf \varphi}_j$, then   
 the matrix $C:=G\Phi$ has entries  $c_{i,j}=\langle {\bf g}_i,{\bf \varphi}_j\rangle$.   
  We denote by ${\bf c}_j$, the $j$-th column of $C$.   It follows from the arithmetic geometric mean inequality 
  for the numbers $\{\|{\bf c}_j\|_{\ell^2}^2\}_{j=1}^m$ that
 \be
 \label{fc}
\prod_{j=1}^m \|{\bf c}_j\|_{\ell^2}^2\le \(\frac{1}{m}\sum_{j=1}^m \|{\bf c}_j\|^2_{\ell^2}\)^m=
\(\frac{1}{m}\sum_{j=1}^m \sum_{i=1}^K\langle {\bf g}_i,\varphi_j\rangle^2\)^m=
\(\frac{1}{m}\sum_{i=1} ^K\|P{\bf g}_i\|_{\ell^2}^2\)^m.
 \ee
 Similarly, 
 
  \be
 \label{lc}
 \prod_{j=m+1}^K \|{\bf c}_j\|_{\ell^2}^2\le \(\frac{1}{K-m} \sum_{j=m+1}^K\|{\bf c}_j\|_{\ell^2}^2\)^{K-m}
 = \(\frac{1}{K-m} \sum_{i=1}^K\|{\bf g}_i-P{\bf g}_i\|_{\ell^2}^2\)^{K-m},
 \ee
 where we have used the fact that $\varphi_j$ is orthogonal to $W$ when $j>m$.
Now, we invoke Hadamard's inequality for the matrix $C$, which says that
\be
(\det C)^2 \leq \prod_{j=1}^K\|{\bf c}_j\|^2_{\ell^2},
\ee
and combine it with relations \eref{fc} and \eref{lc} to obtain
\be
\label{T13}
(\det C)^2\
\le\(\frac{1}{m}\sum_{i=1}^K \|P{\bf g}_i\|_{\ell^2}^2\)^m\( \frac{1}{K-m}\sum_{i=1}^K\|{\bf g}_i-P{\bf g}_i\|_{\ell^2}^2\)^{K-m}.
\ee
 The latter inequality and  the fact that  $|\det C| =|\det G|$ gives \eref{L11}.
$\hfill\Box$
 \nl
 
Let us now see how this lemma is utilized to derive convergence results for the greedy algorithm.  We continue to restrict ourselves to the case of a Hilbert space and the weak greedy algorithm with constant $\gamma$.  Later, we indicate how the results change when $V$ is a general Banach space.  The following theorem, taken from \cite{DPW2},
relates  the errors $\sigma_n(\cK)_V$ to the $n$-widths $d_n(\cK)_V$.

  \begin{theorem}
 \label{T10}  For the weak greedy algorithm with constant $\gamma$ in a Hilbert space  $V$ and for any compact set $\cK$, the following inequalities between $\sigma_n:= \sigma_n(\cK)_V$ and $d_n:=d_n(\cK)_V$
 hold for any $N\ge 0$, $K\ge 1$, and $1\le m< K$, 
 \be
 \label{T1}
 \prod_{i=1}^K\sigma^2_{N+i} \le \gamma^{-2K}\(\frac{K}{m}\)^m\(\frac{K}{K-m}\)^{K-m}
 \sigma_{N+1}^{2m}d_m^{2K-2m}
 \ee

 \end{theorem}
 
 \noindent
 {\bf Proof:}  
  In what follows, we assume that there exists a space $W_m$ 
which achieves the infimum in the definition of the $m$-width of $\cK$, that is, such that
\be
\max_{g\in \cK} \min_{ w\in W_m}\|g-w\|_{V}=d_m.
\label{dmspace}
\ee
If such a space does not exist, we may, for each $\e>0$, find one such that 
$$\max_{g\in \cK} \min_{ w\in W_m}\|g-w\|_{V}\leq d_m+\e$$
 and
modify the proof below by a limiting argument so as to reach the same conclusion.

We consider  the $K\times K$ matrix $G=(g_{i,j})$ which is  formed by the rows and columns of $A$ with indices from $\{N+1,\dots,N+K\}$.  
 Each row ${\bf g}_i$ is the restriction of row  $N+i$  of $A$ to the coordinates $N+1,\dots,N+K$.  
The space $W_m$ determines a sequence space $\bar W_m\subset \ell^2$ such that  $\dist({\bf g}_{i},\bar W_m)_{\ell^2}\leq d_m$ for $i=1,\dots K$. Let $\t W$ be the linear space which is the restriction of $W_m$ to the coordinates $N+1,\dots,N+K$.  Obviously, we have $\dim (\t W)\leq m$. Let $W$ be an  $m$ dimensional space, $W\subset \span\{e_{N+1},\dots,e_{N+K}\}$,  such that $\t W\subset W$ and let $P$ and $\t P$ be the projections in $\R^K$ onto $W$ and $\t W$, respectively.
 Clearly, 
\be
 \label{also}
 \|P{\bf g}_i\|_{\ell^2}\le \|{\bf g}_i\|_{\ell^2}\le \sigma_{N+1},\quad i=1,\ldots,K,
 \ee 
 where we have used Property  {\bf P2} in the last inequality.
 Note that
 \be
 \label{alsoo}
\|{\bf g}_i-P{\bf g}_i\|_{\ell^2}\leq \|{\bf g}_i-\t P{\bf g}_i\|_{\ell^2}=\dist({\bf g}_i,\t W)_{\ell^2}\le 
 \dist({\bf g}_i,W)_{\ell^2} \leq d_m, \quad i=1,\dots , K.
 \ee
 It follows from   Property {\bf P1} that
   \be
 \label{Tha}
 \gamma^K \prod_{i=1}^K\sigma_{N+i}\le  \prod_{i=1}^K|a_{N+i,N+i}|.
 \ee
We now  apply Lemma \ref{L1} for this $G$ and $W$, and use estimates \eref{also}, \eref{alsoo},  and \eref{Tha} to derive \eref{T1}. 
 $\hfill\Box$
 \nl
  
Using this theorem, we now establish convergence results for the weak greedy algorithm,
showing in particular that if $d_n(\cK)_V$ decays with an algebraic or exponential convergence
rates, then a similar rate holds for $\sigma_n(\cK)_V$.

 \begin{cor}
 \label{C1}
 For the weak greedy algorithm  with constant $\gamma$ in a Hilbert space $V$, we have the following:
 \vskip .1in
 \noindent
 {\rm (i)} For any compact set $\cK$, we have
 \be
 \label{C11}
 \sigma_n(\cK)_V\le  \sqrt{2}\gamma^{-1}d_0(\cK)_V^{\frac{m}{n}}\min_{1\le m< n} d_m^{\frac{n-m}{n}}(\cK)_V,\quad n\ge 1.
\ee 
In particular $\sigma_{2n}(\cK)_V\le \gamma^{-1}\sqrt{2d_0(\cK)_Vd_n(\cK)_V}$ for all $n\ge 1$.
\vskip .1in
\noindent
{\rm (ii)} For any $s>0$ and $C_0>0$,
\be
d_{n}(\cK)_V\le C_0(\max\{1,n\})^{-s}, \quad n\ge 0\quad \Rightarrow \quad \sigma_n(\cK)_V\le C_1(\max\{1,n\})^{-s}, \quad n\ge 0,
\ee 
where $C_1:= \gamma^{-2} 2^{4s+1}C_0$. 
\vskip .1in
\noindent
{\rm (iii)} For any $s>0$ and $c_0,C_0>0$, 
\be 
d_{n}(\cK)_V\le C_0e^{-c_0n^{s}}, \quad n\ge 0\quad \Rightarrow \quad \sigma_n(\cK)_V\le \t C_1e^{-c_1n^{s}}, \quad n\ge 0,
\ee 
where $c_1=\frac {c_0} 2 3^{-s}$ and $\t C_1:=C_0\max\{\sqrt{2}\gamma^{-1},e^{c_1}\}$.
\end{cor}
 
 \noindent
 {\bf Proof:}  (i) We take $N=0$, $K=n$ and any $1\le m< n$ in Theorem \ref{T10}. Using the monotonicity of 
 $(\sigma_n)_{n\geq 0}$ and the fact that $\sigma_1\le \sigma_0\le d_0$, we obtain
 \be
 \label{C14}
 \sigma_n^{2n}\le \prod_{j=1}^n \sigma_j^2\le  \gamma^{-2n} \(\frac{n}{m} \)^m \(\frac{n}{n-m}\)^{n-m}d_m^{2n-2m}\d_0^{2m}.
 \ee
 Since $x^{-x}(1-x)^{x-1}\le 2$ for $0< x< 1$, we derive \eref{C11}.
 \vskip .1in
 \noindent
 (ii)   It follows from the monotonicity of $(\sigma_n)_{ n\geq 0}$ and  \eref{T1} for $N=K=n$ and any $1\leq m<n$
 that
 $$
  \sigma_{2n}^{2n}\le  \prod_{j=n+1}^{2n}\sigma_j^2
\le \gamma^{-2n}\(\frac{n}{m}\)^m\(\frac{n}{n-m} \)^{n-m}\sigma_{n}^{2m}d_m^{2n-2m}.
 $$
In the case $n=2k$ and  $m=k$ we have for any positive integer $k$,
 \be
 \label{Cnew}
 \sigma_{4k}\leq \sqrt{2}\gamma^{-1}\sqrt{\sigma_{2k}d_k}.
\ee 
Assuming that $d_n(\cK)_V\le C_0(\max\{1,n\})^{-s}$ for all $n\ge 0$, 
we obtain by induction that for all $j\geq 0$ and $n=2^{j}$,
\be
\sigma_n =\sigma_{2^j}\leq C2^{-sj} \leq n^{-s}, \quad C:=\ 2^{3s+1} \gamma^{-2}C_0.
\ee
Indeed, the above obviously holds for $j=0$ or $1$  since for these values, we have 
\be
\sigma_{2^j}\leq \sigma_0= d_0\le C_0 \leq C2^{-sj}.
\ee
Assuming 
its validity for some $j\geq 1$, we find that
$$
\begin{disarray}{ll}
\sigma_{2^{j+1}} & \leq \sqrt{2}\gamma^{-1}\sqrt{\sigma_{2^j}d_{2^{j-1}}} \\
& \leq \gamma^{-1}  2^{\frac {3s}2} \sqrt {2CC_0} 2^{-s(j+1)}\\
& = \sqrt C\sqrt{2^{3s+1}C_0\gamma^{-2}}2^{-s(j+1)} = C2^{-s(j+1)},
\end{disarray}
$$
where we have used the definition of $C$.
For values $2^j<n<2^{j+1}$, we obtain the general result by writing
\be
\sigma_{n} \leq \sigma_{2^{j}} \leq C2^{-sj} \leq 2^sCn^{-s} = C_1n^{-s}.
\ee
In the case $n=0$, we simply write $\sigma_0=d_0\leq C_0\leq C_1$.
 \vskip .1in
 \noindent
 (iii) Assuming that $d_{n}(\cK)_V\le C_0e^{-c_0n^{s}}$ for all $n\ge 0$, we obtain from (i) for all $n\geq 1$,
 \be
 \label{C17}
 \sigma_{2n+1} \le \sigma_{2n} \le \sqrt{2}\gamma^{-1}\sqrt{d_nd_0}\le \sqrt{2C_0d_0}\gamma^{-1}e^{-\frac{c_0}{2}n^s}\leq \sqrt{2}C_0\gamma^{-1}
 e^{-\frac{c_0}{2}3^{-s}(2n+1)^s}.
 \ee
 This proves 
 \be
 \sigma_n \leq \t C_1e^{-c_1 n^s}, \quad n\geq 2.
 \ee
For the values $n=0$ or $n=1$, we simply write
\be
\sigma_n\leq \sigma_0=d_0 \leq C_0 \leq \t C_1 e^{-c_1} \leq  \t C_1 e^{-c_1n^s},
\ee
which concludes the proof of (iii). \hfill $\Box$

\begin{remark}
\label{remrange}
Inspection of the above proof shows that in items {\rm (ii)} and {\rm (iii)}
of Corollary \ref{C1}, if the same decay rates of $d_n(\cK)_V$ are 
only assumed within a limited range $0\leq n\leq N$, then the same
decay rates of $\sigma_n(\cK)_V$ are achieved for the same rate $0\leq n\leq N$
up to some changes in the expressions of the constants $C_1$, $c_1$ and $\t C_1$.
\end{remark}

\subsection{Convergence analysis of greedy algorithms in a Banach space}
\label{subconvrbb}

We   now turn our attention to  the performance of the  weak  greedy algorithm for a compact set $\cK$ in a general Banach space $V$.    
We use the abbreviation $\sigma_n:= \sigma_n(\cK)_V$ and $d_n:=d_n(\cK)_V$. While the development is quite similar to the case of a 
Hilbert space, there is a slight loss in the comparison between $\sigma_n$ and $d_n$ due to the lack of Hilbert space orthogonality.

As in the Hilbert space case, we  associate with the greedy procedure  a lower triangular matrix $A=(a_{i,j})_{i,j=0}^\infty$  in  the following way.
For each $j=0,1,\dots$, we let $\lambda_j\in V^*$ be  the linear functional of norm $\|\lambda_j\|_{V^*}=1$ that  satisfies
 \be
 \label{functionals}
  {\rm (i)}\  \lambda_j(g)=0, \quad g\in V_j,\quad  {\rm and}Ê\quad {\rm (ii)}\ \lambda_j(g_j)=\dist(g_j,V_j)_V. 
\ee
 The existence of such a functional is a simple consequence of the Hahn-Banach theorem.   We now let $A$ be the matrix
  with  entries 
 $$
 a_{i,j}=\lambda_j(g_i).
 $$
From (ii) of \eref{functionals}, we see that $A$ is lower triangular. Its  diagonal elements $a_{j,j}$ satisfy the inequality
\be
\label{diag}
\gamma \sigma_j\leq  a_{j,j}=\dist(g_j,V_j)_V=\sigma_j,
 \ee
 because of the weak greedy selection property \eref{gae1}.
   Also, each entry $a_{i,j}$ satisfies
 $$
  |a_{i,j}|=|\lambda_j(g_i)| = |\lambda_j(g_i-g)|\leq  \|\lambda_j\|_{V^*}\|g_i-g\|_V=\|g_i-g\|_V, \quad j<i,
 $$
 for every $g\in V_j$, since $\lambda_j(V_j)=0$. Therefore, we have
 \be
 \label{entries}
 |a_{i,j}|\leq   \dist(g_i,V_{j})_V\le \sigma_{j}, \quad j<i.
 \ee

\begin{theorem}
 \label{T20}
 For the weak greedy algorithm with constant $\gamma$ in a Banach space $V$ and for any compact set $\cK$ contained in $V$, we have the following inequalities between $\sigma_n:= \sigma_n(\cK)_V$ and $d_n:=d_n(\cK)_V$:  for 
 any $N\ge 0$, $K\ge 1$, and $1\le m < K$, 
 \be
 \label{T3}
 \prod_{i=1}^K\sigma_{N+i}^2\le 2^K K^{K-m}\gamma^{-2K} \left ( \sum_{i=1}^K \sigma_{N+i}^2\right )^m d^{2K-2m}_m.
  \ee
 \end{theorem}
 
 \noindent
 {\bf Proof:} 
  As in  the proof of Theorem \ref{T10}, we consider  the   $K\times K$ matrix $G$ which is  formed by the rows and columns of 
 $A$ with indices from $\{N+1,\dots,N+K\}$. 
 Let $V_m$ be  a Kolmogorov subspace of $V$ for which $\dist(\cK,V_m)_V= d_m$.   Again, we assume that such a space $V_m$ exists.  Otherwise we modify the proof given below by adding an arbitrary $\e>0$ to $d_m$ and then letting $\e$ tend to zero at the end.
 
 For each $i$, there is an element
  $h_i\in V_m$ such that 
  $$
  \|g_i-h_i\|_V=\dist(g_i,V_m)_V\le d_m,
  $$  
   and therefore
  \be
  \label{approx1}
  |\lambda_j(g_i)-\lambda_j(h_i)|=|\lambda_j(g_i-h_i)|\le \|\lambda_j\|_{V^*}\|g_i-h_i\|_V\leq d_m.
  \ee

  We now consider the collection of  vectors
  $(\lambda_{N+1}(h),  \dots, \lambda_{N+K}(h))$ for all  $h\in V_m$. They span a space $W_m\subset \R^K$ of dimension at most $m$.  We  assume that $\dim(W_m)=m$ (a slight notational adjustment has to be made if $\dim(W_m)<m$ without affecting the final result).   
 It follows from \eref{approx1} that each row ${\bf g}_i$ of $G$  can be approximated by a vector  from  $W_m$ in the $\ell^\infty$ norm to accuracy $d_m$, and therefore    in the $\ell^2$ norm to accuracy  $\sqrt{K}d_m$.
 Let $P$ be the orthogonal projection of $\R^K$ onto $W$.   Hence, we have  
 \be
 \label{also1}
\|{\bf g}_i-P{\bf g}_i\|_{\ell_2}\le \sqrt{K} d_m, \quad i=1,\dots , K.
 \ee
It also follows from \eref{entries} that 
\be
\nonumber
 \|P{\bf g}_i\|_{\ell_2}\le \|{\bf g}_i\|_{\ell_2}\le \left (\sum_{j=1}^i\sigma^2_{N+j}\right )^{1/2},
\ee
and therefore
\be
 \label{also2}
\sum_{i=1}^K \|P{\bf g}_i\|_{\ell_2}^2\le \sum_{i=1}^K\sum_{j=1}^i\sigma^2_{N+j}\leq K \sum_{i=1}^K \sigma^2_{N+i}.
\ee
Next, we apply Lemma \ref{L1} for this $G$ and $W$ and use estimates \eref{diag}, \eref{also1} and \eref{also2} to derive 
\begin{eqnarray}
\nonumber
 \gamma^{2K}\prod_{i=1}^K\sigma_{N+i}^2&\leq&\left (\frac{K}{m}\sum_{i=1}^K \sigma^2_{N+i}\right )^m\left ( \frac{K^2}{K-m}d^2_m\right )^{K-m}\\
\nonumber
&=& K^{K-m}\left (\frac{K}{m}\right )^m\left(\frac{K}{K-m}\right)^{K-m}
\(\sum_{i=1}^K\sigma_{N+i}^2\)^m d^{2(K-m)}_m\\
\nonumber
&\le &2^KK^{K-m}\left(\sum_{i=1}^K\sigma_{N+i}^2\right)^m d^{2(K-m)}_m ,
\end{eqnarray}
 and the proof is complete.
\hfill $\Box$
\nl

In analogy with Corollary \ref{C1}, we can use the above result to establish convergence
theorem for the weak
greedy algorithm in a general Banach space. Since the proof is very similar to that of Corollary \ref{C1},
except that we use \eref{T3} in place of \eref{T1}, we only state the result and refer to \cite{DPW2} for more details.

 \begin{cor}
 \label{C2}
 Suppose that $V$ is a Banach space.  For the weak greedy algorithm  with a constant $\gamma$, applied to a compact set $\cK\subset V$, 
 we have the following:
 \vskip .1in
 \noindent
 {\rm (i)} For any  $n\ge 1$, we have
 \be
 \label{C21}
 \sigma_n(\cK)_V\le  \sqrt{2}\gamma^{-1}\min_{1\le m< n} n^{\frac{n-m}{2n}}\left (\sum_{i=1}^n\sigma_i(\cK)_V^2\right )^{\frac{m}{2n}}d_m(\cK)_V^{\frac{n-m}{n}}.
\ee 
In particular $\sigma_{2n}(\cK)_V\le 2\gamma^{-1} \sqrt{n d_0(\cK)_Vd_n(\cK)_V}$ for all $n\geq 1$.
 \vskip .1in
\noindent
{\rm (ii)} For any  $s>0$, $C_0>0$ and $\e>0$, we have
 \be
 d_{n}(\cK)_V\le C_0(\max\{1,n\})^{-s}, \; n\geq 0  \; \Rightarrow \;\sigma_n(\cK)_V\le C_1(\max\{1,n\})^{-(s-\e-1/2)}, \; n\geq 0, 
 \ee
 where $C_1$ depends on $C_0$, $s$, $\gamma$ and $\e$.
\vskip .1in
\noindent
{\rm (iii)} For any  $s>0$ and $c_0,C_0>0$, we have
\be
d_{n}(\cK)_V\le C_0e^{-c_0n^{s}}, \quad n\geq 0 \quad \Rightarrow\quad    
 \sigma_n(\cK)_V \leq  \t C_1 e^{-c_1n^s}, \quad n\geq 0, 
 \ee
where $c_1$ depends on $s$ and $c_0$, and where $\t C_1$ depends on $C_0$, $\gamma$, $s$, and $c_0$.  
\end{cor}
  
 The statement (ii) in the above corollary shows that there is a loss of $\frac 1 2$ in the algebraic rate
 of decay of $\sigma_n$ compared to that of $d_n$.
 It is natural to ask whether this loss is unavoidable when proving results in a Banach space.  
We next  provide an example which shows that a loss of this type is in general unavoidable.   However,
  there is still a small gap between the above corollary and what the
 example below provide.

Let us begin by considering the space $V=\ell^\infty (\N\cup\{0\})$ equipped with its usual norm.  We consider a monotone non-increasing sequence $x_0\ge x_1\ge x_2\ge \cdots$
of positive real numbers   which converges to zero and we define  
\be
\nonumber
f_j:=x_je_j,\quad j=0,1,\dots,
\ee
where $e_j$ is the Kroenecker sequence with $1$ at position $j$.  We consider the compact set 
\be
\cK:=\{f_0,f_1,\dots \}.
\ee  
From the monotonicity of the $x_j$'s, the greedy algorithm for $\cK$ in $X$ can  
choose the elements from $\cK$ in the natural order $f_0,f_1,\dots$.  Hence, 
\be
\nonumber
\sigma_j=\sigma_j(\cK)_{V}=x_j,\quad j\ge 0.
\ee
We want to give an upper bound for the Kolmogorov width of $\cK$.  For this, we shall use the following
result (see (7.2) of Chapter 14  in  \cite{LGM}) on  Kolmogorov $n$-widths of the $m$-dimensional unit ball $b_1^m$ of $\ell^1$ in 
the $\ell^\infty$ metric, in $\R^m$:
\be
\label{pajor}
d_n(b_1^m)_{V}\le C_0 \(\log_2 (m/n)\)^{1/2} n^{-1/2}, \quad 1\le n\le m/2.
\ee
Let us now  define  the  sequence  $\{x_j\}_{j\geq 0}$ so that in position $2^{k-1}\le j\le 2^{k}-1$ it has the constant value $2^{-ks}$, for 
all $k\geq 0$, where $s>1/2$.  It follows that
\be
\nonumber
\sigma_n(\cK)_{V} \geq c n^{-s}, \quad n\geq 1,
\ee
for some $c>0$.
We now bound the $n$-width of $\cK$ when $n=2^{k+2}$ by constructing a good space $V_n$ 
of dimension at most $n$ for approximating $\cK$.   The space $V_n$ is defined 
as the span of a set $E$ of at most $n$ vectors which we construct as follows.  First, we  place into $E$ all of the vectors, $e_0,e_1,\dots,e_{2^k-1}$.  Next, for each $j=0,1,\dots n$, we  use \eref{pajor}  to choose a    basis for the  space of dimension
$2^{n-j} $ whose vectors are supported on $[2^{k+j},2^{k+j+1}-1]$ and this space approximates in $V$ each of the $f_i$ for 
$i=2^{k+j},\dots,2^{k+j+1}-1$,  to accuracy 
$$C_02^{-(k+j+1)s} \sqrt{2j}2^{-(k-j)/2}\le C_02^{-(k+j)s} \sqrt{j}2^{-(k-j)/2},$$
where we used the fact that $s>1/2$.
We place these basis vectors into $E$ so that
\be
\#(E)\leq 2^k+2^{k+1}-1 \leq n.
\ee
Notice that $|x_i|\le 2^{-2ks}$ for $i\ge 2^{2k}$. This means that for the space $V_n:=\span(E)$, with $n=2^{k+2}$,
\begin{eqnarray}
\nonumber
d_n(\cK)_{V} &\le& \dist (\cK,V_n)_V\le \max  \left\{2^{-2ks}, \max_{1\le j\leq n}C_02^{-(k+j)s} 2^{-(k-j)/2} \sqrt{
j}\right\}\nonumber \\
&=& \max  \left\{2^{-2ks}, C_02^{-k(s+1/2)}\cdot\max_{1\le j\leq k}2^{-j(s-1/2)}\sqrt{j}\right\}
\le C_1 n^{-(s+1/2)}.
\nonumber
\end{eqnarray}
From the monotonicity of $(d_n(\cK)_V)_ {n\geq 0}$, we obtain that
\be 
\nonumber
d_n(\cK)_{V}  \le C_2 n^{-s-1/2},\quad n\geq 1.
\ee
This example shows that the loss of $\frac 12$ which appears in (ii) of Corollary \ref{C2} can in general not be avoided.

\subsection{Space discretization and convergence analysis}
\label{subrbnum}

The   greedy algorithms introduced in the previous section are at this stage only theoretical algorithms 
 because they involve several steps that cannot be implemented numerically.   To describe a numerical version of these algorithms that are applicable to solving parametric PDEs, we place ourselves
 in the following numerical setting.   We assume that we are given a target accuracy $\e>0$ and
we wish to find a space $V_n=\span\{g_1,\dots,g_n\}$ where $n=n(\e)$ such that
\be
\label{goal}
\dist(\cM,V_n)_V:=\max_{v\in\cM} \dist(v,V_n)_V=\sup_{a\in\cA}  \dist(u(a),V_n)_V \le \e,
\ee
and of course we want $n$ to be small. In the reduced basis method, the functions $g_i$ are picked from
the solution manifold $\cM$, or equivalently, are of the form
\be
g_i=u(a^i),
\ee
where $\{a^1,\dots,a^n\}$ are picked from the parameter set $\cA$. Our benchmark is given by the $n$-width of $\cM$.
Namely, we know that as soon as  $d_n(\cM)_V\le \e$ then there is a space of  this dimension $n$ which
satisfies \eref{goal}.   We have seen that the theoretical greedy algorithms  also
give us such a space $V_n$  with provable bounds on performance, namely with rate guarantees on the growth of $n$ with respect to
$\e$ comparable to the $n$-width, as expressed by Corollaries \ref{C1} and \ref{C2}.
However, the greedy algorithm as it stands cannot be implemented numerically for several reasons that we now delineate.  
\nl

 \noindent
 {\bf Issue 1: Computing the greedy selection $g_k$:}   {\it Once the parameter $a^k$ of the $k$-th greedy selection is identified,  the function $g_k:=u(a^k)$ cannot be computed exactly. In practice, it is computed approximately by space discretization in the finite element method
 in the space $V_h$.}
 \nl

 This means that we take 
 \be
 \label{discreteh}
 g_k=u_h(a^k)\in V_h,
 \ee
 and so the spaces $V_n$ are subspaces of $V_h$.
As explained further, this may be viewed as applying the weak greedy algorithm to the {\em approximate
solution manifold} defined as
\be
\cM_h:=\{u_h(a)\; : \; a\in \cA\}.
\ee  
Recall that 
\be
{\rm dist}(\cM_h,\cM)_V=\max_{a\in\cA}\|u(a)-u_h(a)\|_V \leq \e(h). 
\ee
where $\e(h)$ is the accuracy of the numerical solver. In order to reach the goal \eref{goal},
we pick $h$ such that $\e(h)\leq \e/3$.
\nl

 \noindent
 {\bf Issue 2: Search over the manifold $\cM_h$.}  {\it The $k$-th greedy step requires a search over the entire manifold $\cM_h$ to choose the next basis function $g_k$.
Since the manifold is typically an infinite set,  this search has to be discretized. }
\nl

One way to handle this issue is by finding a finite set $\cM_{h,\e}\subset \cM_h$ such that
 each element in $\cM_h$ is at distance at most $\e/3$ from $\cM_{h,\e}$, i.e.
 \be
  \label{discret}
 \sup_{a\in\cA}  \dist(u(a),\cM_{h,\e})_V \le  \e /3.
 \ee
 In practice this discretization is done on the parameter side so that each $v\in\cM_\e$ is of the form $u(a)$, $a\in\cA_\e$, where $\cA_\e$ is a finite subset of $\cA$.  
 If we apply the weak greedy algorithm to  $\cM_{h,\e}$ until we are guaranteed that the resulting space $V_n$ satisfies $\dist(\cM_{h,\e},V_n)_V\le \e/3$, then we are guaranteed that
 the goal \eref{goal} is met since
 \be
 \dist(\cM,V_n)_V\leq \dist(\cM,\cM_{h})_V+\dist(\cM_{h},\cM_{h,\e})_V+\dist(\cM_{h,\e},V_n)_V\leq \e.
 \ee
 
 \vskip .1in
 \noindent{\bf Issue 3:  Computation of $\dist (u_h(a),V_k)_V$ for $a\in\cA$ (or $a\in\cA_\e$).} {\it At each iteration $k$  of the
 greedy algorithm,  we need to compute $\dist(u_h(a),V_k)_V$ to a sufficient accuracy so when selecting $g_k$ based on these computed distances we are certain that the weak greedy criterion \eref{gae1} is satisfied. }
\nl

Here, we want to avoid computing $u_h(a)$ itself since this is too costly
and must be done many times, i.e. for each $a\in\cA_\e$.   Instead, this computation is done
by a surrogate $d(a,V_k)_V$ which is typically evaluated by a residual-based a posteriori analysis
from the Galerkin approximation to $u_h(a)$ from $V_k$. This surrogate satisfies  
\be
 \label{surrogate1}
 \delta d(a,V_k)_V\le \dist(u_h(a),V_k)_V\le \beta d(a,V_k)_V.
 \ee
A practical construction of this surrogate is discussed further in the particular 
case of the elliptic problem \iref{ellip}. It follows that maximizing this surrogate in place of the true error 
amounts in applying the weak greedy algorithm with constant $\gamma:=\frac {\delta}{\beta}$ to
the approximate solution manifold $\cM_h$.
  \nl
   
We can now put together the proposed  solutions to each of the stated numerical issues 1, 2 and 3,
and form the following numerical version of the weak greedy algorithm.
 
 \vskip .1in
 \noindent
 {\bf Numerical Weak Greedy Algorithm:}  {\it  We assume we are given  a numerical tolerance $\e$ and that for each subspace $V_n$ of $V_h$, we have, in hand, a surrogate $d(a,V_n)$ 
 which satisfies \eref{surrogate1} with uniform constants $\delta,\beta$.  We first construct a set $\cA_\e$ of parameters for which the discrete set $\cM_{h,\e}$ satisfies \eref{discret}.  We now run the pure greedy algorithm on the compact set $\cK:=\cM_{h,\e}$ 
however using the surrogate $d(a,V_n)$ in place of $\dist(u(a), V_n)_V$. This means that the new
element $g_n=u_h(a^n)$ is defined by
\be
a^n:={\rm argmax}\{d(a, V_{n-1}) \; : \; a\in \cA_\e\}.
\ee
We  stop the algorithm at the first value $n=n(\e)$ for which 
\be
\max\{d(a, V_{n}) \; : \; a\in \cA_\e\} \leq   \frac \e {3\beta}
\ee
 The output of this perturbed greedy algorithm is our reduced basis space $V_n$.}
\nl

In view of the previous discussion, on issues 1, 2 and 3, the output space satisfies the
goal \iref{goal}. As an immediate consequence of Corollary \ref{C1}, we obtain one
first result on its number of steps $n(\e)$, which uses assumptions on the $n$-width of $\cM_h$.
 
 \begin{theorem}
 \label{Tng}
For the above algorithm, we have:
\vskip .1in
\noindent
{\rm (i)} For any $s>0$ and $C_0>0$,
\be
d_{n}(\cM_h)_V\le C_0(\max\{1,n\})^{-s}, \quad n\ge 0\quad \Rightarrow \quad n(\e) \le \( \frac \e {3\beta C_1}\)^{-1/s}, \quad \e> 0,
\ee 
where $C_1:= \gamma^{-2} 2^{4s+1}C_0$. 
\vskip .1in
\noindent
{\rm (ii)} For any $s>0$ and $c_0,C_0>0$, 
\be 
d_{n}(\cM_h)_V\le C_0e^{-c_0n^{s}}, \quad n\ge 0\quad \Rightarrow \quad n(\e)\le  \(\frac 1 {c_1}\max\Big \{\log\( \frac \e {3\beta \t C_1}\),0\Big\}\)^{1/s}, \quad e>0,
\ee 
where $c_1=\frac {c_0} 2 3^{-s}$ and $\t C_1:=C_0\max\{\sqrt{2}\gamma^{-1},e^{c_1}\}$.
\end{theorem}

\noindent
{\bf Proof:} This is a direct application of items (ii) and (iii) in Corollary \ref{C1}, 
using the fact that $d_n(\cM_{h,\e})_V\leq d_n(\cM_{h})_V$. \hfill $\Box$
\nl

Let us observe that the assumptions in the above theorem
are on the decay of the $n$-widths of $\cM_h$, in contrast to Corollary \ref{C1}
which uses assumptions on the decay of the $n$-widths of $\cM$.
As already explained in \S \ref{compint}, 
the approximate solution map $u_h$ may often be viewed as
the solution map of a 
discrete parametrized problem of the form \iref{genparh} 
with similar properties as the original parametric problem \iref{genpar}.
This allows us to apply the same techniques as in \S 4 in order
to evaluate $d_n(\cM_h)_V$ and justify the validity of the
assumptions in the above corollary for relevant instances of parametric PDEs.

In more general cases, we may be able justify the decay of $d_n(\cM)_V$
but not of $d_n(\cM_h)_V$. This occurs for example if the solver involves a different
finite element space for each instance, such as in adaptive methods. Then, we may still write
\be
d_n(\cM_h)_V\leq d_n(\cM)_V+\e(h).
\ee
This means, for example that if we start from the assumption that $d_n(\cM)_V\leq C_0(\max\{1,n\})^{-s}$,
we need to study the weak greedy algorithm applied to $\cM_{h,\e}$, however under the modified assumption 
\be
d_n(\cM_{h})_V\leq C_0(\max\{1,n\})^{-s}+\e(h).
\ee
We may then separate $n$ between the ranges $\{1,\dots,N(h)\}$ where $\e(h)\leq C_0(\max\{1,n\})^{-s}$
and the larger values of $n$. Then, having fixed $\e(h)=\e/3$ and using Remark \ref{remrange}, we reach a similar conclusion
on the order of magnitude of $n(\e)$ as in Theorem \ref{Tng}. The same holds for exponential rates.

\subsection{Computational cost}
\label{subcost}
 
We now turn to the analysis of the computational cost required 
by the numerical weak greedy algorithm in order to reach the 
accuracy goal \iref{goal}. For simplicity, we restrict our attention to
the regime of algebraic rates, that is described by item (ii) in Theorem \ref{Tng}.
A similar analysis can be carried out for exponential rates.
Here, we only consider linear elliptic problems
expressed in variational form, which are particular cases of those treated in \S 7: find $u\in V$ such that 
\be
B(u,v;a)=L(v),\quad v\in V.
\ee
where where $B(\cdot,\cdot;a)$ and $L$ are continuous sesquilinear
and antilinear forms over $V\times V$ and $\t V$ respectively, 
and where we make the additional assumption that
\be
a\mapsto B(\cdot,\cdot;a),
\ee 
is a continuous {\it linear} map $X$ to ${\frak B}$ the set of continuous sesquilinear forms over $V\times V$.
We work under the following 
symmetric elliptic version of {\bf Assumption AL}.
\nl
\nl
{\bf Assumption ALE:} {\it The parameter set $\cA$ has a complete affine representer 
$(\psi_j)_{j\geq 1}$ and, for all $a\in a(U)$, the sesquilinear form $B(\cdot,\cdot;a)$ satisfies
the coercivity conditions \iref{coerc} and it is symmetric when restricted to real valued functions of $V$.}
\nl

Under such an assumption, the approximate solution $u_h(a)\in V_h$ is defined by the
Galerkin method, that is,
\be
B(u_h(a),v_h;a)=L(v_h),\quad v_h\in V_h,
\ee
and can be computed for any given $a\in a(U)$ by the numerical solver at cost $C_h$.

We turn now to the online cost of the numerical weak greedy algorithm
assuming that we have already computed in the offline stage the reduced basis elements $g_k=u_h(a^k)$,
$k=0,\dots,n-1$, by using the possibly expensive finite element solver for $u_h$.  Given a query $a\in\cA$,  the online stage computes $u_n(a)\in V_n$,  where $V_n$ is the reduced basis space. We recall that $V_n\subset V_h$. We find $u_n(a)$
by the Galerkin method  for $V_n$, that is,
\be
B(u_n(a),v_n;a)=L(v_n),\quad v_n\in V_n,
\ee
This amounts in solving an $n\times n$ linear system, where the unknowns are the coefficients
$\alpha_l(a)$ in the decomposition
\be
u_n(a)=\sum_{l=0}^{n-1} \alpha_l(a) g_l.
\ee
Note that, as opposed to stiffness matrices resulting from the  discretization of PDEs
in a nodal finite element basis, the resulting
stiffness matrix 
\be
\bB_n(a)=(B(g_k,g_l;a))_{k,l=0,\dots,n-1}.
\ee
is generally full.
Using a direct solver,  such as Gauss elimination, the cost of solving this system
is therefore or order
\be
n(\e)^3\sim \e^{-3/s}.
\ee

However, we also need to take into account the cost of assembling the system,
that is, computing the above stiffness matrix which depends on $a$.
Since the data vector ${\bf F}_n:=(L(g_k))_{k=0,\dots,n-1}$ of this system does not depend on $a$,
its computation can be performed during the offline stage. In order to compute the
stiffness matrix, we recall the bilinear forms $\o B$, $B_j$ and $B(\cdot,\cdot,y)$
defined in \S \ref{tayrec}.  If $y\in U$  and 
\be
a=a(y)=\o a+\sum_{j\geq 1} y_j\psi_j,
\ee
the stiffness matrix is  
\be
\bB_n(y)=\o\bB_n+\sum_{j\geq 1} y_j\o\bB_{n,j},
\ee
where
\be
\o\bB_n:=(\o B(g_k,g_l))_{k,l=0,\dots,n-1}\quad {\rm and}\quad  \o\bB_{n,j}:= (B_j(g_k,g_l))_{k,l=0,\dots,n-1},
\ee
are $n\times n$ matrices. Each of these matrices can be computed in the offline stage, however
the infinite sum over $j\geq 1$ needs to be truncated at some prescribed level $J$. In the case where 
$(\|\psi_j\|_X)_{j\geq 1}$ is $\ell^p$ summable, and if the $\psi_j$ are
organized such that the $\|\psi_j\|_X$ are non-increasing with $j$,  then we then know
that the $L^\infty(U,V)$ error in the approximation of the solution map $y\mapsto u_h(y)$
resulting from this truncation is of the order $\cO(J^{-s})$ where $s:=\frac 1 p-1$, 
and therefore the order of accuracy $\e$ can be preserved by taking
\be
J=J(\e)\sim \e^{-1/s}.
\ee
We may thus incorporate such a truncation in the
definition of the approximation map $y\to u_h(y)$ used to handle {\bf Issue 1}.  Note that the choice of $J$ depends only on $\e$ and is independent of $h$.
Therefore, using this $u_h$, the conclusion of (i) in Theorem \ref{Tng} is retained and  the cost of assembling the system is
\be
J(\e) n(\e)^2\sim \e^{-3/s}. 
\ee
Note that, once the coefficents $\alpha_k(a)$ are found, 
computing the finite element representation of  
the solution $u_n(a)=\sum_{k=0}^{n-1} \alpha_k(a) g_k$ has a cost of $n(\e)N_h$.
In conclusion, the total online cost is of the order 
\be
C_{\rm on}(\e) \sim \e^{-3/s} +n(\e)N_h.
\ee

However, note that in some applications, one may only work with the reduced basis representation
$(\alpha_k(a))_{k=0,\dots,n-1}$, without the need to recompute the finite element
representation. This the case for instance when manipulating a quantity of interest
such as a linear scalar functional 
\be
Q(u_n(a))=\sum_{k=0}^{n-1} \alpha_k(a) Q(g_k).
\ee
Having pre-computed the quantities $Q(g_k)$ in the offline stage, the online evaluation of this quantity 
is therefore executed at cost of order $\e^{-3/s}$.

Also note that $u_n(a)$ is not the best approximation of $u(a)$ from $V_n$
in the norm $V$ since it is the Galerkin projection onto $V_n$, however, from  Cea's lemma 
one has
\be
\|u(a)-u_n(a)\|_V \leq \sqrt {\frac R \alpha} \min_{v\in V_n}Ê\|u(a)-v\|_V = \sqrt {\frac R \alpha}{\rm dist}(\cM,V_n)_V,
\ee
where $\alpha$ is the constant in \iref{coerc} and $R:=\max_{a\in\cA}Ê\|B(\cdot,\cdot;a)\|_{\frak B}$. This guarantees
that we reach an error of the prescribed order $\e$ between $u(a)$ and its reduced basis approximation $u_n(a)$.
\nl

We next discuss the offline cost.  The  first step is to find an $\e/3$ net $\cM_{h,\e}$ for $\cM_h$.   We describe this net only through the parameter set $\cA$, namely
as 
\be
\cM_{h,\e}=u_h(\cA_\e),
\ee 
where $\cA_\e$ is a finite subset of $\cA$ such that
 \be
\label{discreta}
 \sup_{a\in\cA}  \dist(a,\cA_{\e})_X \le  \frac {\e}{3C},
 \ee
and  $C$ is a Lipschitz constant for the map $a\mapsto u_h(a)$. 
 The same type of computation as done in \iref{lip} for the particular problem \iref{ellip}
 shows that 
 \be
 C=\frac {\|L\|_W}{\alpha^2},
 \ee
 is an admissible Lipschitz constant. This implies that the resulting $\cM_{h,\e}$ satisfies \iref{discret}.
 Note that we do not need to compute the elements of $\cM_{h,\e}$ but only the parameter values in $\cA_\e$.

 We can bound the cardinality of $\cA_\e$ from results on covering numbers and $n$-widths.  Let us recall that the covering number $N_\delta:=N_\delta(\cA,X)$ is the smallest number of $X$-balls of radius $\delta$ that cover $\cA$.   Let $B(a_i,\delta)$,  
 $i=1,\dots N_\delta$,  be such a covering. 
 Note that the $a_i$ need not be from $\cA$ but this is easily remedied.   
Namely, for any such ball we have $B(a_i,\delta)\cap \cA\neq \emptyset$, 
and so we choose an $a_i\in B(a_i,\delta)\cap \cA$.  Then the balls $B(a_i,2\delta)$ are a covering of $\cA$
with centers from $\cA$.  Therefore, taking $\eta=\frac {2\e}{3C}$, we can find a set
$\cA_\e\subset \cA$ satisfying \eref{discreta} with
\be
\label{boundA'} 
\#(\cA_\e)\le N_{\eta/2}(\cA,X) = N_{\frac {\e}{3C}}(\cA,X).
\ee
The   well-known  Carl's inequality  \cite{P}  gives a bound on the covering numbers $N_\delta(\cA,X)$
in terms of the $n$-widths $d_n(\cA)_X$ of $\cA$.  In our case,
 this inequality gives   
 \be
 \label{carl}
 N_\eta(\cA,L^\infty)\le C_12^{\eta^{-1/s}}, 
 \ee
 where $C_1$ is a constant depending on $s$.   
This gives us the bound
 \be
\label{boundA'} 
\#(\cA_\e)\le C_1 2^{c_1\e^{-1/s}},
\ee
for a constant $c_1$ that also depends on $s$. While it is generally not possible, at least in any reasonable way, to find a minimal set $\cA_\e$, in typical settings we can give a simple description of a set $\cA_\e$ so that \eref{boundA'} is still satisfied for an appropriate constant $C_1$.
For example, whenever we can construct   
a sequence of spaces $W_n$ for which $\dist(\cA,X_n)_{X}=\cO(n^{-s})$, then the proof of Carl's inequality
(see e.g. \cite{LGM}) gives an explicit description of such an $\cA_\e$.   
In particular, under the assumption $(\|\psi_j\|_{X})\in \ell^p$
and with the $\psi_j$ organized in non-increasing $X$ norms, we can take $X_n:=\span\{\psi_1,\dots,\psi_n\}$, for each $n\ge 1$,
and the description of $\cA_\e$ amounts in defining a specific lattice discretization $U_\e$ of $U$ such that $\cA_\e=a(U_\e)$.

Let us now evaluate the cost of the $k$-th step of the numerical weak-greedy algorithm. 
This step includes the computation of the reduced basis element $g_k:=u_h(a^k)$ once $a^k$ has
been chosen, using the possibly expensive finite element solver,
which has cost of order $C_h$.  On the other hand, we must also account for the 
maximization of the surrogate $d(a,V_{k-1})_V$ over the set $\cA_e$. This cost is of order
\be
\#(\cA_\e) s_k,
\ee
where $s_k$ is the cost of computing $d(a,V_{k-1})_V$ for one value of $a$.

We now give a derivation of a possible surrogate and evaluate the cost $s_k$ for this particular surrogate. Since for the reduced basis
solution $u_k(a)\in V_k$, we have
\be
\sqrt{ \frac \alpha R} \|u_h(a)-u_k(a)\|_V \leq {\rm dist}(u_h(a),V_k)_V \leq \|u_h(a)-u_k(a)\|_V 
\label{cea1}
\ee
this surrogate should be an equivalent quantity to $\|u_h(a)-u_k(a)\|_V$.
We introduce the $N_h\times N_h$ stiffness
matrix
\be
\bB_h(y)=\o \bB_h+\sum_{j=1}^J y_j \o \bB_{h,j},
\ee
for the sesquilinear form $B(\cdot,\cdot;y)$ in the finite element basis, where $\o \bB_h$ and $\bB_{h,j}$ are
the corresponding stiffness matrices for $\o B$ and $B_j$. Therefore, the coordinate vector 
$U_h(y)$ of $u_h(y)=u_h(a(y))$ in the finite element basis is the solution of the $N_h\times N_h$ system
\be
\bB_h(y)U_h(y)=F_h,
\ee
where the right side $F_h$ does not depend on $y$. Here $J=J(\e)$ is the truncation level, already introduced 
for the evaluation of the online cost and $u_h$ is defined
as the discrete solution of the trunctated problem. We also introduce the coordinate vectors $G_i$
of the reduced basis elements $u_h(a^i)$ in the finite element basis. Therefore, a reduced basis solution
$u_k(y)=u_k(a(y))$ is represented in the finite element basis by the vector
\be
U_k(y)=\sum_{i=0}^{k-1} \alpha_i(y) G_i.
\ee
We introduce an hilbertian norm $\|\cdot\|_*$ on $\R^{N_h}$, defined in such way that 
\be
\|W_h\|_*:=\|w_h\|_V.
\label{norm}
\ee
whenever $W_h$ is the coordinate vector of $w_h\in V_h$. Note that this would coincide with the euclidean norm if
the finite basis were orthonormal in $V$. This is generally not the case, but nevertheless the computation of this norm
is usually of complexity $N_h$. We may thus write
\be
\|u_h(y)-u_k(y)\|_V = \|U_h(y)-U_k(y)\|_*.
\ee
By \iref{cea1}, it follows that
\be
 \frac {\alpha} {R}\|U_h(y)-U_k(y)\|^2_* \le {\rm dist}(u_h(y),V_k)_V^2 \leq  \|U_h(y)-U_k(y)\|^2_*.
\ee
Our next observation is that since for any $y\in U_\cA$,
\be
 \frac {\alpha}{R}\<\o \bB_h W_h,W_h\>\leq \<\bB_h(y) W_h,W_h\>\leq\frac {R}{\alpha} \<\o \bB_h W_h,W_h\>,
\ee
one has the norm equivalence
\be
\frac {\alpha}{R}\|W_h\|_*\leq  \|\o \bB_h^{-1}\bB_h(y) W_h\|_* \leq \frac {R}{\alpha} \|W_h\|_*.
\ee 
Therefore, we can define a surrogate quantity by
\be
d(y,V_{k-1})_V:=\|\o \bB_h^{-1}\bB_h(y)(U_h(y)-U_k(y))\|_*,
\ee
and obtain the equivalence \iref{surrogate1}Ê with constants
\be
\delta= \(\frac {\alpha}{R}\)^{3/2} \quad {\rm and}\quad \beta=\frac R {\alpha}.
\ee 
This surrogate is computable
since we have
$$
\begin{disarray}{ll}
d(y,V_{k-1})^2_V & =\Big \|\o \bB_h^{-1} F_h-\o \bB_h^{-1}\bB_h(y) U_k(y)\Big\|^2_*\\
& =Ê\Big \|\o \bB_h^{-1} F_h-\o \bB_h^{-1}\(\o \bB_h+\sum_{j=0}^J y_j \o \bB_{h,j}\)\sum_{i=0}^{k-1} \alpha_i(y) G_i\Big \|^2_*\\
\end{disarray}
$$
Developing this square norm, we find that it
the sum of the constant term $\|\o \bB_h^{-1} F_h\|^2_*$
and of a linear combination of the real numbers
$\alpha_i(y)$, $\alpha_i(y)\alpha_{i'}(y)$, $y_j\alpha_i(y)$ and $y_jy_{j'}\alpha_i(y)\alpha_{i'}(y)$
for $i,i'=0,\dots,k-1$ and $j,j'=1,\dots,J$. The coefficients of these linear combinations are
given by the $\<\cdot,\cdot\>_*$ inner products (associated to the $\|\cdot\|_*$ norm)  
between pairs of vectors chosen from
\be
\o \bB_h^{-1} F_h, \quad  G_i, \quad  \o \bB_h^{-1}\o \bB_{h,j} G_i, \quad i=0,\dots,k-1, \; j=1,\dots, J.
\ee 
The precomputation of these vectors
and of their inner product has a cost of order
\be
kJC_h+ k^2J^2N_h.
\ee
Then, the computation of the surrogate $d(y,V_{k-1})_V$ for each $y$ has cost of order $k^2J^2$ for the linear combination
to which we must add the cost of computing the $\alpha_i(y)$, which according to the discussion on the online cost
is of order $k^3$. Therefore
\be
s_k\sim k^2 J^2+k^3.
\ee
In summary, the total cost of step $k$ of the algorithm, without including the precomputations, is 
\be
C_h+\#(\cA_\e) (k^2J^2+k^3)
\ee
so that the total cost up to step $n=n(\e)$ is of order
\be
n(\e)C_{h(\e)}+\#(\cA_\e) n(\e)^3J(\e)^2+\#(\cA_\e)n(\e)^4.
\ee
We need to addÊ the cost of precomputing:
\begin{itemize}
\item[(i)]
the vectors 
$\o \bB_h^{-1} F_h$,  $G_i$ and $\o \bB_h^{-1}\o \bB_{h,j} G_i$,
for  $i=0,\dots,n-1$ and $j=1,\dots, J$ and their $\<\cdot,\cdot\>_*$ inner products.
\item[(ii)]
the matrices $\o \bB_k$ and
$\bB_{k,j}$ for $k=0,\dots,n-1$, which entries are 
given by the euclidean inner product between the vectors $G_k$ and 
the vectors $\o \bB_hG_i$ and $\bB_{h,j}G_i$. 
\end{itemize}
This precomputing cost
is of total order 
\be
n(\e)J(\e)C_{h(\e)}+n(\e)^2J(\e)^2N_{h(\e)}.
\ee
In summary, the total offline cost is of order
\be
C_{\rm off}\sim n(\e)J(\e)C_{h(\e)}+n(\e)^2J(\e)^2N_{h(\e)}+\#(\cA_\e) n(\e)^3J(\e)^2+\#(\cA_\e)n(\e)^4.
\ee
Among these terms, the largest is typically the third one which in our algebraic rate regime is of order
$\e^{-5/s}2^{c_1\e^{-1/s}}$ in view of \iref{boundA'}.

This offline cost is thus potentially extremely large. Note however that it is due to the fact that we are
using a brutal discrete search over $\cA_\e$ for the maximization of the surrogate quantity, so that there is room
for improvement by using more sophisticated optimization strategies. Note also
that in the case of a parametric problem with moderate number $d$ of parameters, the quantity $J(\e)$
can simply be replaced by $d$.

In conclusion, we find that, compared to the polynomial methods discussed in \S 6 and \S 7, 
the reduced basis method suffers from a very high offline cost, especially in high parametric dimension.
This can be compensated by the fact that this method captures the same rate of decay 
as achieved by the optimal $n$-width spaces, so that a prescribed accuracy $\e$
may be achieved with a number $n=n(\e)$ of reduced basis elements 
much smaller than the number of terms in polynomial expansions for the same accuracy,
making the online cost potentially lower.

\end{document}